\documentstyle[11pt]{report}
\input amssym.def
\input amssym
\newcommand{\x}{X_1,\dots ,X_n}
\newcommand{\y}{Y_1,\dots,Y_r}

\newcommand{\xx}{x_1,\dots ,x_n}
\newcommand{\yy}{y_1,\dots ,y_r}
\newcommand{\X}{\rm X}
\newcommand{\Y}{\rm Y}

\newcommand{\fm}{\frak m}
\newcommand{\fn}{\frak n}
\newcommand{\fp}{\frak p}
\newcommand{\g}{e}
\newcommand{\A}{\alpha}

\newcommand{\Kk}{\mbox{\rm $k[(I^e)_c]$} }
\newcommand{\Rr}{\mbox{\rm $R_A (I)$} }

\newcommand{\Rdd}{\mbox{\rm $R_\Delta$}}

\newcommand{\ExsL}{\mbox{\rm $0\to {D_t}\to \cdots \to {D_1} \to
{D_0} \to L \to 0$}\,}

\newcommand{\B}{\Box}

\newcommand{\dg}{\mbox{\rm{deg}} \,}
\newcommand{\dpp}{\mbox{\rm proj.dim} }

\newcommand{\rdim}{\mbox{\rm{rel.dim}} \,}
\newcommand{\Ker}{\mbox{\rm{Ker}} \,}
\newcommand{\Sym}{\mbox{\rm{Sym}} \,}
\newcommand{\Ann}{\mbox{\rm{Ann}} \,}
\newcommand{\rad}{\mbox{\rm{rad}} \,}
\newcommand{\Tor}{\mbox{\rm{Tor}} \,}
\newcommand{\lL}{\rm length}
\newcommand{\length}{\rm length}
\newcommand{\Hom}{\mbox{\rm{Hom}} \,}
\newcommand{\Ext}{\mbox{\rm{Ext}} \,}
\newcommand{\Supp}{\mbox{\rm{Supp}} \,}
\newcommand{\supp}{\mbox{\rm{supp}} \,}

\newcommand{\reg}{\mbox{\rm{reg}} \,}
\newcommand{\Ass}{\mbox{\rm{Ass}} \,}
\newcommand{\uHom}{  {\rm {\underline{Hom}}}   }
\newcommand{\uH}{  {\rm {\underline{H}}}   }
\newcommand{\uExt}{{\rm {\underline{Ext}}}}

\newcommand{\h}{\mbox{\rm{ht}} \,}
\newcommand{\Proj}{\mbox{\rm Proj } }
\newcommand{\Sp}{\mbox{\rm Spec} }
\newcommand{\Spec}{\mbox{\rm Spec} }
\newcommand{\depth}{\mbox{\rm depth} }

\newcommand{\Ll}{\cal L}
\newcommand{\M}{\cal M}
\newcommand{\MD}{{\cal M}_{\Delta}}
\newcommand{\F}{ F_{\frak m}(I)}
\newcommand{\p}{\rm P}
\newcommand{\pq}{ {(p, q)} }
\newcommand{\mq}{ {(m, q)} }
\newcommand{\nn}{\overline n}
\newcommand{\bn}{{\bf n}}
\newcommand{\bm}{{\bf m}}
\newcommand{\bd}{{\bf d}}

\newcommand{\XX}{\underline X}

\newcommand{\dd}{\delta}

\newtheorem{thm}{Theorem}[section]{\bf}{\it}
\newtheorem{prop}[thm]{Proposition}{\bf}{\it}
\newtheorem{lem}[thm]{Lemma}{\bf}{\it}
\newtheorem{cor}[thm]{Corollary}{\bf}{\it}
\newtheorem{rem}[thm]{Remark}{\bf}{\rm}
\newtheorem{defn}[thm]{Definition}{\bf}{\it}
\newtheorem{ex}[thm]{Example}{\bf}{\rm}
{\bf}{\it}

\newcommand{\pf}{\noindent {\bf Proof. }}

\begin{document}

\title{{\huge   {\bf On the diagonals of a Rees algebra}} \\
\vspace{5mm} {\huge Olga Lavila--Vidal}\\
\vspace{20mm} PhD Thesis \\
Universitat de Barcelona (1999)}

\maketitle

\pagenumbering{roman}

\tableofcontents


\chapter*{Introduction}
\typeout{Introduction}

\addcontentsline{toc}{chapter}{Introduction}

\markboth{INTRODUCTION}{INTRODUCTION}

\bigskip
The aim of this work is to study the ring-theoretic properties of
the diagonals of a Rees algebra, which from a geometric point of view
are the homogeneous coordinate rings of embeddings of blow-ups of
projective varieties along a subvariety.
First we are going to introduce the subject and the main problems.
After that we shall review the known results about these
problems, and finally we will give a summary of the contents and
results obtained in this work.

\medskip
Let $A$ be a noetherian graded algebra generated over a field $k$ by
homogeneous elements of degree $1$, that is, $A$ has a presentation
$A= k[{\x}]/K= k[{\xx}]$, where $K$ is a homogeneous ideal of the
polynomial ring $k[{\x}]$ with the usual grading. Given a homogeneous
ideal $I$ of $A$, let $X$ be the projective variety obtained by
blowing-up the projective scheme $Y= \Proj (A)$ along the sheaf of
ideals ${\cal I}= \widetilde{I}$, that is,
$X= {\cal P}roj (\bigoplus_{n \geq 0}{\cal I}^n)$. For a given
 $c \in \Bbb Z$, let us denote by $I_c$ the
$c$-graded component of $I$. If $I$ is generated by forms of degree
less or equal than $d$, then  $(I^e)_c$ corresponds to a complete linear
system on $X$ very ample for $c \geq de+1$ which embeds $X$ in a
projective space $X \cong \Proj (\Kk) \subset \Bbb P^{N-1}_k$, with $N=
\dim_k (I^e)_c$ \cite[Lemma 1.1]{CH}.

\medskip
Our main purpose is to study the arithmetic properties of
the $k$-algebras $\Kk$, where $c,e$ are positive integers and
$I$ is any homogeneous ideal of $A$.
This problem was first started in the work by A. Gimigliano
\cite{G}, A. Geramita and A. Gimigliano \cite{GG}, and
A. Geramita, A. Gimigliano and B. Harbourne \cite{GGH} who treated
similar problems for the rational projective surfaces which
arise as embeddings of blow-ups of a projective plane at a set of
distinct points.

\medskip
Let $k$ be an algebraically closed field and
$s = {d+1 \choose 2}$, $d \geq 2$. In \cite{G} the particular case of
the blow-up of $\Bbb P^2_k$ at a set of  $s$ different
points $P_1, \dots, P_s$ which do not lie on a curve of degree $d-1$
and such that there is no subset of $d$ points on a line (if $d \geq 3$)
is studied in detail.
In this case, the defining ideal $I$ of the set of points is
generated by forms of degree $d$ and the rational maps defined by the
linear systems $I_c$ give embeddings of the blow-up for $c \geq d$.
In the case $c=d$ the surface obtained is called {\it White Surface},
and for $c=d+1$ {\it Room Surface}. It is then shown that White surfaces are
contained in $\Bbb P^d_k$ as surfaces of degree ${ d \choose 2}$ with
defining ideal generated by the maximal minors of a $3 \times d$
matrix of linear forms. In particular, $k[I_d]$ is Cohen-Macaulay and
it has a resolution given by the Eagon-Northcott complex
\cite[Proposition 1.1]{G}. On the other hand, Room Surfaces are
arithmetically Cohen-Macaulay
\cite[Theorem B]{GG} with defining ideal generated by quadrics
\cite[Theorem 1.2]{GG}.

\medskip
This detailed study of White and Room Surfaces is the first
step to consider the following more general case. Let
$P_1, \dots, P_s$ be $s$ distinct points in $\Bbb P^2_k$, with $k$
an algebraically closed field, let $I$ be its defining ideal and $d =
\reg (I)$ the regularity of $I$. Assume that the points do not lie on a
curve of degree $d-1$ and that there is no subset of $d$ points on a
line. Then the linear systems $I_c$ give embeddings of the blowing-up of $\Bbb
P^2_k$ at this set of points for $c \geq d$. The resultant surfaces
are arithmetically Cohen-Macaulay \cite[Theorem B]{GG} and its defining
ideal is defined by quadrics if $c \geq d+1 $ \cite[Theorem 2.1]{GG}.

\medskip
Even more generally, A. Geramita, A. Gimigliano and Y. Pitteloud
\cite{GGP} consider the blow-up of $\Bbb P^n_k$ along an ideal of fat
points, with $k$ an algebraically closed field of characteristic zero.
Given a set of points $P_1, \dots, P_s \in \Bbb P^n_k$, let
${\cal P}_1, \dots, {\cal P}_s \subset k[X_0, \dots, X_n]$ be their
defining ideals, and let us take ideals of the type
$I= {\cal P}_1^{m_1} \cap \dots \cap {\cal P}_s^{m_s}$, with $m_1,
\dots, m_s \in \Bbb Z_{\geq 1}$. Then one may study the
projective varieties obtained by embeddings of the blow-up of $\Bbb
P^n_k$ along $\cal I$ via the linear systems corresponding to the graded
pieces of $I$, whenever these linear systems are very ample.
Let $d=\reg (I)$, and let us assume that there are not $d$
points on a line. Then the linear systems $I_c$ are very ample for $c
\geq d$, and the varieties obtained via these embeddings are
projectively normal \cite[Proposition 2.2]{GGP} and arithmetically
Cohen-Macaulay \cite[Theorem 2.4]{GGP}.

\medskip
A new point of view to treat these questions was introduced
by A. Simis, N.V. Trung and G. Valla in \cite{STV}, and later
followed by A. Conca, J. Herzog, N.V. Trung and G. Valla in \cite{CHTV},
to study the more general problem of the blow-up of a projective
space along an arbitrary subvariety.
If $I$ is a homogeneous ideal of $A$, let us consider the Rees
algebra
$R_A(I) = \bigoplus_{n \geq 0} I^n \cong A[It] \subset A[t]$ of $I$
with the natural bigrading given by
$$R_A(I)_{(i,j)} = (I^j)_i.$$
The crucial point now is that all the algebras $\Kk$ are
subalgebras of the Rees algebra in a natural way.
To describe this relationship we need to introduce the diagonal functor.

\medskip
Given positive integers $c, e$, the $(c,e)$-{\it diagonal} of $\Bbb Z^2$
is the set $${\Delta} := \{ (cs, es) \mid s \in {\Bbb Z}\}.$$
For any bigraded algebra $S= \bigoplus_{(i,j) \in \Bbb Z^2} S_{(i,j)}$,
the {\it diagonal} {\it subalgebra} of $S$ along $\Delta$ is the graded
algebra
$$S_{\Delta} := \bigoplus_{s \in \Bbb Z} S_{(cs, es)}.$$
Similarly we may define the diagonal of a bigraded $S$-module $L$ along
$\Delta$ as the graded
$S_\Delta$-module
$$L_{\Delta} := \bigoplus_{s \in \Bbb Z} L_{(cs, es)}.$$
So we have an exact functor
$$(\;\;)_{\Delta} : M^2(S) \to M^1(S_{\Delta}),$$
where $M^2(S)$, $M^1(S_{\Delta})$ denote the categories of bigraded
$S$-modules and graded $S_{\Delta}$-modules respectively.

\medskip
Now we may give a description of the rings $\Kk$ as diagonals
of the Rees algebra in the following way: By taking $\Delta$ to be the
$(c,e)$-diagonal of $\Bbb Z^2$ with $c \geq de+1$, we have
$$R_A(I)_{\Delta} = \bigoplus_{s \geq 0} (I^{es})_{cs} =
k[(I^e)_c].$$
This observation allows an algebraic approach to study the rings
$\Kk$ via the diagonals of $R_A(I)$. This is the starting point in
\cite{STV} to study the case of homogeneous ideals of the polynomial
ring generated by forms of the same degree, and later in \cite{CHTV} to
study arbitrary homogeneous ideals of the polynomial ring. By
paraphrasing \cite{STV}:
{\sl One is to believe that the algebraic approach via the diagonals of
the Rees algebra may throw further light not only on the study of
embedded rational surfaces obtained by blowing-up a set of points in
$\Bbb P^2_k$ but also of the embedded rational $n$-folds obtained, more
generally, by blowing-up $\Bbb P^n_k$ along some special smooth
subvariety}.
On the other hand, the diagonals of any standard bigraded algebra
defined over a local ring have also been studied by E. Hyry \cite{H} by
using both an algebraic approach and a geometric approach. Finally,
S.D. Cutkosky and J. Herzog \cite{CH} have studied the
diagonals of the Rees algebra of a homogeneous ideal in a general
graded $k$-algebra.

\vspace{7mm}
Next we are going to expose the main results of those works.

\medskip
The main contribution of A. Simis et al. \cite{STV}
to the problems considered by A. Geramita et al. is the algebraic
approach via the diagonal of a bigraded algebra, a notion which
generalizes the Segre product of graded algebras.
Given algebraic varieties $V \subset \Bbb P^{n-1}_k$, $W \subset
\Bbb P^{r-1}_k$ with homogeneous coordinate rings $R_1$, $R_2$, the
image of $V \times W \subset \Bbb P^{n-1}_k \times \Bbb P^{r-1}_k $
under the Segre embedding
$$\Bbb P^{n-1}_k \times \Bbb P^{r-1}_k \hookrightarrow \Bbb
P^{nr-1}_k$$
is a variety with homogeneous coordinate ring the Segre product of
$R_1$ and $R_2$:  $$R_1 \underline \otimes_k R_2 =
\bigoplus_{u \in \Bbb N} (R_1)_u \otimes_k (R_2)_u.$$
Given a standard bigraded $k$-algebra
$R= \bigoplus_{(u,v) \in \Bbb N^2} R_{(u,v)}$, its diagonal
$R_\Delta$ is defined as
$R_\Delta= \bigoplus_{u \in \Bbb N} R_{(u,u)}$ (that is, the
$(1,1)$-diagonal). By considering the tensor product
$R=R_1 \otimes_k R_2$ bigraded by means of
$R_{(u,v)} = (R_1)_u \otimes_k (R_2)_v$, we have that
$R_\Delta = R_1 \underline  \otimes_k R_2$.
Classically $R$ is taken to be the bihomogeneous coordinate ring of a
projective subvariety of $\Bbb P^{n-1}_k \times \Bbb P^{r-1}_k$, and
$R_\Delta$ is then
the homogeneous coordinate ring of its image via the Segre embedding.

\medskip
In the first section of \cite{STV}, a relation between the
presentations, the dimensions and the multiplicities of a standard
bigraded $k$-algebra $R$ and its diagonal $R_\Delta$ is obtained.
The key for proving these results
is the existence of the Hilbert polynomial of a standard bigraded
$k$-algebra and the characterization of its degree, due to
D. Katz et al. \cite{KMV} and M. Herrmann et al. \cite{HHRT} among
others. Similarly to the graded case, one may define in this case
the irrelevant ideal, the irrelevant primes and the biprojective scheme
associated to a standard bigraded $k$-algebra.

\medskip
After that, it is studied the behaviour of the normality and
the Cohen-Macaulay property by taking diagonals.
Since there is a Reynolds operator from $R$ to $R_\Delta$, one
immediately gets that the normality of $R$ will be inherited by its
diagonal $R_\Delta$. With respect to the Cohen-Macaulayness,
the strategy is to reduce the problem to a special situation where
the diagonal subalgebra becomes a Segre product, case in which it is
known a criterion for the Cohen-Macaulayness.

\medskip
These results are then applied to the study of the Rees algebra $R_A(I)$
of a homogeneous ideal  $I \subset A=k[{\x}]$ generated by forms of the
same degree $d$ ({\it equigenerated} ideals). In this situation, the
Rees algebra can be bigraded so that becomes standard by means of
$$R_A(I)_{(i,j)} =(I^j)_{i+dj},$$
and then $R_A(I)_\Delta = k[I_{d+1}]$.
Mainly, two classes of ideals are then considered in detail: For
complete intersection ideals generated by a regular sequence of $r$
forms of degree $d$ it is shown that
$k[I_{d+1}]$ is a Cohen-Macaulay ring if $(r-1)d<n$, while $k[I_{d+1}]$
is not a Cohen-Macaulay ring if $(r-1)d>n$ \cite[Theorem 3.7]{STV};
for straightening closed ideals under some restrictions it is shown
that $k[I_{d+1}]$ is a Cohen-Macaulay ring \cite[Theorem 3.13]{STV}.
This second class of ideals includes for instance the determinantal
ideals generated by the maximal minors of a generic matrix.

\vspace{5mm}
As a natural sequel of the work above, A. Conca et al. study in
\cite{CHTV} the diagonals $R_\Delta$ of a bigraded $k$-algebra
$R$ for $\Delta =(c,e)$, with $c,e$ positive integers. The main problem
considered there is to find suitable conditions on $R$ such that certain
algebraic properties of $R$ are inherited by some diagonal
$R_\Delta$, mostly with respect to the Cohen-Macaulay property and
the Koszul property.
Their goal is to apply the results to the case of a
standard bigraded $k$-algebra or the Rees algebra of any homogeneous
ideal $I$ of $A=k[{\x}]$. In the first case, $R$ has a presentation
as a quotient of a polynomial ring $S=k[{\x}, {\y}]$ endowed with
the grading given by $\deg(X_i)=(1,0)$, $\deg(Y_j)=(0,1)$.
As for the Rees algebra, if $I$ is generated by forms $f_1, \dots, f_r$
of degrees $d_1, \dots, d_r$ respectively, we have a natural bigraded
epimorphism
$$\begin{array}{ccc}
S=k[{\x}, {\y}]  & \longrightarrow & R=R_A(I) \\
X_i & \mapsto & X_i \\
Y_j & \mapsto & f_j t
\end{array}$$
where $\deg(X_i)=(1,0)$, $\deg(Y_j)=(d_j,1)$.
Therefore, by working in the category of bigraded $S$-modules for
$S=k[{\x},{\y}]$ the polynomial ring with
$\deg(X_i)=(1,0)$, $\deg(Y_j)=(d_j,1)$, $d_1, \dots, d_r \geq 0$, one
may study both cases at the same time. Let us
denote by $\cal M$ and $m= {\cal M}_\Delta$ the homogeneous maximal
ideals of $S$ and $S_\Delta$ respectively. Denoting by $d= \max \{d_1,
\dots, d_r \}$, we will consider diagonals $\Delta=(c,e)$ with $c \geq
de+1$.

\medskip
Since the arithmetic properties of a module can be often characterized
in terms of its local cohomology modules, it is of interest to
study the local cohomology of the diagonals $L_\Delta$ of any finitely
generated bigraded $S$-module $L$. This is
done from the bigraded minimal free resolution of $L$ over $S$: Let
$$0 \to D_l \to \dots \to D_0 \to L \to 0$$
with $D_p= \bigoplus_{(a,b) \in \Omega_p} S(a,b)$ be the bigraded
minimal free resolution of $L$ over $S$. By taking
diagonals one gets a graded resolution of $L_\Delta$
$$ 0 \to (D_l)_\Delta \to \dots \to (D_0)_\Delta \to L_\Delta \to 0,$$
with $(D_p)_\Delta = \bigoplus_{(a,b) \in \Omega_p} S(a,b)_\Delta$.
The first step is then the computation of the local cohomology modules
of the $S_\Delta$-modules $S(a,b)_\Delta$, which is done in the frame of
a more general study about the local cohomology of the Segre product of
two bigraded $k$-algebras. In particular, it is obtained a criterion for the
Cohen-Macaulay property of $S(a,b)_\Delta$ by means of $a$, $b$ and
$\Delta$. We say that the resolution of $L$ is good if every module
$(D_p)_\Delta$ is Cohen-Macaulay for large diagonals $\Delta$. Then it
is stated the following theorem:

\medskip
\noindent
{\bf Theorem} \cite[Theorem 3.6, Lemma 3.8]{CHTV}
{\it Assume $n \geq r$. For any finitely generated bigraded $S$-module
$L$, there exists a canonical morphism
$$\varphi^q_L : H^q_m(L_\Delta) \rightarrow H^{q+1}_{\cal M}(L)_\Delta,
\; \forall q \geq 0$$ such that
\begin{enumerate}
\item $\varphi^q_L$ is an isomorphism for $q >n$.
\item $\varphi^q_L$ is a quasi-isomorphism for $q \geq 0$.
\item If $L$ has a good resolution, $\varphi^q_L$ is an isomorphism for
large diagonals.
\end{enumerate}}

\medskip
As a corollary one gets necessary and sufficient conditions for
the existence of Cohen-Macaulay or Buchsbaum diagonals $L_\Delta$ of $L$
in terms of the graded pieces of the local cohomology modules of $L$.

\medskip
Given a standard bigraded $k$-algebra $R$,
one may define the graded $k$-subalgebras
${\cal R}_1 = \bigoplus_{i \in \Bbb N} R_{(i,0)}$,
${\cal R}_2 = \bigoplus_{j \in \Bbb N} R_{(0,j)}$.
The following result gives a criterion for the Cohen-Macaulay property
of the diagonals of $R$ by means of ${\cal R}_1$ and ${\cal R}_2$.
Namely,

\medskip
\noindent
{\bf Theorem} \cite[Theorem 3.11]{CHTV}
{\it Let $R$ be a standard bigraded Cohen-Macaulay $k$-algebra. If
the shifts in the resolutions of ${\cal R}_1$ and ${\cal R}_2$
are greater than $-n$ and $-r$ respectively, then $R_\Delta$
is Cohen-Macaulay for large $\Delta$ .}

\medskip
In particular, they get the following corollary:

\medskip
\noindent
{\bf Corollary} \cite[Corollary 3.12]{CHTV}
{\it Let $R$ be a standard bigraded Cohen-Macaulay $k$-algebra. If
${\cal R}_1$, ${\cal R}_2$ are Cohen-Macaulay with
$a({\cal R}_1), a({\cal R}_2)<0$, then $R_\Delta$ is
Cohen-Macaulay for large $\Delta$ .}

\medskip
This result applied to Rees algebras of equigenerated ideals gives
a criterion for the Cohen-Macaulay property of their diagonals.

\medskip
Furthermore, the study done in \cite{STV} for the $(1,1)$-diagonal of
the Rees algebra of an equigenerated complete intersection ideal is
completed and extended to any complete intersection ideal and
any diagonal, by determining exactly
which are the Cohen-Macaulay diagonals. This is the only case where non
equigenerated ideals are considered.

\medskip
\noindent
{\bf Theorem} \cite[Theorem 4.6]{CHTV}
{\it Let $I \subset A=k[{\x}]$ be a homogeneous complete intersection
ideal minimally generated by $r$ forms of degrees $d_1, \dots, d_r$.
Set $u= \sum_{j=1}^r d_j$. For $c \geq de+1$,
$\Kk$ is a Cohen-Macaulay ring if and only if $c> d(e-1)+u-n$.}

\medskip
About the Cohen-Macaulay property of the diagonals of a Rees
algebra is conjectured the following fact:

\medskip
\noindent
{\bf Conjecture}
{\it Let $I \subset A=k[{\x}]$ be a homogeneous ideal.
If $R_A(I)$ is a Cohen-Macaulay ring, then there exists a diagonal
$\Delta$ such that $R_A(I)_\Delta$ is a Cohen-Macaulay ring.}

\medskip
With respect to the Gorenstein property, there is just one statement
referred to the diagonals of the Rees algebra of a homogeneous ideal
generated by a regular sequence of length $2$.

\medskip
\noindent
{\bf Proposition} \cite[Corollary 4.7]{CHTV}
{\it Let $I \subset A=k[{\x}]$ be a homogeneous complete intersection
ideal minimally generated by two forms of degree $d_1 \leq d_2$.
If $n \geq d_2+1$, $k[I_n]$ is a Gorenstein ring with $a$-invariant
$-1$.}

\medskip
Finally, it is shown that large diagonals of the Rees algebra are
always Koszul:

\medskip
\noindent
{\bf Theorem} \cite[Corollary 6.9]{CHTV}
{\it
Let $I \subset A=k[{\x}]$ be a homogeneous ideal generated by forms of
degree $\leq d$. Then there exist integers $a, b$ such that
$k[(I^e)_{c+de}]$ is Koszul for all $c \geq a$ and $e \geq b$.}

\vspace{5mm}
Under a slightly different setting, E. Hyry \cite{H} is concerned with
comparing the Cohen-Macaulay property of the biRees algebra $R_A(I,J)$
with the Cohen-Macaulay property of the Rees algebra $R_A(IJ)$, where
$I,J \subset A$ are ideals of positive height in a local ring. To this
end, he studies the $\Delta=(1,1)$-diagonal
of any standard bigraded ring $R$ defined over a local ring.
The main result
\cite[Theorem 2.5]{H} gives necessary and sufficient conditions for the
Cohen-Macaulayness of a standard bigraded ring $R$ with negative
$a$-invariants by means of the local cohomology of the modules
$R(p,0)_\Delta$ and $R(0,p)_\Delta$ ($p \in \Bbb N$). In particular,
it provides sufficient conditions on $R$ so that the
Cohen-Macaulay property is carried from $R$ to $R_\Delta$:

\medskip
\noindent
{\bf Theorem}
{\it Let $R$ be a standard bigraded ring defined over a local ring.
Suppose that $\dim {\cal R}_1, \dim {\cal R}_2 < \dim R$ and
$a^1(R), a^2(R)<0$. If $R$ is Cohen-Macaulay, then so is $R_\Delta$ for
$\Delta=(1,1)$.}

\vspace{5mm}
Now let $A$ be a noetherian graded $k$-algebra generated in degree $1$
and let $I \subset A$ be a homogeneous ideal. The general problem of
studying the embeddings of the blow-up
$X= {\cal P}roj (\bigoplus_{n  \geq 0} {\cal I}^n)$ of the projective
scheme $Y= \Proj (A)$ along the sheaf of ideals ${\cal I}=
\widetilde{I}$ given by the graded pieces of $I$ is treated by S.D.
Cutkosky and J. Herzog \cite{CH}. They are mainly concerned with
the existence of an integer $f$ such
that $\Kk$ is Cohen-Macaulay for all $e>0$ and $c \geq ef$. The first
example considered is the blow-up of a smooth projective variety
$Y$ along a regular ideal in a field of characteristic zero, where
the Kodaira Vanishing Theorem can be used to prove:

\medskip
\noindent
{\bf Theorem} \cite[Theorem 1.6]{CH}
{\it Suppose that $k$ has characteristic zero, $A$ is Cohen-Macaulay,
$Y$ is smooth, $I$ is equidimensional and $\Proj (A/I)$ is smooth.
 Then there exists a positive integer $f$ such that {\Kk} is
Cohen-Macaulay for all $e>0$ and $c \geq ef$.}

\medskip
Let $\pi: X \rightarrow Y$ be the blow-up morphism,
$E= {\cal P}roj (\bigoplus_{n  \geq 0} {\cal I}^n/{\cal I}^{n+1})$,
and $w_E$ its dualizing sheaf.
The main result they obtain is the following general criterion:

\medskip
\noindent
{\bf Theorem} \cite[Theorem 4.1]{CH}
{\it Suppose that $I \subset A$ is a homogeneous ideal such that
$I \not \subset \fp$, $\forall \fp \in \Ass(A)$, $A$ is
Cohen-Macaulay and $X$ is a Cohen-Macaulay scheme.
Suppose that $\pi_* {\cal O}_E (m) =  {\cal I} ^m / {\cal I} ^{m+1} $
for $m \geq 0$,
$R^i \pi _*  {\cal O} _E (m)$=$0$  for $i>0$ and $m \geq 0$,
$R^i \pi _*  {w} _E (m)$=$0$  for $i>0$ and $m \geq 2$. Then
there exists a positive integer $f$ such that {\Kk}
is Cohen-Macaulay for $e>0$ and $c \geq ef$.}

\medskip
This result is applied there to the following classes of ideals:

\medskip
\noindent
{\bf Corollary} \cite[Corollary 4.2]{CH}
{\it Let $I \subset A$ be a homogeneous ideal such that $I \not
\subset \fp$, $\forall \fp \in \Ass(A)$, $A$ is Cohen-Macaulay and
$I_\fp$ is a complete intersection ideal for any $\fp \in \Proj(A)$.
Then there exists a positive integer $f$ such that
{\Kk} is Cohen-Macaulay for $e>0$, $c \geq ef$.}

\medskip
\noindent
{\bf Corollary} \cite[Corollary 4.4]{CH}
{\it Let $I \subset A$ be a homogeneous ideal such that $I \not
\subset \fp$, $\forall \fp \in \Ass(A)$, $A$ is Cohen-Macaulay and
$I_{(\fp)}$ is strongly Cohen-Macaulay with
$\mu(I_{(\fp)}) \leq \h (\fp)$ for any prime ideal
$\fp \in \Proj(A)$ containing $I$. Then there exists a positive integer
$f$ such that {\Kk} is Cohen-Macaulay for $e>0$, $c \geq ef$.}

\bigskip
\bigskip
As a somehow unexpected by-product, the methods used to study the
diagonals of a Rees algebra also allow to study the regularity of the
powers of an ideal and their asymptotic properties. These problems have
been previously handled by using other techniques. Let $A= k[{\x}]$ be a
polynomial ring with the usual grading and let $I \subset A$ be a
homogeneous ideal.
I. Swanson  \cite{S} has shown that there exists an integer $B$
such that $\reg (I^e) \leq Be$, $\forall e$. The problem is to make $B$
explicit. In some particular cases, such $B$ was already known.
A. Geramita, A. Gimigliano and Y. Pitteloud \cite{GGP} and
K. Chandler \cite{C} had proved that for ideals with
dim$(A/I)= 1$, $\reg(I^e) \leq \reg(I) \, e$.
On the other hand,
R. Sj\"{o}gren \cite{Sj} had given another kind of bound: If $I$
is an ideal generated by forms of degree $\leq d$ with
dim$(A/I) \leq 1$,
$\reg (I^e) < (n-1)de$. Also A. Bertram, L. Ein and R. Lazarsfeld
\cite{BEL} have given a bound for the regularity of the powers of an
ideal in terms of the degrees of its generators: If $I$ is
the ideal of a smooth complex subvariety $X$ of ${\Bbb P}_{\Bbb
C}^{n-1}$ of codimension $c$  generated by forms of degrees
$d_1 \geq d_2 \geq \dots \geq d_r$ , then
$$ H^i ({\Bbb P}_{\Bbb C}^{n-1}, {\cal I}^e (k)) = 0, \,\,
\,\forall i \geq 1,\, \forall k \geq ed_1 + d_2 + \dots +d_c -(n-1).$$

\medskip
Let $(A, \fm, k)$ be a local ring and let $I \subset A$ be an ideal.
Concerning the asymptotic properties of the powers of $I$, a classical
well known result of M. Brodmann \cite{Bro} says that
$\depth \,A/I^j$ takes a constant asymptotic value $C$ for $j>>0$, and
moreover $C \leq \dim A-l(I)$. This value $C$ was determined by
D. Eisenbud and C. Huneke \cite{EH} for ideals under some restrictions:
If $I$ is an ideal of height greater than zero and $G_A(I)$
is Cohen-Macaulay, then $\inf \{\depth \, A/I^j\} = \dim A - l(I)$, and
if $\depth \, A/I^s= \inf \{\depth \, A/I^j\}$, then
$\depth \, A/I^{s+1}= \depth \, A/I^{s}$.
Finally, V. Kodiyalam  \cite{Ko1} has shown that for any fixed
nonnegative integer $p$ and all sufficiently large $j$, the $p$-th
Betti number $\beta_p ^A(I^j) =\dim_k \Tor^A_p(I^j, k)$ and the
$p$-th Bass number $\mu^p_A(I^j) =\dim_k \Ext_A^p(k, I^j)$ are
polynomials in $j$ of degree $\leq l(I)-1$.

\bigskip
\bigskip

Now we are going to set and motivate the concrete problems and
questions considered in this dissertation.

\medskip
The restriction to Rees algebras of equigenerated ideals done by A.
Simis et al. \cite{STV} is due to the fact that in this case
the Rees algebra can be endowed with a bigrading so that it becomes
standard. For standard bigraded algebras one may define its
biprojective scheme (see \cite{STV}, \cite{H}) and there are also known
results about its Hilbert polynomial (see \cite{HHRT}, \cite{KMV}).
If $I$ is an ideal generated by forms $f_1, \dots, f_r$ of degrees
$d_1, \dots, d_r$ respectively, the Rees algebra of $I$ has a
presentation as a quotient of $S=k[{\x},{\y}]$ bigraded by setting
$\deg(X_i)=(1,0)$, $\deg(Y_j)=(d_j,1)$ which is non standard. Our first
problem will be to extend the definitions and known results on bigraded modules
over standard bigraded algebras to the category of bigraded $S$-modules.

\medskip
Several arithmetic properties of a ring such as the
Cohen-Macaulayness
and the Gorenstein property can be characterized by means
of its local cohomology modules. This is the reason why it is
interesting and useful to study when the local cohomology modules and
the diagonal functor commute,
case in which we may conclude that certain arithmetic properties of the
Rees algebra are inherited by its diagonals.
The shifts $(a,b)$ which arise in the bigraded minimal free
resolution of the Rees algebra $R_A(I)$ over the polynomial ring $S$
play an essential role in this problem as it was seen in \cite{CHTV}.
We will study and bound these shifts by relating them to
the local cohomology of the Rees algebra.
After that, we will focus on the obstructions for the local
cohomology modules and the
diagonal functor to commute.

\medskip
Once we have done all those preliminaries, our main purpose will be
to study the Cohen-Macaulayness of the rings $\Kk$. We will consider
different questions such as the existence and the determination of the
diagonals $(c,e)$ for which $\Kk$ is Cohen-Macaulay, problems treated in
\cite{STV}, \cite{CHTV} and \cite{CH}.
Similarly, our next goal will be to study the Gorenstein property of
the $k$-algebras $\Kk$. This has been only done in a very particular
case in \cite{CHTV}.

\medskip
Some of the criteria we will obtain for the Cohen-Macaulayness of the
$k$-algebras $\Kk$ are in terms of the local cohomology modules of the
powers of the ideal $I$. This will lead us to study the $a$-invariants
of the powers of a homogeneous ideal.
We will then show how the bigrading defined in the Rees algebra can
be used to study the $a$-invariants and the asymptotic properties of the
powers of an ideal.

\medskip
Summarizing, the main problems we have considered in this work are:

\begin{itemize}
\item[(\bf 1)]
To extend the definitions and results about the biprojective scheme
and the Hilbert polynomial of finitely generated bigraded modules
defined over standard bigraded $k$-algebras to finitely generated
bigraded $S$-modules, for $S=k[{\x},{\y}]$ the polynomial ring
bigraded by
$\deg(X_i)=(1,0)$,
$\deg(Y_j)=(d_j,1)$, $d_1, \dots, d_r \geq 0$.

\item[(\bf 2)]
To relate the shifts in the bigraded minimal free resolution of any
finitely generated bigraded $S$-module to its $a$-invariants.

\item[(\bf 3)]
To study the local cohomology modules of the diagonals of
any finitely generated bigraded $S$-module.

\item[(\bf 4)]
To study the Cohen-Macaulay property of the rings $\Kk$.

\item[(\bf 5)]
To study the Gorenstein property of the rings $\Kk$.

\item[(\bf 6)]
To study the $a$-invariants of the powers of a homogeneous ideal.

\item[(\bf 7)]
To study the asymptotic properties of the powers of a homogeneous ideal.

\end{itemize}

\bigskip
\bigskip
Now we are ready to describe the results obtained in this work.

\medskip
    In {\bf Chapter 1} we introduce the notations and definitions we
will need throughout this work. We begin the chapter by defining
the category of multigraded modules over a multigraded ring, and by
recalling some well-known results about multigraded local
cohomology and the canonical module mainly following M. Herrmann, E.
Hyry and J. Ribbe \cite{HHR} and S. Goto and K. Watanabe \cite{GW1}.
Then we define the multigraded $a$-invariants
of a module and we study the relationship between these $a$-invariants
and the shifts of its multigraded minimal free resolution. We will
obtain a formula which extends \cite[Example 3.6.15]{BH}, where it was
proved for Cohen-Macaulay modules in the graded case.
This result will be a very useful device used all along this work.
To precise it, let $S$ be a $d$-dimensional $\Bbb N^r$-graded
Cohen-Macaulay $k$-algebra with homogeneous maximal ideal $\cal M$ and
let $M$ be a finitely generated $r$-graded $S$-module of dimension
$m$ and depth $\rho$. For each $i=0, \dots, m$,
we may associate to the $i$-th local cohomology module
of $M$ its multigraded $a_i$-invariant
${\bf a}_i (M) = (a_i^1(M), \dots, a_i^r(M))$, where
$$a_i^j(M)= \max \, \{ n \mid \exists \bn= (n^1, \dots, n^r) \in \Bbb
Z^r  \; s. t. \;  \uH ^i_{\cal M}(M)_{\bn} \not =0, n^j=n \}$$
if $\uH^i_{\cal M}(M) \not =0 $ and $a_i^j (M)= -\infty$ otherwise.
Notice that ${\bf a}_m(M)$ coincides with the usual $a$-invariant, and
so we will denote by
${\bf a}(M) = (a^1(M), \dots, a^r(M))={\bf a}_m(M)$.
Finally, the multigraded $a_*$-invariant of
$M$ is ${\bf a}_*(M) = (a_*^1(M), \dots, a_*^r(M))$, where
$a_*^j(M) = \max_{i=0, \dots, m} \{a_i^j(M) \}$.

\medskip
On the other hand, we may consider the $r$-graded minimal free
resolution of $M$ over $S$. Suppose that this resolution is finite:
$${0\to {D_l}\to \dots \to {D_1} \to {D_0} \to {M} \to 0},$$
with ${D_p}=\bigoplus_q S(a^1_{pq},...,a^r_{pq})$.
For every $p \in \{ 0, \dots, l \}$, $j \in \{ 1, \dots, r \}$,
let us denote by

\vspace{2mm}

\hspace{15mm} $t_p^j(M) = \max_q \{ -a^j_{pq}\}, $

\vspace{2mm}

\hspace{15mm} $t_*^j(M) = \max_{p, q} \{ -a^j_{pq}\} =
\max_p   t_p^j(M),$

\vspace{2mm}

\hspace{15mm} ${\bf t}_* (M) = (t_*^1(M), \dots, t_*^r(M)).$

\vspace{2mm}

Moreover, given a permutation $\sigma$ of the set $\{ 1,...,r \}$,
let us consider $\leq_{\sigma}$ the order in $\Bbb Z^r$ defined by:
$ (u_1,...,u_r) \leq_{\sigma} (v_1,...,v_r)$
iff
$(u_{\sigma(1)},...,u_{\sigma(r)}) \leq_{lex}
 (v_{\sigma(1)},...,v_{\sigma(r)})$,
where $\leq_{lex}$ is the lexicographic order.
Set  ${\bf M}_{p}^{\sigma} = \max_{\leq_{\sigma}} \{\
(-a^1_{pq},...,-a^r_{pq}) \}$.
 Then we can relate the shifts and the $a$-invariants of $M$ in the
following way:

\medskip
\noindent
{\bf Theorem 1} [Theorem 1.3.4]
{\it
For every  $j = 1, \dots, r$,
\begin{itemize}
\item[(i)]
$a_{d-p}^j(M) \leq t_p^j(M) + a^j(S)$,
for $p= d-m, \dots, d- \rho$.
\item[(ii)] Assume that for some $p$ there exists $\sigma$ s.t.
$\sigma(1)=j$ and $M_p^{\sigma}>_{\sigma} M_{p+1}^{\sigma}$.
 Then
 $a_{d-p}^j(M) = t^j_p(M) +a^j(S)$.
\item[(iii)] $a_*^j(M) = t_*^j(M)+ a^j(S)$. That is,
    ${\bf a}_*(M) =  {\bf t}_*(M) + {\bf a}(S)$.
\end{itemize}}

\medskip
 After that, we extend the definition and some of the results about the
multiprojective scheme associated to a standard $r$-graded ring given
by E. Hyry \cite{H} and M. Herrmann et al. \cite{HHRT} to rings endowed
with a more general grading, which will also include the Rees algebra of
a homogeneous ideal.
Let $S$ be a noetherian $\Bbb N^r$-graded ring generated over
$S_0$ by homogeneous elements
 $x_{11}, \dots, x_{1k_1}, \dots, x_{r1}, \dots, x_{rk_r}$ in degrees
$\dg (x_{ij}) = (d_{ij}^1, \dots, d_{ij}^{i-1},1,0, \dots, 0)$, with
 $d_{ij}^l \geq 0$.
For every $j = 1, \dots, r$, let $I_j$ be the ideal of $S$
generated by the homogeneous components of $S$ of degree ${\bf n}=(n_1,
\dots, n_r)$
such that $n_j >0, n_{j+1}=\dots =n_r=0$. The irrelevant ideal of $S$ is
$S_+ = I_1 \cdots I_r$. We may associate to $S$ the $r$-projective
scheme $\Proj^r(S)$ which as a set contains all the homogeneous prime
ideals $P \subset S$ such that $S_+ \not \subset P$. The relevant
dimension of $S$ is
 $$\rdim S = \cases {r-1 & if $\Proj^r(S) = \emptyset$ \cr
\max \, \{\, \dim S/P \mid P \in \Proj^r (S) \} & if $\Proj^r(S) \not =
\emptyset$ \cr} .$$

\medskip
It can be proved that
$\dim \Proj^r(S) = \rdim S -r$ by arguing as in \cite[Lemma
1.2]{H} where the standard $r$-graded case was considered. This
result jointly with the isomorphism of schemes
$\Proj^r(S) \cong \Proj (S_\Delta)$ that we have for certain diagonals
allows to compute the dimension of $S_\Delta$ whenever $S_0$ is
artinian,
by extending \cite[Proposition 2.3]{STV} where this dimension was
determined for the $(1,1)$-diagonal of a standard bigraded $k$-algebra
by different methods.

\medskip
Finally, we extend to the category of $r$-graded modules defined
over the $r$-graded $k$-algebras introduced before the basic
results concerning Hilbert functions and Hilbert polynomials.
Some of them have been established in the standard $r$-graded case
in \cite{HHRT} and \cite{KMV}.

\medskip
   In {\bf Chapter 2} we are concerned with the diagonal functor
in the category of bigraded $S$-modules, where $S$ is the polynomial
ring $S=k[{\x},{\y}]$ bigraded by setting
$\deg X_i=(1,0)$, $\deg Y_j =(d_j,1)$, $d_1, \dots, d_r \geq 0$.
In the first section, we compare the local
cohomology modules of a finitely generated bigraded $S$-module $L$ with
the local cohomology modules of its diagonals.
In particular, we will prove the main results in \cite{CHTV} by a
different and somewhat easier approach. In addition, this approach will
provide more detailed information about several problems related to the
behaviour of the local cohomology when taking diagonals.
Set $d= \max \{d_1, \dots, d_r \}$, and let
$\Delta = (c,e)$ be a diagonal with $c \geq de+1$. Let us consider the
following subalgebras of $S$:
$S_1= k[{\x}]$, $S_2= k[{\y}]$, with homogeneous maximal ideals
${\fm}_1= (X_1, ..., X_n)$ and ${\fm}_2 = ({\y}) $. Let ${\M}_1$,
${\M}_2$ be the ideals of $S$ generated by ${\fm}_1$, ${\fm}_2$
respectively, and let $\M$ be the homogeneous maximal ideal of $S$.
Then:

\medskip
\noindent
{\bf Proposition 2} [Proposition 2.1.3]
{\it Let $L$ be a finitely generated bigraded $S$-module. There
exists a natural exact sequence
$$ ... \to H^q_{\M} (L)_{\Delta}
\to H^q_{{\M}_1} ( L )_{\Delta}
 \oplus H^q_{{\M}_2} ( L )_{\Delta}
\to H^q_{\MD}(L_{\Delta})
\stackrel {\varphi^q_L} \longrightarrow
 H^{q+1}_{\cal M}(L)_{\Delta}  \to ...$$
}

In the rest of the section, we study the obstructions for ${\varphi^q_L}$
to be an isomorphism. Firstly, we relate this question to the vanishing
of the local cohomology with respect to ${\cal M}_1$ and ${\cal M}_2$ of
the modules $S(a,b)$ which arise in the bigraded minimal free
resolution of $L$ over $S$. This allows us as said to recover the main
results in \cite{CHTV}. After that, we study the vanishing of the local
cohomology modules of $L$ with respect to ${\cal M}_1$ and ${\cal M}_2$
by themselves.

\medskip
In Section 2.2 we will focus on standard bigraded
$k$-algebras.
Given a standard bigraded $k$-algebra $R$, let us consider the graded
subalgebras ${\cal R}_1 = \bigoplus _{i \in \Bbb N} R_{(i,0)}$,
${\cal R}_2 = \bigoplus _{j \in \Bbb N} R_{(0,j)}$.
By using Theorem 1, we obtain a characterization for $R$
to have a good resolution in terms of the $a_*$-invariants of
${\cal R}_1$ and ${\cal R}_2$ which, in particular, provides a
criterion for the Cohen-Macaulay property of its diagonals.
We also find necessary and sufficient conditions on the local cohomology
of ${\cal R}_1$ and ${\cal R}_2$ for the existence of Cohen-Macaulay
diagonals of a Cohen-Macaulay standard bigraded $k$-algebra $R$. This
result extends
\cite[Corollary 3.12]{CHTV}.

\medskip
\noindent
{\bf Proposition 3} [Proposition 2.2.7]
{\it Let $R$ be a standard bigraded Cohen-Macaulay $k$-algebra of
relevant dimension $\delta$. There exists $\Delta$ such that
$R_\Delta$ is Cohen-Macaulay if and only if
$H^q_{{\frak m}_1}({\cal R}_1)_0=
H^q_{{\frak m}_2}({\cal R}_2)_0= 0$ for any $q< \delta -1$.}

\medskip
Now let us consider a standard bigraded ring $R$ defined over a
local ring with $a^1(R), a^2(R) <0$. In  \cite[Theorem 2.5]{H} it is
shown that if $R$ is Cohen-Macaulay then the $\Delta=(1,1)$-diagonal of
$R$ has also this property. This result can be extended to any diagonal
of a standard bigraded $k$-algebra:

\medskip
\noindent
{\bf Proposition 4} [Proposition 2.2.6]
{\it Let $R$ be a standard bigraded Cohen-Macaulay $k$-algebra with
$a^1(R), a^2(R) <0$.
Then $R_\Delta$ is Cohen-Macaulay for any diagonal $\Delta$.}

\medskip
At the end of the chapter, we apply the results about bigraded
$k$-algebras to the Rees algebra of a homogeneous ideal.
Let $A$ be a noetherian graded $k$-algebra generated in degree $1$ of
dimension $\nn$ and let $\fm$ be the homogeneous maximal ideal of $A$.
Given a homogeneous ideal $I$ of $A$, the Rees algebra $R=R_A(I)$ of $I$
is bigraded by $R_A(I)_{(i,j)}= (I^j)_i$. If $I$ is generated in degree
$\leq d$, for any diagonal $\Delta =(c,e)$ with $c \geq de+1$ we have:
$$R_A(I)_\Delta = k[(I^e)_c].$$
The diagonals $\Kk$ are graded $k$-algebras of dimension $\nn$
if no associated prime of $A$  contains $I$. In the sequel we will
always assume such hypothesis. We can relate the local cohomology
modules of the $k$-algebras $\Kk$ and those of the powers of $I$.
Denoting by $m$ the homogeneous maximal ideal of $\Kk$, we have:

\medskip
\noindent
{\bf Proposition 5} [Corollary 2.3.5]
{\it For any $c \geq de+1$, $e>a_*^2(R)$, $s>0$, we have
isomorphisms
 $$H^q_{m}(\Kk)_s \cong H^q_{\fm}(I^{es})_{cs}, \forall q \geq 0.$$}

\medskip
In the particular case where $A= k[{\x}]$, A. Conca et al. \cite{CHTV}
conjectured that if the Rees algebra
of a homogeneous ideal $I$ of $A$ is Cohen-Macaulay, then there
exists a Cohen-Macaulay diagonal. The results proved for standard
bigraded $k$-algebras
provide an affirmative answer for equigenerated homogeneous ideals.
In fact, we can give a full answer to this conjecture.

\medskip
\noindent
{\bf Theorem 6} [Theorem 2.3.12]
{\it Let $I$ be a homogeneous ideal of the polynomial ring $A=k[{\x}]$.
If $R_A(I)$ is a Cohen-Macaulay ring, then $R_A(I)$ has a good
resolution. In particular, $\Kk$ is Cohen-Macaulay for $c \gg e \gg 0$.}

\medskip
Furthermore, we obtain sufficient and necessary conditions on the
ring $A$ for the existence of Cohen-Macaulay diagonals of a
Rees algebra $R_A(I)$ with this property. Namely,

\medskip
\noindent
{\bf Theorem 7} [Theorem 2.3.13]
{\it If $R_A(I)$ is Cohen-Macaulay, then the following are
equivalent:
\begin{enumerate}
 \item There exist $c, e$ such that $\Kk$ is Cohen-Macaulay.
 \item $H^i_{\frak m}(A)_0 = 0$ for $i < \overline n$.
\end{enumerate}}

\medskip
  In {\bf Chapter 3} we study in detail the Cohen-Macaulay
property of the rings $\Kk$. We consider the
problem of the existence of Cohen-Macaulay diagonals of the Rees
algebra. Once studied this problem, we will try to determine the diagonals
with this property. The following isomorphisms will play an important
role:

\medskip
\noindent
{\bf Proposition 8} [Proposition 3.1.2]
{\it Let $X$ be the blow-up of $\Proj (A)$ along ${\cal I}=
\widetilde{I}$, where $I$ is a homogeneous ideal of $A$ generated by
forms of degree $\leq d$. For any $c \geq de+1$, there are isomorphisms
of schemes
$$X \cong \Proj^2(R_A(I)) \cong \Proj (\Kk).$$}

First of all, these isomorphisms will be used to give a criterion
for the existence of diagonals $\Kk$ which are generalized
Cohen-Macaulay modules, thereby solving a conjecture of \cite{CHTV}.

\medskip
\noindent
{\bf Proposition 9} [Proposition 3.2.6]
{\it The following are equivalent:
\begin{enumerate}
\item $H^i_{\cal M} (R_A(I))_{(p,q)}$=$0$ for $ i < \nn+1$, $p
\ll q \ll 0$.
\item $\Kk$ is a generalized Cohen-Macaulay module for $c \gg e \gg 0$.
\item There exist $c,e$ such that $\Kk$ is generalized
Cohen-Macaulay.
\item $\Kk$ is a Buchsbaum ring for $c \gg e \gg 0$.
\item There exist $c,e$ such that $\Kk$ is a Buchsbaum ring.
\item There exist $q_0$, $t$ such that
 $H^i_{\cal M} (R_A(I))_{(p,q)} = 0$ for $i < \overline n+1$, $q <
q_0$ and $p < dq +t$.

\end{enumerate}}

\medskip
After that, we use Proposition 8 to give necessary and sufficient
conditions for a Rees algebra to have Cohen-Macaulay diagonals.
Namely,

\medskip
\noindent
{\bf Theorem 10} [Theorem 3.2.3, Corollary 3.2.5]
{\it The following are equivalent:
\begin{enumerate}
\item  There exist $c,e$ such that $\Kk$ is a Cohen-Macaulay ring.
\item
 \begin{itemize}
 \item[(1)]  There exist $q_0, t \in \Bbb Z$ such that
 $H^i_{\cal M} (R_A(I))_{(p,q)}= 0$ for all $i < \overline n+1$,
 $q < q_0$ and $p < dq + t$.
 \item[(2)]
 $ H^{i}_{R_A(I)_+} (R_A(I)) _{(0,0)} = 0$ for all $i < \overline n$.
 \end{itemize}
\item
\begin{itemize}
\item[(1)]  $X$ is an equidimensional Cohen-Macaulay scheme.
\item[(2)]
${\mit \Gamma} ({\X}, \cal O _{\X})$ = $k$,
  $H^i({\X}, \cal O _{\X})$ = $0$ for $0 < i < \overline n -1$.
\end{itemize}
\end{enumerate}

In this case, ${\Kk}$ is a Cohen-Macaulay ring for $c \gg e \gg 0$.}

\medskip
By using this theorem, we can exhibit some general situations
in which we can ensure the existence of Cohen-Macaulay coordinate rings
for $X$. For instance,

\medskip
\noindent
{\bf Proposition 11} [Proposition 3.3.3]
{\it  Let $X$ be the blow-up of $\Bbb P^{n-1}_k$ along a closed
subscheme, where $k$ has $char k = 0$. Assume that $X$ is smooth or
with rational singularities. Then $X$ is arithmetically Cohen-Macaulay.}

\medskip
Our next goal in the chapter will be to determine the Cohen-Macaulay
diagonals once we know its existence. This is a difficult problem which
has been completely solved only for complete intersection ideals in the
polynomial ring \cite[Theorem 4.6]{CHTV}. For equigenerated ideals,
we can give a criterion for the Cohen-Macaulayness of a diagonal in
terms of the local cohomology modules of the powers of the ideal
by just assuming that the Rees algebra is Cohen-Macaulay. Namely,

\medskip
\noindent
{\bf Proposition 12} [Proposition 3.4.1]
{\it Let $I \subset A$ be an ideal generated by forms of degree $d$
whose Rees algebra is Cohen-Macaulay. For any $c \geq de+1$, $\Kk$ is
Cohen-Macaulay if and only if
\begin{enumerate}
\item $H^i_{\frak m}(A)_0 = 0$ for $i < \nn$.
\item $H^i_{\frak m}(I^{es})_{cs} = 0$ for $i < \nn$, $s>0$.
\end{enumerate}}

\medskip
For arbitrary homogeneous ideals, we can also prove a
criterion for the Cohen-Macaulayness of a diagonal by means of the
local cohomology of the powers of the ideal and the local cohomology of
the graded pieces of the canonical module of the Rees algebra.
Let us denote by $K= K_{R_A(I)}= \bigoplus_{(i,j)} K_{(i,j)}$ the
canonical module of the Rees algebra, and for each $e \in \Bbb Z$,
let us consider the graded $A$-module
$K^e = \bigoplus_i K_{(i,e)}$. Then we have:

\medskip
\noindent
{\bf Theorem 13} [Theorem 3.4.3]
{\it Let $I$ be a homogeneous ideal of $A$ generated by forms of
degree $\leq d$ whose Rees algebra is Cohen-Macaulay. For any
$c \geq de+1$,
$\Kk$ is Cohen-Macaulay if and only if
\begin{enumerate}
\item $H^i_{\frak m}(A)_0 = 0$ for $i < \nn$.
\item $H^i_{\frak m}(I^{es})_{cs} = 0$ for $i < \nn$, $s>0$.
\item $H^{\nn -i+1}_{\frak m}(K^{es})_{cs} = 0$ for $1 \leq i < \nn$,
$s>0$. \end{enumerate}}

\medskip
If the form ring is quasi-Gorenstein we can express the criterion
above only in terms of the local cohomology of the powers of the ideal.

\medskip
\noindent
{\bf Theorem 14} [Corollary 3.4.4]
{\it Let $I$ be a homogeneous ideal of $A$ generated by forms of
degree $\leq d$. Assume that $R_A(I)$ is Cohen-Macaulay,
$G_A(I)$ is quasi-Gorenstein. Set $a= -a^2(G_A(I))$, $b = -a(A)$. For any
$c \geq de+1$,
$\Kk$ is Cohen-Macaulay if and only if
\begin{enumerate}
\item $H^i_{\frak m}(A)_0 = 0$ for $i < \nn$.
\item $H^i_{\frak m}(I^{es})_{cs} = 0$ for $i < \nn$, $s>0$.
\item $H^{i}_{\frak m}(I^{es-a+1})_{cs-b} = 0$ for $1<i \leq \nn$, $s>0$.
\end{enumerate}
}

\medskip
We can use Theorem 14 to determine exactly the Cohen-Macaulay
diagonals of the Rees algebra of a complete intersection ideal in
any Cohen-Macaulay ring. In particular, we get a new proof of
\cite[Theorem 4.6]{CHTV} where the case $A= k[{\x}]$ was studied.

\medskip
These criteria will be also applied in the Chapter 5, once we have
studied in detail the local cohomology modules of the powers of several
families of ideals, such as equimultiple ideals or strongly
Cohen-Macaulay ideals.

\medskip
Furthermore, the results and methods used up to now allow us to show
the behaviour of the $a_*$-invariant of the powers of a homogeneous
ideal. The following statement has been obtained independently by
S.D. Cutkosky, J. Herzog and N. V. Trung \cite{CHT} and
V. Kodiyalam \cite{Ko2} by different methods.

\medskip
\noindent
{\bf Theorem 15} [Theorem 3.4.6]
{\it Let $L$ be a finitely generated bigraded $S$-module. Then there
exists $\alpha$ such that
$$a_*(L^e) \leq de+ \alpha, \; \forall e .$$}

After that, we use the bound on the shifts of the bigraded minimal free
resolution of the Rees algebra obtained in Theorem 1 to determine a
family of Cohen-Macaulay diagonals of a Cohen-Macaulay Rees algebra.

\medskip
\noindent
{\bf Theorem 16} [Theorem 3.4.12]
{\it Let $I$ be a homogeneous ideal of $A$ generated by $r$ forms of
degree $d_1 \leq \dots \leq d_r=d$. Assume that $H^i_\fm(A)_0
=0$ for $i < \nn$. Set $u = \sum_{j=1}^r d_j$. If the Rees algebra
is Cohen-Macaulay, then
\begin{enumerate}
\item $k[(I^e)_c]$ is Cohen-Macaulay for
  $c> \max \{d(e-1)+u+a(A), d(e-1)+u-d_1(r-1) \}$.
\item If $I$ is generated by forms in degree $d$,
$\Kk$ is Cohen-Macaulay for $c> d(e -1 +l(I))+ a(A) $.
\end{enumerate}}

\medskip
   Our results can be also applied to study the embedddings of the
blow-up of a projective space along an ideal $I$ of fat points via
the linear systems $(I^e)_c$ whenever these linear systems are very
ample, slightly extending \cite[Theorem 2.4]{GGP} where only the linear
systems $I_c$ were considered.

\medskip
\noindent
{\bf Theorem 17} [Theorem 3.4.15]
{\it Let $I \subset A= k[{\x}]$ be an ideal of fat points, with $k$ a
field of characteristic zero. Then
\begin{enumerate}
\item $\Kk$ is Cohen-Macaulay if and only if
 $H^{i}_{\fm}(I^{es})_{cs} = 0$ for $s>0$, $i<n$.
\item For $c > \reg(I) e$, $\Kk$ is Cohen-Macaulay with
 $a (\Kk) < 0$.  In particular, $\reg (\Kk) < n-1$.
\end{enumerate}}

\medskip
The chapter finishes by studying sufficient conditions for the
existence of a positive integer $f$ such that $\Kk$ is a
Cohen-Macaulay ring for all $c \geq ef$ and $e >0$, a question that
has been treated by S.D. Cutkosky and J. Herzog. Our main result,
which improves \cite[Corollaries 4.2, 4.3 and 4.4]{CH}, is the
following:

\medskip
\noindent
{\bf Theorem 18} [Theorem 3.5.3]
{\it  Let $I$ be a homogeneous ideal of an equidimensional
$k$-algebra $A$ such that $R_{A_{\frak p}} ( I_ {\frak p})$ is
Cohen-Macaulay for any prime ideal ${\frak p} \in \Proj(A)$. Assume
that $H^i_{\frak m}(A)_0 = 0 $ for $i < \nn$.
Then there exists an integer $\alpha$ such that $\Kk$ is
Cohen-Macaulay for all $c \geq de+ \alpha$ and $e>0$.}

\medskip
The aim of {\bf Chapter 4} is to study the Gorenstein property
of the $k$-algebras $\Kk$. About the Cohen-Macaulay property, we have
already proved that if there exists a Cohen-Macaulay diagonal
then there are infinitely many with this property.
We show that the behaviour of the Gorenstein property is
totally different. For instance, by considering the polynomial
ring $S=k[{\x},{\y}]$ with
$\deg X_i=(1,0)$, $\deg Y_j =(d_j,1)$, $d_1, \dots, d_r \geq 0$,
we have that $S_\Delta$ is Cohen-Macaulay for any diagonal $\Delta$ but
there is just a finite set of Gorenstein diagonals.

\medskip
\noindent
{\bf Proposition 19} [Proposition 4.1.1]
 {\it ${S_\Delta}$ is Gorenstein if and only if
  $\frac{r}{e}  =\frac{n+u}{c} = l \in {\Bbb Z}$.
 Then $a(S_\Delta)= -l$.}

\medskip
To determine the rings $\Kk$ which are Gorenstein,
we will compare the canonical module of the Rees algebra with the
canonical module of each diagonal.
For complete intersection ideals of the polynomial ring, it was proved
in \cite[Proposition 4.5]{CHTV} that the canonical module and the
diagonal functor commute. This result can be extended to more general
situations.

\medskip
\noindent
{\bf Proposition 20} [Proposition 4.1.4 and Remark 4.1.5]
{\it Let $A=k[{\x}]$ be the polynomial ring, $n \geq 2$, and let $I$
be a homogeneous ideal of $A$ with $\mu(I) \geq 2$.
\begin{enumerate}
\item If $\mu(I) \leq n$, $K_{R_\Delta} \cong (K_R)_\Delta$.
\item If $I$ is equigenerated and $R$ is Cohen-Macaulay,
$K_{R_\Delta} \cong (K_R)_\Delta$.
\end{enumerate}}

\medskip
Although this isomorphism can be extended to a more general class
of rings, we will restrict our attention to the above two cases. This
will suffice to study the rational surfaces obtained by blowing-up the
projective plane at a set of points.

\medskip
Next we study the behaviour of the Gorenstein property of the Rees
algebra when we take diagonals.
If the Rees algebra is Gorenstein then the form ring is also Gorenstein.
Under this assumption on the form ring, which is less restrictive, we
can determine exactly for which $c,e$ the algebra $\Kk$ is
quasi-Gorenstein. Namely,

\medskip
\noindent
{\bf Theorem 21} [Theorem 4.1.9]
{\it Let $I \subset A= k[{\x}]$ be a homogeneous ideal with $1< \h (I) <n$
whose form ring $G_A(I)$ is Gorenstein. Set $a = -a^2(G_A(I))$. Then
$\Kk$ is a quasi-Gorenstein ring if and only if
$\frac{n}{c}=\frac{a-1}{e} =l_0 \in \Bbb Z$.
In this case, $a(\Kk) =-l_0$.}

\medskip
For homogeneous non principal ideals $I$ of height $1$, the ring $\Kk$
is never
Gorenstein. If $I$ has height $n $, then the diagonals determined in the
theorem are always Gorenstein, but the converse is not true. As a
corollary of this result we can solve the problem of determining
completely the Gorenstein diagonals for complete intersection ideals or
determinantal ideals generated by the maximal minors of a generic
matrix.

\medskip
\noindent
{\bf Corollary 22} [Corollary 4.1.12]
{\it Let $I \subset A= k[{\x}]$ be a complete intersection
ideal minimally generated by $r$ forms of degree
$d_1 \leq \dots \leq d_r=d$, with $r<n$. For any $c \geq de+1$,
 $\Kk$ is a Gorenstein ring if and only if
  $\frac{n}{c}= \frac{r-1}{e} =l_0 \in \Bbb Z$.
 In this case, $a(\Kk) = -l_0$.}

\medskip
\noindent
{\bf Corollary 23} [Example 4.1.13]
{\it  Let  ${\bf X} = (X_{ij})$ denote a matrix of indeterminates, with
$1 \leq i \leq n, 1 \leq j \leq m$ and $m \leq n$.
Let $I \subset A= k[ \bf X]$ denote the ideal generated by the
maximal minors of $\bf X$, where $k$ is a field.
Then:
\begin{enumerate}
\item
 If $ m < n$, then $\Kk$ is Gorenstein if and only if
 $\frac{nm}{c}= \frac{n-m}{e} \in \Bbb Z$.
\item
If $m = n$, then $\Delta= (n(n+1), 1)$ is the only Gorenstein diagonal.
\end{enumerate}}

\medskip
We have shown that if the form ring is Gorenstein there is just
a finite set of Gorenstein diagonals. This fact also holds under the
general assumptions of the chapter. Namely,

\medskip
\noindent
{\bf Proposition 24} [Proposition 4.2.1]
{\it There is a finite set of diagonals $\Delta=(c,e)$ such that
$\Kk$ is quasi-Gorenstein.}

\medskip
If the Rees algebra is  Cohen-Macaulay, then we can bound the diagonals
$\Delta=(c,e)$ for which $\Kk$ is Gorenstein.

\medskip
\noindent
{\bf Proposition 25} [Proposition 4.2.2]
{\it Assume that $\h (I) \geq 2$ and $\Rr$ is Cohen-Macaulay.
Let $a= -a^2(G_A(I))$.
If $\Kk$ is quasi-Gorenstein, then $e \leq a-1$ and $c \leq
n$. Moreover, if $\dim (A/I) >0$ then $\lceil \frac{a}{e} \rceil-1 =
\frac{n}{c} = l \in \Bbb Z$. In particular, if
$a=1$ there are no diagonals $(c,e)$ such that $\Kk$ is
quasi-Gorenstein.
}

\medskip
Finally, we show that in some cases the existence of a diagonal
$(c,e)$ such that $\Kk$ is quasi-Gorenstein forces the form
ring to be Gorenstein. It may be seen as a converse of
Theorem 21 for those cases.

\medskip
\noindent
{\bf Theorem 26} [Theorem 4.2.3]
{\it  Assume that $\Rr$ is Cohen-Macaulay, $\h(I) \geq 2$,
$l(I)<n$ and $I$ is equigenerated. If there exists a diagonal
$(c,e)$ such that $\Kk$ is quasi-Gorenstein then $G_A(I)$ is Gorenstein.}

\medskip
We finish the chapter by applying the previous results to recover
the fact that the Del Pezzo sestic surface in ${\Bbb P}^6$ is
the only Room surface which is Gorenstein.

\medskip
In {\bf Chapter 5} we study the $a$-invariant and the regularity of
any finitely generated bigraded $S$-module $L$, for
$S=k[{\x},{\y}]$ the polynomial ring with $\deg X_i=(1,0)$, $\deg Y_j
=(0,1)$. This class of modules includes for instance any standard
bigraded $k$-algebra $R$.

\medskip
Given a finitely generated bigraded $S$-module $L$, let us
consider the bigraded minimal free resolution of $L$ over $S$
$$0 \to D_l \to \dots \to D_0 \to L \to 0,$$
with $D_p =\bigoplus_{(a,b) \in \Omega_p} S(a,b)$.
The bigraded regularity of $L$ is
${\rm \bf reg}(L) = (\reg_1 L, \reg_2L)$, where
 $$\reg_1 L = \max_p \{-a-p \mid (a,b) \in \Omega_p \},$$
$$\reg_2 L = \max_p \{-b-p \mid (a,b) \in \Omega_p \}.$$

\medskip
For each $e \in \Bbb Z$, we may define the graded $S_1$-module
 $L^e = \bigoplus_{i \in \Bbb Z} L_{(i,e)}$
and the graded $S_2$-module
 $L_e = \bigoplus_{j \in \Bbb Z} L_{(e,j)}$.
Our first result gives a new description of the $a_*$-invariant
${\bf a}_*(L)$ of $L$ and the regularity ${\rm \bf reg} (L)$ of $L$ in
terms of the $a_*$-invariants and the regularities of the graded
modules $L^e$ and $L_e$. Namely,

\medskip
\noindent
{\bf Theorem 27} [Theorem 5.1.1, Theorem 5.1.2]
{\it Let $L$ be a finitely generated bigraded $S$-module. Then:
\begin{enumerate}
\item $a_*^1(L) = \max_e \{a_*(L^e) \} =
                  \max_e \{a_*(L^e) \mid e \leq a_*^2(L) + r \}$.
\item $a_*^2(L) = \max_e \{a_*(L_e) \} =
                  \max_e \{a_*(L_e) \mid e \leq a_*^1(L) + n \}$.
\item $\reg_1 L = \max_e \{\reg (L^e) \} =
                 \max_e \{\reg (L^e) \mid e \leq a_*^2(L) + r \}$.
\item $\reg_2 L = \max_e \{ \reg (L_e) \}=
                  \max_e \{ \reg (L_e) \mid e \leq a_*^1(L) + n \}$.
\end{enumerate}}

\medskip
This result will be used to study the $a_*$-invariant and the
regularity of the powers of a homogeneous ideal $I$ in the polynomial
ring $A= k[{\x}]$. According to Theorem 15, there exists an integer
$\alpha$ such that $a_*(I^e) \leq de+ \alpha$, $\forall e$. The first
aim is to determine such an $\alpha$ explicitly, and this will be done
for any equigenerated ideal by means of a suitable $a$-invariant of
the Rees algebra. For a homogeneous ideal $I$, we will denote by $R$,
$G$ and $F$ the Rees algebra of $I$, its form ring and the fiber cone
respectively. If $I$ is an ideal generated by forms in degree
$d$, let us denote by $R^\varphi$ the Rees algebra endowed
with the bigrading $[R^\varphi]_{(i,j)} =(I^j)_{i+dj}$. Then we have

\medskip
\noindent
{\bf Theorem 28} [Theorem 5.2.1]
{\it Let $I$ be a homogeneous ideal of $A$ generated by forms in degree $d$.
Set $l = l(I)$. Then
\begin{enumerate}
\item $a_*^1 ( R^\varphi ) = \max_{e} \{\,a_*(I^e)-de \} = \max
\, \{\,a_*(I^e)-de  \mid  e \leq a_*^2(R) + l\}.$
\item $\reg_1 ( R^\varphi ) = \max_{e} \{\, \reg(I^e)-de \} = \max
\, \{\, \reg (I^e)-de  \mid  e \leq a_*^2(R) + l\}.$
\end{enumerate}}

\medskip
Therefore, we need to study $a_*^1(R^\varphi)$ to get concrete
bounds for the $a_*$-invariant of the powers of several families of
ideals. If the Rees algebra is Cohen-Macaulay we have

\medskip
\noindent
{\bf Proposition 29} [Proposition 5.2.5]
{\it Let $I$ be a homogeneous ideal generated by forms in degree $d$ whose
Rees algebra is Cohen-Macaulay. Set $l= l(I)$. Then
$$ -n + d(-a^2(G)-1) \leq  \max_{e \geq 0} \{ a_* (I^e) -de \}  \leq -n
+ d (l-1).$$}

The $a_*$-invariants of the powers of a complete intersection ideal are
well-known, and in this case the inequalities above are sharp.
Next we compute explicitly  $a^1_*(R^\varphi)=
\max_{e \geq 0} \{a_*(I^e)-de \}$ for other families of ideals.
First we consider equimultiple ideals.

\medskip
\noindent
{\bf Proposition 30} [Proposition 5.2.8]
{\it Let $I$ be an equimultiple ideal equigenerated in degree $d$ and
set $h= {\rm ht}(I)$. If the Rees algebra is Cohen-Macaulay,
\begin{enumerate}
\item  $a (I^e/I^{e+1}) = de + a(A/I)$. In particular,
$a^1(G^\varphi) = a(A/I)$.
\item   $a_{n-h+1} (I^e) = d (e-1) + a(A/I) $. In particular,
$a^1(R^\varphi) = a(A/I)-d$.
\end{enumerate}
}

\medskip
For ideals whose form ring is Gorenstein we can also
compute explicitly $\max_{e \geq 0} \{a_*(I^e)-de \}$, and then
we get that the lower bound given by Proposition 29 is sharp.

\medskip
\noindent
{\bf Proposition 31} [Proposition 5.2.9]
{\it Let $I$ be a homogeneous ideal equigenerated in degree $d$
whose form ring is Gorenstein. Set $l = l(I)$. Then
\begin{enumerate}
\item $\max_{e \geq 0} \{ a_*(I^e) -de \}=  d(- a^2(G) -1) -n$.
\item For $e> a^2(G) -a(F)$,
$\depth (A/I^e) = n-l$ and
 $a_*(I^e)= a_{n-l} (A/I^e) =d(e-a^2(G)-1)-n$.
\end{enumerate}
}

\medskip
For instance, we may apply this result to determinantal ideals generated
by the maximal minors of a generic matrix as well as to
strongly Cohen-Macaulay ideals satisfying condition $({\cal F}_1)$.

\medskip
The computation of the $a_*$-invariants of the powers of these
families of ideals is then applied to determine the Cohen-Macaulay
diagonals of a Rees algebra. For equimultiple ideals, we
have

\medskip
\noindent
{\bf Proposition 32} [Proposition 5.2.20]
{\it Let $I$ be an equimultiple ideal generated in degree $d$ whose Rees
algebra is Cohen-Macaulay. For any $c \geq de +1$,
$\Kk$ is Cohen-Macaulay if and only if $c> d(e-1) + a(A/I)$.
}

\medskip
For strongly Cohen-Macaulay ideals, we have

\medskip
\noindent
{\bf Proposition 33} [Proposition 5.2.21]
{\it
Let $I$ be a strongly Cohen-Macaulay ideal
such that $\mu(I_{\frak p}) \leq \h(\frak p)$ for any prime ideal
$\frak p \supseteq I$.
Assume that $I$ is minimally generated by forms of degree $d=d_1
\geq  \dots \geq d_r$, and let $h = \h (I)$.
For $c >  d(e-1)+d_1 + \dots + d_h -n $, $\Kk$ is Cohen-Macaulay.
}

\medskip
If the Rees algebra is Cohen-Macaulay, we have proved the existence of
an integer $\alpha$ such that $\Kk$ is a Cohen-Macaulay ring for any
$c >de+ \alpha$ and $e>0$ by Theorem 16.
For equigenerated ideals we had $\alpha= d(l-1)$ as an upper bound.
We can determine the best $\alpha$.

\medskip
\noindent
{\bf Proposition 34} [Proposition 5.2.15, Corollary 5.2.16]
{\it  Let $I$ be an ideal in the polynomial ring $A= k[{\x}]$ generated by
forms in degree $d$ whose Rees algebra is Cohen-Macaulay.
Set $l = l(I)$. For $\alpha \geq 0$, the following are equivalent
\begin{itemize}
 \item[(i)]   $\Kk $ is CM for $c > de + \alpha$.
 \item[(ii)]   $a_i (I^e) \leq d e +\alpha$, $\forall i$, $\forall e$.
 \item[(iii)]   $a_i (I^e) \leq d e +\alpha$, $\forall i$, $\forall e
  \leq  l-1$ .
 \item[(iv)]   $H^{n +1}_{\cal M} (R_A(I))_{(p,q)} = 0$, \,
  $\forall p >dq + \alpha$, that is,
   $ \alpha \geq a^1(R^\varphi)$.
 \item[(v)]   The minimal bigraded free resolution of $R_A(I)$ is
  good for any diagonal $\Delta=(c,e)$ such that $c >de + \alpha$.
\end{itemize}
If the form ring is Gorenstein, these conditions are equivalent to
\begin{itemize}
\item[(vi)] $\alpha \geq d(-a^2(G)-1)-n$.
\end{itemize}}

\medskip
Up to now, we have used Theorem 27 to bound the $a_*$-invariants of
the powers of an ideal, which has been applied to study the
Cohen-Macaulayness of the diagonals. In the last section, we use
this theorem to prove a bigraded version of the Bayer-Stillman theorem
which  characterizes the bigraded regularity of a homogeneous ideal of
$S$ by means of generic homogeneous forms. Next, similarly to the
graded case, we define the generic initial ideal ${\rm \bf gin} I$ of a
homogeneous ideal $I$ of $S$ and we establish its basic properties.
In particular, we may use the Bayer-Stillman theorem to compute
the regularity of a Borel-fix ideal in $S$ when $k$ has
characteristic zero.
For $j=1, 2$,  let us denote by $\delta_j(I)$
the maximum of the $j$-th component of the degrees in a minimal
homogeneous system of generators of $I$. Then we have

\medskip
\noindent
{\bf Proposition 35} [Proposition 5.3.10]
{\it Let $I \subset S$ be a Borel-fix ideal. If $char k =0$, then
$$\reg_1(I) = \delta_1(I),$$
$$\reg_2(I) = \delta_2(I).$$}

\medskip
This result has been also proved by A. Aramova et al. \cite{ADK} by
different methods.
In the graded case, D. Bayer and M. Stillman  \cite{BS1} also
proved the existence of an order in the polynomial ring
$A=k[{\x}]$ (the reverse lexicographic order) such that
$\reg I = \reg ({\rm gin} I)$ for any homogeneous ideal $I$ of $A$.
We finish the chapter by showing that the analogous bigraded
result does not hold because we can find a homogeneous ideal $I$ of
$S$ such that for any order
${\rm \bf reg}(I) \not = {\rm \bf reg}({\rm \bf gin} I)$.

\medskip
In  {\bf Chapter 6} we study the asymptotic properties
of the powers of a homogeneous ideal $I$ in the polynomial ring
$A=k[{\x}]$. We will show how the bigraded structure of the Rees algebra
provides information about the Hilbert polynomials, the
Hilbert series and the graded minimal free resolutions of the powers of
$I$. This grading of the Rees algebra will be also useful
to study the mixed multiplicities of the Rees algebra and the form
ring of an equigenerated ideal.

\medskip
\noindent
{\bf Theorem 36} [Theorem 6.1.1]
{\it Let $I$ be a homogeneous ideal of $A$. Set $c= a_*^2(R_A(I))$, $h
= \h (I)$. Then there are polynomials $e_0(j), \dots, e_{n-h-1}(j)$
with integer values such that for all $j \geq c+1$
$$P_{A/I^j}(s) = \sum_{k=0}^{n-h-1} (-1)^{n-h-1-k}e_{n-h-1-k}(j) {s+k
\choose k}.$$
Furthermore, $\deg e_{n-h-1-k}(j) \leq n-k-1$ for all k. }

\medskip
In particular, this result says that a finite set of Hilbert polynomials
of the powers of an ideal allows to
compute the Hilbert polynomials of its Rees algebra and its
form ring, without needing an explicit presentation of these bigraded
algebras. For equigenerated ideals, we may also compute the
multiplicities of their Rees algebras and form rings.

\medskip
\noindent
{\bf Corollary 37} [Corollary 6.1.8]
{\it Let $I$ be a homogeneous ideal in $A$. Let  $c=a_*^2(R_A(I))$,
$h = \h (I)$. Then the Hilbert polynomials of $I^j$ for $c + 1 \leq j
\leq c+n$ determine
\begin{enumerate}
\item The polynomials $e_{n-h-1-k}(j)$ for $k = 0, \dots, n-h-1$.
\item The Hilbert polynomials of $A/I^j$ for $j> c + n$.
\item The Hilbert polynomial of $R_A(I)$ and the
Hilbert polynomial of $G_A(I)$.
\item If $I$ is equigenerated and not $\frak m$-primary, the
mixed multiplicities of $R_A(I)$ and $G_A(I)$.
\end{enumerate}}

\medskip
A similar result can be proved for the Hilbert series of the powers of
a homogeneous ideal. Namely,

\medskip
\noindent
{\bf Proposition 38} [Theorem 6.2.1, Proposition 6.2.7]
{\it Let $I$ be a homogeneous ideal. Set $r = \mu (I)$,
$l=l(I)$, $c =a_*^2(R_A(I))$. Then:
\begin{enumerate}
\item The Hilbert series of $I^j$ for $ j \leq c+r$ determine
the Hilbert series of $I^j$ for $j > c +r$.
\item If $I$ is an equigenerated ideal, the Hilbert series of $I^j$
for $c + 1 \leq j \leq c + l$  determine the Hilbert series of
$I^j$ for $j > c + l$.
\end{enumerate}
}

\medskip
Next we study the behaviour of the projective dimension of
the powers of an ideal. As a by-product, we recover the classic
result of M. P. Brodmann \cite{Bro} which says that the depth of the
powers of an ideal becomes constant asymptotically, and a result of
D. Eisenbud and C. Huneke \cite{EH} which precises this asymptotic value
under some restrictions. Moreover, for ideals whose form ring is
Gorenstein we may determine exactly the powers of the ideal for which
the projective dimension takes the asymptotic value. Namely,

\medskip
\noindent
{\bf Proposition 39} [Proposition 6.3.2]
{\it Let $I$ be a homogeneous ideal in $A$ and set $l=l(I)$. If $G$ is
Gorenstein, ${\dpp}_A (I^j) \leq l-1$ for all $j $, and
${\dpp}_A I^j = l - 1$ if and only if
 $j > a^2(G) - a(F)$.}

\medskip
Finally, we show that the graded minimal free resolutions of the powers
of an ideal also have a uniform behaviour. For equigenerated
ideals, we can prove that the shifts which arise in the minimal
resolutions are linear functions asymptotically and the Betti numbers
are polynomial functions asymptotically. More explicitly,

\medskip
\noindent
{\bf Proposition 40} [Proposition 6.3.6]
{\it Let $I$ be a homogeneous ideal generated in degree $d$. Set $l =
l(I)$, $s = n - \depth_{(\frak m R)} (R)$. Then there is a finite set of
integers $\{\alpha_{pi} \mid 0 \leq p \leq s, 1 \leq i \leq k_p \}$
and polynomials
$\{Q_{\alpha_{pi}}(j) : 0 \leq p \leq s, 1 \leq i \leq k_p \}$
of degree $\leq l-1$ such that the graded minimal free
resolution of $I^j$ for $j$ large enough is
$$ 0 \to D_s^j \to \dots \to D_0^j \to I^j \to 0 \,\, ,$$
\noindent
with $D_p^j = \bigoplus_i A(-\alpha_{pi} -dj)^{\beta_{pi}^j}$
and ${\beta_{pi}^j}= Q_{\alpha_{pi}} (j)$.}

\medskip
From this result, we may deduce that a finite number of the graded
minimal free resolutions of the powers of an ideal determine the
rest of them. This finite set of resolutions can be found for ideals
with a very particular behaviour. For instance, we get

\medskip
\noindent
{\bf Proposition 41} [Proposition 6.3.10]
{\it Let $I$ be an equigenerated homogeneous ideal, and $b =
a_*^2(R_A(I)) + l(I)$. If the graded minimal free resolutions of $I,
I^2, \dots, I^{b}$ are linear, then the graded minimal free resolutions
of $I^j$ are also linear for any $j$. Furthermore, the minimal free
resolutions of $I, I^2, \dots, I^b$ determine the minimal graded free
resolutions of $I^j$ for any $j$.}

\vspace{30mm}
Some parts of this work have already appeared published in:

\medskip
\noindent
- O. Lavila--Vidal, {\em On the Cohen-Macaulay property of diagonal
subalgebras of the Rees algebra}, manuscripta math. 95 (1998),
47--58.

\medskip
\noindent
- O. Lavila--Vidal, S. Zarzuela, {\em On the Gorenstein property
of the diagonals of the Rees algebra},
Collect. Math. 49, 2-3 (1998), 383--397.


\chapter*{$\;\;\;$}

\quad

\vspace{50mm} It is a pleasure to thank my thesis advisor Santiago
Zarzuela for his help, guidance and support throughout these
years.

I would also thank the late Pr.Dr. Manfred Herrmann and his
Seminar for their hospitality while I was visiting the
Mathematisches Institut der Universit\"{a}t zu K\"{o}ln.

I also wish to thank my fellows at the Departament d'\`{A}lgebra i
Geometria and at the Seminari d'\`{A}lgebra Commutativa, especially
Maria Alberich.

Finally, let me thank my friends, my family and Eduard for their
love, patience and encouragement.

\newpage

\pagenumbering{arabic}
\chapter{Multigraded rings}
\typeout{Multigraded rings}

\bigskip
In this chapter we collect some basic definitions and facts of the
theory of multigraded rings which we will need in the next
chapters. We also state the multigraded versions of some
well-known results in the category of graded rings. Rings are
always assumed to be noetherian.

\bigskip
\section{Multigraded rings and modules}
\markboth{CHAPTER I. MULTIGRADED RINGS}{MULTIGRADED RINGS AND
MODULES}

\medskip
The general theory of multigraded rings and modules is analogous
to that of graded rings and modules. We first recall some basic
definitions. The main sources are \cite{BH}, \cite{HHR} and
\cite{GW1}.

\medskip
We use the following multi-index notation. For $\bn =(n^1, \dots,
n^r) \in \Bbb Z^r$, we set $\vert \bn \vert = n^1+\dots+n^r$, and
for $\bn, \bm \in \Bbb Z^r$, we define their sum $\bn + \bm =
(n^1+m^1, \dots, n^r+m^r)$, and we set $\bn < \bm$ ($\bn \leq
\bm$) if
 $n^i < m^i$ ($n^i \leq m^i$) for every $i$.

\medskip
  A $\Bbb Z^r$-graded ring (or $r$-graded ring) is a ring $S$ endowed with a direct sum decomposition
$S = \bigoplus_{\bn \in \Bbb Z^r} S_{\bn}$, such that $S_{\bn}
S_{\bf m} \subset S_{\bn + \bf m}$ for all $\bn, {\bf m} \in \Bbb
Z^r$.  An $r$-graded $S$-module is an $S$-module $M$ endowed with
a decomposition $M = \bigoplus_{\bn \in \Bbb Z^r} M_{\bn}$, such
that $S_{\bn} M_{\bf m} \subset M_{\bn + \bf m}$ for all $\bn,
{\bf m} \in \Bbb Z^r$. We shall call $M_{\bn}$ the homogeneous
component of $M$ of degree $\bn$.  An element $x \in M$ is
homogeneous of degree $\bn$ if $x \in M_{\bn}$. The degree of $x$
is then denoted by $\dg  x$. For any $r$-graded $S$-module $M$, we
define the support of $M$ to be the set $\supp M = \{{\bf n} \in
\Bbb Z ^r \mid M_{\bf n} \not =0 \}$.

\medskip
   For a given $r$-graded ring $S$, we may consider the category of
    $r$-graded $S$-modules $M^r(S)$. Its objects are the $r$-graded
$S$-modules, and a morphism $f : M \to N$ in $M^r(S)$ is an
$S$-module morphism such that $f(M_{\bn}) \subset N_{\bn}$ for all
$\bn \in \Bbb Z^r$.

\medskip
    Given an $r$-graded $S$-module $M$, an $r$-graded submodule is a
submodule $N \subset M$ such that $N = \bigoplus_{\bn \in \Bbb
Z^r} N \cap M_{\bn}$, equivalently, $N$ is generated by
homogeneous elements. The $r$-graded submodules of $S$ are called
homogeneous ideals. For an arbitrary ideal $I$ of $S$, the
homogeneous ideal $I^*$ is defined to be the ideal generated by
all the homogeneous elements of $I$.

\medskip
   As a first example of $r$-graded ring we have the
polynomial ring $S=A[{\x}]$ defined over an arbitrary ring $A$.
For every choice of elements $\bd_1,\dots,\bd_n \in \Bbb Z^r$, we
have a unique $r$-grading on $S$ such that $\dg X_i = \bd _i$ and
$\dg  a = 0$ for all $a \in A$.

\medskip
    For an $r$-graded $S$-module $M$ and ${\bf k} \in \Bbb Z^r$, then
$M({\bf k})$ denotes the $S$-module $M$ with the grading given by
$M({\bf k})_{\bn} = M_{{\bf k} + \bn} $.

\medskip
If $M, N$ are $r$-graded $S$-modules, we denote by
$\uHom_S(M,N)_0$ the abelian group of all the homomorphisms of
$r$-graded $S$-modules from $M$ into $N$. We set $\uHom_S(M,N)=
\bigoplus_{\bn \in \Bbb Z^r} \uHom_S(M,N(\bn))_0$. Note that
$\uHom_S(M,N)_{\bf k}$ is nothing but the abelian group of
$S$-module homomorphisms $f: M \to N$ such that $f(M_{\bn})
\subset N_{\bn + \bf k}$ for all $\bn \in \Bbb Z^r$. The derived
functors of $\uHom_S(\;,\;)$ are $\uExt_S^i(\;,\;)$, with $i \in
\Bbb N$.

\bigskip
\section{Multigraded cohomology}
\label{B} \markboth{CHAPTER I. MULTIGRADED RINGS} {MULTIGRADED
COHOMOLOGY}

\medskip
Next we are going to introduce the local cohomology functor in the
category of multigraded modules, mainly following \cite{HHR}. The
basic results are the multigraded version of the Local Duality
Theorem and the good behaviour of the local cohomology modules
under a change of grading.

\medskip
From now on in this chapter, we assume that $S = \bigoplus _{{\bn}
\in \Bbb N^r} S_\bn$ is an $r$-graded ring defined over a local
ring $S_0 =A$. Then $S$ has a unique homogeneous maximal ideal $
{\cal M} = {\frak m} \oplus ( \bigoplus_{{\bn} \not = 0}
S_{\bn})$, where $\frak m$ is the maximal ideal of $A$. Set $ d =
\dim S$.

\medskip
 If $I \subset S$ is a homogeneous ideal and $M$ is an
$r$-graded $S$-module, we denote by
 $\uH^0_{I}(M) =  \Gamma_{I}(M) =
 \{ x \in M : I^k x = 0 \;{\rm for} \; {\rm some} \; k \geq 0 \}$.
Note that $\uH^0_I(M)$ is an $r$-graded submodule of $M$. The
local cohomology functors ${\uH^i_{I}} (\;)$ are the right derived
functors of $\Gamma_{I}(\;)$ in the category of $r$-graded
$S$-modules. If no confusion, we will usually denote them by
$H^i_{I} (\;).$

\medskip
 An $r$-graded $S$-module $K_S$ is called a canonical
module of $S$ if
 $$ K_S \otimes_A \widehat{A}  \cong
  \uHom_S (\uH^d_{\cal M}(S), {\underline {\rm E}}_S(k)) \,,$$
where $k$ is the residue field of $A$ and ${\underline {\rm
E}}_S(k)$ is the injective envelope of $k$ in the category of
$r$-graded $S$-modules. The injective envelope ${\underline {\rm
E}}_S(k)$ of $k$ is $\uHom_A (S, {\rm E}_A(k))$, where $A$ is
thought as an $r$-graded ring concentrated in degree 0, and both
$S$ and ${\rm E}_A(k)$ are considered as $r$-graded $A$-modules.
Therefore, we have
$$  K_S \otimes_A \widehat{A} \cong
 \uHom_A (\uH^d_{\cal M}(S), {\rm E}_A(k))
= \bigoplus_{\bn \in \Bbb Z^r} \Hom_A ([\uH^d_{\cal M}(S)]_{-\bn},
E_A(k)).$$ If a canonical module exists, it is finitely generated
and unique up to an isomorphism. In the particular case where
$A=k$ is a field, the canonical module of $S$ exists and
$$K_S \cong \uHom_k(\uH^d_{\cal M}(S), k).$$
The next results are the extension to the $r$-graded case of two
of the main properties of the canonical module, well-known for the
graded case (see \cite[Theorem 2.2.2]{GW2}).

\begin{thm}
\label{B2} (Local Duality) Let $S$ be an $r$-graded ring defined
over a complete local ring $A$.
 Let ${\cal M}$ be the homogeneous maximal ideal of $S$.
Then $S$ is Cohen-Macaulay if and only if every finitely generated
$r$-graded $S$-module $M$ satisfies
$${\uHom}_S ( \uH^i_{\cal M}(M), {\underline {\rm E}}_S(k))
\cong {\uExt}^{d-i}_S (M, K_S) \;, \, i = 0, \dots, d.$$
\end{thm}

\begin{cor}
\label{B3} Let $S$ be a Cohen-Macaulay $r$-graded ring with
canonical module $K_S$. Let $T$ be an $r$-graded ring defined over
a local ring such that there exists a finite $r$-graded ring
morphism $S \to T$. Then $T$ has canonical module
$$K_T = {\uExt}^e_S (T, K_S) \,,$$
where $e= \dim S - \dim T$.
\end{cor}

\medskip
Often we are going to consider the ring $S$ endowed with a
different grading obtained in the following way: given a group
morphism
 $\varphi: {\Bbb Z}^r \rightarrow {\Bbb Z}^q$
such that $\varphi ( $supp $S) \subset {\Bbb N}^q$, we can define
the ${\Bbb N}^ q$-graded ring
 $$S^\varphi := \bigoplus_{{\bf m} \in {\Bbb N}^q}
(\bigoplus_{\varphi({\bf n})= {\bf m}} S_{\bf n}).$$ Similarly,
given an $r$-graded $S$-module $M$, we may define the $q$-graded
$S^\varphi$-module $M^\varphi$ as
 $$M^\varphi:= \bigoplus_{{\bf m} \in {\Bbb Z}^q}
(\bigoplus_{\varphi({\bf n})={\bf m}} M_{\bf n}).$$ Then  $(\;
\;)^\varphi :
 M^r(S)  \rightarrow M^q(S^\varphi )$ is an exact functor.
 By considering $\varphi_j :  {\Bbb Z}^r  \rightarrow  {\Bbb Z} $ the
projection on the $j$-component, that is, $\varphi_j({\bf
n})=n^j$, we denote by $S_j= S^{\varphi_j}$ and by $M_j=
M^{\varphi_j}$. Note that $S_j$ is just the ring $S$ graded by the
$j$-th partial degree.

\medskip
 The next lemma shows that the local cohomology modules
behave well under a change of grading.

\begin{lem}
\label{B1} \cite[Lemma 1.1]{HHR} Let $S$ be an $r$-graded ring
defined over a local ring.
 Let ${\cal M}$ be the homogeneous maximal ideal of $S$.
Let $\varphi : {\Bbb Z}^r \rightarrow {\Bbb Z}^q$ be a morphism
such that $\varphi ( \supp S) \subset {\Bbb N}^q$. For every
$r$-graded $S$-module $L$, we have
         $${\underline H^i_{\cal M} (L)}^\varphi = \underline H^i_{{\cal
M}^\varphi}(L^\varphi),  \;\forall i .$$
\end{lem}

\bigskip
\section{Multigraded a-invariants}
\label{B} \markboth{CHAPTER I. MULTIGRADED RINGS} {MULTIGRADED
$a$-INVARIANTS}

\medskip
We begin this section by extending the definition of the
$a$-invariants of a graded module to the multigraded case. After
that, under some mild assumptions, we relate the multigraded
$a$-invariants of a multigraded module to the shifts which appear
in its multigraded minimal free resolution. This result will be
essential in the next chapters. In the graded case, a similar
result can be found in \cite[Example 3.6.15]{BH} for
Cohen-Macaulay modules.

\medskip
Let $S$ be a $d$-dimensional $\Bbb N^r$-graded ring defined over a
local ring. For each $i= 0,..., d$, the multigraded
$a_i$-invariant of $S$ is ${\bf a}_i (S) = (a_i^1(S),...,
a_i^r(S))$, where
$$a_i^j (S)= \max \,\{ m \in {\Bbb Z} \mid \exists
     {\bf n} \in {\Bbb Z}^r :
     \varphi_j ({\bf n}) = m ,
      [\uH^i_{\cal M} (S)]_{\bf n} \not =0\}$$
\noindent if $\uH^i_{\cal M}(S) \not =0$ and $a_i^j(S)=-\infty$
otherwise. We will denote by ${\bf a}(S) = {\bf a}_d (S)$. Note
that by Lemma \ref{B1}
 $$a_i^j (S) =  \max \,\{ m \in {\Bbb Z}  \mid [\uH^i_{{\cal M}_j}
(S_j)]_m \not =0 \} = a_i(S_j).$$ Following N.V. Trung \cite{T2},
the multigraded $a_*$-invariant of $S$ is defined as ${\bf a}_*(S)
= (a_*^1(S),..., a_*^r(S))$, where $a_*^j (S) = \max \{ a_0^j (S),
\dots, a_d^j(S) \}$. Similarly, for any finitely generated
$r$-graded $S$-module $M$ we may define the a-invariants ${\bf
a}_i(M)$ of $M$ and the $a_*$-invariant ${\bf a}_*(M)$ of $M$.

\medskip
Observe that if there exists $K_S$ the canonical module of $S$ ,
then
 $$a^j(S) = a^j_d(S) = -\min \, \{ m \in \Bbb Z \mid \exists
         {\bf n} \in {\Bbb Z^r} :
         \varphi_j ({\bf n}) = m ,
         [K_S]_{\bf n} \not = 0 \}.$$

\medskip
If $S$ has a canonical module $K_S$, $S$ is said to be
quasi-Gorenstein if there exists an $r$-graded isomorphism $K_S
\cong S({\bf a}(S))$, and Gorenstein if in addition $S$ is
Cohen-Macaulay.

\medskip
From now on in this section we assume that $S$ is a noetherian
$\Bbb N^r$-graded algebra defined over a field $k$, and let $\cal
M$ be its homogeneous maximal ideal. Our main purpose is then to
compute the multigraded a-invariants of a finitely generated
$r$-graded $S$-module $M$ from an $r$-graded minimal finite free
resolution of $M$ over $S$, whenever it exists and $S$ is
Cohen-Macaulay. To begin with, let us consider
$${\dots \to {D_t}\to \dots \to {D_1} \to {D_0} \to 0}$$
an exact sequence of finitely generated $r$-graded $S$-modules
such that Im${(D_{p+1}) \subset {\cal M} D_p}$, for all $p \geq
0$. Let us denote by  $ \{\,{\bf v}_{pq} \}$ the set of degree
vectors of a minimal homogeneous system of generators of ${D_p}$.
Note that this set is uniquely determined because it can be
obtained as the homogeneous components of the vector space $D_p
\otimes_S k$ which are not zero. We set ${\bf
m}_p=\min_{\leq_{lex}} \{\,{\bf v}_{pq} \}$ and ${\bf M}_p=\max_{
\leq_{lex}} \{\, {\bf v}_{pq} \}$, where $\leq_{lex}$ is the
lexicographic order. Let us denote by  $ n^j_p=\min_q\{ v^j_{pq}
\}$,
 $ t^j_p =\max_q\{ v^j_{pq} \} $, where
 ${\bf v}_{pq}= (v^1_{pq}, \dots, v^r_{pq})$, and
 ${\bf n}_p=(n^1_p,...,n^r_p)$,
${\bf t}_p=(t^1_p,...,t^r_p)$. Let us also consider $\leq$ the
partial order in $\Bbb Z^r$ defined coefficientwise. Then we have

\begin{lem}
\label{B4}

\begin{enumerate}
\item   ${\bf n}_p \leq {\bf n}_{p+1}$.
\item  ${\bf m}_p <_{lex} {\bf m}_{p+1} $
\end{enumerate}
\end{lem}

{\pf}   Let $C_p= {\rm Coker}(D_{p+1} \rightarrow D_p)$, $\forall
p \geq 1$. Then there are short exact sequences
  $$ 0 \rightarrow C_{p+2} \rightarrow D_{p+1} \rightarrow C_{p+1}
     \rightarrow 0 , \; \; \forall p \geq 0.$$
Applying the functor $- \bigotimes_S k $, we get exact sequences
  $$ C_{p+2}/ {\cal M} C_{p+2}  \rightarrow
  D_{p+1}/ {\cal M} D_{p+1}
\rightarrow C_{p+1}/ {\cal M} C_{p+1}
 \rightarrow 0 , \; \; \forall p \geq 0.$$
Since $C_{p+2} \subset {\cal M} D_{p+1}$, then the first map is
the zero morphism. Therefore we get isomorphisms
 $$D_{p+1}/ {\cal M} D_{p+1} \stackrel{\cong} \longrightarrow
C_{p+1}/ {\cal M} C_{p+1} .$$ Let us denote by $\{\, e_{pq}\}$ a
minimal homogeneous system of generators of $D_p$ with
deg$(e_{pq})= {\bf v}_{pq}$, and let $f$ be the map from $D_{p+1}$
to $D_p$. From the isomorphism it follows that
 $f(e_{p+1,q})  \not =0$, for all $q$. Now let us fix $q$. We can write
 $f(e_{p+1,q}) =\sum_{l} \lambda_l e_{pl}$,
where $\lambda_l $ are homogeneous elements of ${\cal M}$. Set
$\dg (\lambda_l) = (\lambda^1_l,...,\lambda^r_l) \in \Bbb N^r $
and note that $\dg(\lambda_l) \not =  0$ if $\lambda_l \not =0$.
Looking at the $j$-th component of the degree, we get
 $v^j_{p+1,q} \geq  \min_l\{\ v^j_{pl} \} =n^j_p$,
and so $n^j_{p+1} \geq n^j_{p}$ for all $j$.

To obtain $(ii)$, it is enough to prove that ${\bf v}_{p+1,q}
>_{lex} {\bf m}_p$ for all $q$. We have already shown that
$v^1_{p+1,q} \geq \min_l\{\ v^1_{pl}\} = m^1_p $. If $v^1_{p+1,q}
> m^1_p$, we are done. Otherwise, $v^1_{p+1,q} = m^1_p$ and  so
$\lambda^1 _l =0 $  for each $l$ such that $\lambda _l   \not =0$.
Then we have $v^2_{p+1,q} \geq  \min_l\{\ v^2_{pl} \mid v^1_{pl} =
m^1_p \}=  m^2_p $. By repeating this argument, we get the result
since there exist $l,j$ such that $\lambda_l^j >0$. $\B$

\medskip
Let $S$ be a $d$-dimensional $r$-graded Cohen-Macaulay
$k$-algebra. Assume that $M$ is a finitely generated $r$-graded
$S$-module with a finite minimal $r$-graded free resolution over
$S$
$${0\to {D_l}\to \dots \to {D_1} \to {D_0} \to {M} \to 0},$$
with ${D_p}=\bigoplus_q S(a^1_{pq},...,a^r_{pq})$. Set $m = \dim
M$, $\rho= \depth \,M$. Note that $l=d-\rho$ by the graded
Auslander-Buchsbaum formula. Next we are going to study the shifts
which appear in this resolution.

\medskip
Note that, with the notation introduced before,

\vspace{3mm}

\hspace{25mm} $n^j_p = \min_q \{-a^j_{pq} \},$

\vspace{1.5mm}

\hspace{25mm} $t^j_p =\max_q\{ - a^j_{pq} \} $,

\vspace{1.5mm}

\hspace{25mm} ${\bf m}_p = \min_{\leq_{lex}}
\{\,(-a^1_{pq},...,-a^r_{pq})\}$,

\vspace{1.5mm}

\hspace{25mm} ${\bf M}_p=\max_{ \leq_{lex}} \{\,
(-a^1_{pq},...,-a^r_{pq})\}$.

\vspace{3mm} \noindent We will also denote by $t_p^j(M)=t_p^j$,
$t_*^j(M)= \max \{t_0^j, \dots, t_l^j\}$, ${\bf t}_*(M)=
(t_*^1(M), \dots, t_*^r(M)) $. From Lemma \ref{B4}, we have ${\bf
n}_p \leq {\bf n}_{p+1}$, ${\bf m}_p <_{lex} {\bf m}_{p+1}$.
Furthermore,

\begin{lem}
\label{B5}
\begin{enumerate}
\item ${\bf M}_{0} <_{lex} {\bf M}_1 <_{lex} \cdots <_{lex}
{\bf M}_{d-m-1} <_{lex}  {\bf M}_{d-m}$.
\item ${\bf t}_{0} \leq {\bf t}_1 \leq \cdots \leq {\bf t}_{d-m-1}
\leq  {\bf t}_{d-m}$.
\end{enumerate}
\end{lem}

{\pf}   Let $K_S$ be the canonical module of $S$. Note that it
exists because $S$ is a finitely generated $k$-algebra. By setting
$C_p= {\rm Coker}(D_{p+1} \rightarrow D_p)$ for $p \geq 0$, we get
short exact sequences
    $$ 0 \to C_{p+1} \to D_{p} \to C_{p} \to 0,  $$
\noindent for $0 \leq p \leq l-1$, where $C_0= M$, $C_l = D_l$.
For any $p<d-m-1$, we have
$$\uExt^1_S(C_p, K_S) \cong
\uExt^2_S(C_{p-1}, K_S) \cong  \cdots
 \cong \uExt^{p+1}_S (M, K_S) =0$$
by Theorem \ref{B2}. Therefore, by applying the functor $(\,\,)^*
= \uHom_S (\,\,,K_S)$ to the sequences above for $p \leq d-m-1$,
we get exact sequences
 $$ 0 \to C_{p}^* \to
      D_{p}^* \to C_{p+1}^* \to 0, \; \; {\rm for} \; p \leq d-m-2, $$
 $$  0 \to C_{d-m-1}^* \to D_{d-m-1}^* \to C_{d-m}^* \to
     H_{\M}^{m}(M)^\vee \to 0,$$
where $(\;)^\vee = \uHom_k (\;, k)$. By gluing these exact
sequences, we get the $r$-graded exact sequence
    $$ 0 \to D_0^* \to \dots \to
     D_{d-m-1}^* \to
     C_{d-m}^* \to H_{\M}^m (M)^\vee \to 0.
     $$
\noindent Observe that $D_p^* = \bigoplus_q
K_S(-a^1_{pq},...,-a^r_{pq})$. One can also check that
Im$(D_{p}^*)   \subset  {\cal M}{D_{p+1}^*}$ for all $p \leq
d-m-2$.

 Let $\{\,{\bf b}_1,...,{\bf b}_k \}$ be the set of degree vectors of a
minimal homogeneous system of generators of $K_S$. If we denote by
${\bf a}_{pq}=(a^1_{pq},...,a^r_{pq})$, then the vectors ${\bf
a}_{pq}+{\bf b}_i$ are the degrees of a minimal homogeneous system
of generators of $D^*_p$. For $p \leq d-m-1$, let us consider

\vspace{3mm}

\hspace{15mm} $ {\bf \widetilde{m}}_p= \min_{\leq_{lex}} \{\,
 {\bf a}_{pq}+ {\bf b}_i \} = -{\bf M}_p +
   \min_{\leq_{lex}} \{\,{\bf b}_i \} \;,$

\vspace{3mm}

\hspace{15mm} $\widetilde{n}^j_p= \min_{q,i} \{\,a^j_{pq}+b^j_i
\}= -t^j_p   + \min_i \{\,b^j_i \} \;.$

\vspace{3mm} \noindent According to Lemma \ref{B4}, we have
$\widetilde{n}^j_{p+1}  \leq \widetilde{n}^j_p$ and ${\bf
\widetilde{m}}_{p+1} <_{lex} {\bf \widetilde{m}}_p$, so

\vspace{3mm}

\hspace{15mm} ${\bf t}_{0} \leq {\bf t}_1 \leq \cdots \leq {\bf
t}_{d-m-2} \leq  {\bf t}_{d-m-1}$

\vspace{1.5mm}

\hspace{15mm} ${\bf M}_{0} <_{lex} {\bf M}_1 <_{lex} \cdots
<_{lex} {\bf M}_{d-m-2} <_{lex}  {\bf M}_{d-m-1}.$

\vspace{3mm}

Next we want to show that ${\bf M}_{d-m} >_{lex} {\bf M}_{d-m-1}$.
To this end, let us study the morphism $ D_{d-m} \to D_{d-m-1}$
and for that denote by $\nu : C_{d-m} \to D_{d-m-1}$. Assume that
there is an element $u$ in the basis of $D_{d-m-1}$ of degree
${\bf M}_{d-m-1} \geq_{lex} {\bf M}_{d-m}$. If $g$ is a
homogeneous minimal generator of $C_{d-m}$, then $g$ has trivial
terms in $u$: Otherwise, we would have that ${\bf M}_{d-m-1}<
_{lex} \deg g \leq_{lex} {\bf M}_{d-m}$ because $C_{d-m} \subset
{\cal M} D_{d-m-1}$. Let ${\bf b}= \min_{<_{lex}} \{\,{\bf
b}_1,...,{\bf b}_k \}$, and let us take $c \in [K_S]_{\bf b}$, $c
\not =0$. Let $w: D_{d-m-1} \rightarrow K_S$ defined by $w(u)=c$,
$w(v)=0$ for any $v \not = u$ homogeneous element in the basis of
$D_{d-m-1}$. Then $\nu^* : D_{d-m-1}^* \to C_{d-m}^*$  satisfies
$\nu^*(w) = 0$, hence $\nu^*$ is not a monomorphism in degree $
\deg w = \deg w(u) - \deg (u)= {\bf b}- {\bf M}_{d-m-1}$.
Therefore $[C_{d-m-1}^*]_{ {\bf b}-{\bf M}_{d-m-1}} \not =0$, and
then $[D_{d-m-2}^*]_{ {\bf b} -{\bf M}_{d-m-1}} \not =0$, so there
exists a shift  ${\bf a}=(a^1, \dots, a^r)$ in $D_{d-m-2}$ such
that $ - {\bf a} \geq_{lex} {\bf M}_{d-m-1}$. So we obtain ${\bf
M}_{d-m-2} \geq_{lex} {\bf M}_{d-m-1}$ which is a contradiction.

    Furthermore, note that the first component of ${\bf M}_p$ is
$t^1_p$. Therefore, we have $t_{d-m-1}^1 \leq t_{d-m}^1$ since
${\bf M}_{d-m-1}<_{lex} {\bf M}_{d-m}$. The inequalities
$t_{d-m-1}^j \leq t_{d-m}^j$ for $j=2, \dots, r$ follow directly
from the next remark. $\B$

\begin{rem}
\label{B6} {\rm Given a permutation $\sigma$ of $\{1, \dots, r
\}$, we may define $\leq_{\sigma}$ to be the order in $\Bbb Z^r$
defined by
$$ (u_1,...,u_r) \leq_{\sigma} (v_1,...,v_r) \iff
(u_{\sigma(1)},...,u_{\sigma(r)}) \leq_{lex}
(v_{\sigma(1)},...,v_{\sigma(r)}).$$ Then Lemmas \ref{B4} and
\ref{B5} also hold if we define
$${\bf m}_{p}^{\sigma} = \min_{{\leq}_{\sigma}}
\{\,(-a^1_{pq},...,-a^r_{pq}) \},$$
$${\bf M}_{p}^{\sigma} = \max_{{\leq}_{\sigma}}
\{\,(-a^1_{pq},...,-a^r_{pq}) \}. $$

}
\end{rem}

\medskip
The following result gives a formula for the multigraded
$a_*$-invariant of $M$ by means of the shifts which arise in its
resolution over $S$ (see \cite[Example 3.6.15]{BH} for the case of
a $\Bbb Z$-graded Cohen-Macaulay module).

\begin{thm}
\label{B77} For each $j = 1, \dots, r$, we have
\begin{enumerate}
\item $a_{d-p}^j(M) \leq t_p^j(M) +a^j(S)$, for $d-m \leq p \leq
d-\rho$.
\item Assume that for some $p$ there exists $\sigma$
s.t. $\sigma(1)=j$ and ${\bf M}_{p}^\sigma >_{\sigma} {\bf
M}_{p+1}^\sigma$. Then $a_{d-p}^j(M) = t_p^j(M) +a^j(S)$.
\item
$a_*^j(M) = t_*^j(M)+ a^j(S)$. That is, ${\bf a}_*(M) = {\bf
t}_*(M) + {\bf a}(S)$.
\end{enumerate}

\end{thm}

{\pf} From the minimal $r$-graded free resolution of $M$ over $S$
    $$ 0 \to D_l \to \dots \to D_0 \to M \to 0 ,$$
\noindent by setting $C_p= {\rm Coker}(D_{p+1} \rightarrow D_p)$,
we have short exact sequences
    $$ 0 \to C_{p+1} \to D_{p} \to C_{p} \to 0,  $$
\noindent for $0 \leq p \leq l-1$. By Theorem \ref{B2}, if we
apply the functor $(\,\,)^* = \uHom_S (\,\,,K_S)$ to the sequences
above we get exact sequences

\vspace{3mm}

\hspace{8mm} $ (1) \;\; 0 \to D_0^* \to \dots \to D_{d-m-1}^* \to
     C_{d-m}^* \to H_{\M}^m (M)^\vee \to 0, $

\vspace{2mm}

\noindent and

\vspace{2mm}

\hspace{8mm} $ (2) \;\; 0 \to C_{p}^* \to
   D_{p}^* \to C_{p+1}^* \to 0, \; \; {\rm for} \;p \leq d-m-2, $

\vspace{3mm}

\hspace{8mm} $ (3) \;\;
   0 \to C_{p-1}^* \to D_{p-1}^* \to C_{p}^* \to
     H_{\M}^{d-p}(M)^\vee \to 0,\;\; {\rm for} \;  p \geq d-m \,,$

\vspace{3mm} \noindent where $(\;)^\vee = \uHom_k (\;, k)$. Note
that for $d-m \leq p \leq d-\rho$ we have monomorphisms
$$0 \to C_{p}^* \to D_{p}^*=
\bigoplus_q K_S (-a_{pq}^1, \dots, -a_{pq}^r),$$

\noindent and so $[C_{p}^*]_{-\bf i} = 0$ for any $\bf i$ such
that $i^1 > t_p^1 + a^1(S)$. Now from the epimorphisms
$$C_{p}^* \to H_{\M}^{d-p}(M)^\vee \to 0,$$

\noindent we get $H_{\M}^{d-p}(M)_{\bf i} = 0$ if $i^1  > t_p^1 +
a^1(S)$, and therefore $a_{d-p}^1(M) \leq t_p^1(M)+a^1(S)$. This
proves $(i)$ for the case $j=1$.

Assume now that there exists $p$ with ${\bf M}_p >_{lex}{\bf
M}_{p+1}$ (then $p \geq d-m$ according to Lemma \ref{B5}). Let
${\bf b}=(b^1, \dots, b^r)$ be the minimum with respect to the
lexicographic order such that $[K_S]_{\bf b} \not =0$. Note that
$b^1=-a^1(S)$. Let ${\bf i} = {\bf M}_p - {\bf b}$. Since
$[D_{p+1}^*]_{- {\bf i}}=0$ because ${\bf M}_p >_{lex} {\bf
M}_{p+1}$, we have $[C_{p+1}^*]_{-{\bf i}}=0$ by $(3)$, and so
$[C_{p}^*]_{-{\bf i}} = [D_{p}^*]_{-{\bf i}}$ also by $(3)$. Then,
denoting by $f: D_p \to D_{p-1}$, we get an exact sequence
$$[D_{p-1}^*]_{-{\bf i}} \stackrel{f^*}{\rightarrow}
[D_{p}^*]_{-{\bf i}} \rightarrow [H_{\M}^{d-p}(M)]_{{\bf i}} \to
0.$$ \noindent Let $e_1, \dots, e_s$ be the elements of the
canonical basis of $D_p$ with degree ${\bf M}_p$, and let $v_1,
\dots, v_m$ be the canonical basis of $D_{p-1}$. Since $f^*(
D_{p-1}^*) \subset {\cal M} D_{p}^*$ and $[D_{p}^*]_{-{\bf i}} =
[K_S]_{{\bf b}} e_1^* \oplus \dots \oplus [K_S]_{{\bf b}} e_s^*$,
we have that $ [K_S]_{{\bf b}} e_1^* \not \subset {\rm Im} f^*$.
In particular, $f^*$ is not an epimorphism, and so
$[H_{\M}^{d-p}(M)]_{{\bf i}} \not = 0$. Therefore, $a_{d-p}^1(M)
\geq i^1 =M_p^1 -b^1 = t_p^1(M) +a^1(S)$. This proves $(ii)$ for
the case $j=1$, $\sigma= Id$.

Let $p$ be the greatest integer such that ${\bf M}_{p}=
\max_{\leq_{lex}} \{ {\bf M}_{0},\dots, {\bf M}_{l} \}$. Then,
${\bf M}_{p+1}<_{lex} {\bf M}_p$, so $a_{d-p}^1(M)= t_p^1(M)
+a^1(S)$ by $(ii)$. Therefore, $a_*^1(M)= t_*^1(M) +a^1(S)$ and we
have $(iii)$ for $j=1$. The proof of the statement for $j=
2,...,r$ follows from Remark \ref{B6}. $\B$

\bigskip
\section{Scheme associated to a multigraded ring}
\label{U} \markboth{CHAPTER I. MULTIGRADED RINGS} {SCHEME
ASSOCIATED TO A MULTIGRADED RING}

\medskip
Let $S$ be a noetherian $\Bbb N^r$-graded ring. We call $S$
standard if $S$ may be generated over $S_0$ by elements in degrees
$(1, 0,\dots,0),\dots, (0,\dots,0,1)$. Similarly to the graded
case, we may associate to such a ring a multigraded scheme in a
natural way (see \cite{H}). Our purpose is to extend this
construction to a more general class of multigraded rings, which
will recover the standard case as well as the Rees algebra of any
homogeneous ideal in a graded $k$-algebra.

\medskip
     Let $S$ be a noetherian $\Bbb N^r$-graded ring finitely generated
over $S_0$ by homogeneous elements $x_{11}, \dots, x_{1k_1},
\dots, x_{r1}, \dots, x_{rk_r}$ of degrees $\dg (x_{ij}) =
(d_{ij}^1, \dots, d_{ij}^{i-1},1,0, \dots, 0)$, where $d_{ij}^l $
are non-negative integers, and set $d_i^l = \max_j \{ d_{ij}^l
\}$. This class of rings includes for instance any standard $\Bbb
N^r$-graded ring by taking $d_{ij}^l =0$. For every $i = 1, \dots,
r$, let $I_i$ be the ideal of $S$ generated by the homogeneous
components of $S$ of degree ${\bf n}=(n_1, \dots, n_r)$ such that
$n_i >0, n_{i+1}=\dots =n_r=0$. Then we define the irrelevant
ideal of $S$ as $S_+ = I_1 \cdots I_r$. We are going to associate
a scheme to $S$ in the following way. A homogeneous prime ideal
$P$ of $S$ is said to be relevant if $P$ does not contain $S_+$.
Then we define the set $\Proj^r(S)$ to be the set of all relevant
homogeneous prime ideals $P$. It is easy to check that $\dim S/P
\geq r$ for any relevant prime ideal (see the proof of Lemma
\ref{U1}). Following \cite{STV} (where the standard bigraded case
was studied), we define the relevant dimension of $S$ as
 $$\rdim S = \cases {r-1 & if $\Proj^r(S) = \emptyset$ \cr
\max \{\, \dim S/P \mid P \in \Proj^r (S) \} & if $\Proj^r(S) \not
= \emptyset$ \cr} .$$ If $I$ is a homogeneous ideal of $S$, we
define the subset $V_+(I):= \{ P \in \Proj^r(S) \mid I \subset P
\}$. We can define a topology on $\Proj^r (S)$ by taking as closed
subsets the subsets of the form $V_+(I)$. Next, to define a sheaf
 of rings $\cal O$ in $\Proj^r(S)$, we first consider for each $P \in
\Proj^r (S)$ the homogeneous localization by $P$
$$ S_{(P)} = \{\, \frac{a}{s} \mid s \not \in P, \,\,   a,s \in S_{\bf
n}, \,\, {\bf n} \in {\Bbb Z^r} \, \}.  $$ For any open subset $U
\subset \Proj^r (S)$, we define ${\cal O} (U)$ to be the set of
functions $t: U \to \bigsqcup_{P \in U} S_{(P)}$ such that for
each $P \in U$, $t(P) \in S_{(P)}$ and $t$ is locally a quotient
of elements of $S$. Then, $\cal O$ is a sheaf of rings. We call
$\Proj^r(S)$ the $r$-projective scheme associated to $S$. Defining
for any homogeneous $f \in S_+$ the set $D_+(f) = \{P \in
\Proj^r(S) \mid f \not \in P \}$ we have an open cover of
$\Proj^r(S)$, and for each such open set we have an isomorphism of
locally ringed spaces $$(D_+(f), {\cal O} \vert D_+(f)) \cong \Sp
(S_{(f)}).$$ Moreover, ${\cal O}_P \cong S_{(P)}$ for any relevant
prime ideal $P$, hence $\Proj^r(S)$ is a scheme in a natural way.
This construction extends the usual one given in the standard case
(see \cite{H}).

\medskip
The next lemma computes the dimension of $\Proj^r(S)$. Its proof
follows the same arguments as in \cite[Lemma 1.2]{H}, but we
include it for completeness.

\begin{lem}
\label{U1} $\dim \Proj^r(S) = \rdim S -r$.
\end{lem}

{\pf} We may assume that $\Proj^r (S) \not = \emptyset$ (otherwise
the result is trivial). Let $P \in \Proj^r(S)$ be a closed point.
Since the projection $\Proj^r(S) \to \Spec(S_0)$ is proper, we
have that $P_0 = P \cap S_0$ is a closed point of $\Spec (S_0)$,
so $(S/P)_0= S_0/P_0$ is a field. Let us denote by $T= S/P$, and
note that $\dim \Proj^r(T)= 0$. For $j = 1, \dots, r$, let $J_j$
be the ideal of $T$ generated by the homogeneous components of $T$
of degree $\bf n$ such that $n_j >0, n_{j+1}=\dots =n_r=0$. We
have a maximal chain of homogeneous prime ideals
$$0 \subset J_r \subset J_{r-1}+J_r \subset \dots \subset J_1+ \dots
+J_r,$$ so $\dim T=r$ because $T$ is a catenary ring. On the other
hand, for a given minimal prime $Q_0 \in \Proj^r(S)$, we have a
chain of homogeneous prime ideals of type $Q_0 \subset \dots
\subset Q_s \subset \dots \subset Q_{s+r}$, with $Q_s$ a closed
point of $\Proj^r(S)$. Therefore,

\vspace{2mm}

\hspace{10mm} $\dim \Proj^r(S) = \sup \, \{ \h Q : Q \in
\Proj^r(S) \}$

\vspace{2mm}

\hspace{34mm} $=\sup \, \{ \dim S/Q : Q \in \Proj^r(S) \} -r$

\vspace{2mm}

\hspace{34mm} $= \rdim S-r. \B$

\medskip
Next we are going to define the diagonal functor. Given $e_1$,...,
$e_r$ positive integers, the set
$$\Delta:= \{ (e_1 s,..., e_r s) \mid s \in \Bbb Z \}$$ is
called the $(e_1,..., e_r)$-diagonal of  ${\Bbb Z^r}$. We may then
define the diagonal of $S$ along $\Delta$ as the graded ring
 $S_{\Delta} := \bigoplus_{s \in \Bbb Z}
S_{(e_1 s,..., e_r s)} $. Similarly, given an $r$-graded
$S$-module $M$ we define the diagonal of $M$ along $\Delta$ as the
graded $S_{\Delta}$-module $M_{\Delta}:= \bigoplus_{s \in \Bbb Z}
M_{(e_1 s,..., e_r s)}$. Then we have an exact functor
$$( \,\,)_{\Delta}: M^r(S) \rightarrow  M^1(S_{\Delta}),$$
called diagonal functor.

\medskip
Let us denote by $X = \Proj^r (S)$, and for each $\Delta$, let
$X_\Delta = \Proj (S_\Delta)$. By considering diagonals $\Delta =
(e_1, \dots,e_r)$ such that $e_r>0$, $e_{r-1}> d_r^{r-1} e_{r}$,
$\dots$, $e_1 > d_2^1 e_2 + \dots+ d_r^1 e_r$, then the sheaf of
ideals ${\cal L}= ( S_{(e_1,\dots, e_r)} ) \; {\cal O}_X$ defines
an isomorphism $X \stackrel{\cong} \longrightarrow X_\Delta$. In
particular, this isomorphism allows us to compute the dimension of
$S_\Delta$, extending \cite[Proposition 2.3]{STV} where this
dimension was computed for bigraded standard $k$-algebras.

\begin{lem}
\label{U2} Assume that $S_0$ is an artinian local ring. Then
 $\dim S_{\Delta} = {\rdim} \,S -r +1$, for any
$\Delta = (e_1, \dots,e_r)$ with $e_r>0$, $e_{r-1}> d_r^{r-1}
e_{r}, \dots, e_1 > d_2^1 e_2 + \dots+ d_r^1 e_r$.
\end{lem}

{\pf} From the isomorphism $X \cong X_\Delta$, we have that $\rdim
S_\Delta = \rdim S-r+1$ by Lemma \ref{U1}. Moreover, since $S_0$
is artinian, any minimal prime ideal of $S_\Delta$ is relevant,
and so $\rdim S_\Delta =\dim S_\Delta$. $\B$

\medskip
Classically,  $S$ is the multihomogeneous coordinate ring of a
multiprojective variety $V$ contained in some multiprojective
space ${\Bbb P}^{n_1}_k \times \dots \times {\Bbb P}^{n_r}_k$. By
taking the $(1, \dots, 1)$-diagonal, $S_{\Delta}$ is then the
homogeneous coordinate ring of the image of $V$ via the Segre
embedding ${\Bbb P}^{n_1}_k \times \dots
 \times {\Bbb P}^{n_r}_k \rightarrow  \Bbb P^N_k$, where $N = (n_1+1)
\dots (n_r+1) -1$.

\bigskip
\section{Hilbert polynomial of multigraded modules}
\label{F} \markboth{CHAPTER I. MULTIGRADED RINGS} {HILBERT
POLYNOMIAL OF MULTIGRADED MODULES}

\medskip
Let $S= \bigoplus_{\bn \in \Bbb N^r} S_{\bn}$ be an $r$-graded
ring defined over an artinian local ring $S_0=A$. If $S$ is
standard, then we have that the Hilbert function of any finitely
generated $r$-graded $S$-module $L$, $H(L, \bn)=
\length_A(L_{\bn})$, is a polynomial function; that is, there
exists a polynomial
 $P_L(t_1, \dots, t_r) \in \Bbb Q [t_1, \dots,t_r]$, called Hilbert
polynomial of $L$, such that for any ${\bf n} \gg 0$, $P_L(n_1,
\dots, n_r) = {\rm length}_A (L_{\bn})$ (see \cite{HHRT},
\cite{KMV}). In this section we are going to extend the existence
of such a polynomial for the larger class of finitely generated
$r$-graded modules defined over the multigraded rings introduced
in Section 1.4. Furthermore, we will state a formula for the
difference between the Hilbert polynomial and the Hilbert function
of any finitely generated $r$-graded module analogous to the one
known in the graded case.

\medskip
Let $S$ be a noetherian $\Bbb N^r$-graded ring generated over
$S_0=A$ by homogeneous elements $x_{11}, \dots, x_{1k_1}, \dots,
x_{r1}, \dots, x_{rk_r}$ in degrees $\dg (x_{ij}) = (d_{ij}^1,
\dots, d_{ij}^{i-1},1,0, \dots, 0)$, where $d_{ij}^l \geq 0$. Set
$d^l_i = \max _j \{ d_{ij}^l \}$.

\medskip
Given a finitely generated $r$-graded $S$-module $L$, let us
define its homogeneous support as ${\rm Supp}_+(L) = \{\, P \in
\Proj^r (S) \mid L_{P} \not = 0 \}$. Note that ${\rm Supp}_+(L) =
V_+({\rm Ann} \, L)$ is a closed subset of $\Proj^r(S)$. We define
the relevant dimension of $L$ as
 $$\rdim L = \cases {r-1 & if $\Supp_+(L) = \emptyset$ \cr
\max \{\, \dim S/P \mid P \in \Supp_+(L) \} & if $\Supp_+(L) \not
= \emptyset$ \cr} .$$ One can check that $\rdim L = \dim \Supp_+ L
+ r $.

\medskip
From now on in this section we will assume that $A$ is an artinian
local ring. Given a finitely generated $r$-graded $S$-module $L$,
its homogeneous components $L_{\bn}$ are finitely generated
$A$-modules, and hence have finite length. The numerical function
$H(L, \; \cdot \;) : \Bbb Z^r \rightarrow \Bbb Z$ with $H(L, \bn)
= \lL_{A} ( L_{\bn} )$ is the Hilbert function of $L$. Next result
shows the existence of the Hilbert polynomial for any finitely
generated $r$-graded $S$-module.

\begin{prop}
\label{F0}
   Let $L$ be a finitely generated $r$-graded $S$-module of
relevant dimension
    $\delta$. Then there exists a polynomial  $P_L(t_1, \dots, t_r) \in
\Bbb Q[t_1, \dots, t_r]$ of total degree $\delta -r $ such that
$H(L, i_1, \dots, i_r)= P_L(i_1, \dots, i_r)$ for
 $i_1 \gg d_2^1 i_2 + \dots + d_r^1 i_r$, \dots, $i_{r-1} \gg d_r^{r-1}
i_r$, $i_r \gg 0$. Moreover,
$$P_L (t_1, \dots, t_r) = \sum_{\vert \bn \vert \leq \delta -r}
a_{\bn} {t_1-d_2^1 t_2 - \dots -d_r^1 t_r \choose n_1} \dots
{t_{r-1}-d_r^{r-1} t_r \choose n_{r-1}} {t_{r} \choose n_{r}}
\,,$$ where $a_{\bn} \in \Bbb Z$, $a_{\bn} \geq 0$ if $\vert \bn
\vert = \delta-r$.
 \end{prop}

{\pf} Given a finitely generated $r$-graded $S$-module $L$, first
note that there is a chain
$$ 0= L_0 \subset L_1 \subset \dots \subset L_s=L$$
of $r$-graded submodules of $L$ such that for each $i \geq 1$,
$L_i/L_{i-1} \cong (S/P_i) ({\bm}_i)$, where $P_i \in \Supp L$ is
a homogeneous prime ideal and ${\bm}_i \in \Bbb Z^r$. Indeed, we
may assume $L \not =0$. Choose $P_1 \in \Ass L$. Then $P_1$ is a
homogeneous prime ideal, and there exists an $r$-graded submodule
$L_1 \subset L$ such that $L_1 \cong (S/P_1) ({\bm}_1)$. If $L_1
\not = L$, by repeating the procedure with $L/L_1$ we get an
$r$-graded submodule $L_2 \subset L$ such that $L_2/ L_1 \cong
(S/P_2) ({\bm}_2)$. Since $L$ is noetherian, this process finishes
after a finite number of steps. From this chain, we obtain
$$H(L, \bn) = \sum_{i=1}^s H(S/P_i, \bn + \bm _i).$$
So it is enough to prove the result for the rings $T=S/P$, with
$P$ a homogeneous prime ideal. To this end, we will reduce the
problem to the standard case where the result is already known.

  Set $B= T_0$. Let us consider $\overline T \subset T$ the
$B$-algebra generated by the homogeneous elements of $T$ of degree
$(\g_1, \dots, \g_r)$ such that
$$
 \begin{array}{llll}
\g_{r-1} \geq d_r^{r-1} \g_{r} \\
\g_{r-2} \geq d_{r-1}^{r-2} \g_{r-1}+ d_r^{r-2} \g_{r} \\
\dots \dots \\
\g_1 \geq d_2^1 \g_2 + \dots + d_r^1 \g_r .
\end{array}
$$

\noindent Then one has $\overline T_{\bn} = T_{\bn}$ for each $\bn
\in \Bbb N^r$ satisfying the inequalities before. Let us consider
the morphism
    $$  \begin{array}{lll}
\psi : \hspace{5mm} \Bbb Z^r & \longrightarrow & \hspace{5mm} \Bbb
Z^r
\\
 (x_1, \dots, x_r)   & \mapsto &
(x_1-d_2^1 x_2- \cdots- d_r^1 x_r, \dots, x_{r-1}-d_r^{r-1} x_r,
x_r)
\\
\end{array} $$
Note that $\psi (\supp \overline T) \subset \Bbb N^r$, so
$\overline T^{\psi}$ is again a $\Bbb N^r$-graded ring.
Furthermore, we have $\rdim \overline T^\psi = \rdim \overline T =
\rdim T = \delta$. If $\overline T^{\psi}$ is standard, by
\cite[Theorem 4.1]{HHRT} there exists a polynomial $Q(t_1, \dots,
t_r) \in \Bbb Q[t_1, \dots, t_r]$ of total degree $\delta -r$
$$Q (t_1, \dots, t_r) = \sum_{\vert \bn \vert \leq \delta -r}
a_{\bn} {t_1 \choose n_1} \dots {t_{r} \choose n_{r}} \,,$$ with
$a_{\bn} \in \Bbb Z$, $a_{\bn} \geq 0$ if $\vert \bn \vert =
\delta -r$ such that for ${\bf i} \gg 0$
$$Q (i_1, \dots, i_r)= \length_B
 [\,\overline T^\psi \,]_{(i_1, \dots, i_r)}.$$
Then, by defining $P (t_1, \dots, t_r)= Q(t_1-d_2^1 t_2 - \dots
d_r^1 t_r, \dots, t_{r-1}-d_r^{r-1} t_r, t_r) $, let us observe
that for $i_1 \gg d_2^1 i_2 + \dots + d_r^1 i_r, \dots, i_{r-1}
\gg d_r^{r-1} i_{r}, i_r \gg 0$, we have

\vspace{3mm}

\hspace{10 mm} $ P(i_1, \dots, i_r) =
 \length_B \; [\overline T^\psi]_{(i_1- d_2^1 i_2 - \dots
-d_r^1 i_r, \dots, i_{r-1}- d_r^{r-1} i_r, i_r)} $

\vspace{3mm}

\hspace{32mm} $= \length_B \; [\overline T]_{(i_1, \dots, i_r)} $

\vspace{3mm}

\hspace{32mm} $= \length_A \; [T]_{(i_1, \dots, i_r)}, $

\vspace{3mm} \noindent so we get the statement.

     Therefore we only have to prove that $\overline T^{\psi}$ is
standard or, equivalently, that $\overline T$ can be generated
over $B$ by homogeneous elements in degrees $(\g_1, \dots, \g_r)$
such that $\g_{i+1}= \dots =\g_r = 0$, $\g_i =1$, $\g_{i-1} =
d_i^{i-1} \g_{i}$, $\g_{i-2} = d_{i-1}^{i-2} \g_{i-1} + d_i^{i-2}
\g_i$,  $\dots$, $\g_1 = d_2^1 \g_2 + \dots + d_i^1 \g_i$. Assume
that $T$ is generated over $B$ by homogeneous elements $z_{11},
\dots, z_{1k_1}, \dots, z_{r1}, \dots, z_{rk_r}$ in degrees $\dg
(z_{ij}) = (d_{ij}^1, \dots, d_{ij}^{i-1},1,0, \dots, 0)$. Let us
take a homogeneous element $z$ in $\overline T$, with $\dg z=
(\A_1, \dots, \A_r)$. Let $j$ be such that $\A_j \not =0$,
$\A_{j+1}= \dots = \A_r=0$ ($j$ is 0 if $z \in B$). We are going
to prove by induction on $j$ that $z$ can be generated over $B$ by
the homogeneous elements whose degrees satisfy the equalities
before. If $j =0$, there is nothing to prove. If $j=1$, then $\dg
\, z = (\A_1, 0, \dots, 0)$ and we can write $z$ as a linear
combination with coefficients in $B$ of products of $\A_1$
elements among $z_{11}, \dots, z_{1k_1}$, so the result is
trivial. Assume now that $j >1$. By forgetting the first component
of the degree, we have by induction hypothesis that $z$ can be
written as a sum of terms of the type $\lambda w_{1} \dots w_l$
with $\lambda \in B[z_{11}, \dots, z_{1k_1}]$, and the degree of
the elements $w_i$ satisfying the $r-1$ first equalities. Set $\dg
w_j = ( s^1_j, \dots, s^r_j)$,  $\dg \lambda= (s, 0, \dots, 0)$.
We will finish if we prove that
$$\alpha_1 \geq \sum_{j=1}^l (d_2^1 s^2_j + \dots + d_r^1
s^r_j).$$ But note that $\A_1 \geq d_2^1 \A_2 + \dots + d_r^1 \A_r
 =  \sum_{j=1}^l d_2^1 s^2_j + \dots + d_r^1 s^r_j.$
$\B$

\medskip
Our next aim will be to study for a given finitely generated
$r$-graded $S$-module $L$, the $A$-modules $H^i_{S_+}(L)_{\bn}$
for $i \geq 0$, $\bn \in \Bbb Z^r$. We need two previous lemmas.

\begin{lem}
\label{F12}
   Let $L$ be a finitely generated $r$-graded $S$-module such that
$(S_+)^u L = 0$ for an integer $u$. Then there exits $\bm =(m_1,
\dots, m_r) \in \Bbb Z^r$ such that
$$L_{\bn} =0 ,$$
for $\bn =(n_1, \dots, n_r)$ such that $n_1 > d_2^1 n_2 + \dots +
d_r^1 n_r + m_1$, \dots,
 $n_{r-1} > d_r^{r-1} n_r+ m_{r-1}$, $n_r > m_r$.
\end{lem}

{\pf} Since there exists $u \in \Bbb Z$ such that $(S_+)^u L = 0$,
we have $\Supp_+(L) =V_+(\Ann L) \subset V_+(S_+) = \emptyset$, so
$\rdim L= r-1$. Then the result follows from  Proposition
\ref{F0}. $\Box$

\begin{lem}
\label{F1} (Homogeneous Prime Avoidance) Let $P_1, \dots, P_m \in
\Proj ^r(S)$. If $I$ is any homogeneous ideal of $S$ such that $I
\not \subset P_i$ for $i = 1,\dots, m$, then there is a
homogeneous element $a$ such that $a \in I$, $a \not \in P_1 \cup
\dots \cup P_m$.
 \end{lem}

{\pf}
 We may assume that $P_j \not \subset P_i$ for $i \not = j$, so for a
given $i$, we have that for any $j \not = i$ there exists a
homogeneous element
 $p_{ij} \in P_j$, $p_{ij} \not \in P_i$. Then
  $p_i = \prod_{j \not = i} p_{ij}$ satisfies that
$p_i \not \in P_i$, but
 $p_i \in P_j$ for all $j \not = i$.
  Next we may take homogeneous elements $a_i \in I $, $a_i \not \in
P_i$ for $i=1, \dots, m$. Set
 $\deg a_i p_i =(\alpha_{i1}, \dots, \alpha_{ir})$.
Since $S_+ \not \subset P_i$, there exists an element of the type
$x_{1 j_1} \ldots x_{r j_r} \not \in P_i$. So multiplying each
$a_i p_i$ by a power of the corresponding $x_{r j_r}$ we can
assume that $\alpha_{1r} = \dots = \alpha_{mr} = \alpha_r$. Then,
multiplying by suitable powers of each $x_{r-1,j_{r-1}}$ we may
also assume that $\alpha_{1,r-1} = \dots = \alpha_{m,r-1} =
\alpha_{r-1}$. By repeating this procedure as many times as
necessary, we can assume at the end that $\deg(a_1 p_1)= \dots=
\deg(a_m p_m ) = (\alpha_{1}, \dots, \alpha_{r})$. Now
      $a= a_1 p_1 + \dots +a_m p_m$
is homogeneous  and $a \in I$, $a \not \in P_1 \cup \dots \cup
P_m$. $\B$

\medskip
Now we are ready to prove that if $L$ is a finitely generated
$r$-graded $S$-module, then the $A$-modules $H^i_{S_+}(L)_{\bn}$
are finitely generated for all $\bn \in \Bbb Z^r$, $i \geq 0$, and
vanish for all sufficiently large $\bn$. Here, the artinian
assumption on $A$ is not necessary. In the graded case, this is a
classical result due to J.P. Serre.

\begin{prop}
\label{F11}
   Let $L$ be a finitely generated $r$-graded $S$-module. Then
\begin{enumerate}
\item For all $i \geq 0$, $\bn \in \Bbb Z^r$, the $A$-module
$H^i_{S_+}(L)_{\bn}$ is finitely generated.
\item There exits $\bm =(m_1, \dots, m_r) \in \Bbb Z^r$ such that
$H^i_{S_+}(L)_{\bn} =0 $ for all $i \geq 0$, $\bn= (n_1, \dots,
n_r)$ such that $n_1 > d_2^1 n_2 + \dots + d_r^1 n_r + m_1$,
\dots,
 $n_{r-1} > d_r^{r-1} n_r+ m_{r-1}$, $n_r > m_r$.
\end{enumerate}
\end{prop}

{\pf} We will follow the same lines as the proof of the graded
version in \cite[Proposition 15.1.5]{BS3}. We will prove by
induction on $i$ that $H^i_{S_+}(L)_{\bn}$ is a finitely generated
$A$-module for all $\bn \in \Bbb Z^r$, and that it is zero for all
sufficiently large values of $\bn$. This proves the statement
because $H^i_{S_+}(L) =0$ for all $i$ greater than the minimal
number of generators of $S_+$.

Assume $i=0$. Then $H^0_{S_+}(L)$ is a finitely generated
$r$-graded $S$-module since it is a submodule of $L$, and so
$H^0_{S_+}(L)_\bn $ is a finitely generated $A$-module and there
exists $u \in \Bbb N$ such that $(S_+)^u H^0_{S_+}(L) =0 $. Then,
according to Lemma \ref{F12} there exists $\bm \in \Bbb Z^r$ such
that $H^0_{S_+}(L)_{\bn} = 0$ for $n_1>d^1_2 n_2 + \dots +d^1_r
n_r + m_1$, $\dots$, $n_{r-1} > d_r^{r-1}n_r+ m_{r-1}$, $n_r >
m_r$.

Now let us assume $i >0$. From the $r$-graded isomorphism
$H^i_{S_+}(L) \cong H^i_{S_+} (L/ H^0_{S_+}(L))$, we may replace
$L$ by $L/ H^0_{S_+}(L)$ and then assume that $H^0_{S_+}(L)=0$.
Then $S_+ \not \subset P$ for all $P \in {\rm Ass}(L)$, and so by
the Prime Avoidance Lemma there exists a homogeneous element $x
\in S_+$ of degree ${\bf k}=(k_1, \dots, k_r)$ such that $x \not
\in P$ for all $P \in {\rm Ass}(L)$. Looking at the proof of the
Prime Avoidance Lemma, notice that we can choose $x$ such that
$\bf k$ satisfies $k_1 >d_2^1 k_2 + \cdots+ d_r^1 k_r$, $\dots$,
$k_{r-1}> d_r^{r-1} k_r$. Then we get an $r$-graded exact sequence
$$0 \to L(-{\bf k})
\stackrel{\cdot  x} \longrightarrow L \to L/xL \to 0,$$ which
induces for all $\bn \in \Bbb Z^r$ the exact sequence of
$A$-modules
$$ H^{i-1}_{S_+}(L/ xL)_\bn \to H^i_{S_+}(L)_{\bn -{\bf k}}
\stackrel{\cdot x} \longrightarrow H^{i}_{S_+}(L)_\bn \;\;.$$ By
the induction hypothesis, there exists ${\bf {\overline m}} \in
\Bbb Z^r$ such that $H^{i-1}_{S_+}(L/x L)_{\bn} =0 $ for all $\bn
=(n_1, \dots, n_r)$ such that $n_1 > d_2^1 n_2 + \dots + d_r^1 n_r
+ \overline m_1$, \dots,
 $n_{r-1} > d_r^{r-1} n_r+ \overline m_{r-1}$, $n_r > \overline m_r$.
Now let $\bn$ verifying these inequalities. Then note that for any
$s \geq 1$, ${\bf n}+s{\bf k}$ also satifies them, and so we have
exact sequences
$$ 0 \to H^i_{S_+}(L)_{\bn -{\bf k}}
\stackrel{\cdot x^s} \longrightarrow
H^{i}_{S_+}(L)_{{\bn}+(s-1){\bf k}} \;\;.$$ Since $H^i_{S_+}(L)$
is $S_+$-torsion and $x \in S_+$, we have $ H^i_{S_+}(L)_{\bn
-{\bf k}} = 0$. Therefore, by taking $\bm = {\bf {\overline m}} -
{\bf k}$, we obtain $H^i_{S_+}(L)_{\bn} =0 $ for all $\bn$ such
that $n_1 > d_2^1 n_2 + \dots + d_r^1 n_r + m_1$, \dots, $n_{r-1}
> d_r^{r-1} n_r+ m_{r-1}$, $n_r > m_r$.

Now let us fix $\bn \in \Bbb Z^r$. If $n_1 > d_2^1 n_2 + \dots +
d_r^1 n_r + m_1$, \dots, $n_{r-1} > d_r^{r-1} n_r+ m_{r-1}$, $n_r
> m_r$, we have that $H^i_{S_+}(L)_{\bn} =0 $, and so it is a
finitely generated $A$-module. Otherwise, let us take $y \in S_+$
such that $y \not \in \bigcup_{P \in {\rm Ass}(L)} P$ with degree
${\bf l}=(l_1, \dots, l_r)$ such that $\bn + \bf l$ satisfies the
previous inequalities (we can find such a $y$ by the Prime
Avoidance Lemma). Then we have the graded exact sequence
$$ H^{i-1}_{S_+}(L/ yL)_{\bn+ {\bf l}} \to  H^i_{S_+}(L)_{\bn}
\stackrel{ \cdot y} \longrightarrow H^{i}_{S_+}(L)_{\bn+ {\bf l}}
=0, $$ and from the induction hypothesis we also have that
$H^{i-1}_{S_+}(L/y L)_{\bn+{\bf l}} $ is a finitely generated
$A$-module, and so $H^i_{S_+}(L)_{\bn}$. $\B$

\medskip
We have already shown that the Hilbert function of any finitely
generated $r$-graded $S$-module is a polynomial function for large
$\bn$. Our next result precises the difference between the Hilbert
function and the Hilbert polynomial for any $\bn$.

\begin{prop}
\label{F2}
 Let $L$ be a finitely generated $r$-graded $S$-module.
Then for all $\bn \in \Bbb Z^r$
 $$H(L,\bn)- P_L(\bn) = \sum_q (-1)^q \, \lL_A( H^q_{S_+} (L)_{\bn}).$$
\end{prop}

{\pf} We will follow the proof of the graded version from
\cite[Theorem 4.3.5]{BH}. For an arbitrary finitely generated
$r$-graded $S$-module $L$, let us define the series

\vspace{3mm}

\hspace{15mm} $H'_L (u_1, \dots, u_r) = \sum_{\bn \in  \Bbb Z^r} (
H(L, \bn) - P_L (\bn)) u^{\bn}$

\vspace{3mm}

\hspace{15mm} $H''_L (u_1, \dots, u_r) = \sum_{\bn \in  \Bbb Z^r}
(\, \sum_q (-1)^q \, \lL_A( H^q_{S_+} (L)_{\bn}) \,) u^{\bn}. $

\vspace{3mm}

We will prove the statement by induction on $\delta = \rdim L$. If
$\delta =r-1$, then ${\rm Supp}_+ L = \emptyset$, and so there
exists $m$ such that $S_+ ^m \subset {\rm Ann} (L)$. Therefore
$H^0_{S_+}(L)= L$, and hence the result is trivial. Assume now
$\delta \geq r$, and let us consider $\overline L = L/
H^0_{S_+}(L)$. Since $H^0_{S_+}(L)$ is a finitely generated
$r$-graded $S$-module which is vanished by some power of $S_+$,
there are integers $i_1, \dots, i_r$ such that
 $H^0_{S_+}(L)_{\bn} = 0$ for $n_1>d^1_2 n_2 + \dots +d^1_r n_r
+ i_1$, $\dots$, $n_{r-1} > d_r^{r-1}n_r+i_{r-1}$, $n_r > i_r$. So
we have $P_L({\bf t}) = P_{\,\overline L}({\bf t})$, and it is
enough to prove the result for $\overline L$ because then, for all
$\bn = (n_1, \dots, n_r)$

\vspace{3mm}

\hspace{6mm} $H(L, \bn)- P_L( \bn) = H(\overline L, \bn)+ \lL_A
(H^0_{S_+}(L)_{\bn})- P_{\overline L}( \bn) $

\vspace{3mm}

\hspace{35mm} $= \sum_q (-1)^q \lL_A ( H^q_{S_+} (\,\overline
L\,)_{\bn}) + \lL_A \,(H^0_{S_+}(L)_{\bn}) $

\vspace{3mm}

\hspace{35mm} $= \sum_q (-1)^q \lL_A \,( H^q_{S_+} (L)_{\bn}).$

\vspace{3mm} \noindent So let us assume $H^0_{S_+}(L)=0$. Then
$S_+ \not \subset P$ for all $P \in {\rm Ass}(L)$, and so by the
Prime Avoidance Lemma there exists a homogeneous element $x \in
S_+$ of degree ${\bf k} =(k_1,\dots, k_r)$ such that $x \not \in
P$ for all $P \in {\rm Ass} (L)$. Then we have the $r$-graded
exact sequence
 $$0 \to L(-{\bf k}) \to L \to L/xL \to 0,$$
with $\rdim L/xL < \rdim L$. Note that $H(L/x L, \bn) = H(L, \bn)
- H(L, \bn - {\bf k})$ for all $\bn$, and so $P_{L/x L}({\bf t}) =
P_L ({\bf t}) - P_L( {\bf t} - {\bf k})$. We conclude $H'_{L/ x L}
( {\bf u}) = (1- {\bf u}^{\bf k}) H'_L ({\bf u})$. From the long
exact sequence of local cohomology, we also get
 $H''_{L/ x L} ( {\bf u}) = (1- {\bf u}^{\bf k}) H''_L ({\bf u})$.
By the induction hypothesis, we have
 $H'_{L/ x L} ( {\bf u}) = H''_{L/ x L} ( {\bf u})$ and so
 $H'_{L} ( {\bf u}) = H''_{L} ( {\bf u})$. $\B$

\newpage
\newpage
\chapter*{$\;\;\;$}

\newpage

\medskip

\newpage

\chapter{The diagonals of a bigraded module}
\typeout{The diagonals of a bigraded module}


\markboth{CHAPTER II. THE DIAGONALS OF A BIGRADED MODULE} {THE
DIAGONALS OF A BIGRADED MODULE}

\medskip
 Throughout this chapter we will study in more detail the diagonal
functor in the category of bigraded $S$-modules, where
$S=k[{\x},{\y}]$ is the polynomial ring in $n+r$ variables with
the bigrading given by deg$(X_i)= (1, 0)$, deg$(Y_j)= (d_j,1)$,
and $d_1, \dots, d_r \geq 0$. This category includes any standard
bigraded $k$-algebra, by taking $d_1= \dots=d_r=0$, as well as the
Rees ring and the form ring of a homogeneous ideal in a graded
$k$-algebra, when those rings are endowed with an appropiate
bigrading (see Section 2.3).

\medskip
For a given $c, e$ positive integers, let $\Delta$ be the
$(c,e)$-diagonal of  ${\Bbb Z}^2$. Our purpose is to study the
exact functor $(\; \;)_\Delta : M^2(S) \rightarrow M^1(S_\Delta)$
(see Chapter 1, Section 4). We are mainly interested in studying
how the arithmetic properties of a bigraded $S$-module $L$ and its
diagonals $L_\Delta$ are related. Most of these properties, like
the Cohen-Macaulayness or the Gorenstein property, can be
characterized by means of the local cohomology modules. So it
would be very useful to relate the local cohomology modules of $L$
with the local cohomology modules of its diagonals. This has been
done by A. Conca et al. in \cite{CHTV} from the study of the
bigraded minimal free resolution of $L$ over $S$, after developing
a theory of generalized Segre products of bigraded algebras. In
Section 2.1 we are going to present their results by a different
and somewhat easier approach. In addition, this approach will
provide more detailed information about several problems
concerning to the behaviour of the local cohomology when taking
diagonals.

\medskip
In Section 2.2 we focus our study on standard bigraded
$k$-algebras. For a such $k$-algebra $R$, let ${\cal R}_1 =
\bigoplus_{i \in \Bbb N} R_{(i,0)}$, ${\cal R}_2 =  \bigoplus_{j
\in \Bbb N} R_{(0,j)}$. In this case, we give a characterization
for $R$ to have a good resolution in terms of the $a_*$-invariants
of ${\cal R}_1$ and ${\cal R}_2$ which, in particular, provides a
criterion for the Cohen-Macaulayness of its diagonals. We also
find necessary and sufficient conditions on the local cohomology
of ${\cal R}_1$ and ${\cal R}_2$ for the existence of
Cohen-Macaulay diagonals of $R$, whenever $R$ is Cohen-Macaulay.

\medskip
Given a homogeneous ideal $I$ in a graded $k$-algebra $A$, the
Rees algebra $R_A(I)= \bigoplus_{n \geq 0} I^n$ of $I$ can be
endowed with the bigrading $R_A(I)_{(i,j)} =(I^j)_i$. The last
section of the chapter is devoted to study the diagonals of the
Rees algebra. In the case where $A$ is the polynomial ring, we
will show that if the Rees algebra is Cohen-Macaulay then there
exists some diagonal with this property, thus proving a conjecture
stated in \cite{CHTV}. Furthermore, we will give necessary and
sufficient conditions on the ring $A$ for the existence of a
Cohen-Macaulay diagonal of a Cohen-Macaulay Rees algebra.

\bigskip
\section{The diagonal functor on the category of bigraded modules}
\markboth{CHAPTER II. THE DIAGONALS OF A BIGRADED MODULE} {THE
DIAGONAL FUNCTOR} \label{W}

\medskip
Let $S=k[{\x},{\y}]$ be the polynomial ring in $n+r$ variables
over a field $k$ with the bigrading given by deg$(X_i)= (1, 0)$,
deg$(Y_j)= (d_j, 1)$, where $d_1, \dots, d_r \geq 0$. Set $d =
\max \{d_1, \dots, d_r \}$, $u = \sum_{j=1}^r d_j$. Let us denote
by $\cal M$ the homogeneous maximal ideal of $S$. Note that the
irrelevant ideal $S_+$ of $S$ is the ideal generated by the
products $X_i Y_j$, for $i=1, \dots, n$, $j=1, \dots,r$.

\medskip
Given $c, e$ positive integers, let $\Delta$ be the
$(c,e)$-diagonal of ${\Bbb Z}^2$. For any bigraded $S$-module $L$,
let us recall that the
 diagonal of $L$ along $\Delta$ is defined as
$L_\Delta := \bigoplus_{s \in \Bbb Z} L_{(cs,es)}$, which is a
graded module over the graded ring $S_\Delta := \bigoplus_{s \geq
0} S_{(cs,es)}$. Our first lemma computes the dimension of the
diagonals of a finitely generated bigraded $S$-module.

\begin{lem}
\label{W0} Let $L$ be a finitely generated bigraded $S$-module.
For $\Delta = (c,e)$ with $c \geq de+1$,
 $\dim L_{\Delta} = {\rdim} \,L -1$.
\end{lem}

{\pf} The proof follows the same lines as the one given for the
bigraded standard case by A. Simis et al. in \cite[Proposition
2.3]{STV}. Set $\delta = \rdim L$. According to Proposition
\ref{F0}, there is a polynomial $P(s,t) \in \Bbb Q[s,t]$ of total
degree $\delta -2$ of the type
$$ P(s,t) = \sum_{k+l \leq \delta -2} a_{kl} { s -dt \choose k}
{t \choose l} \, ,$$ with $a_{kl} \geq 0$ for any $k, l$ verifying
$k+l= \delta -2$ such that for $i \gg dj $, $j \gg 0$, $P(i,j)=
\dim_k L_{(i,j)}$. For any $c \geq de +1$, let us consider the
polynomial $Q(u)= P(cu, eu) \in \Bbb Q[u]$. Then $Q(u) = \dim_k
L_{(cu,eu)} = \dim_k (L_\Delta)_u$  for $u$ large enough and deg
$Q (u) = \delta -2$. Therefore $\dim L_\Delta = \delta -1$. $\B$

\medskip
From now on in the chapter we will always consider diagonals
$\Delta = (c, e)$ with $c \geq de+1$. The next two propositions
are inspired in some results and techniques used by E. Hyry in
\cite{H}. The first one shows how the local cohomology modules of
$L$ with respect to $S_+$ are related to the local cohomology
modules of $L_{\Delta}$ with respect to $\MD$.

\begin{prop}
\label{W1} Let $L$ be a finitely generated bigraded $S$-module.
Then there are graded isomorphisms
 $$H^q_{S_+}(L)_{\Delta} \cong H^q_{{\cal M}_{\Delta}}(L_{\Delta})
\;\;, \forall q.$$
\end{prop}

{\pf} Let $\cal N$ be the ideal of $S$ generated by ${\cal
M}_\Delta$. Observe that $\sqrt {S_+} = \sqrt{ \cal N}$, so we
immediately get a bigraded isomorphism $H^q_{S_+}(L) \cong
H^q_{{\cal N}}(L)$, $\forall q \geq 0$. Denoting by $g_1, \dots,
g_s$ a $k$-basis of $S_{(c,e)}$, we have that $\cal N$ can be
generated by $g_1, \dots, g_s$. So we may compute the local
cohomology modules of $L$ with respect to $\cal N$ from the
C\v{e}ch complex built up from these elements
$$ {\bf C^{\cdot}} : 0 \to C^0 \to C^1 \to \dots \to C^s \to 0,$$
$$C^t= \bigoplus_{1 \leq i_1 < i_2 < \dots < i_t \leq s}
L_{g_{i_1} g_{i_2} \dots g_{i_t}},$$ with the differentation $d^t
: C^t \to C^{t+1}$ defined on the component
$$L_{g_{i_1} g_{i_2} \dots g_{i_t}}
\longrightarrow L_{g_{j_1} g_{j_2} \dots g_{j_t} g_{j_{t+1}}  }$$
to be the homomorphism $(-1)^{m-1} {\rm nat} : L_{g_{i_1} g_{i_2}
\dots g_{i_t}} \to ( L_{g_{i_1} g_{i_2} \dots g_{i_t}})_{g_{j_m}}$
 if $\{ i_1, \dots, i_t \} = \{ j_1, \dots,
\widehat{j_m}, \dots, j_{t+1} \}$, and $0$ otherwise. Then
$H^q_{\cal N}(L) \cong H^q ({\bf C^{\cdot}})$. We can also consider the
C\v{e}ch complex associated to $L_\Delta$ built up from $g_1,
\dots, g_s$
$$ {\bf D^{\cdot}} : 0 \to D^0 \to D^1 \to \dots \to D^s \to 0,$$
$$D^t= \bigoplus_{1 \leq i_1 < i_2 < \dots < i_t \leq s}
(L_\Delta)_{g_{i_1} g_{i_2} \dots g_{i_t}},$$ with the
differentation $\delta^t : D^t \to D^{t+1}$ defined on the
component $$(L_\Delta)_{g_{i_1} g_{i_2} \dots g_{i_t}}
\longrightarrow (L_\Delta)_{g_{j_1} g_{j_2} \dots g_{j_t}
g_{j_{t+1}}}$$ to be the homomorphism $(-1)^{m-1} {\rm nat}:
(L_\Delta)_{g_{i_1} g_{i_2} \dots g_{i_t}} \to (
(L_\Delta)_{g_{i_1} g_{i_2} \dots g_{i_t}})_{g_{j_m}}$
 if $\{ i_1, \dots, i_t \} = \{ j_1, \dots,
\widehat{j_m}, \dots, j_{t+1} \}$, and $0$ otherwise. Then
$H^q_{{\cal M}_{\Delta}}(L_\Delta) \cong H^q ({\bf D^{\cdot}})$. Note
that $d^t \vert D^t = \delta ^t$, so we have $(\Ker d^t)_\Delta =
\Ker \delta^t$, $({\rm Im} \, d^t)_\Delta = {\rm Im} \, \delta^t$.
Therefore we may conclude $H^q_{\cal N}(L) _\Delta \cong
H^q_{{\cal M}_{\Delta}}(L_\Delta)$. $\B$

\medskip
Now, let $S_1, S_2$ be the bigraded subalgebras of $S$ defined by
$S_1= k[{\x}]$, $S_2= k[{\y}]$, and note that the ideals
 ${\fm}_1= (X_1, ..., X_n)$ and ${\fm}_2 = ({\y}) $ are the
homogeneous maximal ideals of $S_1$ and $S_2$ respectively. Then
let us define ${\M}_1$ to be the ideal of $S$ generated by
${\fm}_1$ and ${\M}_2$ to be the ideal of $S$ generated by
${\fm}_2$. Note that ${\M}_1 +{\M}_2 = \M$  and ${\M}_1 \cap
{\M}_2 = S_{+}$. Therefore we have

\begin{prop}
\label{W2} Let $L$ be a finitely generated bigraded $S$-module.
There is a natural graded exact sequence
$$ ... \to H^q_{\M} (L)_{\Delta}
\to H^q_{{\M}_1} ( L )_{\Delta}
 \oplus H^q_{{\M}_2} ( L )_{\Delta}
\to H^q_{\MD}(L_{\Delta}) \stackrel {\varphi^q_L} \longrightarrow
 H^{q+1}_{\cal M}(L)_{\Delta}  \to ...$$
\end{prop}

{\pf}
   We get the result by applying the diagonal functor to the
Mayer-Vietoris sequence associated to ${\M}_1$, ${\M}_2$ and by
then using Proposition \ref{W1}. $\B$

\medskip
As a first consequence we may recover the following result by A.
Conca et al. in \cite{CHTV}.

\begin{cor}
\label{W3} \cite[Theorem 3.6]{CHTV} Let $L$ be a finitely
generated bigraded $S$-module. For all $q \geq 0$, there exists a
canonical graded homomorphism
$$\varphi^q_L : H^q_{{\cal M}_{\Delta}} (L_{\Delta}) \rightarrow
H^{q+1}_{\cal M} (L)_{\Delta}, $$ which is an isomorphism for $q
>\max \{n, r\}$.
\end{cor}

{\pf} Since ${\cal M}_1$ is generated by $n$ elements, we have
that $H^q_{{\M}_1}( L)=0$ for any $q>n$. Similarly, $H^q_{{\M}_2}
( L )=0$ for any $q >r$. Now, the corollary follows from
Proposition \ref{W2}. $\B$

\medskip
Moreover, let us also notice that Proposition \ref{W2} precises
the obstructions for $\varphi^q_L$ to be isomorphism. Denote by
$[\varphi^q_L]_s : H^q_{{\cal M}_{\Delta}} (L_{\Delta})_s
\rightarrow H^{q+1}_{\cal M} (L)_{(cs,es)} $ the component of
degree $s$ of the map $\varphi^q_L$. Then we have

\begin{cor}
\label{W31} Let $L$ be a finitely generated bigraded $S$-module.
For a given $s \in \Bbb Z$, the following are equivalent
\begin{itemize}
\item[(i)]
$[\varphi^q_L]_s$ is an isomorphism, for all $q \geq 0$.
\item[(ii)]
$H^q_{{\M}_1} ( L )_{(cs,es)} = H^q_{{\M}_2} ( L )_{(cs,es)} = 0$,
for all $q \geq 0$.
\end{itemize}
In particular, $\varphi^q_L$ is an isomorphism for all $q \geq 0$
if and only if $H^q_{{\M}_1} ( L )_\Delta = H^q_{{\M}_2} ( L
)_\Delta = 0$ for all $q \geq 0$.
\end{cor}

\medskip
Therefore, the obstructions for the maps $\varphi^q_L$ to be
isomorphisms are located in the vanishing of the local cohomology
modules with respect to ${\M}_1$ and ${\M}_2$. So our next goal
will be to study these local cohomology modules. For that, let us
consider
$${\ExsL}$$
the ${\Bbb Z^2}$-graded minimal free resolution of ${L}$ over $S$.
For every $p$, $D_p$ is a finite direct sum of $S$-modules of the
type $S(a,b)$. If we apply the diagonal functor to this
resolution, we get a resolution of $L_{\Delta}$ by means of the
modules $S(a,b)_\Delta$. Let us begin by studying the local
cohomology modules of the bigraded $S$-modules obtained by
shifting $S$ with degree $(a, b)$.

\medskip
First, let us fix some notations. For $ \alpha = (\alpha_1, \dots,
\alpha_n) \in \Bbb Z^n$, and $ \beta = (\beta_1, \dots, \beta_r)
\in \Bbb Z^r$, we write $X^{\alpha}$ for the monomial
$X_1^{\alpha_1} \cdots X_n^{\alpha_n}$ and $Y^{\beta}$ for the
monomial $Y_1^{\beta_1} \cdots Y_r^{\beta_r}$. Note that
deg$(X^{\alpha})= (\sum_{i=1}^{n} {\alpha_i}, 0)$,
deg$(Y^{\beta})=( \sum_{i=1}^{r} d_i \beta_i, \sum_{i=1}^r
\beta_i)$. We will write $\alpha < 0$ (or $\alpha \geq 0$) if all
the components of $\alpha$ satisfy this condition, and the same
for $\beta$. Then we have:

\begin{prop}
\label{W60} Let $a, b \in {\Bbb Z}$.
\begin{itemize}
\item[(i)]
 $$H^q_{{\M}_1} (S(a,b)) =
 \cases{ 0  & if  $q \not = n$ \cr
  \big( \bigoplus_{\alpha <0, \beta \geq 0} k X^{\alpha}  Y^{\beta}
\big) (a,b) & if $q=n$ \cr}
            $$

\item[(ii)]
 $$ H^q_{{\M}_2} (S(a,b))=
 \cases{0  & if  $q \not = r$ \cr
 \big( \bigoplus_{\alpha \geq 0, \beta < 0} k X^{\alpha}  Y^{\beta}
\big) (a,b) & if $q=r$ \cr}
            $$

\end{itemize}
\end{prop}

{\pf}
 Since $S(a,b)$ is a free $S_1$-module with basis the monomials in the
variables ${\y}$, we have that
 $H^q_{{\M}_1}(S(a,b)) = 0$ for all $q \not = n$, and
 $H^n_{{\M}_1}(S(a,b)) =
\big( \bigoplus_{\beta \geq 0} H^n_{{\frak m}_1}(S_1) Y^{\beta}
\big) (a,b) = \big( \bigoplus_{\alpha <0, \beta \geq 0} k
X^{\alpha}  Y^{\beta} \big) (a,b)$. By taking into account that
$S(a,b)$ is a free $S_2$-module with basis the monomials in the
variables ${\x}$, we also get $H^q_{{\M}_2}(S(a,b)) = 0$ for all
$q \not = r$, and $H^r_{{\M}_2} (S(a,b))= \big( \bigoplus_{\alpha
\geq 0} H^r_{{\frak m}_2}(S_2) X^{\alpha} \big) (a,b) = \big(
\bigoplus_{\alpha \geq 0, \beta < 0} k X^{\alpha}  Y^{\beta} \big)
(a,b)$. $\B$

\begin{cor}
\label{W6} Let $a, b \in {\Bbb Z}$.
\begin{itemize}
\item[(i)]
 $$
  {\supp} (\,H^q_{{\M}_1} (S(a,b))_\Delta \,) =
 \cases{ {\emptyset}  & if  $q \not = n$ \cr
 \{\, s  \in {\Bbb Z} \mid
\frac{-b}{e}  \le  s \le
 \frac{bd-a-n}{c-ed} \} & if $q=n$ \cr}
            $$

\item[(ii)]
 $$
  {\supp} (\,H^q_{{\M}_2} (S(a,b))_\Delta \,) =
 \cases{ {\emptyset}  & if  $q \not = r$ \cr
 \{\,s \in {\Bbb Z} \mid \frac{(b+r)d-u-a}{c-ed} \le  s \le
        \frac{-b-r}{e} \} & if $q=r$ \cr}
            $$

\end{itemize}
\end{cor}

{\pf} From Proposition \ref{W60}, a straightforward computation
gives the support by taking into account that a monomial
$X^{\alpha} Y^{\beta}$ in $H^n_{{\M}_1} (S(a,b))$ has degree $(p,
q)$ with $p =\sum_{i=1}^{n} {\alpha_i} + \sum_{j=1}^{r} d_j
\beta_j -a$ and $ q = \sum_{j=1}^r \beta_j-b$. Similarly one gets
$(ii)$. $\B$

\medskip
 For a real number $x$, let us denote by $[x] =\max \{ n \in \Bbb Z
\mid n \leq x \}$ the integral part of $x$. The following
corollary gives necessary and sufficient numerical conditions for
$S(a,b)_\Delta$ to be Cohen-Macaulay in terms of the diagonal
$\Delta$ and the shift $(a,b)$. In particular, notice that
$S_\Delta$ is Cohen-Macaulay for any $\Delta$.

\begin{cor}
\label{W666} \cite[Proposition 3.4]{CHTV} Assume $n, r \geq 2$.
For any $a, b \in \Bbb Z$, $S(a,b)_ \Delta$ is a Cohen-Macaulay
$S_\Delta$-module if and only if $ [\frac{bd-a-n}{c-ed}] <
\frac{-b}{e} $ and $ [\frac{-b-r}{e}] <
 \frac{(b+r)d-u-a}{c-ed}$.
\end{cor}

{\pf} Since $S$ is a domain, we have that $\rdim S(a,b) = \rdim S
= \dim S=n+r$, and so $\dim S(a,b)_\Delta= n+r-1$ by Lemma
\ref{W0}. Therefore, $S(a,b)_\Delta$ is Cohen-Macaulay if and only
if $H^q_{{\M}_\Delta} (S(a,b)_\Delta) =0$ for any $q<n+r-1$. By
Proposition \ref{W2}, note that for $q<n+r-1$ we have that
$$H^q_{{\M}_\Delta} (S(a,b)_\Delta) \cong
H^q_{{\cal M}_1} (S(a,b))_\Delta \oplus H^q_{{\cal M}_2}
(S(a,b))_\Delta.$$ Since $n,r \geq 2$, we get $n+r-2 \geq n, r$,
and then the result follows from Corollary \ref{W6}. $\B$

\begin{rem}
\label{W06} {\rm Note that if $n=r=1$, the proof above shows that
$S(a,b)_\Delta$ is always Cohen-Macaulay. In the case where $n
\geq 2$, $r=1$, we get that $S(a,b)_\Delta$ is Cohen-Macaulay if
and only if $ [\frac{-b-r}{e}] <
 \frac{(b+r)d-u-a}{c-ed}$, while if $n=1$, $r \geq 2$,
$S(a,b)_\Delta$ is Cohen-Macaulay if and only if $
[\frac{bd-a-n}{c-ed}] < \frac{-b}{e} $. }
\end{rem}

\medskip
For simplicity, from now on we will assume $n,r \geq 2$. Now, let
$L$ be a finitely generated bigraded $S$-module. For any $p \geq
0$, let us denote by $\Omega_{p,L}$ the set of shifts $(a,b)$
which appear in the place $p$ of its bigraded minimal free
resolution, and $\Omega_L$ the union of all these sets. Often we
will write $\Omega_p$, $\Omega$ if there is not danger of
confusion with respect to the module $L$. The next result relates
the local cohomology of the diagonals $L_\Delta$ of $L$ to the
local cohomology of the diagonals $S(a,b)_\Delta$ of the modules
$S(a,b)$ which arise in its minimal free resolution.

\begin{prop}
\label{W61} Let $L$ be a finitely generated bigraded $S$-module.
Then
\begin{itemize}
\item[(i)]
 If $H^q_{{\M}_1} (L)_{(cs,es)} \not = 0$, then there exists a shift
 $(a,b) \in \Omega_{n-q,L}$ such that
$H^n_{{\M}_1} (S(a,b))_{(cs,es)} \not =0$, and so $\frac{-b}{e}
\le  s \le  \frac{bd-a-n}{c-ed} $.
\item[(ii)]
 If $H^q_{{\M}_2} (L)_{(cs,es)} \not = 0$, then there exists a shift
 $(a,b) \in \Omega_{r-q,L}$ such that
$H^r_{{\M}_2} (S(a,b))_{(cs,es)} \not =0$, and so
 $ \frac{(b+r)d-u-a}{c-ed} \le  s \le \frac{-b-r}{e}$.
\end{itemize}
\end{prop}

{\pf} To prove $(i)$, let $0 \to D_t \to \dots \to D_0 \to L \to
0$ be the bigraded minimal free resolution of $L$ over $S$. By
considering $C_p= {\rm Coker}(D_{p+1} \rightarrow D_p)$ for $p
\geq 0$, this yields the short exact sequences
  $$ 0 \rightarrow C_{p+1} \rightarrow D_{p} \rightarrow C_{p}
     \rightarrow 0 , \; \; \forall p \geq 0.$$
If $H^q_{{\M}_1}(L) \not  = 0$, then $q \leq n$ because ${\cal
M}_1$ is generated by $n$ elements. In the case $q = n$, from the
short exact sequence $ 0 \rightarrow C_{1} \rightarrow D_{0}
\rightarrow L \rightarrow 0 $, we obtain a  bigraded epimorphism
$H^n_{{\M}_1}(D_0) \rightarrow  H^n_{{\M}_1}(L) $. Therefore, if
$H^n_{{\M}_1} (L)_{(cs,es)} \not = 0$ then $H^n_{{\M}_1}
(D_0)_{(cs,es)} \not = 0$, so by Corollary \ref{W6} there exists a
shift $(a,b) \in {\Omega}_0$ such that $\frac{-b}{e}  \le  s \le
\frac{bd-a-n}{c-ed} $. If $q<n$, since $H^v_{{\M}_1} (D_p) =0$ for
any $v \not = n$, we have bigraded isomorphisms
$$H^q_{{\M}_1}(L) \cong H^{q+1}_{{\M}_1}(C_1) \cong
H^{q+2}_{{\M}_1}(C_2) \cong \dots \cong
H^{n-1}_{{\M}_1}(C_{n-q-1}),$$ a bigraded monomorphism
$$0 \to H^{n-1}_{{\M}_1}(C_{n-q-1}) \to H^{n}_{{\M}_1}(C_{n-q}) ,$$
 and a bigraded epimorphism
$$ H^{n}_{{\M}_1}(D_{n-q}) \to H^{n}_{{\M}_1}(C_{n-q}) \to 0.$$
Therefore, $H^q_{{\M}_1} (L)_{(cs,es)} \not = 0$ implies $
H^{n}_{{\M}_1}(D_{n-q})_{(cs,es)} \not = 0$, and we are done by
Corollary \ref{W6}. Similarly one can prove $(ii)$. $\B$

\begin{rem}
\label{W10} {\rm Given a finitely generated bigraded $S$-module
$L$, for each diagonal $\Delta =(c,e)$ let us consider the sets of
integers $$X_p^\Delta= \bigcup _{(a,b) \in \Omega_{p,L}} \supp
(H^n_{{\M}_1}(S(a,b))_\Delta), $$
$$Y_p^\Delta = \bigcup _{(a,b) \in \Omega_{p,L}}
\supp (H^r_{{\M}_2}(S(a,b))_\Delta), $$ where $X_p ^\Delta = Y_p
^\Delta = \emptyset$ if $p<0$. Let $X^\Delta= \bigcup_{p}
X^\Delta_p$, $Y^\Delta= \bigcup_{p} Y^\Delta_p$. Then, Proposition
\ref{W61} jointly with Proposition \ref{W2} says that if $s \not
\in X^{\Delta}_{n-q} \cup Y_{r-q}^{\Delta}$, then
$[\varphi^q_L]_s$ is a monomorphism and $[\varphi^{q-1}_L]_s$ is
an epimorphism.  In particular, for an integer $s \not \in
X^{\Delta} \cup Y^{\Delta}$ then $[\varphi^q_L]_s $ is an
isomorphism for any $q$. (In fact, it is enough to define
$X^\Delta= \bigcup_{p \le n} X^\Delta_p$ and $Y^\Delta= \bigcup_{p
\le r} Y^\Delta_p$). }
\end{rem}

\medskip
Note that the set of integers $s$ satisfying that $\frac{-b}{e}
\le  s \le  \frac{bd-a-n}{c-ed} $ and $ \frac{(b+r)d-u-a}{c-ed}
\le  s \le \frac{-b-r}{e}$ is empty or it only contains the
integer $0$ for suitable $c,e$, that is, the sets $X^{\Delta}$ and
$Y^{\Delta}$ are contained in $\{0\}$ for those $\Delta=(c,e)$. So
we immediately get:

\begin{cor}
\label{W63} Let $L$ be a finitely generated bigraded $S$-module.
There exist positive integers $e_0, \alpha$ such that for any
$e>e_0$, $c > d e + \alpha$, we have isomorphisms $[\varphi^q_L]_s
: H^q_{{\cal M}_{\Delta}} (L_{\Delta})_s \rightarrow H^{q+1}_{\cal
M} (L)_{(cs,es)}$ for all $q \geq 0$ and $s \not =0$.
\end{cor}

{\pf} It is enough to take $e_0 \geq \max \{b, -b-r : (a,b) \in
\Omega_L \}$ and $\alpha \geq \max \{bd-a-n, u+a-(b+r)d : (a,b)
\in \Omega_L \}$. $\B$

\medskip
A similar result has been obtained in \cite[Lemma 3.8]{CHTV}.
Therefore, we have that for diagonals large enough the only
obstruction for the map $\varphi^q_L$ to be an isomorphism is
located in the component of degree $0$.

\begin{defn}
\label{W8} {\rm Let $L$ be a finitely generated bigraded
$S$-module and let
$$ \begin{array}{l}  {\ExsL} \end{array}  $$
be the bigraded minimal free resolution of $L$ over $S$. Let
$\Delta$ be a diagonal. We say that the resolution is good for
$\Delta$ if all the modules $(D_p)_{\Delta}$ are Cohen-Macaulay,
that is, $X^\Delta= Y^\Delta= \emptyset$. We say that the
resolution is good if there exists $\Delta$ such that the
resolution is good for $\Delta$.}
\end{defn}

\medskip
From Remark \ref{W10}, we immediately get that if $L$ has a good
resolution for $\Delta$ then the corresponding maps $\varphi^q_L $
are isomorphisms.

\begin{cor}
\label{W64} Let $L$ be a finitely generated bigraded $S$-module
whose resolution is good for $\Delta$. Then we have graded
isomorphisms
$$\varphi^q_L : H^q_{{\cal M}_{\Delta}} (L_{\Delta})
\rightarrow H^{q+1}_{\cal M} (L)_{\Delta}, \; \forall q \geq 0.$$
\end{cor}

\medskip
Our next goal is to study the existence of diagonals $\Delta$ for
which the bigraded minimal free resolution of $L$ is good for
$\Delta$. To this end, following \cite{CHTV} we define:

\begin{defn}
\label{W7} {\rm We say that a property holds for $c \gg 0$
relatively to $e \gg 0$ if there exists $e_0$ such that for all $e
> e_0$ there exists a positive integer $c(e)$ (depending on $e$)
such that this property holds for all $(c,e)$ with $c > c(e)$. We
will often write $c \gg e \gg 0$. In fact, in the statements we
will prove we could replace the condition $c \gg e\gg 0$ by the
stronger one that there exist positive integers $e_0, \alpha$ such
that the property holds for $e>e_0$, $c>de+\alpha$. For
simplicity, we will keep the notation and definition of $c \gg e
\gg 0$ from \cite{CHTV}. }
\end{defn}

\medskip
Next result provides necessary and sufficient numerical conditions
for the Cohen-Macaulayness of $S(a,b)_\Delta$ for $c \gg e \gg 0$.
Namely,

\begin{prop}
\label{W9} \cite[Corollary 3.5]{CHTV} Let $a, b \in \Bbb Z$. Then
$S(a,b)_ \Delta$ is a Cohen-Macaulay module  for $c \gg e \gg 0$
if and only if $a, b$ satisfy one of the following conditions:
\begin{enumerate}
\item $b \leq -r$ and $(b+r)d-u-a>0$,
\item $-r<b<0$,
\item $b \geq 0$ and $bd-a-n<0$.
\end{enumerate}
\end{prop}

{\pf} From Proposition \ref{W2}, we have that $S(a,b)_\Delta$ is
Cohen-Macaulay for $c \gg e \gg 0$ if and only if $0 \not \in
\supp (H^n_{{\M}_1}(S(a,b))_\Delta) \cup \supp
(H^r_{{\M}_2}(S(a,b))_\Delta) $ for $c \gg e \gg 0$. Then the
result follows from Corollary \ref{W6}. $\B$

\medskip
Notice that, for a given diagonal $\Delta$, we have that
$S(a,b)_\Delta$ is Cohen-Macaulay if and only if the corresponding
maps $\varphi^q_{S(a,b)}$ are isomorphisms for all $q \geq 0$. On
the other hand, from the proof of Proposition \ref{W9}, observe
that if $(a,b)$ does not satisfy any of the conditions above then
$S(a,b)_\Delta$ is never Cohen-Macaulay. Therefore, we can not
hope to extend Corollary \ref{W63} to the component of degree $0$
of the maps $\varphi^q_L$. In fact, the proof of Proposition
\ref{W2} shows that $[\varphi^q_L]_0$ does not depend on the
diagonal $\Delta$.

\medskip
Furthermore, note that the proof of Proposition \ref{W9} also
shows that if there exists $\Delta$ such that $S(a,b)_\Delta$ is
Cohen-Macaulay then $S(a,b)_\Delta$ is Cohen-Macaulay for $c \gg e
\gg 0$. Therefore, if a finitely generated bigraded $S$-module $L$
has a good resolution, then the resolution of $L$ is good for
diagonals $\Delta=(c,e)$ with $c \gg e \gg 0$.

\medskip
Up to now we have related the vanishing of the local cohomology
with respect to ${\M}_1$ and ${\M}_2$ of a bigraded $S$-module $L$
with the vanishing of the local cohomology with respect to
${\M}_1$ and ${\M}_2$ of the modules $S(a,b)$ which arise in the
bigraded minimal free resolution of $L$ over $S$. This study has
led us to get sufficient conditions on the shifts $(a,b)$ in order
to $\varphi^q_L$ to be isomorphisms. In the rest of the section we
shall deal with the computation of the local cohomology modules of
a bigraded $S$-module $L$ with respect to the ideals ${\M}_1$ and
${\M}_2$ by themselves.

\medskip
In Corollary \ref{W64} we have given sufficient conditions on the
shifts in $\Omega_L$ to get that the maps $\varphi^q _L$ are
isomorphisms for large diagonals. Next we give necessary and
sufficient conditions for the maps $\varphi^q_L$ to be
isomorphisms in terms of the local cohomology modules of $L$ with
respect to ${\M}_1$ and ${\M}_2$. Namely,

\begin{prop}
\label{W633} Let $L$ be a finitely generated bigraded $S$-module.
Then the following are equivalent:
\begin{enumerate}
\item There exists $\Delta$ such that $\varphi^q_L$ is an
 isomorphism for all $q \geq 0$.
\item For large diagonals $\Delta$, $\varphi^q_L$ is an isomorphism for
 all $q \geq 0$.
\item $H^q_{{\M}_1}(L)_{(0,0)}= H^q_{{\M}_2}(L)_{(0,0)}=0$ for all
 $q \geq 0$.
\end{enumerate}
\end{prop}

\medskip
For an integer $e$ and a bigraded $S$-module $L$, let us define
the graded $S_1$-module $L^e = \oplus_{i \in \Bbb Z} L_{(i,e)}$.
Then we have an exact functor $(\;)^e :M^2(S) \rightarrow
M^1(S_1)$. The bigraded initial degree of a bigraded $S$-module
$L$ is defined by ${\rm {\bf indeg}}(L)= ({\rm indeg}_1(L),{\rm
indeg}_2(L))$, where
$${\rm indeg}_1(L) = \min
\{ i \mid \exists j \; s.t. \; L_{(i,j)} \not = 0 \},$$
$${\rm indeg}_2(L) = \min
\{ j \mid \exists i \; s.t. \; L_{(i,j)} \not = 0 \}. $$

\begin{prop}
\label{W11} Let $L$ be a finitely generated bigraded $S$-module.
Then:
\begin{enumerate}
\item
$H^q_{{\M}_1} ( L )_{(i,j)} = H^q_{{\frak m}_1}(L^j)_i$ . In
particular, $H^q_{{\M}_1} ( L )_{(i,j)}= 0$ for $i> a_q(L^j)$ or
$j < {\rm indeg}_2(L)$.
\item
$H^q_{{\M}_2} ( L )_{(i,j)} = 0$ for $j> a^2_*(L)$.
\end{enumerate}
\end{prop}

{\pf}
 As $S_1$-module, $L$ is the direct sum of the modules
$L^e = \oplus_i L_{(i,e)}$. Since ${\M}_1$ is the ideal of $S$
generated by ${\frak m}_1 =(\x)$, we have that $H^q_{{\cal
M}_1}(L) = \oplus_j H^q_{{\frak m}_1}(L^j)$, and so we get $(i)$.

Now let ${\ExsL}$ be the bigraded minimal free resolution of ${L}$
over $S$, where $D_p=\bigoplus_{(a,b) \in \Omega_p} S(a,b)$. By
taking short exact sequences as in Proposition \ref{W61}, it is
just enough to prove that if $j>a_*^2(L)$ then $H^q_{{\M}_2}
(S(a,b))_{(i,j)} = 0$ for any $(a,b) \in \Omega_L$ and $q \geq 0$.
The case $q \not = r$ is trivial. From Proposition \ref{W60}, we
may deduce that $H^r_{{\M}_2} (S(a,b))_{(i,j)} = 0$ for $j >
-b-r$. This finishes the proof because, according to Theorem
\ref{B77}, $a^2_*(L) \geq -b-r$ for any $(a,b) \in \Omega_L$. $\B$

\medskip
In the particular case $d_1 = \dots = d_r = d$, $S$ can be thought
as a standard bigraded $k$-algebra by a change of grading. If we
consider the morphism $\varphi(p, q) = (p-dq, q)$, observe that
$\varphi( \supp S) \subset \Bbb N^2$, so $S^\varphi$ is a $\Bbb
N^2$-graded ring with $[S^\varphi]_{(p, q)} = S_{(p+dq, q)}$.
Noting that $\deg(X_i)=(1,0)$ for $i=1,\dots,n$, and
$\deg(Y_j)=(0,1)$ for $j=1, \dots, r$ as elements of $S^\varphi$,
we have that $S^\varphi$ is standard. For a bigraded $S$-module
$L$, let us recall that the $S^\varphi$-module $L^\varphi$ is the
$S$-module $L$ with the grading defined by $[L^\varphi]_{(p, q)} =
L_{(p+dq, q)}$.

\medskip
Furthermore, in this case, given an integer $e $ we can define an
exact functor $(\;)_e :M^2(S) \rightarrow M^1(S_2)$ in the
following way: For any bigraded $S$-module $L$, we define $L_e$ to
be the graded $S_2$-module $L_e = \bigoplus_{j \in \Bbb Z}
L_{(e+dj,j)}$. Then we have

\begin{prop}
\label{W32} Assume that $d_1= \dots = d_r=d$. For any finitely
generated bigraded $S$-module $L$, we have
\begin{enumerate}
\item
$H^q_{{\M}_1} ( L )_{(i,j)} = 0$ for $i> dj + a^1_*(L^\varphi)$.
\item
$H^q_{{\M}_2} ( L )_{(i,j)} = H^q_{{\frak m}_2}(L_{i-dj})_j$ . In
particular, $H^q_{{\M}_2} ( L )_{(i,j)}= 0$ for $j>
a_q(L_{i-dj})$.
\end{enumerate}
\end{prop}

{\pf} Let ${\ExsL}$ be the bigraded minimal free resolution of
${L}$ over $S$. Observe that

\vspace{3mm}

\hspace{20mm} $S(a,b)^{\varphi} =  \bigoplus_{(i,j)}
S(a,b)_{(i+dj,j)}$

\vspace{2mm}

\hspace{35mm} $=\bigoplus_{(i,j)} S_{(a+i+dj,b+j)}$

\vspace{2mm}

\hspace{35mm} $= \bigoplus_{(i,j)} S_{(a-db+i+d(b+j),b+j)}$

\vspace{2mm}

\hspace{35mm} $=S^{\varphi}(a-db,b),$

\vspace{1.5mm}

\noindent so in particular $S(a,b)^{\varphi}$ is a free
$S^\varphi$-module. Therefore, by applying the exact functor
$(\;)^\varphi$ to the resolution of $L$ we get that
$$ 0 \to D_{t}^{\varphi} \to \dots \to D_0^{\varphi} \to L^{\varphi}\to 0$$
is a bigraded minimal free resolution of $L^{\varphi}$ over
$S^{\varphi}$. Since $a^1(S^\varphi)=-n$, from Theorem \ref{B77}
it follows that
 $$a_*^1 (L^\varphi) = \max \,\{\,db-a \mid (a,b) \in {\Omega}_L
  \} - n.$$
If $i, j$ are such that $i > dj + a_*^1(L^\varphi)$, then  we have
that $i > dj + db -a-n$ for any shift $(a,b) \in {\Omega}_L$, and
so from Proposition \ref{W60} we have that $H^q_{{\M}_1}
(S(a,b))_{(i,j)} = 0$ for any $q \geq 0$. By taking short exact
sequences as in Proposition \ref{W11}, we then obtain
$H^q_{{\M}_1} (L)_{(i,j)} =0$ for $q \geq 0$, $i > dj +
a_*^1(L^\varphi)$.

To prove $(ii)$, note that since $d_1 = \dots = d_r=d$ we may
decompose $L$ as the direct sum of the $S_2$-modules $L_i$. Then,
by using that ${\M}_2$ is the ideal of $S$ generated by ${\frak
m}_2 =(\y)$,  we obtain $H^q_{{\cal M}_2}(L) = \bigoplus_i
H^q_{{\frak m}_2} (L_i)$. Noting that $\dg (Y_1)= \dots = \dg
(Y_r)=(d, 1)$, we finally get $H^q_{{\cal M}_2}(L)_{(i,j)} =
H^q_{{\frak m}_2}(L_{i-dj})_j$. $\B$

\bigskip
\section{Case study: Standard bigraded $k$-algebras}
\label{Y} \markboth{CHAPTER II. THE DIAGONALS OF A BIGRADED
MODULE} {STANDARD BIGRADED $k$-ALGEBRAS}

\medskip
Our aim in this section is to particularize and improve for
standard bigraded $k$-algebras several results proved in Section
2.1. So let $R$ be a standard bigraded $k$-algebra generated by
homogeneous elements ${\xx}$, ${\yy}$ in degrees $\dg(x_i)
=(1,0)$, $i= 1, \dots, n$, $\dg (y_j)=(0,1)$, $j=1, \dots, r$. By
taking the polynomial ring $S=k[{\x},{\y}]$ with the bigrading
given by deg$(X_i)= (1, 0)$, deg$(Y_j)= (0, 1)$, we have that $R$
is a finitely generated bigraded $S$-module in a natural way.

\medskip
In this case, denote by ${\cal R}_1 = R^0 = \bigoplus_{i \in \Bbb
N} R_{(i,0)}$, ${\cal R}_2 = R_0 = \bigoplus_{j \in \Bbb N}
R_{(0,j)}$. Observe that ${\cal R}_1$ and ${\cal R}_2$ are graded
$k$-algebras, and denote by $\fm_1$ and $\fm_2$ their homogeneous
maximal ideals. Given $e \in \Bbb Z$, we may define the graded
${\cal R}_1$-module $R^e= \oplus_{i \in \Bbb Z} R_{(i,e)}$ and the
graded ${\cal R}_2$-module $R_e= \oplus_{j \in \Bbb Z} R_{(e,j)}$.
By a straightforward application of Proposition \ref{W2}, we get

\begin{prop}
\label{Y2} There is a natural graded exact sequence
$$ ... \to H^q_{\M} (R)_{\Delta}
\to H^q_{{\M}_1} ( R )_{\Delta}
 \oplus H^q_{{\M}_2} ( R )_{\Delta}
\to H^q_{\MD}(R_{\Delta}) \stackrel {\varphi^q_R} \longrightarrow
 H^{q+1}_{\cal M}(R)_{\Delta}  \to ...$$
In particular, given $s \in \Bbb Z$, the following are equivalent:
\begin{itemize}
\item[(i)] $[\varphi^q_{R}]_s$ is isomorphism, $\forall q \geq 0$.
\item[(ii)] $H^q_{\fm_1} (R^{es})_{cs} = 0$ and
 $H^q_{\fm_2}(R_{cs})_{es}=0$,  $\forall q \geq 0$.
\end{itemize}
\end{prop}

{\pf} It follows from Proposition \ref{W2}, Proposition \ref{W11}
and Proposition \ref{W32}. $\B$

\medskip
As a direct consequence of Proposition \ref{Y2} we have:

\begin{cor}
\label{Y4} $\varphi^q_R$ is an isomorphism for $q> \max \{ \dim
{{\cal R}_1}, \dim {{\cal R}_2} \}$.
\end{cor}

{\pf} Set $d_1= \dim {{\cal R}_1}$, $d_2= \dim {{\cal R}_2}$. It
is enough to prove that $H^q_{\fm_1} (R^{e}) = H^q_{\fm_2} (R_{e})
= 0$ for any $e \in \Bbb Z$, $q> \max \{d_1, d_2 \}$. But note
that $R^e$ is a graded ${\cal R}_1$-module, so $H^q_{\fm_1}
(R^{e}) =0$ for $q>d_1$. Similarly, $H^q_{\fm_2} (R_{e}) =0$ for
$q>d_2$ and we are done. $\B$

\medskip
From Proposition \ref{Y2} we can also determine a set of integers
$s$, depending on the diagonal $\Delta$, for which
$[\varphi^q_R]_s$ is an isomorphism for all $q \geq 0$. More
explicitly,

\begin{cor}
\label{Y22}
\begin{enumerate}
\item $[\varphi^q_R]_s$ is isomorphism for $s < 0$.
\item $[\varphi^q_R]_s$ is isomorphism for
 $s > \max \{\, a_*^1(R)/c,  a_*^2(R)/e \}$.
In particular,
 $a_*(R_\Delta) \leq \max \{\, a_*^1(R)/c,  a_*^2(R)/e \}$.
\end{enumerate}
\end{cor}

{\pf} It is a direct consequence of Proposition \ref{Y2},
Proposition \ref{W11} and Proposition \ref{W32}. $\B$

\medskip
We have shown that $[\varphi^q_R]_s$ is an isomorphism for any $s
< 0$. Moreover, note that if $c>a^1_*(R)$ and $e>a^2_*(R)$, then
$[\varphi^q_R]_s$ is an isomorphism for any $s>0$. We may ensure
that $[\varphi^q_R]_0$ is an isomorphism for any $q$ if $R$ has a
good resolution. Next we study the existence of a such resolution.
The following result provides a useful characterization for a
standard bigraded $k$-algebra $R$ to have a good resolution by
means of the $a_*$-invariants of ${\cal R}_1$ and ${\cal R}_2$.
Namely,

\begin{prop}
\label{Y1} The following are equivalent:
\begin{itemize}
\item[(i)] $R$ has a good resolution.
\item[(ii)] $a_*({\cal R}_1) <0$, $a_*({\cal R}_2) < 0$.
\end{itemize}
\end{prop}

{\pf} Let us consider
$$ {\bf D_{\bf .}} : \,\, 0 \to D_t \to \cdots \to D_1 \to D_0 =S \to R
\to 0 $$ the bigraded minimal free resolution of $R$ over $S$,
where $D_p = \bigoplus_{(a,b) \in {\Omega}_p} S(a, b)$. Note that
$a,b$ are non-positive integers. So, by Proposition \ref{W9} the
resolution is good if and only if all the shifts $(a, b)$ satisfy
one of the three following conditions
\begin{enumerate}
\item $-r < b < 0$.
\item $ b= 0$ and $-n< a$.
\item $b \leq -r$ and $a<0$.
\end{enumerate}
It is not hard to check that these conditions are equivalent to
that for any shift $(a,0) \in {\Omega}_R$ we have $a >-n$, and for
any shift $(0, b) \in {\Omega}_R$ we have $b>-r$. Observe that
 $$ S(a, b)^0 = \bigoplus_{j} S_{(a+ j, b)} =
 \cases{ 0  & if $b < 0$ \cr
  S_1 (a)  &  if $b= 0$ \cr}
 $$
So by applying the functor $(\;\;)^0$ to the resolution of $R$ we
obtain a graded free resolution of ${\cal R}_1$ over $S_1$:
$$ {\bf F_{\bf .}} : \,\,0 \to F_t \to \cdots \to F_1 \to F_0 =S_1   \to
{{\cal R}_1} \to 0 \;,$$ where $F_p = (D_p)^0 = \bigoplus_{a \in
\gamma_p} S_1(a)$, and $\gamma_p = \{ a \in \Bbb Z : (a, 0) \in
\Omega_p \}$ ($F_p =0$ if $\gamma_p = \emptyset$). Furthermore, we
have that Im$(F_p) \subset  \fm_1 F_{p-1}$ for all $p= 1,..., t$.
Hence this resolution is in fact the graded minimal free
resolution of ${{\cal R}_1}$ over $S_1$. Then we can use Theorem
\ref{B77} to compute $a_*({\cal R}_1)$ :

\vspace{3mm}

\hspace{10 mm} $a_*({\cal R}_1) =
    \max \{-a \mid a \in \cup_p \gamma_p \} + a(S_1)  = $

\vspace{3mm}

\hspace{23 mm}   $ =  \max \{ -a \mid (a, 0) \in \Omega_R \} - n.$

\vspace{3mm}

\noindent Therefore, any shift $(a, 0) \in {\Omega}_R$ satisfies
$a>-n$ if and only if $a_* ({\cal R}_1) < 0$. Similarly, any shift
$(0,b) \in {\Omega}_R$ satisfies $b >-r$ if and only if $a_*
({\cal R}_2)< 0$. $\Box$

\medskip
As an immediate consequence we get a criterion for the existence
of Cohen-Macaulay diagonals of a standard bigraded $k$-algebra
which extends \cite[Corollary 3.12]{CHTV}. More explicitly,

\begin{cor}
\label{Y11} Let $R$ be a standard bigraded $k$-algebra with
$a_*({\cal R}_1) <0$, $a_*({\cal R}_2) < 0$. Then $\depth \,
R_\Delta \geq \depth \,R-1$ for large $\Delta$. In particular, if
$R$ is Cohen-Macaulay, then so $R_\Delta$ for large $\Delta$.
\end{cor}

\medskip
For a standard bigraded ring $R$ defined over a local ring with
$a^1(R), a^2(R) <0$, it has been shown in \cite[Theorem 2.5]{H}
that if $R$ is Cohen-Macaulay, then its $(1,1)$-diagonal inherits
this property. This result can be extended to any diagonal of a
standard bigraded $k$-algebra.

\begin{prop}
\label{Y12} Let $R$ be a standard bigraded Cohen-Macaulay
$k$-algebra with $a^1(R), a^2(R) < 0$. Then $R_\Delta$ is
Cohen-Macaulay for any diagonal $\Delta$.
\end{prop}

{\pf} The bigraded standard $k$-algebra $R$ has a presentation as
a quotient of the polynomial ring $S=k[{\x},{\y}]$ bigraded by
$\deg(X_i)=(1,0)$, $\deg(Y_j)=(0,1)$. According to Theorem
\ref{B77},  for any shift $(a,b) \in \Omega_R$ we have that
$$0 \leq -a \leq a^1(R)-a^1(S) < n$$
$$0 \leq -b \leq a^2(R)-a^2(S) < r.$$
Then note that for any diagonal $\Delta=(c,e)$ with $c,e>0$,
$$X^\Delta= \bigcup_{(a,b) \in \Omega_R} \{ \; s \in \Bbb Z \mid
\frac{-b}{e}  \le  s \le \frac{-a-n}{c} \; \} = \emptyset $$
$$Y^\Delta= \bigcup_{(a,b) \in \Omega_R} \{ \; s \in \Bbb Z \mid
 \frac{-a}{c} \le  s \le \frac{-b-r}{e} \; \}= \emptyset ,$$
so the resolution is good for any $\Delta$. Now, by Corollary
\ref{W64} we have $H^q_{{\cal M}_{\Delta}} (R_{\Delta}) \cong
H^{q+1}_{\cal M} (R)_{\Delta} =0$ for $q< \dim R - 1$, so we are
done. $\B$

\medskip
We finish this section by giving necessary and sufficient
conditions on the local cohomology of ${\cal R}_1$ and ${\cal
R}_2$ for the existence of a Cohen-Macaulay diagonal of a standard
bigraded Cohen-Macaulay $k$-algebra $R$. Namely,

\begin{prop}
\label{Y01} Let $R$ be a standard bigraded Cohen-Macaulay
$k$-algebra of relevant dimension $\delta$. Then there exists
$\Delta$ such that $R_\Delta$ is Cohen-Macaulay if and only if
$H^q_{{\fm}_1}({\cal R}_1)_{0}= H^q_{{\fm}_2}({\cal R}_2)_{0}= 0$
for any $q< \delta -1$.
\end{prop}

{\pf} According to Lemma \ref{W0}, we have that $\dim R_\Delta=
\delta-1$ for any $\Delta$. By taking into account Corollary
\ref{W63}, there exists $\Delta$ such that $R_\Delta$ is
Cohen-Macaulay if and only if there exists $\Delta$  such that
$H^q_{\MD}(R_\Delta)_0=0$ for any $q<\delta-1$. But from
Proposition $\ref{Y2}$, for any $q<\delta-1$ we have
$$H^q_{\MD}(R_{\Delta})_0
\cong H^q_{{\fm}_1} ({\cal R}_1 )_{0}
 \oplus H^q_{{\fm}_2} ( {\cal R}_2 )_{0}.$$
This finishes the proof. $\B$



\bigskip
\section{Case study: Rees algebras}
\markboth{CHAPTER II. THE DIAGONALS OF A BIGRADED MODULE} {REES
ALGEBRAS} \label{ED}

\medskip
Let $A$ be a noetherian graded algebra generated in degree 1 over
a field $k$. Then $A$ has a presentation $A= k[{\x}]/ K =
k[{\xx}]$, where $K$ is a homogeneous ideal of the polynomial ring
$k[{\x}]$ with the usual grading. Let $\fm$ be the graded maximal
ideal of $A$. For a homogeneous ideal $I$ of $A$, let us consider
the Rees algebra
$$R= R_A(I) = \bigoplus_{n \geq 0}  I^n t^n \subset A[t]$$
of $I$ endowed with the natural bigrading given by
$$R_{(i,j)} =(I^j)_i,$$
introduced by A. Simis et al. in \cite{STV}. If $I$ is generated
by forms $f_1, \dots, f_r$ in degrees $d_1, \dots, d_r$
respectively, note that $R$ is a $k$-algebra finitely generated by
${\xx}, f_1 t,...., f_r t$ with deg$(x_i)= (1,0)$, deg$(f_j
t)=(d_j,1)$, and that it has a unique homogeneous maximal ideal
${\cal M}=({\xx}, f_1 t,...., f_r t)$. By considering the
polynomial ring $S=k[{\x},{\y}]$ with the grading determined by
setting deg$(X_i)= (1, 0)$ and deg$(Y_j)=(d_j, 1)$, we have a
bigraded epimorphism:
    $$  \begin{array}{lll}
     S   & \longrightarrow  &   R  \\
        X_i    & \longmapsto  &  x_i  \\
        Y_j    & \longmapsto   & f_j t
        \end{array} $$
so that $R$ has a natural structure as finitely generated bigraded
$S$-module. Set $d= \max \{d_1, \dots, d_r\}$. For any $c \geq
de+1$, the $\Delta=(c,e)$-diagonal of the Rees algebra is
$$R_A(I)_\Delta = \bigoplus_{s \geq 0} (I^{es})_{cs} = k[(I^e)_c].$$
Note that $\Kk$ is a graded $k$-algebra with a unique homogeneous
maximal ideal $m= \MD$. The interest of these algebras $\Kk$ is,
as we will show in the next chapter, that for any $c \geq de+1$
there is a projective embedding of the blow-up $X$ of $\Proj (A)$
along $\widetilde{I}$ so that $X \cong \Proj (\Kk)$.

\medskip
Set $\overline n = \dim A$. The next lemma computes the dimension
of the rings $\Kk$ extending \cite[Lemma 1.3]{CHTV} where the case
$A=k[{\x}]$ was studied.

\begin{lem}
\label{ED1} Assume $I \not \subset \frak p$, for all $\frak p \in
\Ass (A)$. Then $\dim {\Kk} = \overline n$ for all $c \geq de+1$.
\end{lem}

{\pf}
 Since $I$ is not contained in any associated prime of $A$, we have that
any associated prime ideal of the Rees algebra $R$ is relevant. So
$\rdim \,R = \dim R$. Furthermore, $\dim R = \nn +1$ by
\cite[Exercise 4.4.12]{BH}. So we may conclude $\dim {\Kk} = \nn$
by Lemma \ref{W0}.    $\Box$

\medskip
From now on we will always assume that $I \not \subset \frak p$,
for all $\frak p \in \Ass (A)$. The following result relates the
local cohomology of the graded $k$-algebras $\Kk$ to the local
cohomology of the Rees algebra. By setting $\frak n= (It)=(f_1t,
\dots, f_rt) \subset k[It]= k[f_1t, \dots, f_rt]$, we have

\begin{prop}
\label{ED2} Let $I$ be an ideal of $A$ generated by forms of
degree $\leq d$. For any diagonal $\Delta =(c,e)$ with $c \geq
de+1$, there is a natural graded exact sequence
$$ ... \to H^q_{\M} (R)_{\Delta}
\to H^q_{\fm R} ( R )_{\Delta} \oplus H^q_{\fn R} ( R )_{\Delta}
\to H^q_{m}( \Kk) \stackrel {\varphi^q_R} \longrightarrow
 H^{q+1}_{\cal M}(R)_{\Delta}  \to ...$$
\end{prop}

{\pf} It is clear that $H^q_{{\M}_1} ( R ) = H^q_{{\frak m} R}(R)$
and $H^q_{{\M}_2} ( R ) = H^q_{{\frak n} R}(R)$. Then the result
follows immediately by applying Proposition \ref{W2} to the Rees
algebra $R$ of $I$.$\B$

\medskip
In Corollary \ref{W3} we proved that the maps $\varphi^q_L$ become
isomorphisms for $q > \max \{n,r \}$. This bound was refined for
standard bigraded $k$-algebras in Corollary \ref{Y4}. Next we want
to consider the case of the Rees algebras. To this end, we are
going to study the vanishing of the local cohomology modules of
$R$ with respect to $\fn R$.  For any ideal $I$ of $A$, the fiber
cone of $I$ is defined as the graded $k$-algebra $F=F_\fm (I)=
\bigoplus_{n \geq 0}I^n/ \fm I^n$. The analytic spread $l(I)$ of
$I$ is then the dimension of the fiber cone, that is, $l(I)= \dim
F$. Note that if $I$ is generated by forms of the same degree $d$,
the fiber cone is nothing but $\F =k[I_d]$. The following lemma
shows the known result that the local cohomology modules of $R$
with respect to ${\fn} R$ vanish in order $q >l(I)$, but not in
order $l(I)$. We include the proof for the sake of completeness.

\begin{lem}
\label{ED5} Let $I$ be a homogeneous ideal of $A$. Set $l=l(I)$.
Then $H^q_{{\fn} R}(R)=0 \;, \forall q>l$ and $H^l_{{\fn} R}(R)
\not = 0.$
\end{lem}

{\pf} We may assume that the field $k$ is infinite. Then there
exists an ideal $J \subset I$ generated by $l=l(I)$ elements of
$A$ such that $I^m = J I^{m-1}$ for $m \gg 0$, that is, there
exists a reduction $J$ of $I$ generated by $l$ elements (see
\cite[Proposition 4.6.8]{BH}). Note that $(It) R$ and $(Jt) R$ are
ideals with the same radical, so $H^q_{{\fn} R}(R)= H^q_{I R}(R)=
H^q_{J R}(R)= 0, \forall q > l$. Moreover, from the presentation
$R \rightarrow R/ \fm R= F_{\fm}(I)$ we get the epimorphism
$$H^l_{{\fn} R}(R) \rightarrow H^l_{{\fn} R}(F)\not = 0,$$
so $H^l_{{\fn} R}(R) \not =0$. $\B$

\medskip
As a consequence, we get:

\begin{cor}
\label{ED22} Let $I$ be an ideal of $A$ generated by forms in
degree $\leq d$. For any $\Delta$, we have a graded epimorphism
$$ H^{\nn}_{m}( \Kk)
\stackrel {\varphi^{\nn}_R} \longrightarrow
 H^{\nn +1}_{\cal M}(R)_{\Delta}.$$
\end{cor}

\medskip
From Proposition \ref{ED2} we may also deduce that for diagonals
large enough the positive components of the local cohomology of
the diagonals of the Rees algebra coincide with the positive
components of the local cohomology of the powers of the ideal.
Namely,

\begin{cor}
\label{ED4} Let $I$ be an ideal generated by forms of degree $\leq
d$. For any $c \geq de+1$, $e>a_*^2(R)$, $s>0$, there are
isomorphisms
$$H^q_{m} (\Kk)_s \cong H^{q}_{\fm} (I^{es})_{cs} \;\;,
\forall q \geq 0.$$
\end{cor}

{\pf} Let $c,e$ be integers such that $c \geq de+1$, $e>a_*^2(R)$.
For any $s>0$, we have that $H^q_{{\M}_1} ( R )_{(cs,es)} =
H^q_{{\frak m}}(I^{es})_{cs}$ and $H^q_{{\M}_2} ( R )_{(cs,es)} =
0$ by Proposition \ref{W11}. Thus from Proposition \ref{ED2} we
get the isomorphisms $H^q_{m} (\Kk)_s \cong H^{q}_{\fm}
(I^{es})_{cs}$.$\B$

%

\medskip
In the case where $I$ is generated by forms in the same degree $d$
(that is, $I$ is equigenerated) the Rees algebra is a standard
bigraded $k$-algebra by setting
$$R_A(I)_{(i,j)}= (I^j)_{i+dj}.$$
Then we may apply the results in Section 2.2 to these Rees
algebras. From Lemma \ref{Y1},  we get a useful characterization
for the Rees algebra to have a good resolution by means of the
$a_*$-invariants of the ring $A$ and the fiber cone of $I$.
Namely,

\begin{prop}
\label{ED7} Let $I$ be an ideal of $A$ generated by forms in
degree $d$. The following are equivalent:
\begin{enumerate}
\item   The Rees algebra $R_A (I)$ has a good resolution.
\item   $a_*(A) < 0$,
        $a_*(F_{\frak m}(I)) < 0$.
\end{enumerate}
\end{prop}

{\pf} We have already noted that the Rees algebra is a standard
bigraded ring by means of $R_{(i,j)}= (I^j)_{i+dj}$. With this
grading, notice that  ${\cal R}_1=A$ and ${\cal R}_2= k [I_d]=
F_{\fm}(I)$. Then the result follows from Lemma \ref{Y1}. $\B$


\medskip
To apply Proposition \ref{ED7}, we need to know the
$a_*$-invariant of the fiber cone. The next two lemmas bound it.
The first one gives a lower bound by means of the reduction number
of $I$ (compare with \cite{T}, \cite{Sc2}), while the second one
gives an upper bound by means of the $a_*$-invariant of the Rees
algebra.

\begin{lem}
\label{ED9}
  Let $(A, \frak m)$ be a local noetherian ring with an infinite
residue field. Let $I \subset \frak m$ be an arbitrary ideal of
$A$, $J$ a minimal reduction of $I$ and $l$ the analytic spread of
$I$. Then
$$a_l ({\F}) + l \leq  r_J (I)  \leq \max \{ a_i ({\F}) + i \}
= \reg ({\F}).$$

\end{lem}

{\pf} Let $a_1,..., a_l$ be a minimal system of generators of $J$.
For $a \in  I$, denote by $a^0$ the class of $a$ in $I/\fm I$.
Then $a_1^0,..., a_l^0 $ are a (homogeneous) system of parameters
of ${\F}$ (see \cite[Proposition 10.17]{HIO}). According to
\cite[Lemma 45.1]{HIO}, we have
   $$a_l ({\F}) + l \leq
\max \{ n \mid \bigg[ \frac{{\F}}{(a_1^0,..., a_l^0)} \bigg]_n
\not =0 \}
    \leq \max \{ a_i ({\F}) + i \}$$
\noindent On the other side,

\vspace{2mm}

\hspace{10mm} $r_J (I) = \min \{n \mid I^{n+1} = J I^n \} $

\vspace{2mm}

\hspace{20mm} $=\min \{\; n \mid  \frac {I^{n+1}} {\fm I^{n+1}} =
 \frac { J I^n + \fm I^{n+1}} {\fm I^{n+1}}   \; \} $

\vspace{2mm}

\hspace{20mm} $= \min \{ \; n \mid \bigg[ \frac {{\F}}
{(a_1^0,...,a_l^0)} \bigg]_{n+1} = 0 \; \}$

\vspace{2mm}

\hspace{20mm}  $= \max \{ \; n \mid \bigg[ \frac {{\F}}
{(a_1^0,..., a_l^0)} \bigg]_n \not =0 \; \}. $

\vspace{2mm}

\noindent
 This concludes the lemma. $\B$

\medskip
The next lemma bounds the $a_*$-invariant of the fiber cone by
means of the $a_*$-invariant of the Rees algebra. Namely,

\begin{lem}
\label{ED10} Let $I$ be an equigenerated homogeneous ideal of $A$.
Then $a_*({\F}) \leq a_*^2 (R_A(I))$.
\end{lem}

{\pf} Let
$$ {\bf D_.} : \,\, 0 \to D_t \to \cdots \to D_1 \to D_0 =S  \to
R=R_A(I) \to 0 $$ be the bigraded minimal free resolution of the
Rees algebra $R$ over $S$, where $D_p = \bigoplus_{(a,b) \in
{\Omega}_p} S(a, b)$. Note that $a,b$ are non-positive integers
with $a \leq db$. Therefore, we have that
 $$
 S(a, b)_{0} = \bigoplus_{j} S_{(a+ dj, b+j)} =
 \cases{ 0  & if $a < db$ \cr
  S_2 (b)  &  if $a = db$ \cr}
 $$
Then by applying the functor $(\;)_0$ to the resolution ${\bf
D_.}$ we get a graded free resolution of $R_0= {\F}$ over $S_2$:
$$ {\bf F_.} : \,\,0 \to F_t \to \cdots \to F_1 \to F_0 =S_2 \to F
\to 0  \;,$$ where $F_p = (D_p)_{0} = \bigoplus_{b \in \gamma_p}
S_2 (b)$, and $\gamma_p =\{ b \in \Bbb Z : (db,b) \in \Omega_p
\}$. Moreover, for any $p=1, \dots, t$ we have that Im$(F_p)
\subset  {\fm}_2 F_{p-1}$, so $ {\bf F_.}$ is in fact the graded
minimal free resolution of $\F$ over $S_2$. Then by Theorem
\ref{B77} we have:

\vspace{3mm}

\hspace{10mm} $  a_*(F)  =
   \max \{-b \mid b \in \cup_p \gamma_p \} + a(S_2)$

\vspace{2mm}

\hspace{21mm} $ \leq \max \{-b \mid (a,b) \in {\Omega}_R \} + a
(S_2) $

\vspace{2mm}

\hspace{21mm} $ = a_*^2 (R). \; \; $ $\B$

\medskip
Now we are ready to exhibit some families of ideals such that the
diagonal functor and the local cohomology functor commute whenever
we take diagonals large enough.

\begin{ex}
\label{ED11} {\rm Let $I$ be an equigenerated ideal in a ring $A$
with $a_*(A) < 0$ (for instance, we may take $A=k[{\x}]$). Set
$r(I)$ the reduction number of $I$ and assume that ${\F}$ is
Cohen-Macaulay with negative a-invariant. Note that $a(F)<0$ is
equivalent to $r(I) <l(I)$ by Lemma \ref{ED9}. This class of
ideals includes:
\begin{enumerate}
\item ideals $I$ with reduction number $r(I) = 0$ (for instance,
complete intersection ideals and ideals of linear type).
\item  $\frak m$-primary ideals with $r(I) \leq 1 < l(I)$ in a
Cohen-Macaulay ring $A$ (see \cite{HS}) .
\item  equimultiple ideals with $r(I) \leq 1 < l(I)$ in a
Cohen-Macaulay ring
 $A$ (see \cite{Sh}).
\item generically complete intersection ideals with ad$(I) = 1$,
$r(I)\leq 1 < l(I)$ in a Cohen-Macaulay ring A (see \cite{CZ}).
\end{enumerate}

For these families of ideals, we have that the Rees algebra has a
good resolution according to Lemma \ref{ED7}. Then we have graded
isomorphisms
$$ H^q_{m} (\Kk) \cong H^{q+1}_{\cal M} (R)_{\Delta},$$
for any $q \geq 0$ and $c \gg e \gg 0$ by Corollary \ref{W64}.
Therefore, we have that for large diagonals $\depth (\Kk) \geq
\depth (R)-1$. In particular, if the Rees algebra is
Cohen-Macaulay then its large diagonals will be also
Cohen-Macaulay. }
\end{ex}

\medskip
\begin{rem}
\label{ED14} {\rm Recall that the form ring $G_A(I)$ of an ideal
$I$ in $A$ is
$$G=G_A(I)= \bigoplus_{n \geq 0}{I^n}/{I^{n+1}} =R_A(I)/ IR_A(I).$$
If $I$ is a homogeneous ideal, the form ring has a natural
bigrading by means of $G_{(i,j)}=(I^j/I^{j+1})_i$. We can get for
the form ring similar results to the ones obtained for the Rees
algebra. For instance, for an equigenerated ideal $I$ we have that
$G_A(I)$ has a good resolution if and only if $a_*(A/I) < 0$,
$a_*(F_{\frak m}(I)) < 0$  and it holds $a_*({\F}) \leq a_*^2
(G_A(I))$. }
\end{rem}

\begin{rem}
\label{ED13}
{\rm For an equigenerated ideal $I$ of $A$, note that we can
recover several relationships between $r (I)$, $l(I)$ and $a^2(G)$
proved with more generality in \cite{AHT}. By applying the
diagonal functor to the minimal bigraded free resolution of
$R_A(I)$ or $G_A(I)$, we obtain the minimal graded free resolution
of ${\F} = k[I_d]$, and so $ a_*({\F}) \leq  a_*^2(R_A(I))$ and $
a_*({\F})  \leq  a_*^2(G_A(I))$. Now, according to Lemma
\ref{ED9}, given $J$ an arbitrary minimal reduction of $I$ we have
$r_J(I) - l(I) \leq a_*^2(R_A(I))$ and the same formula for
$G_A(I)$. In particular, if $R_A(I)$ is CM we get $r_J(I) \leq l
(I) + a^2(R_A(I)) \leq l(I) -1$. We can also obtain that if
$R_A(I)$ is CM then
 ${\rm reltype} (I) \leq \mu (I)-1$.
}
\end{rem}

\medskip
Once we have studied the equigenerated case, in the rest of the
section we consider the Rees algebra of an arbitrary homogeneous
ideal $I$ of $A$. In the case where $A$ is the polynomial ring it
was conjectured in \cite{CHTV} that if the Rees algebra is
Cohen-Macaulay there exist large diagonals which are
Cohen-Macaulay. Note that for equigenerated ideals this follows
from Proposition \ref{ED7} and Lemma \ref{ED10}. Next we give a
full answer to this conjecture.

\begin{thm}
\label{ED120} Let $I$ be a homogeneous ideal of the polynomial
ring $A=k[{\x}]$. If $R_A(I)$ is Cohen-Macaulay, then $R_A(I)$ has
a good resolution. In particular, $\Kk$ is Cohen-Macaulay for $c
\gg e \gg 0$.
\end{thm}

{\pf} Let us consider the bigraded minimal free resolution of $R$
over $S$:
$$ 0 \to D_{r-1} \to \dots \to D_1 \to D_0=S \to R \to 0.$$
The first morphism in this resolution maps each $X_i$ to $X_i$, so
we immediately get that any shift $(a,b) \in \Omega_1$ satisfies
$b<0$. Then, by Lemma \ref{B4} jointly with Remark \ref{B6}, we
get that for any $p \geq 1$ and any shift $(a,b) \in \Omega_p$ we
have $b<0$. On the other hand, since $R$ is Cohen-Macaulay,
$a_*^2(R)=a^2(R)$. Now, by Lemma \ref{B1} we have $a^2(R)=a(R_2)$,
where $R_2 = \bigoplus_j (\bigoplus_i (I^j)_i) =\bigoplus_j I^j$
is nothing but the Rees algebra with the usual ${\Bbb Z}$-
grading, and so $a(R_2)=-1$. By then applying Theorem \ref{B77},
for any $(a,b) \in \Omega_R$ we have
$$-b \leq a^2(S)-a_*^2(R)=r-1 <r,$$
so by Proposition \ref{W9} the Rees algebra has a good resolution.
Then the result follows from Corollary \ref{W64}. $\B$

\medskip
Assume that a Rees algebra is Cohen-Macaulay. Is there any
Cohen-Macaulay diagonal?  We know that this holds if $A$ is the
polynomial ring. The next result provides the obstruction for the
existence of such diagonals in the general case. Note that the
obstruction depends only on the ring $A$.

\begin{thm}
\label{ED121} If  $R_A(I)$ is Cohen-Macaulay, then the following
are equivalent:
\begin{itemize}
     \item[(i)] There exist $c, e$ such that $\Kk$
                 is Cohen-Macaulay.
     \item[(ii)] $H^i_{\frak m}(A)_0 = 0$ for all $i < \overline n$.
\end{itemize}
\end{thm}

{\pf} If $R$ is Cohen-Macaulay, then $a_*^2(R) = a^2(R)=-1$, so by
Proposition \ref{W11} we have $H^q_{{\fn} R}(R)_{(0,0)}=0$ for all
$q$. Then, according to Proposition \ref{ED2} and Proposition
\ref{W11}, $H^q_m(\Kk)_0 = H^q_{{\fm} R}(R)_{(0,0)}=
H^q_{\fm}(A)_0$. Now the statement follows from Corollary
\ref{W63}. $\B$

\chapter{Cohen-Macaulay coordinate rings of blow-up schemes}
\typeout{Cohen-Macaulay coordinate rings of blow-up schemes}

\bigskip
\markboth{CHAPTER III. COHEN-MACAULAY BLOW-UP SCHEMES}
{INTRODUCTION}

\medskip
After introducing in the previous chapters the basic tools we will
need along this work, we are now ready to study in detail the
Cohen-Macaulayness of the coordinate rings of blow-ups of
projective varieties. Let $k$ be a field, and let $Y$ be a closed
subscheme of ${\Bbb P}^{n-1}_k$ with coordinate ring $A=
k[{\x}]/K$, where $K \subset k[{\x}]$ is a homogeneous ideal.
Given $I \subset A$ a homogeneous ideal, let denote by ${\cal
I}=\widetilde{I}$ the sheaf associated to $I$ in $Y = \Proj (A)$.
Let $X$ be the projective scheme obtained by blowing up $Y$ along
$\cal I$, that is, $X= {\cal P}roj (\bigoplus_{n \geq 0} {\cal
I}^n)$. If $I$ is generated by forms of degree $\leq d$, then
$(I^e)_c$ corresponds to a complete linear system on $X$ very
ample for $c \geq de+1$ which gives a projective embedding of $X$
so that $X \cong \Proj(\Kk) \subset \Bbb P^{N-1}_k$, where
$N=\dim_k(I^e)_c$ (see \cite[Lemma 1.1]{CH}).

\medskip
For a given homogeneous ideal $I\subset A$, we can consider the
Rees algebra ${\Rr}= \oplus_{j \geq 0} I^j$ of $I$ endowed with
the natural bigrading $R_A(I)_{(i,j)} =(I^j)_i$. By taking
diagonals $\Delta =(c,e)$ with $c \geq de+1$, we have that
$R_A(I)_\Delta= k[(I^e)_c]$. In Chapter 2 we used this fact to
study the existence of algebras $\Kk$ which are Cohen-Macaulay in
the case where the Rees algebra also has this property (see
Theorems \ref{ED120} and \ref{ED121}). Our aim in this chapter is
to get some general criteria for the existence of (at least) one
coordinate ring $\Kk$ with the Cohen-Macaulay property. In Section
3.2 we will give sufficient and necessary conditions to ensure
this existence by means of the local cohomology of $\Rr$ and the
sheaf cohomology $H^i({\X}, \cal O _{\X})$. This result will be
applied in Section 3.3 to exhibit several situations where we can
ensure the existence of Cohen-Macaulay coordinate rings for $X$.
We also give a criterion for the existence of Buchsbaum coordinate
rings, proving in particular a conjecture stated by A. Conca et
al. \cite{CHTV}.

\medskip
Once we have studied the existence of Cohen-Macaulay diagonals of
a Rees algebra, in Section 3.4 our aim will be to precise these
diagonals. This is a difficult problem which has only been
completely solved for complete intersection ideals in the
polynomial ring \cite[Theorem 4.6]{CHTV}. We will give several
criteria to decide if a given diagonal is Cohen-Macaulay, which
will allow us to recover and  extend the result on complete
intersection ideals as well as to determine the Cohen-Macaulay
diagonals for new families of ideals. If the Rees algebra is
Cohen-Macaulay, we can also determine a family of Cohen-Macaulay
diagonals. The section finishes by studying the coordinate rings
of the embeddings of the blow-up of a projective space along an
ideal of fat points.

\medskip
   The last section is devoted to study sufficient conditions for the
existence of a constant $f$ ensuring that {\Kk} is Cohen-Macaulay
for any $c \geq ef$ and $e>0$, a question that has been treated by
S.D. Cutkosky and J. Herzog in \cite{CH}. The main result shows
that this holds for homogeneous ideals in a Cohen-Macaulay ring
$A$ whose Rees algebra is Cohen-Macaulay at any $\fp \in \Proj
(A)$.

\bigskip
\section{The blow-up of a projective variety}
\label{G} \markboth{CHAPTER III. COHEN-MACAULAY BLOW-UP SCHEMES}
{THE BLOW-UP OF A PROJECTIVE VARIETY}

\medskip
    From now on in this chapter we will have the following assumptions.
Let $k$ be a field and $A$ a noetherian graded $k$-algebra
generated in degree 1. Then $A$ has a presentation $A= k[{\x}]/ K
= k[{\xx}]$, where $K$ is a homogeneous ideal in the polynomial
ring $k[{\x}]$ with the usual grading. We will denote by $\fm$ the
graded maximal ideal of $A$. Let $Y$ be the projective scheme $
{\rm Proj}(A) \subset \Bbb P_k ^ {n-1}$. Let $I$ be a homogeneous
ideal not contained in any associated prime ideal of $A$, and let
${\cal I}$ be the sheaf associated to $I$ in $Y$. Then $\cal I$
can be blown up to produce the projective scheme $X ={\cal P}roj
(\bigoplus_{n \geq 0} {\cal I} \, ^n) $ together with a natural
morphism $\pi : X \to Y$. Let us recall the construction of the
${\cal P}roj$ of a sheaf of graded algebras $\cal R$ over a scheme
$Y$ (see \cite[Chapter II, Section 7]{Ha}). For each open affine
subset $U=\Spec(B)$ of $Y$, let ${\cal R}(U)$ be the graded
$B$-algebra ${\Gamma}(U,{\cal R} \vert U)$. Then we can consider
${\rm Proj}({\cal R}(U))$ and its natural morphism $\pi \vert U:
{\rm Proj}({\cal R}(U)) \rightarrow U$. These schemes can be glued
to obtain the scheme ${\cal P}roj({\cal R})$ with the morphism
$\pi:{\cal P}roj({\cal R}) \rightarrow Y$ such that for each open
affine $U \subset Y$, $\pi^{-1}(U) \cong {\rm Proj}({\cal R}(U))$.

\medskip
Assume that $I$ is generated by forms $f_1, \dots, f_r$ in degrees
$d_1, \dots, d_r$ respectively. Let $d= \max \{d_1, \dots, d_r\}$.
For any $c \geq d+1$, let us consider the invertible sheaf of
ideals ${\cal L}= {\cal I}(c) {\cal O}_X$. We are going to show
that $\Ll$ defines a morphism of $X$ in a projective space
$\varphi: X \to \Bbb P^{N-1}_k$ which is a closed immersion so
that $X \cong \Proj (k[I_c])$. Since the blow-up of $Y$ along
${\cal I}^e$ is isomorphic to $X$, we will also have $X \cong
\Proj (\Kk)$ for any $c \geq de+1$. For that, we are going to
follow the proof of \cite[Lemma 1.1]{CH}. First of all, notice
that we have an affine cover of $X$ by considering the set
 $\{ U_{ij} \mid 1 \leq i \leq n, 1 \leq j \leq r \}$,
where $U_{ij}= \Spec (R_{ij})$, and
 $$R_{ij}=\bigg( k \bigg[
 \frac{X_1}{X_i}, \dots, \frac{X_n}{X_i} \bigg] /K_i \bigg)
\bigg[ \frac{f_1 X_i^{d_j-d_1}}{f_j}, \dots, \frac{f_r
X_i^{d_j-d_r}}{f_j} \bigg] \;.$$

\noindent Furthermore, $\Gamma (U_{ij},{\cal I}(c) \, {\cal O}_X)
= f_j x_i^{c-d_j} R_{ij}$. Since $f_j x_i^{c-d_j} \in I_c$ and
$I_c \subset \Gamma(X, {\cal I}(c) \, {\cal O}_X)$, we have that
${\Ll}= (I_c) \, {\cal O}_X$.

\medskip
$I_c$ is a $k$-vector space generated by the elements $s$ of the
type $s = f_j x_1^{l_1} \dots x_n^{l_n}$ with degree $c$, that is,
such that $d_j + l_1 + \dots+ l_n =c$. By considering $X_s = \{ P
\in X \mid s_P \not \in {\fm}_P {\Ll}_P \}$ with $s \in I_c$, we
have an open covering of $X$. Since $c>d$, there exists some $i$
with $l_i >0$, so denoting by $u = (\frac{x_1}{x_i})^{l_1} \dots
(\frac{x_n}{x_i})^{l_n}$ we have that
$$X_s = \Spec ((R_{ij})_{u}) = \Spec \bigg (R_{ij}
\bigg[ \bigg( \frac{x_i}{x_1} \bigg)^{l_1} \dots \bigg(
\frac{x_i}{x_n} \bigg)^{l_n} \bigg] \bigg) $$ is an open affine.

\medskip
 Set $N= \dim_k I_c$. Let $\Bbb P^{N-1}_k = \Proj (k[\{Z_s\}_{s \in
\Lambda}])$, where $\Lambda$ is a $k$-basis of $I_c$, and $V_s =
D_+(Z_s) \subset \Bbb P^{N-1}_k$. The $k$-linear maps defined by
$$\begin{array}{ccc}
\Gamma(V_s,{\cal O}_{V_s})= k[T_{\overline s} : \overline{s}\not =
s] \longrightarrow
\Gamma(X_s, {\cal O}_{X_s}) \\
\hspace{25 mm} T_{\overline s} \hspace{5mm} \mapsto \hspace{5mm}
\frac{\overline s}{s}
\end{array}$$
are epimorphisms which define morphisms of schemes $X_s
\rightarrow V_s$. By gluing them, we get a closed immersion
$\varphi : X \rightarrow \Bbb P^{N-1}_k$ so that $X \cong
\Proj(k[I_c])$ by \cite[Proposition II.7.2]{Ha}.

\medskip
Let ${\cal L} = \widetilde{I} {\cal O}_X$, ${\cal M }= \pi^* {\cal
O}_Y (1)$, so that $(I^e)_c {\cal O}_X = {\cal L}^e \otimes {\cal
M}^c$. Then, the classical Serre's exact sequence allows to relate
the local cohomology of the rings $\Kk$ and the global cohomology.

\begin{rem}
\label{G1} \cite[Lemma 1.2]{CH} {\rm There is an exact graded
sequence
   $$ 0 \to H^0_m(\Kk)   \to \Kk \to
   \bigoplus_{s \in \Bbb Z} {\mit \Gamma} ({\X}, {\cal L}^{es}
\otimes {\cal M}^{cs}) \to H^1_m (\Kk) \to 0$$
    and isomorphisms
$$H^i_m (\Kk) \cong \bigoplus_{s \in \Bbb Z} H^{i-1}
      ({\X}, {\cal L}^{es} \otimes {\cal M}^{cs})$$ for $i>1$.
In particular, looking at the homogeneous component of degree 0 we
get the exact sequence
   $$ 0 \to H^0_m(\Kk)_0   \to k \to
   {\mit \Gamma} ({\X}, {\cal O}_{\X})
    \to H^1_m (\Kk)_0  \to 0$$
    and isomorphisms  $H^i_m (\Kk)_0 \cong H^{i-1} ({\X}, {\cal O}_{\X})$
for $i>1$.}
 \end{rem}

\medskip
For a homogenous ideal $I$ of $A$, let us consider the Rees
algebra $R= R_A(I) = \bigoplus_{n \geq 0}  I^n t^n \subset A[t]$
of $I$ endowed with the natural bigrading given by $R_{(i,j)} =
(I^j)_i$. Then, by taking a diagonal $\Delta =(c,e)$ with $c \geq
de +1$, we have that $R_A(I)_\Delta = \bigoplus_{s \geq 0}
(I^{es})_{cs} = k[(I^e)_c]$. The natural inclusion $\Kk = R_\Delta
\hookrightarrow R$, gives the isomorphism of schemes
$\Proj^2(R_A(I)) \cong \Proj (\Kk)$. Summarizing, we have:

\begin{prop}
\label{D10} Let $X$ be the blow-up of $Y= \Proj (A)$ along ${\cal
I}= \widetilde{I}$, where $I$ is a homogeneous ideal of $A$
generated by forms of degree $\leq d$. For any $c \geq de+1$, we
have isomorphisms of schemes
$$ X \cong \Proj^2 (R_A(I)) \cong \Proj (\Kk).$$
\end{prop}

\medskip
In Chapter 2, Section 3, we have followed an algebraic approach to
study the local cohomology modules of the rings $\Kk$ in terms of
the local cohomology of the Rees algebra and the diagonal functor.
Next we give a new approach to these modules by using sheaf
cohomology.

\medskip
Notice that from Remark \ref{G1} we may determine the local
cohomology modules of the $k$-algebras $\Kk$ by means of the
cohomology modules $H^{i}({\X}, {\cal L}^{es} \otimes {\cal
M}^{cs})$. On the other hand, we can get some information about
these modules from the Leray spectral sequence
$$ E_2^{i,j}= H^i ({\Y}, R^j \pi_* ({\cal L}^{es} \otimes {\cal
M}^{cs})) \Rightarrow H^{i+j}({\X}, {\cal L}^{es} \otimes {\cal
M}^{cs}),$$ and the vanishing of the higher direct-image sheaves $
R^j \pi_* ({\cal L}^{es} \otimes {\cal M}^{cs})$. First of all,
let us study the vanishing of $ R^j \pi_* ({\cal L}^{es})$.

\begin{thm}
\label{G2} Set $e_0= \max \{\, a_*(R_{A_{\fp}}(I_{\fp})): \fp \in
\Proj(A) \,\}$. For any $e > e_0$, $j>0$, $R^j \pi_* {\cal L}^e =
0$ and $\pi_* {\cal L}^e = \widetilde {I^e}$.
\end{thm}

{\pf} Let us denote by
     $A_i= A_{(x_i)}$,
     $I_i= I_{(x_i)}$ and
     $R_i =R_{A_i}(I_i) = A_i [I_i t]$.
Note that by defining  $Y_i= Y - V_+(x_i)= D_+(x_i) \cong
\Spec(A_i)$, we have that $\{Y_i : 1 \leq i \leq n\}$ is an open
affine cover of $Y$. Then, given $j>0, e>0$, $R^{j} \pi_* {\Ll}^e
= 0$ if and only if $(R^{j} \pi_* {\Ll}^e) \mid Y_i = 0$ for all
$i$. Denoting by $X_i = \pi^{-1} Y_i = \Proj(R_i)$, by
\cite[Corollary III.8.2 and Proposition III.8.5]{Ha} we have that
for $j >0$
$$R^j \pi_* ({\Ll}^e) \mid Y_i =
 R^j \pi_* ({\Ll}^e \mid X_i) =
  H^j (X_i, {\Ll}^e \mid {{\X}_i}) \,\,\widetilde{} =
 H^{j+1}_{(R_i)_+} ( (I_i)^e R_i)_0 \,\,\widetilde{} \,\,\,.$$

 From the graded exact sequence
 $$ 0 \to  (I_i)^e R_i (-e) \to R_i \to \bigoplus_{q<e}(I_i)^q \to 0
\;,$$
 it follows that
 $H^{j+1}_{(R_i)_+} ( (I_i)^e R_i)_0 =
 H^{j+1}_{(R_i)_+} ( (I_i)^e R_i (-e))_e =
 H^{j+1}_{(R_i)_+} (R_i)_e$.
 Similarly, $\pi_* {\cal L}^e = \widetilde{I^e}$ if
 $H^{0}_{(R_i)_+} (R_i)_e  = H^{1}_{(R_i)_+} (R_i)_e = 0$ for all $i$.
Therefore, we have reduced the problem to prove that
$H^{j}_{(R_i)_+} (R_i)_e = 0$ for all $i, j$ if $e >e_0$.

Set $\overline R_i =R_{A_{x_i}}(I_{x_i})$. We can think $\overline
R_i$ as a ${\Bbb Z}$-graded ring with
deg$(\frac{\;x_j\;}{x_i^m})=1-m$, deg$(\frac{\;f_j t
\;}{x_i^m})=d_j-m$. Note that with this grading we have
$[{\overline R}_i]_0 = R_i$ and $\frac{x_i}{1}$ is an invertible
element in $\overline R_i$ of degree 1. Then we may define the
graded isomorphism
$$\begin{array}{lll}
R_i [T, T^{-1}] & \stackrel{\psi}\longrightarrow  &  \overline R_i \\
\;\;\;\;\;\;     T    &  \mapsto   &   \frac{x_i}{1}  \; ,\\
 \end{array} $$
where $\psi \vert R_i = id$ and deg$(T)=1$. Since $R_i
\hookrightarrow \overline R_i$ is a flat morphism, we have that
$$H^{j}_{(\overline R_i)_{+}} (\overline R_i)=
H^{j}_{(R_i)_{+}} (R_i) \otimes_{R_i} {\overline R_i} =
H^{j}_{(R_i)_{+}} (R_i) [T, T^{-1}], $$ so that $H^{j}_{(\overline
R_i)_+} (\overline R_i)_e= H^{j}_{(R_i)_+} (R_i)_e \, [T,
T^{-1}]$. Therefore, it suffices to prove that $H^{j}_{(\overline
R_i)_+} (\overline R_i)_e = 0$ for all $i, j$ if $e >e_0$.

 Given a homogeneous prime $\frak q \in \Spec (A_{x_i})$, we have that
$\frak q = \fp A_{x_i}$ with $\fp \in \Proj (A)$. Localizing
$\overline R_i$ at $\frak q$, we have that $ \overline R_i \otimes
(A_{x_i})_{\frak q} = R_{A_\fp}(I_{\fp})$. Denoting by $B
=R_{A_\fp}(I_{\fp})$, note that $B$ is a standard graded ring
whose homogeneous component of degree $0$ is the local ring
$A_\fp$. So $B$ has a unique homogeneous maximal ideal $\frak n$,
with $\frak n = \fp A_\fp  \oplus B_{+}$. Since $H^j_{\frak
n}(B)_e = 0$ for all $j \geq 0$ and $e>e_0$, according to
\cite[Lemma 2.3]{H} we also have $H^j_{B_+} (B)_e = 0$ for all $j
\geq 0$ and $e > e_0$. Therefore,
$$[H^{j}_{(\overline R_i)_+} (\overline R_i)_e]_{\frak q} =
[H^{j}_{(\overline R_i)_+} (\overline R_i)_{\frak q}]_e =
H^{j}_{B_+} (B)_e =0.$$ Hence $(H^j_{(\overline R_i)_+} (\overline
R_i)_e)_\frak q = 0$ for any homogeneous ideal $\frak q \in
\Spec(A_{x_i})$, and we conclude $H^j_{(\overline R_i)_+}
(\overline R_i)_e= 0$ for $j \geq 0$ and $e>e_0$. $\Box$

\medskip
\begin{cor}
\label{G02} Assume that $R_{A_{\fp}}(I_{\fp})$ is Cohen-Macaulay
for any $ \fp \in \Proj(A)$. Then, for any $e \geq 0$, $j>0$, $R^j
\pi_* {\cal L}^e = 0$ and $\pi_* {\cal L}^e = \widetilde {I^e}$.
\end{cor}

\medskip
 Given a homogeneous ideal $I$ of $A$, let us denote by
$I^* = \{ f \in A \mid \frak m^k f \subset I \,\,{\rm for}\,{\rm
some}\,k \}$ the saturation of $I$. Note that $H^0_{\frak m} (A/I)
= I^* /I$. Next we use Theorem \ref{G2} to relate the local
cohomology modules of the $k$-algebras $\Kk$ and the local
cohomology of the powers of the ideal (compare with Corollary
\ref{ED4}).

\begin{cor}
\label{G3} Let $I$ be a homogeneous ideal of $A$, and $e_0= \max
\{\, a_*(R_{A_{\fp}}(I_{\fp})): \fp \in \Proj(A) \,\}$. For any $c
\geq de+1$, $e>e_0$, $s>0$, there is an exact sequence
$$0 \to H^0_{m} (\Kk)_s
\to (I^{es})_{cs} \to (I^{es})^*_{cs} \to H^1_{m}(\Kk)_s \to 0$$
and isomorphisms $H^i_{m} (\Kk)_s \cong H^i_{\frak
m}(I^{es})_{cs}$ for $i>1$.
\end{cor}

{\pf} By the Leray spectral sequence we have
$$ H^i ({\Y}, R^j \pi_* ({\cal L}^{es} \otimes {\cal M}^{cs}))
\Rightarrow H^{i+j}({\X}, {\cal L}^{es} \otimes {\cal M}^{cs}).$$
On the other hand, by the Projection formula \cite[Exercise
III.8.3]{Ha} and Theorem \ref{G2}, we get that for $e >e_0$,
$s>0$,
$$ \pi_* ({\cal L}^{es} \otimes {\cal M}^{cs}) =
 \pi_* ({\cal L}^{es}) \otimes {\cal O}_{\Y}(cs) =
\widetilde{I^{es}} (cs)$$
 $$R^j \pi_* ({\cal L}^{es} \otimes {\cal M}^{cs}) =
 R^j \pi_* ({\cal L}^{es}) \otimes {\cal O}_{\Y}(cs) = 0,  \,\,
 {\rm for} \,\, {\rm all} \, j>0.$$
  Therefore we may conclude that for any $i \geq 1$,
$H^i({\X},{\cal L}^{es} \otimes {\cal M}^{cs}) =
 H^i({\Y}, \widetilde{I^{es}}(cs)) =
 H^{i+1}_{\frak m}(I^{es})_{cs}$, and
$\Gamma( X, {\Ll}^{es} \otimes {\M}^{cs}) =
 \Gamma( Y, \widetilde{I^{es}} (cs)) =
 (I^{es})^*_{cs}$. Now the result follows from Remark \ref{G1}.
$\Box$

\medskip
We have shown that the positive components of the local cohomology
modules of the rings $\Kk$ are closely related to the local
cohomology modules of the powers of the ideal $I^e$. Next we want
to study the negative components. In the case where $X$ is
Cohen-Macaulay, we will express them by means of the local
cohomology of the canonical module of the Rees algebra. Recall
that the canonical module of the Rees algebra is defined in the
category of bigraded $S$-modules, so let us write $K_R=
\bigoplus_{(i,j)} K_{(i,j)}$. For any integer $e$, we denote by
$K^e= (K_R)^e =\bigoplus_i K_{(i,e)}$. Then we have

\begin{prop}
\label{GG01} Assume that $X$ is an equidimensional Cohen-Macaulay
scheme. For any $c \geq de+1$, $e>a_*^2(K_R)$, $s> 0$, $1 \leq i <
\nn $,
$$H^i_{m}(\Kk)_{-s} \cong H^{\nn -i+1}_\fm (K^{es})_{cs}.$$
\end{prop}

{\pf} By Serre's duality and Remark \ref{G1}
 $$H^i_m(\Kk)_{-s} \cong
  H^{i-1}(X, {\cal L}^{-es} \otimes {\cal M}^{-cs}) \cong
   H^{\nn-i} (X, w_X \otimes {\cal L}^{es} \otimes {\cal M}^{cs}).$$

\noindent Then, by taking  $c \geq de+1$, $e>a^2_*(K_R)$, $i <
\nn$ we get

\vspace{3mm}

\hspace{10mm} $H^i_m (\Kk)_{-s} =  H^{\nn
+1-i}_{R_+}(K_R)_{(cs,es)}$

\vspace{3mm}

\hspace{37mm} $=H^{\nn +1-i}_{\fm}(K^{es})_{cs}, \; \; {\rm by} \;
{\rm Proposition} \; \ref{W11}$. $\B$

\medskip
\begin{rem}
\label{GG3} {\rm   According to \cite[Theorem 2.1]{HHK} we can
also express the negative components of the local cohomology of
the rings $\Kk$ by means of the local cohomology of their
canonical modules whenever $X$ is Cohen-Macaulay equidimensional.
In this case, there are also isomorphisms $H^i_{m}(\Kk)_{-s} \cong
H^{\nn + 1 -i}_{m}(K_{R_\Delta})_s$ for any $s> 0$ and  $1 \leq i
< \nn $.

But notice that $[K_{R_\Delta}]_j \cong [(K_R)_\Delta]_j$ for any
$j \gg 0$, so we immediately get
$$H^{i}_{m}(K_{R_\Delta})_s \cong
H^{i}_{m}((K_R)_\Delta)_s \cong H^i_{R_+}(K_R)_{(cs,es)} \cong
H^{i}_\fm (K^{es})_{cs}$$ for any $s \geq 0$, $i>1$.

}
\end{rem}

\bigskip
\section{Existence of Cohen-Macaulay coordinate rings}
\label{J} \markboth{CHAPTER III. COHEN-MACAULAY BLOW-UP SCHEMES}
{EXISTENCE OF COHEN-MACAULAY COORDINATE RINGS}

\medskip
   Our aim in this section is to find necessary and sufficient
conditions for the existence of integers $c,e$, with $c \geq
de+1$, such that the ring $\Kk$ is Cohen-Macaulay. Before proving
our main result, we need two previous lemmas. The first one may be
seen as a Nakayama's Lemma adapted to our situation, and in fact
it is just Lemma \ref{F12} for the case $r=2$.

\begin{lem}
\label{J21}
         Let $L$ be a finitely generated bigraded $R$-module and
$m$ an integer such that $R_+^m  L$ =$0$. Then,
 there exist integers $q_0, t$ such that $L_{(p,q)}$=$0$ for all $p> dq
+ t$, $q > q_0$.
\end{lem}

\medskip
The second lemma provides restrictions on the local cohomology
modules of the Rees algebra whenever $X$ is Cohen-Macaulay.

\begin{lem}
\label{J22} If $X$ is Cohen-Macaulay equidimensional, then there
are integers $q_0$, $t$ such that
   $H^i_{\cal M} (R_A(I))_{(p,q)} = 0$ for all $i < \overline
     n+1$, $q < q_0$ and $p < dq + t$.
\end{lem}

{\pf} Let ${\p} \in X$. Then $R_+ \not \subset {\p}$ and so there
exist $i,j$ such that $x_i \not \in {\p}$, $f_j t \not \in {\p}$.
Denote by $R_{<{\p}>}$ = $T^{-1}R$, where $T$ is the
multiplicative system consisting of all homogeneous elements of
$R$ which are not in ${\p}$. Note that $R_{(\p)} =
[R_{<{\p}>}]_{(0,0)}$. Furthermore, $\frac{x_i}{1}$ and $\frac{f_j
t}{x_i^{d_j}}$ are invertible elements in $R_{<{\p}>}$ with
$$ \deg \frac{x_i}{1}= (1,0),
 \hspace{2mm} \deg \, \, \frac{\;\;f_j t\;\;}{x_i^{d_j}} \, = (0,1). $$
Then we may define a bigraded isomorphism $\psi$:
 $$  \begin{array} {lll}
 R_{({\p})} [U,U^{-1},V,V^{-1}]
   & \stackrel{\psi} \longrightarrow  & R_{<{\p}>}  \\
  \hspace{25mm}  U   &      \longmapsto  &   \frac{x_i}{1} \\
  \hspace{25mm}  V   &    \longmapsto  &  \frac{\;\;f_j t\;\;}{x_i^{d_j}} \\
    \end{array} $$
where $\psi \vert R_{({\p})} = id$, and deg$(U)$ = $(1,0)$,
deg$(V)$ = $(0,1)$. Since $X$ is CM,  ${\cal O}_{{\X},{\p}} =
R_{({\p})}$ is CM and so $R_{<{\p}>}$ too. Then, localizing at
 $ {\p}  R_{<{\p} >}$, we have that $R_{\p}$ is CM.

 Now let ${\p} \in \Spec(R)$
and denote by ${\p}^*$ the ideal generated by the homogeneous
elements of ${\p}$. By \cite[Corollary 1.2.4]{GW2}, $R_{\p}$ is CM
if and only if $R_{{\p}^*}$ is CM, so we have that $R_{\p}$ is CM
for any prime ideal ${\p}$ such that $R_+ \not \subset {\p}$.
Localizing the Rees algebra $R$ at the homogeneous maximal ideal
${\cal M}$ we then have that $R_{\cal M}$ is a generalized
Cohen-Macaulay module with respect to $R_+ R_{\cal M}$ \cite[Lemma
43.2]{HIO}. Therefore there exists $m \geq 0$ such that $R_+^m
H^i_{\cal M} (R)$ = $0$ for all $i < \overline n + 1$. From the
presentation of $R$ as a quotient of the polynomial ring
$S=k[{\x}, {\y}]$, by Theorem \ref{B2} we get
   $$ H^i_{\cal M} (R)= \uExt^{n+r-i}_S (R, K_S)^{\vee }, $$
and so $R_{+}^{m} \uExt^{n+r-i}_S (R, K_S) = 0$ for $i <
{\overline n}+1$. By Lemma $\ref{J21}$ we then
 obtain that there exist integers $q_1, t_1$ such that
$\uExt^{n+r-i}_S (R, K_S)_{(p,q)} = 0$ for all $q > q_1$, $p > dq
+ t_1$ and $i < \overline n + 1$. The proof finishes by dualizing
again.   $\B$

\medskip
Now we may formulate the main result of this section.

\begin{thm}
\label{J1}
   The following are equivalent:
\begin{itemize}
\item[(i)]   There exist $c,e$ such that
$\Kk$ is Cohen-Macaulay.
\item[(ii)]
    \begin{itemize}
     \item[(1)]  There exist integers $q_0$, $t$ such that
            $H^i_{\cal M} (R_A(I))_{(p,q)}= 0$ for all $i < \overline
n+1$, $q < q_0$ and $p < dq + t$.
     \item[(2)] ${\mit \Gamma} ({\X}, \cal O _{\X})$ = $k$ and
              $H^i({\X}, \cal O _{\X})$=$0$ for $0 < i < \overline n -1$.
     \end{itemize}

\end{itemize}

In this case, ${\Kk}$ is Cohen-Macaulay for $c \gg 0$ relatively
to $e \gg 0$.
\end{thm}

{\pf}  If $(i)$ is satisfied, then the scheme $X = \Proj^2(R)
\cong \Proj (\Kk)$ is CM and equidimensional, so we have $(1)$ of
$(ii)$ by Lemma \ref{J22}. Furthermore, $H^{i}_{m}(\Kk)_0 = 0$ for
any $i < \overline n$ and then by using Remark \ref{G1} we get
$(2)$ of $(ii)$.

    Assume now that $(ii)$ is satisfied. We want to find a diagonal
$\Delta$ such that $R_{\Delta} =\Kk$ is CM. By Remark \ref{G1} and
$(2)$ of $(ii)$, we have that $H^{i}_{m}(R_{\Delta})_0 = 0$ for
any diagonal $\Delta$ and $i < \overline n$. On the other hand,
since $H^{i}_{\cal M} (R)$ are artinian modules there exists $p_1$
such that $H^{i}_{\cal M} (R)_{(p,q)} = 0$ for all $i$ and $p
>p_1$. Furthermore, by Corollary \ref{W63} there are positive
integers $e_0, \alpha$  such that for $e >e_0$, $c >de+\alpha$ we
have
  $$ H^{i}_{m}(R_{\Delta})_{j} \cong
  H^{i+1}_{\cal M}(R)_{(cj, ej)} , \, \forall i, \forall j \not = 0. $$
Now, let us consider $q_0,t$ given by $(1)$ of $(ii)$. Note that
we can assume that $q_0, t$ are negative. Then, by taking
diagonals $\Delta =(c,e)$ with $e > \max \{ e_0, -q_0 \}$, $c >
 \max \{ de+ \alpha, p_1, de-t \}$, we have that
$H^{i}_{m}(R_{\Delta})_{j} = 0$ for all $j$ and $i < \overline n$,
and
 therefore ${\Kk}$ are CM for all these $c, e$. $\Box$

\medskip
\begin{rem}
\label{J26} {\rm  Assume that $(A, {\frak m})$ is a noetherian
local ring and
  let $I \subset \frak m$, $I \not = 0$ be an ideal. Denote by $X =
\Proj (R_A(I))$ the blow-up of $\Spec (A)$ along $I$. Then, it was
proved by J. Lipman \cite[Theorem 4.1]{L} that there exists a
positive integer $e$ such that $R_A(I^e)$ is Cohen-Macaulay if and
only if $X$ is
  Cohen-Macaulay, ${\mit \Gamma}({\X}, {\cal O}_{\X}) = A$ and
  $H^i ({\X}, {\cal O}_{\X}) = 0$ for all $i > 0$. The following
corollary may be seen as a projective version of this result.}
\end{rem}

\begin{cor}
\label{J27}
     The following are equivalent:
  \begin{itemize}
 \item[(i)]  There exist $c,e$ such that $\Kk$ is
Cohen-Macaulay.
 \item[(ii)]   $X$ is a Cohen-Macaulay equidimensional scheme,
 ${\mit \Gamma}({\X}, {\cal O}_{\X}) = k$ and
  $H^i ({\X}, {\cal O}_{\X}) = 0$ for all $0 <i < \overline n -1$.
 \item[(iii)]   $X$ is Cohen-Macaulay equidimensional and
     $ H^{i}_{R_+} (R) _{(0,0)} = 0$  for  $i < \overline n$.
  \end{itemize}
\end{cor}

{\pf}  It is enough to note that we have an exact bigraded
sequence
$$0 \to H^0_{R_+} (R)  \to  R  \to
\bigoplus_ {(p, q)} {\mit \Gamma} ({\X}, {\cal O}_{\X} (p, q))
      \to  H^1_{R_+}(R)  \to  0 ,$$
 and  isomorphisms
 $H^{i+1}_{R_+} (R) \cong
 \bigoplus_{(p, q)} H^i ({\X}, {\cal O}_{\X}
(p, q))$  for $i >0$.             $\Box$

\medskip
We can also give sufficient and necessary conditions for the
existence of
 generalized Cohen-Macaulay or Buchsbaum diagonals of the Rees algebra,
in particular proving a conjecture of A. Conca et al. in
\cite{CHTV}.

\begin{prop}
\label{J28} The following are equivalent:
\begin{enumerate}
\item $H^i_{\cal M} (R_A(I))_{(p,q)} = 0$ for $ i < \nn+1$ and $p
\ll 0$ relatively to $q \ll 0$.
\item $\Kk$ is a generalized
Cohen-Macaulay module for $c \gg e \gg 0$.
\item
There exist $c,e$ such that $\Kk$ is generalized Cohen-Macaulay.
\item $\Kk$ is a Buchsbaum ring for $c \gg e \gg 0$.
\item There exist $c,e$ such that $\Kk$ is a Buchsbaum ring.
\item There exist integers $q_0$, $t$ such that $H^i_{\cal M}
(R_A(I))_{(p,q)} = 0$ for $i < \overline n+1$, $q < q_0$ and $p <
dq + t$.

\end{enumerate}
\end{prop}

{\pf}
 $(i) \Rightarrow (ii)$   Assume that $(i)$ is satisfied. By
Corollary \ref{W63}, we get $H^{i}_{m}(R_{\Delta})_s = 0$ for $c
\gg 0$ relatively to $e \gg 0$, $s \not = 0$ and $i < \overline
n$. So $\Kk$ is generalized CM for $c\gg e \gg 0$.

 $(ii) \Rightarrow (iii)$   Obvious.

 $(iii) \Rightarrow (vi)$   Let $\Delta$ be a diagonal such that
$R_{\Delta}$ is  generalized CM. Then $(R_{\Delta})_{\frak p}$ is
CM with $\dim (R_{\Delta})_{\frak p}+ \dim (R_{\Delta} /{\frak
p})= \dim R_{\Delta} $ for any ${\frak p} \in \Proj (R_{\Delta})$
by
 \cite[Lemma 43.3]{HIO}, and so
$X \cong \Proj(R_{\Delta})$ is CM and equidimensional. By using
Lemma \ref{J22} one obtains $(vi)$.

 $(vi) \Rightarrow (i)$   Obvious.

The implications $(i) \Rightarrow (iv)  \Rightarrow (v)
\Rightarrow (vi)$ may be proved similarly. $\Box$

\bigskip
\section{Applications}
\label{R} \markboth{CHAPTER III. COHEN-MACAULAY BLOW-UP SCHEMES}
{APPLICATIONS}

\medskip
In this section we show several situations in which we can ensure
the existence of Cohen-Macaulay coordinate rings for the blow-up
scheme by using Theorem \ref{J1}. First lemma provides  sufficient
conditions to have ${\mit \Gamma} ({\X}, \cal O _{\X})$ = $k$ and
$H^i({\X}, {\cal O}_{\X}) = 0$ for all $0 < i < \overline n-1$.

\medskip
\begin{lem}
\label{R2} Assume $a_*^2 (R) < 0$, $a_*(A) <0$. Then
  ${\mit \Gamma} ({\X}, {\cal O}_{\X})$ = $k$ and
  $H^i({\X}, {\cal O}_{\X}) = 0$ for $0 < i < \overline n - 1$.

\end{lem}

{\pf} Note that for any $i$, we have $H^{i}_ {\cal M} (R)_{(0,0)}=
H^{i}_{{\cal M}_2}(R)_{(0,0)}=0$ and $H^{i}_ {{\cal M}_1}
(R)_{(0,0)}= H^{i}_{\fm}(A)_{0}=0$ by Proposition \ref{W11}. Then,
from the Mayer-Vietoris exact sequence associated to ${\cal M}_1$
and ${\cal M}_2$, we get $H^{i}_ {R_+} (R)_{(0,0)}= 0$ for any
$i$, so we are done. $\Box$

\medskip
As an immediate consequence we get:

\begin{cor}
\label{R60} Suppose that $a_*^2(R) <0$, $a_*(A) < 0$. If $X$ is
Cohen-Macaulay equidimensional, then $\Kk$ is Cohen-Macaulay for
$c \gg e \gg 0$.
\end{cor}

\medskip
     It is known that there are smooth projective varieties
with no arithmetically Cohen-Macaulay embeddings (see for instance
\cite[Theorem 3.4]{M}). Next we exhibit a situation where this
implication is true.

\begin{prop}
\label{R5}
   Let $X$ be the blow-up of $\Bbb P^{n-1}_k$ along a closed
subscheme, where $k$ has $char k= 0$. Assume that $X$ is smooth or
with rational singularities. Then $X$ is arithmetically
Cohen-Macaulay.
\end{prop}

{\pf}  Let $\pi : X \rightarrow \Bbb P^{n-1}_k$
 be the blow-up morphism. From $\cite{K}$, we have that
 $\pi_* {\cal O}_{\X}$ = ${\cal O}_{\Bbb P^{n-1}_k}$ and
  $R^{j} \pi_* {\cal O}_{\X} = 0$
   for all $j > 0$. This implies that the Leray spectral sequence
$$E_2^{i,j} = H^{i}(\Bbb P^{n-1}_k, R^{j} \pi_* {\cal O}_{\X})
    \Longrightarrow
    H^{i+j}({\X}, {\cal O}_{\X})  $$
degenerates. Therefore we have
  ${\mit \Gamma} ({\X}, {\cal O}_{\X}) =
  {\mit \Gamma} ({\Bbb P^{n-1}_k}, {\cal O}_{\Bbb P^{n-1}_k})$ =
  $k$ and
  $H^i({\X}, {\cal O}_{\X}) =
   H^i ({\Bbb P^{n-1}_k}, {\cal O}_{\Bbb P^{n-1}_k}) = 0$ for all $i >
0$. Then the result follows from Corollary \ref{J27}.  $\Box$

\medskip
Assume that $A$ is Cohen-Macaulay. S.D. Cutkosky and J. Herzog
proved in \cite{CH} that the Rees algebra has Cohen-Macaulay
diagonals for locally complete intersection ideals and for ideals
whose homogeneous localizations are strongly Cohen-Macaulay
satisfying condition $({\cal F}_1)$. In the first case, observe
that  $R_{A_{\frak p}} ( I_ {\frak p})$ is Cohen-Macaulay for any
${\frak p} \in \Proj(A)$, while in the second one $R_{A_{(\frak
p)}} ( I_ {(\frak p)})$ is Cohen-Macaulay for any ${\frak p} \in
\Proj(A)$. Next we want to study those examples.

\medskip

\begin{prop}
\label{R8}
     The following are equivalent:
\begin{itemize}
\item [(i)]  $R_{A_{x_i}}(I_{x_i})$ is Cohen-Macaulay for all $1 \leq i
              \leq n$.
\item [(ii)]  $R_{A_{(x_i)}}(I_{(x_i)})$ is Cohen-Macaulay for all $1
             \leq i \leq n$.
\item[(iii)]  $R_{A_{\frak p}} ( I_ {\frak p})$ is Cohen-Macaulay
 for all ${\frak p} \in \Proj(A)$.
\item[(iv)]  $R_{A_{(\frak p)}} ( I_ {(\frak p)})$ is Cohen-Macaulay
 for all ${\frak p} \in \Proj(A)$.
\end{itemize}
\end{prop}

{\pf}
Set  $\overline R_i =R_{A_{x_i}}(I_{x_i})$, $R_i
=R_{A_{(x_i)}}(I_{(x_i)})$. We have already shown in the proof of
Theorem \ref{G2} that there exists an isomorphism $\overline R_i
\cong R_i [T, T^{-1}]$. Therefore, $R_i$ is CM if and only if
$\overline R_i$ is CM, and so the two first conditions are
equivalent.

 Now let us prove $(i) \Longleftrightarrow (iii)$. First assume
$(i)$, and for any prime ideal $\frak p \in \Proj (A)$ let us take
$x_i \not \in \frak p$. Note that $R_{A_{\frak p}} ( I_ {\frak p})
= \overline R_i \otimes_{A_{x_i}} (A_{x_i})_{\frak p}$, and so
$R_{A_{\frak p}} (I_{\frak p})$ is Cohen-Macaulay. Now assume
$(iii)$, and let us think $\overline R_i$ as a bigraded ring.
Then, to prove $(i)$, it is enough to show that for any
homogeneous prime ideal $Q \in \Spec (\overline R_i)$, we have
that $(\overline R_i)_Q$ is CM. Given such a $Q$, denote by
${\frak q} A_{x_i}= Q \cap A_{x_i}$, where ${\frak q} \in
\Spec(A)$ is a homogeneous prime which does not contain $x_i$,
that is, $\frak q \in \Spec(A_{(x_i)}) \subset \Proj(A)$. Then we
have $(\overline R_i)_Q = (R_{A_{\frak q}} (I_{\frak q}))_Q$, and
so $(\overline R_i)_Q$ is CM.

Finally, let us prove $(ii) \Longleftrightarrow (iv)$. Given
$\frak p \in \Proj (A)$ and $x_i \not \in \frak p$, let $\frak q =
\fp A_{x_i} \cap A_{(x_i)} \in \Spec (A_{(x_i)})$. Then,
$(A_{(x_i)})_{\frak q} = A_{(\fp)}$ and so $R_{A_{(\frak p)}}
(I_{(\frak p)}) = R_i \otimes_{A_{(x_i)}} (A_{(x_i)})_{\frak q}$.
Therefore, $(ii)$ implies $(iv)$. Now let us assume $(iv)$. Given
any homogeneous prime ideal $Q \in \Spec (R_i)$, let $\frak q = Q
\cap A_{(x_i)} \in \Spec(A_{(x_i)}) \subset \Proj(A)$, and let
$\fp \in \Proj (A)$ such that $ \fp A_{x_i} = \frak q [x_i,
x_i^{-1}]$. Since $(R_i)_Q = (R_{A_{({\fp})}} (I_{({\fp})}) )_Q$ ,
we have that $(R_i)_Q$ is Cohen-Macaulay, and so $R_i$ is CM.
$\Box$

\medskip
Now we can prove that the Rees algebra of a homogeneous ideal $I$
in a Cohen-Macaulay ring $A$ satisfying any of the equivalent
conditions above has Cohen-Macaulay diagonals. More generally, we
have:

\begin{thm}
\label{R7}
 Assume that $A$ is equidimensional and
 $R_{A_{\frak p}} ( I_ {\frak p})$ is Cohen-Macaulay
 for all ${\frak p} \in \Proj(A)$.
 Then {\Kk} is Cohen-Macaulay for $c \gg 0$ relatively to $e \gg 0$ if
and only if $H^i_{\frak m} (A)_0 = 0 $ for all $i < \nn$.
\end{thm}

{\pf} Given $P \in X$, let us denote by $\frak p = P \cap A \in
\Proj(A)$. Then $R_A (I)_{P} = (R_{A_{\frak p}}(I_{\frak p}))_{P}$
is CM and so $X$ is CM. Then, by Corollary \ref{G02} and the Leray
spectral sequence
$$E_2^{i,j} = H^{i}(Y, R^{j} \pi_* {\cal O}_{\X})
 \Longrightarrow  H^{i+j}({\X}, {\cal O}_{\X}) , $$
we get $H^j({\X}, \cal O _{\X})$ = $H^j({\Y}, \cal O _{\Y})$ =
$H^{j+1}_{\fm} (A)_0$ for $0 <j < \overline n -1$, and the exact
sequence  $0 \to H^0_{\fm}(A)_0 \to k \to {\mit \Gamma} ({\X},
{\cal O}_{\X}) = {\mit \Gamma} ({\Y}, {\cal O}_{\Y}) \to
H^1_{\fm}(A)_0 \to 0$, so we get the statement. $\Box$


\medskip
Denote by $E$ the exceptional divisor of the blow-up and by $w_E$
its dualizing sheaf. The last result of the section shows that
weaker assumptions on \cite[Lemma 2.1]{CH} are enough to ensure
that the rings $\Kk$ are Cohen-Macaulay for $c \gg e \gg 0$.

\begin{prop}
\label{R88} Suppose that $A$ is Cohen-Macaulay, $X$ is a
Cohen-Macaulay scheme, $\pi_* {\cal O}_E (m) =  \widetilde{I} ^m /
\widetilde{I} ^{m+1} $ for $m \geq 0$ and  $R^i \pi _*  {\cal O}
_E (m)$=$0$  for $i>0$ and $m \geq 0$. Then {\Kk} is
Cohen-Macaulay for $c >>0$ relatively to $e >>0$.
\end{prop}

 {\pf}
  $R^{i} \pi_* {\cal O}_{\X} = 0$  for $i>0$ and
 $\pi_* {\cal O}_{\X}$ = ${\cal O}_{\Y}$ by
 \cite[Lemma 2.1]{CH}. Then, from the Leray spectral sequence, we obtain
 $H^{i}({\X}, {\cal O}_{\X}) = H^{i}({\Y}, {\cal O}_{\Y}) =
 H^{i+1}_{\frak m} (A)_0 = 0$ for $0 < i < \overline n - 1$ and
 ${\mit \Gamma} ({\X}, {\cal O} _{\X})$ =
 $ {\mit \Gamma} ({\Y},{\cal O} _{\Y})$
 = $k$. Now, the proposition follows from Corollary \ref{J27}. $\Box$

\bigskip
\section{Cohen-Macaulay diagonals}
\label{S} \markboth{CHAPTER III. COHEN-MACAULAY BLOW-UP SCHEMES}
{COHEN--MACAULAY DIAGONALS}

\medskip
      Once we have studied the problem of the existence of
Cohen-Macaulay diagonals of a Rees algebra, now we would like to
study in more detail which diagonals are Cohen-Macaulay. This
question has been totally answered only for complete intersection
ideals in the polynomial ring \cite[Theorem 4.6]{CHTV}. Our
approach to this problem will give us criteria to decide if a
diagonal is Cohen-Macaulay, which will allow us to recover and
extend the result in \cite{CHTV} to any Cohen-Macaulay ring as
well as to precise the Cohen-Macaulay diagonals for new families
of ideals.

\medskip
The first criterion gives necessary and sufficient conditions for
a diagonal of a Cohen-Macaulay Rees algebra to have this property
in the case where $I$ is equigenerated. Namely,

\begin{prop}
\label{S611}
 Let $I \subset A$ be a homogeneous ideal generated by forms of degree
$d$ whose Rees algebra is Cohen-Macaulay. For any $c \geq de+1$,
$\Kk$ is Cohen-Macaulay if and only if
\begin{enumerate}
\item $H^i_{\frak m}(A)_0 = 0$, for $i < \nn$.
\item $H^i_{\frak m}(I^{es})_{cs} = 0$, for $i < \nn$, $s>0$.
\end{enumerate}
\end{prop}

{\pf} First, recall that the assumptions on the local cohomology
of $A$ are necessary and sufficient conditions for the existence
of Cohen-Macaulay diagonals (Theorem \ref{ED121}). Then, for any
$c \geq de+1$ and $i<\nn$, we have $H^i_{m}(\Kk)_0 = 0$ by Theorem
\ref{J1} and Remark \ref{G1}.

     On the other hand, by applying Proposition \ref{W11} and
Proposition \ref{W32}, for any $s<0$ we have:
$$H^q_{{\M}_1} ( R )_{(cs,es)} = H^q_{{\frak m}_1}(R^{es})_{cs} =0 $$
$$H^q_{{\M}_2} ( R )_{(cs,es)} = H^q_{{\frak m}_2}(R_{cs-des})_{es}=0 $$
because $R^{es}=0$ and $R_{cs-des}=0$. Therefore, for any diagonal
and any $i< \nn$, $s<0$, we get $H^i_{m}({\Kk})_{s} = 0$ according
to Proposition \ref{W2}. The statement, then, follows from
Corollary \ref{ED4}. $\B$

\medskip
   We may apply Proposition \ref{S611} to study in detail the following
example considered by L. Robbiano and G. Valla in \cite{RV}.

\begin{cor}
\label{S14} Let $\{L_{ij} \}$  be a set of $d \times (d+1)$
homogeneous linear forms of a polynomial ring $A=k[{\x}]$, $i= 1,
\dots, d$; $j = 1, \dots, d+1$, and let $M$ be the matrix
$(L_{ij})$. Let $I_t(M)$ be the ideal generated by the $t \times
t$ minors of $M$ and assume that $\h (I_t(M)) \geq d-t+2$ for $1
\leq t \leq d$. Denoting by $I = I_d(M)$, then $k[(I^e)_c]$ is
Cohen-Macaulay for any $c \geq de+1$.
\end{cor}

{\pf}
  The ideal $I$ is generated by $d+1$ forms of degree $d$, and the
Rees algebra has a presentation of the form
        $$R_A(I)= k[{\x},Y_1, \dots, Y_{d+1}]/(\phi_1, \dots, \phi_d),$$
where the polynomial ring $S = k[{\x},Y_1, \dots, Y_{d+1}]$ is
bigraded by $\deg (X_i)=(1,0)$, deg$(Y_j) =(d,1)$, and $\phi_1,
\dots, \phi_d$ is a regular sequence in $S$ with
deg$(\phi_l)=(d+1,1)$ (see the proof of \cite[Theorem 5.11]{RV}).
Then we have a bigraded minimal free resolution of the Rees
algebra $R_A(I)$ as $S$-module given by the Koszul complex
associated to $\phi_1, \dots, \phi_d$:
      $$ 0 \to F_d \to \cdots \to F_1 \to F_0=S \to R_A(I) \to 0,$$
with $F_p = S(-(d+1)p,-p)^{d \choose p}$. By applying the functor
$(\;)^e$ to this resolution, we have a graded free resolution of
$I^e$ over $A$:
$$ 0 \to F_p^e \to \cdots \to F_1^e \to F_0^e= S^e \to I^e \to
0,$$ with $p = \min \{ e, d \}$, $F_p^e= A(-p-de)^{\rho_p^e}$ for
certain $\rho_p^e \in \Bbb Z$, $\rho_p^e \not = 0$. The minimal
graded free resolution of $I^e$ is then obtained by picking out
some terms, but in any case it must have length $p$ because the
Hilbert series of $A/I^e$ is given by (\cite[Example 6.1]{RV})
 $$H_{A/I^e}(z) = \frac{1-\sum_{j=0}^d (-1)^j{d \choose j}
{d+e-j \choose e-j} z^{de+j}}{(1-z)^n}$$ (note that $z^{p+de}$
appears in the numerator). So by Theorem \ref{B77} we can compute
the $a_*$-invariant of $I^e$ and we get
    $$ a_*(I^e)
 = \cases{de + e -n & if $e<d$ \cr
de + d -n & if $e \geq d$. \cr}
 $$
On the other hand, since $n \geq \h (I_1(M)) \geq d+1$, we have
that $ d \leq n-1$, and so $a_*(I^e) < de$. Therefore, for any $c
\geq de+1$, $s \geq 1$, we have that $H^i_{\frak m}(I^{es})_{cs}=0
$ for all $i$. So $\Kk$ is Cohen-Macaulay by Proposition
\ref{S611}. Furthermore, note that $a(\Kk)<0$. $\Box$

\medskip
For arbitrary homogeneous ideals, we can also get a criterion for
the Cohen-Macaulayness of the diagonals by means of the local
cohomology of the powers of the ideal and the local cohomology of
the graded pieces of the canonical module of the Rees algebra.
More explicitly,

\begin{thm}
\label{S5} Let $I$ be a homogeneous ideal in $A$ generated by
forms of degree $\leq d$ whose Rees algebra is Cohen-Macaulay. For
any $c \geq de+1$, $\Kk$ is Cohen-Macaulay if and only if
\begin{enumerate}
\item $H^i_{\frak m}(A)_0 = 0$ for $i < \nn$.
\item $H^i_{\frak m}(I^{es})_{cs} = 0$ for $i < \nn$, $s>0$.
\item $H^{\nn -i+1}_{\frak m}(K^{es})_{cs} = 0$ for $1 \leq i < \nn$,
$s>0$.
\end{enumerate}
\end{thm}

{\pf} As in the proof of Proposition \ref{S611}, the assumptions
on the local cohomology of $A$ are necessary and sufficient
conditions for the existence of Cohen-Macaulay diagonals. Then we
have that $H^i_{m}(\Kk)_0 = 0$ for $i < \nn$.

  Since $R$ is Cohen-Macaulay, we have that $K_R$ is Cohen-Macaulay
with $a_*^2(K_R)=0$. Therefore, for any $s>0$, $1 \leq i < \nn$,
$H^i_m(\Kk)_{-s} =H^{\nn-i+1}_{\fm} (K^{es})_{cs}$ by Proposition
\ref{GG01}. Moreover, note that $H^0_m(\Kk)_{-s} =0$ for any
$s>0$. Then the statement follows from Corollary \ref{ED4}. $\B$

\medskip
      Let us denote by $G=G_A(I)= \bigoplus_{n \geq 0} I^n/I^{n+1}$ the
form ring of $I$ with the natural bigrading as a quotient of the
Rees algebra. For ideals whose form ring is quasi-Gorenstein, we
may get necessary and sufficient conditions for a diagonal to be
Cohen-Macaulay only in terms of the powers of the ideal.

\begin{cor}
\label{S61}
 Let $I$ be a homogeneous ideal in $A$ generated by forms of degree
$\leq d$. Assume that the Rees algebra is Cohen-Macaulay and the
form ring is quasi-Gorenstein. Let $a= -a^2(G_A(I))$, $b=-a(A)$.
For any $c \geq de+1$, $\Kk$ is Cohen-Macaulay if and only if
\begin{enumerate}
\item $H^i_{\frak m}(A)_0 = 0$, for  $i < \nn$.
\item $H^i_{\frak m}(I^{es})_{cs} = 0$, for $i < \nn$, $s>0$.
\item $H^i_{\frak m}(I^{es-a+1})_{cs-b} = 0$, for $ 1 < i \leq \nn$,
$s>0$.
\end{enumerate}
\end{cor}

{\pf} Under these assumptions $K_R$ has the expected form, that
is, there is a bigraded isomorphism
$$K_{R_A(I)} \cong \bigoplus_{(l,m),\, m \geq 1} [I^{m-a+1}]_{l-b}$$
(see Corollary \ref{K4} for more details about the isomorphism).
Then, for any $s > 0$ we have $K^{es} \cong I^{es-a+1}(-b)$, and
now the result follows from Theorem \ref{S5}. $\B$

\medskip
We can use Corollary \ref{S61} to precise the Cohen-Macaulay
diagonals for a complete intersection ideal of a Cohen-Macaulay
ring. In particular, this gives a new proof of \cite[Theorem
4.6]{CHTV} where the case $A= k[{\x}]$ was studied.

\begin{prop}
\label{S4} Let $I$ be a complete intersection ideal of a
Cohen-Macaulay ring $A$ minimally generated by $r$ forms of
degrees $d_1, \dots, d_r$.  Set $u = \sum_{i=1}^r d_i$. For any $c
\geq de +1$, $\Kk$ is Cohen-Macaulay if and only if $c > (e-1)d +
u + a(A)$.
\end{prop}

{\pf} From the bigraded isomorphism $G_A(I) \cong A/I [Y_1, \dots,
Y_r]$, with $\deg (Y_j)=(d_j,1)$, it is easy to prove by induction
on $e$ that the $a_*$-invariant of $A/I^e$ is:
$$a_*(A/I^e) = a(A/I^e) = (e-1)d+ u + a(A).$$
On the other hand,  we also have that $H^{\nn-r}_{\frak m}
(A/I^e)_s \not = 0$, for all
 $s \leq a(A/I^e)= (e-1)d+u+a(A)$ (see Lemma \ref{T55}).

Let $ \Delta = (c,e)$ be a diagonal with $c \geq de+1$. Since
$a^2(G)=-\h (I)= -r$, by Corollary \ref{S61} we have that $\Kk$ is
Cohen-Macaulay if and only if $cs> a(A/I^{es})$ and $cs + a(A) >
a(A/I^{es-r+1})$ for all $s>0$. The first condition is equivalent
to $(c-de)s >  u -d +a(A)$ for all $s>0$, that is, $c-de >  u -d
+a(A)$. The other one is equivalent to $(c-de)s >  u -dr$ for all
$s>0$, and this always holds because $u-dr \leq 0$. $\B$

\medskip
Until now we have given criteria to decide if a diagonal $\Kk$ is
Cohen-Macaulay once we know the local cohomology of the powers of
$I$, and the local cohomology of the graded pieces of the
canonical module of the Rees algebra. We will apply them in
Chapter 5, Section 2, after computing the local cohomology of the
powers of certain families of ideals.

\medskip
The following result shows the behaviour of the $a_*$-invariant
for the graded pieces of any finitely generated bigraded
$S$-module, so in particular for the powers of an ideal and the
pieces of the canonical module by applying it to the Rees algebra
and its canonical module respectively. This fact has been also
obtained independently by S.D. Cutkosky, J. Herzog and N. V. Trung
\cite{CHT} and V. Kodiyalam \cite{Ko2} by different methods (see
Chapter 5 for more details).

\begin{thm}
\label{S8} Let $L$ be a finitely generated bigraded $S$-module.
Then there exists $\alpha$ such that for any $e$
$$a_*(L^e) \leq de+ \alpha.$$
\end{thm}

{\pf} Let $e_0 = a_*^2(L)$. By Proposition \ref{W11}, $H^i_{{\cal
M}_2} (L)_{(c,e)}=0$ for $i \geq 0$, $e>e_0$. Then, by Proposition
\ref{W2} and Proposition \ref{W11}, we have that for any $c \geq
de+1$, $e>e_0$, $i \geq 0$, there are isomorphisms
$$H^i_{m} (L_{\Delta})_1 \cong
H^i_{{\cal M}_1} (L)_{(c,e)} \cong H^i_{\frak m}(L^{e})_{c}.$$ On
the other hand, from Corollary \ref{W63} there exist positive
integers $e_1$, $\alpha_1$ such that $H^i_m(L_\Delta)_s \cong
H^{i+1}_{\cal M}(L)_{(cs,es)}$ for $s \not = 0$, $e> e_1$, $c >de
+ \alpha_1$. Therefore, we have $H^i_{\frak m}(L^{e})_{c} = 0$ for
 $e>\max \,\{e_0, e_1 \}$, $c>de+\alpha_1$, $i \geq 0$. This proves the
statement. $\Box$

\medskip
Next we will show how to obtain a family of Cohen-Macaulay
diagonals from the bound on the shifts in the bigraded minimal
free resolution of the Rees algebra given by Theorem \ref{B77}. To
begin with, let us study the bigraded a-invariant of the Rees
algebra.

\begin{lem}
\label{S1}
\begin{enumerate}
\item  $a^1(R) \leq a(A)$.
\item  If $R$ is Cohen-Macaulay and $a^2(G)<-1$, then
       $a^1(R)= a(A)$.
\end{enumerate}
\end{lem}

{\pf} By setting $R_{++} = \bigoplus_{j>0} R_{(i,j)}$, we have the
following bigraded exact sequences:
   $$ 0  \rightarrow  R_{++} \rightarrow R \rightarrow A
   \rightarrow 0$$
   $$ 0 \rightarrow  R_{++}(0,1) \rightarrow R \rightarrow G
   \rightarrow 0 .$$
For each $(i,j)$, we get exact sequences:
   $$ ... \rightarrow
     H^{\nn}_{\cal M} (A)_{(i,j)}
   \rightarrow
      H^{\nn +1}_{\cal M}(R_{++})_{(i,j)}
   \rightarrow
     H^{\nn +1}_{\cal M}(R)_{(i,j)}
   \rightarrow
    0 \hspace{5mm}(1)$$
   $$ ... \rightarrow
   H^{\nn}_{\cal M}(G)_{(i,j)}
   \rightarrow
   H^{\nn +1}_{\cal M}(R_{++})_{(i,j+1)}
   \rightarrow
   H^{\nn +1}_{\cal M}(R)_{(i,j)}
   \rightarrow
    0 \hspace{5mm}(2)$$
Note that $A_{(i,j)} =0$ if $j \not =0$ and so
    $H^{\nn}_{\cal M}(A)_{(i,j)}=0$ if $j \not = 0$.

We want to determine $a^1(R)= \max \{i \mid \exists j: H^{\nn
+1}_{\cal M}(R)_{(i,j)} \not =0 \}$. Suppose $H^{\nn+1}_{\cal
M}(R)_{(i,j)} \not =0$. Since $a^2(R)=-1$, we have $j \leq -1$.
Then, from (2), we get
   $H^{\nn +1}_{\cal M}(R_{++})_{(i,j+1)}  \not =0$. If $j+1<0$, from
(1) we obtain  $ H^{\nn+1}_{\cal M}(R)_{(i,j+1)} \cong
H^{\nn+1}_{\cal M}(R_{++})_{(i,j+1)} \not =0$. By repeating this
argument, we obtain $ H^{\nn+1}_{\cal M}(R_{++})_{(i,0)} \not =0$
and, since $H^{\nn+1}_{\cal M}(R)_{(i,0)} =0$, from (1) we get
$H^{\nn}_{\fm} (A)_i= H^{\nn}_{\cal M}(A)_{(i,0)}  \not =0$. Thus
$i \leq a(A)$, and then it follows that $a^1(R) \leq a(A)$.

Assume now that $R$ is Cohen-Macaulay and $a^2(G)<-1$. From (2),
we have that if $H^{\nn + 1}_{\cal M}(R)_{(a(A),-1)} = 0$ then
$H^{\nn +1}_{\cal M}(R_{++})_{(a(A),0)} = 0$. Since $R$ is
Cohen-Macaulay, from (1) we get $H^{\nn}_{\fm}(A)_{a(A)} =
H^{\nn}_{\cal M}(A)_{(a(A),0)}=0$, which is a contradiction. $\B$

\medskip
\begin{rem}
\label{S11}
  {\rm  Note that in the proof of the Lemma \ref{S1} $(ii)$
      it is enough to assume
      $H^{\nn}_{\cal M}(G)_{(a(A),-1)}=0$ and $H^{\nn}_{\cal
M}(R)_{(a(A),0)} =0$.}

\end{rem}

\begin{rem}
\label{S12}
  {\rm Let us consider the group morphism $\psi : \Bbb Z^2   \rightarrow
\Bbb Z$ defined by $\psi(i,j) = i+j $. By Lemma \ref{B1}, $
H^{\nn+1}_{\cal M}(R^\psi)_l =\bigoplus_{i+j=l} H^{\nn +1}_{\cal
M}(R)_{(i,j)}$. Then, by applying Lemma \ref{S1} we get $a(R^\psi)
\leq a(A)-1$. If $R$ is Cohen-Macaulay and $a^2(G)<-1$, we have
proved $H^{\nn+1}_{\cal M} (R)_{(a(A),-1)} \not = 0$ and so
                        $a(R^\psi)= a(A)-1$.}
\end{rem}

\medskip
We can use the upper bound for the bigraded a-invariant of the
Rees algebra found in Lemma \ref{S1} to get bounds for the shifts
$(a,b)$ in its resolution. Namely,

\begin{lem}
\label{S13} Let $I$ be an ideal of $A$ generated by $r$ forms in
degrees $d_1 \leq...\leq d_r$ whose Rees algebra is
Cohen-Macaulay. Set $u=\sum_{j=1}^r d_j$.
 Let
$$0\to {D_{m}}\to \dots \to {D_1} \to
{D_0=S} \to R_A (I) \to 0$$
 be the minimal bigraded
 free resolution of ${\Rr}$ over $S$.
 Given $p \geq 1$ and $(a,b) \in {\Omega}_p$, we have
 \begin{enumerate}
 \item  ${a \leq 0}$, ${b \leq 0}$ , ${a \leq d_1 b}$.
 \item  ${-a-b \leq u+ \nn +a(A)+ p}$.
 \item  ${-a \leq u+ \nn + a(A)+ p-(r-1)}$.
        In particular, ${-a \leq u+ n + a(A)}$.
 \item  ${-b< r}$.
 \end{enumerate}
\end{lem}

{\pf} It is clear that ${a \leq 0}$, ${b \leq 0}$ , ${a \leq d_1
b}$. Also note that $m = \dpp_S R= n+r-\nn-1$. To prove $(ii)$,
let us consider the morphism $\psi : \Bbb Z^2   \rightarrow   \Bbb
Z$ defined by $\psi(i,j) = i+j $, and note that $S(a,b)^\psi =
S^\psi (a+b)$. Applying the functor $(\; \;)^\psi$ to the
resolution, we get a ${\Bbb Z}$-graded minimal free resolution of
$R^\psi$ over $S^\psi$. Moreover $a(S^\psi)= -n-u-r$ and
$a(R^\psi) \leq a(A)-1$ (see Remark \ref{S12}).
 Given $(a,b) \in {\Omega}_p$, from Theorem \ref{B77} we get:

\vspace{3mm}

\hspace{10mm} ${-a-b} \leq
  \max\{\, -\alpha-\beta \mid (\alpha,\beta) \in {\Omega}_p \} \leq $

\vspace{3mm}

\hspace{22mm}  $ \leq  \max\{\, -\alpha-\beta \mid (\alpha,\beta)
\in {\Omega}_{m} \} +p-m =$

\vspace{3mm}

\hspace{22mm} $= a(R^\psi)-a(S^\psi)+ p-m  \leq  u+ a(A)+ \nn + p.
$

\vspace{3mm}

To prove $(iii)$, observe that by Theorem \ref{B77} we have

\vspace{3mm}

\hspace{13mm} $ {-a} \leq \max \{-\alpha \mid (\alpha, \beta) \in
{\Omega}_p \}  \leq $

\vspace{3mm}

\hspace{19mm}  $ \leq  \max\{\, -\alpha \mid (\alpha,\beta) \in
{\Omega}_m \} +p-m =$

\vspace{3mm}

\hspace{19mm}  $= a^1(R)-a^1(S)+ p - m  \leq u + a(A) +\nn -r +1 +
p. $

\vspace{3mm}

Finally, by using Theorem \ref{B77} we also obtain:

\vspace{3mm}

\hspace{13mm} $ {-b} \leq
       \max\{\, -\beta \mid (\alpha,\beta) \in {\Omega}_p \} \leq $

\vspace{3mm}

\hspace{18mm} $ \leq
       \max\{\, -\beta \mid (\alpha,\beta) \in {\Omega}_m \} =$

\vspace{3mm}

\hspace{18mm} $ =a^2(R)-a^2(S) = -1 +r \;,$

\vspace{3mm}

\noindent so $(iv)$ is  proved. $\B$

\begin{rem}
\label{S144} {\rm When $I$ is a complete intersection ideal of the
polynomial ring $A=k[{\x}]$, all the shifts in the resolution may
be explicitly computed. In fact, by the Eagon-Northcott complex
the shifts $(a,b) \in {\Omega}_p$ are of the type:
         $$ a=-d_{j_1}- ... -d_{j_{p+1}}, \hspace{2mm} b=-m $$
where $1 \leq j_1 \leq ... \leq j_{p+1} \leq r$, $1 \leq m \leq p$
(see \cite[Lemma 4.1]{CHTV}). Note that $b$ takes all the values
between $-r$ and $0$ and the bounds of Lemma \ref{S13} $(ii)$,
$(iii)$ are sharp for $p=r-1$.}
\end{rem}

\medskip
       Now we are ready to determine a family of diagonals of the
Rees algebra with the Cohen-Macaulay property when the Rees
algebra is Cohen-Macaulay. Namely,

\begin{thm}
\label{S2} Let $I \subset A$ be a homogeneous ideal generated by
$r$ forms of degrees $d_1 \leq \dots \leq d_r=d$. Assume that
$H^i_\fm(A)_0 =0$ for all $i < \nn$. Set $u = \sum_{j=1}^r d_j$.
If the Rees algebra is Cohen-Macaulay, then
\begin{enumerate}
\item $k[(I^e)_c]$ is Cohen-Macaulay for
      $c> \max \{d(e-1)+u+a(A), d(e-1)+u-d_1(r-1) \}$.
\item If $I$ is equigenerated by forms of degree $d$, $\Kk$ is
      Cohen-Macaulay for $c> d(e -1 +l)+ a(A) $.
\end{enumerate}
\end{thm}

{\pf} We have already shown that the assumptions on the local
cohomology of $A$ imply that $H^i_{m}(\Kk)_0 = 0$ for $i < \nn$.
Now, let us consider the bigraded minimal free resolution of $R$
over $S$:
$$0 \to D_{m} \to \dots \to D_0= S \to R \to 0 \;,$$
where $D_p = \bigoplus_{(a,b) \in {\Omega}_p} S(a,b)$. From Remark
\ref{W10}, recall that if we define

\vspace{3mm}

\hspace{15mm}$X^\Delta = \bigcup_{(a,b) \in \Omega_R} \{ \; s \in
\Bbb Z
 \mid  \frac{-b}{e} \leq s \leq \frac{bd-a-n}{c-ed} \; \},$

\vspace{4mm}

\hspace{15mm}$Y^\Delta = \bigcup_{(a,b) \in \Omega_R} \{ \; s \in
\Bbb Z \mid  \frac{(b+r)d-u-a}{c-ed} \leq s \leq  \frac{-b-r}{e}
\; \} \,,$

\vspace{3mm}

\noindent then we have  $H^i_m (\Kk)_s = H^{i+1}_{\cal
M}(R)_{(cs,es)} = 0$ for $i< \nn$, $s \not \in X^\Delta \cup
Y^\Delta$. Therefore, $\Kk$ is Cohen-Macaulay for any diagonal
$\Delta=(c,e)$ such that $X^\Delta \cup Y^\Delta \subset \{0\}$.
Since $b \leq 0$, any  $s \in X^\Delta$ satisfies $s \geq 0$. If
$b \leq -1$, then $bd-a-n \leq -d +u+a(A)$ by Lemma \ref{S13}. If
$b=0$, then note that $[{\bf D}.]_0$ is a graded minimal free
resolution of $A$ over $S_1$ with $[D_p]_0 = \bigoplus_{(a,0) \in
\Omega_p} S_1(a)$, so $bd-a-n = -a-n \leq a(A)$ by Theorem
\ref{B77}. Therefore, by taking $c> (e-1)d+ u +a(A)$, we have
$\frac{bd-a-n}{c-ed} <1 $ and so $X^\Delta \subset \{0\}$. On the
other hand, any shift $(a,b) \in \Omega_R$ satisfies $b>-r$ by
Lemma \ref{S13}, so if $s \in Y^\Delta$ then $s \leq -1$. By
taking $c> d(e-1)+u-d_1(r-1)$, one can check
$\frac{(b+r)d-u-a}{c-ed} > -1$, so $Y^\Delta = \emptyset$. This
proves $(i)$.

Now, let us assume that $I$ is generated in degree $d$. From the
proof of Proposition \ref{S611}, we have that $H^i_m(\Kk)_s = 0$
for $i < \nn$, $s<0$. So it is just enough to study the positive
components of these local cohomology modules. Tensoriazing by
$k(T)$ we may assume that the field $k$ is infinite. Then, since
the fiber cone $F_{\fm}(I)$ of $I$ is a $k$-algebra generated by
homogeneous elements in degree $(d,1)$, there exists a minimal
reduction $J$ of $I$ generated by $l$ forms of degree $d$. Now, by
considering the polynomial ring $S=k[{\x},Y_1, \dots,Y_l]$, we
have a natural epimorhism $S \rightarrow R_A(J)$. Then $R_A(J)$ is
a finitely generated bigraded $S$-module, and so $R_A(I)$ because
it is a finitely generated $R_A(J)$-module. Then we may consider
the bigraded minimal free resolution of $R_A(I)$ over $S$, and it
suffices to check that the sets $X^\Delta$ and $Y^\Delta$
associated to this resolution do not have positive integers for
$c> d(e -1 +l)+ a(A) $.  $\B$

\medskip
In the case where $A=k[{\x}]$ we can improve the bounds slightly.
More explicitly,

\begin{thm}
\label{S3} Let $I \subset A= k[{\x}]$ be a homogeneous ideal
generated by $r$ forms of degrees $d_1 \leq \dots \leq d_r$.
Assume that the Rees algebra ${\Rr}$ is Cohen-Macaulay. Then, by
defining
$$\alpha= \min \, \{d(e-1)+u-n, e(u-n)\} \;,$$
$$\beta= \min \, \{d(e-1)+u-d_1(r-1), e(u-d_1)\} \;,$$
we have that $k[(I^e)_c]$ is Cohen-Macaulay for all $c > \max
\{de, \alpha, \beta \} $.
\end{thm}

{\pf} Note that the first homomorphism in the resolution of the
Rees algebra is:
 $$ \begin{array}{lll}
 D_0 = S & \longrightarrow  & R_A(I) \\
   \;\; X_i  &\mapsto & X_i  \\
   \;\; Y_j  & \mapsto  & f_j t  \hspace{2mm},\\
\end{array} $$
so any shift $(a,b) \in {\Omega}_p$, with $p \geq 1$, satisfies
$b< 0$. Note that if $\frac{bd-a-n}{c-ed} < \frac{-b}{e}$ for all
$(a,b) \in \Omega_p$, $p \geq 1$, then $X^\Delta$ is empty. This
condition is equivalent to $e(-a-n)<-bc$. Since $e(-a-n) \leq
e(-n-u)$ by Lemma \ref{S13} and $-bc \geq c$, it suffices to take
$c>e(u-n)$ to get this condition. Similarly, if $c>e(u-d_1)$ then
$Y^\Delta = \emptyset$ and we are done. $\B$

\medskip
\begin{rem}
\label{S003} {\rm  With the notation above, note that $\alpha=
e(u-n)$ if and only if $u-d-n<0$, and $\beta=e(u-d_1)$ if and only
if $u-d-d_1<0$ and $e>\frac{u-d-d_1(r-1)}{u-d_1-d}$. For instance,
if $u-d< n$ then $\Kk$ is Cohen-Macaulay for all $c > \max \{de,
\beta \} $.}
\end{rem}

\medskip
We finish this section with an application of Corollary \ref{G02}
to the study of the $(n-1)$-folds obtained from $\Bbb P^{n-1}_k$
by blowing-up a finite set of distinct points. Let $P_1, \dots,
P_s \in {\Bbb P}^{n-1}_k$ be distinct points, and for each $i=1,
\dots, s$, denote by ${\cal P}_i \subset A=k[{\x}]$ the
homogeneous prime ideal which corresponds to $P_i$. Let us
consider the ideal of fat points $I= {\cal P}_1^{m_1} \cap \dots
\cap {\cal P}_s^{m_s}$, with $m_1, \dots, m_s \in \Bbb Z_{\geq
1}$. Next we study the embeddings of the blow-up of $\Bbb
P^{n-1}_k$ along $\cal I$ via the linear systems $(I^e)_c$,
whenever these linear systems are very ample, slightly extending
\cite[Theorem 2.4]{GGP} where only the divisors $(I_c)$ were
considered.

\begin{thm}
\label{S7} Let $I \subset A= k[{\x}]$ be an ideal of fat points,
where $k$ is a field with characteristic 0. Then:
\begin{enumerate}
\item
$\Kk$ is Cohen-Macaulay if and only if $H^{i}_{\fm}(I^{es})_{cs} =
0$ for any $s>0$, $i<n$.
\item For $c > \reg(I) \, e$, $\Kk$ is
Cohen-Macaulay with $a (\Kk) < 0$.  In particular, $\reg (\Kk) <
n-1$.
\end{enumerate}
\end{thm}

{\pf}
   Let $X$ be the blow-up of the projective space $\Bbb P^{n-1}_k$
along $\cal I$. Assume that $I$ is generated by forms in degree
$\leq d$. Then we have shown that ${\cal L}^{e} \otimes {\M}^{c}$
is very ample if $c>de$. Therefore, for any $s<0$, $i<n-1$,
$c>de$, we have that $H^{i}(X, {\cal L}^{es} \otimes {\M}^{cs}) =
0$ by the Kodaira vanishing theorem (see for instance \cite[Remark
III.7.5]{Ha}). Then, $H^i_{m}(\Kk)_s = 0$ for $i <n$, $s<0$ by
Remark \ref{G1}.

On the other hand, from Proposition \ref{R5} we get ${\mit \Gamma}
({\X}, \cal O _{\X})$ = $k$ and $H^i({\X}, {\cal O}_{\X}) = 0$ for
all $ i >0$. Then, according to Remark \ref{G1}, we have
$H^i_{m}(\Kk)_0 = 0$ for any $i$.

Finally, note that for a given $\frak p \in \Proj (A)$ we have:
$$I_\frak p =
\cases{ A_\frak p & if $\frak p \not \in \{ {\cal P}_1, \dots,
{\cal P}_s \}$  \cr {\cal P}_i^{m_i} A_{{\cal P}_i} & if  $\frak p
={\cal P}_i$  .\cr}$$ In both cases, $R_{A_\frak p}(I_{\frak p})$
is Cohen-Macaulay. So, according to Corollary \ref{G3}, we have
that for any $s>0$ and $c \geq de+1$ there is an exact sequence
$$0 \to H^0_{m} (\Kk)_s
\to (I^{es})_{cs} \to (I^{es})^*_{cs} \to H^1_{m}(\Kk)_s \to 0$$
and isomorphisms $H^i_{m} (\Kk)_s \cong H^i_{\frak
m}(I^{es})_{cs}$ for $i>1$. Therefore, we immediately get $(i)$.
From \cite[Theorem 1.1]{GGP} or \cite[Theorem 6]{C}, we have
$a_*(I^e) \leq {\rm reg}(I^e) \leq e \; {\rm reg}(I)$. Then, by
taking $c> \reg(I) e$, we have $H^i_{m}(\Kk)_s = 0$  for any
$s>0$. So $\Kk$ is Cohen-Macaulay with $a(\Kk) < 0$. $\B$


%

\bigskip
\section{Linear bounds}
\label{UU} \markboth{CHAPTER III. COHEN-MACAULAY BLOW-UP SCHEMES}
{LINEAR BOUNDS}

\medskip
S.D. Cutkosky and J. Herzog  \cite{CH} studied sufficient
conditions for the existence of a constant $f$ satisfying that the
rings {\Kk} are Cohen-Macaulay for all $c \geq ef$ and $e>0$, that
is, for the existence of a linear bound on $c$ and $e$ ensuring
that $\Kk$ is Cohen-Macaulay. Note that, according to Theorem
\ref{S2}, this holds for any homogeneous ideal in a Cohen-Macaulay
ring whose Rees algebra is Cohen-Macaulay. Our first purpose is to
show that this also holds under the weaker assumption that
$R_{A_\fp}(I_\fp)$ is Cohen-Macaulay for any $\fp \in \Proj(A)$.
This would recover for instance locally complete intersection
ideals.

\medskip
Let $K=K_R =\bigoplus_{(i,j)} K_{(i,j)}$ be the canonical module
of $R=R_A(I)$, and let $K^e$ be the graded $A$-module $K^e=
\bigoplus_i K_{(i,e)}$. Then we have

\begin{thm}
\label{UU8} Assume that $R_{A_{\frak p}} ( I_ {\frak p})$ is
Cohen-Macaulay for all ${\frak p} \in \Proj(A)$. Then $\pi_* (w_X
\otimes {\cal L}^e) = \widetilde {K^e}$ and $R^j \pi_* (w_X
\otimes {\cal L}^e) = 0$ for $e >0$, $j>0$.
\end{thm}

{\pf} Let $A_i =A_{(x_i)}$, $I_i= I_{(x_i)}$, $R_i =A_i[I_i t]$
and $K_i =K_R \otimes R_i$. Let us consider the affine cover $\{
Y_i : 1 \leq i \leq n \}$ of $Y$, where $Y_i= Y- V_+ (x_i) \cong
\Spec(A_i)$. Denote by $X_i = \pi^{-1} Y_i = \Proj (R_i)$. Then,
for a given $j$ and $e>0$ we have that $R^j \pi_* (w_X \otimes
{\cal L}^e) = 0$ if and only if $R^j \pi_* (w_X \otimes {\cal
L}^e) \mid Y_i = 0$ for all $1 \leq i \leq n$. Furthermore, we
have a diagram
$$\begin{array}{ccccccccc}
 X_i = {\rm Proj }(R_i) & \hookrightarrow & X= \Proj^2(R) \\
  \pi' \downarrow &  \;    & \pi \downarrow   \\
Y_i = \Spec(A_i) & \hookrightarrow &  Y={\rm Proj}(A)
\end{array} $$
Now, by Corollary III.8.2 and Proposition III.8.5 of \cite{Ha},
for any $e > 0$ and $j > 0$ we have
$$R^j \pi_* (w_X \otimes {\cal L}^e) \mid Y_i =
R^j \pi'_* ((w_X \otimes {\cal L}^e) \mid X_i) = H^j (X_i, (w_X
\otimes {\cal L}^e) \mid {\X}_i) \,\,\widetilde \,\,.$$ Since $
(w_X \otimes {\cal L}^e) \mid X_i = \widetilde {K_i(e)}$, we have
reduced the problem to show that $H^{j+1}_{(R_i)_+} (K_{i})_ e=0$.
Similarly, $\pi_* (w_X \otimes {\cal L}^e) = \widetilde{K^e}$ if
$H^{0}_{(R_i)_+} (K_{i})_ e = H^{1}_{(R_i)_+} (K_{i})_ e =0$.

 Denote by $\overline R_i = R_{A_{x_i}}(I_{x_i})$. Tensoriazing
by $\overline R_i$, we have
 $$H^{j}_{(\overline R_i)_+} (K \otimes \overline R_{i})_ e=
 H^{j}_{(R_i)_+} (K_{i})_ e [T, T^{-1}],$$
so it is enough to show that $H^{j}_{(\overline R_i)_+} (K \otimes
\overline R_{i})_ e= 0$ for any $i,j$ and $e >0$. Let $\frak q \in
\Spec(A_{x_i})$ be a homogeneous prime, and let $\fp \in \Proj
(A)$ be such that $\frak q= \fp A_{x_i}$. Denote by $B=
R_{A_\fp}(I_{\fp})$. Then
$$[H^{j}_{(\overline R_i)_+} (K \otimes \overline R_{i})_ e]_\frak q =
[H^{j}_{(\overline R_i)_+} (K \otimes \overline R_{i})_\frak q]_e
= [H^{j}_{B_+} (K \otimes_A A_\fp)]_e.  $$ By taking into account
that $B$ is Cohen-Macaulay, standard arguments allow to check that
$K \otimes_A A_\fp = K_B$ or $K \otimes_A A_\fp = 0$. In any case,
we have that $[H^{j}_{B_+} (K \otimes_A A_\fp)]_e=0$ for any $j$
and $e>0$, so we are done. $\Box$

\medskip
From this result we can obtain a simple criterion for having a
linear bound for the Cohen-Macaulay property. First, let us notice
the following interesting fact.

\begin{prop}
\label{UU80} Let $A$ be an equidimensional ring. Assume
$R_{A_{\frak p}} ( I_ {\frak p})$ is Cohen-Macaulay for all
${\frak p} \in \Proj(A)$. Then, for any $c \geq de+1$ and $e>0$:
\begin{enumerate}
\item For $s>0$, there is an exact sequence
$$ 0 \to H^0_m(\Kk)_s  \to (I^{es})_{cs} \to
(I^{es})^*_{cs} \to H^1_m (\Kk)_s \to 0$$ and isomorphisms $H^i_m
(\Kk)_s \cong H^{i}_{\fm}(I^{es})_{cs}$ for $i>1$.
\item For $s>0$, $1 \leq i \leq \nn - 1$,
$H^i_m (\Kk)_{-s} \cong H^{\nn-i+1}_{\fm}(K^{es})_{cs}$.
\end{enumerate}
\end{prop}

{\pf} The first part of the statement follows directly from
Corollary \ref{G3}. To prove $(ii)$, let $s>0$, $i \geq 1$. Then

\vspace{3mm}

\hspace{1mm} $H^i_{\fm}(\Kk)_{-s} \cong H^{i-1}(X, {\cal L}^{-es}
\otimes {\cal M}^{-cs})$  \,
 by Remark \ref{G1}

\vspace{3mm}

\hspace{15mm}
 $= H^{\nn-i}(X , w_X \otimes {\cal L}^{es} \otimes {\cal M}^{cs})$ \,
 by Serre's duality

\vspace{3mm}

\hspace{15mm}
 $= H^{\nn-i}(Y ,
\pi_* (w_X \otimes {\cal L}^{es}) \otimes {\cal M}^{cs})$ \,
 by Theorem \ref{UU8}

\vspace{3mm}

\hspace{15mm} $= H^{\nn-i}(Y, \widetilde{K^{es}}(cs))$ \, by
Theorem \ref{UU8}

\vspace{3mm}

\hspace{15mm} $= H^{\nn-i+1}_{\fm}(K^{es})_{cs}.$ $\B$

\begin{thm}
\label{UU81} Assume that $A$ is an equidimensional ring with
$H^i_{\fm}(A)_0=0$ for $i<\nn$. If $I$ is a homogeneous ideal of
$A$ such that $R_{A_{\frak p}} ( I_{\frak p})$ is Cohen-Macaulay
for all ${\frak p} \in \Proj(A)$, then there exists $\alpha$ such
that $\Kk$ is Cohen-Macaulay for $c \geq de+\alpha$, $e>0$.
\end{thm}

{\pf} From Proposition \ref{R7} we have that $\Kk$ is
Cohen-Macaulay for $c \gg e \gg 0$. So, in particular, by Theorem
\ref{J1} and Remark \ref{G1} we have that $H^i_{m}(\Kk)_0=0$ for
any $c \geq de+1$, $i< \nn$. On the other hand, according to
Theorem \ref{S8}, there exists $\alpha>0 $ such that $a_*(I^e) <
de + \alpha$, $a_*(K^e) < de + \alpha$, for all $e$. Then, $\Kk$
is Cohen-Macaulay for any $c \geq de+ \alpha$ by Proposition
\ref{UU80}. $\Box$

\medskip
In particular, we can recover Corollary 4.2 and Corollary 4.4 in
\cite{CH}. Furthermore, note that the bound has been improved
slightly.

\begin{cor}
\label{UU2} Let $I$ be a locally complete intersection ideal in a
Cohen-Macaulay ring $A$. Then there exists $\alpha$ such that such
that $\Kk$ is Cohen-Macaulay for any $c \geq de+\alpha$ and $e>0$.
\end{cor}

\begin{cor}
\label{UU22} Let $I$ be a homogeneous ideal such that $I_{(\fp)}$
is a strongly Cohen-Macaulay ideal with
 $\mu(I_{(\fp)}) \leq \h (\fp)$ for any prime ideal $\fp \supseteq I$.
Then there exists $\alpha$ such that such that $\Kk$ is
Cohen-Macaulay for any $c \geq de+\alpha$ and $e>0$.
\end{cor}

\medskip
We can also characterize the existence of linear bounds for the
Cohen-Macaulay property of the rings $\Kk$ by means of the local
cohomology modules of the Rees algebra and its canonical module.
Namely,

\begin{prop}
\label{UU82} Assume that there exist $c,e$ such that $\Kk$ is
Cohen-Macaulay. Then the following are equivalent
\begin{enumerate}
\item There exists $f$ such that $\Kk$ is Cohen-Macaulay
for $c \geq ef$, $e>0$.
\item  There exists $f$ such that $H^i_{R_+}(R)_{(c,e)} = 0$,
$H^{\nn-i+1}_{R_+}(K_R)_{(c,e)} = 0$, for $i<\nn$, $ c\geq ef$ and
$e>0$.
\item  There exists $f$ such that $H^i_{{\M}_2}(R)_{(c,e)} = 0$,
$H^{\nn-i+1}_{{\M}_2}(K_R)_{(c,e)} = 0$, for $i<\nn$, $ c \geq ef$
and $e>0$.
\end{enumerate}
\end{prop}

{\pf} From Proposition \ref{W1}, we have  $H^i_m (\Kk)_s =
H^i_{R_+}(R)_{(cs, es)}$ for any $s>0$. Moreover, since $X$ is
Cohen-Macaulay and equidimensional we also have that $H^i_m
(\Kk)_{-s} = H^{\nn-i+1}_{R_+} (K_R)_{(cs, es)}$ for $1 \leq i <
\nn$ and $s>0$ according to Proposition \ref{GG01}. Therefore two
first conditions are equivalent.

 To prove $(ii) \iff (iii)$, first we will show that there exists
$\overline f$ such that for all $i$, $e>0$, $c \geq e \overline f$
it holds
$$H^i_{\M}(R)_{(c,e)}= H^i_{{\M}_1}(R)_{(c,e)} = 0,$$
$$H^i_{\M}(K_R)_{(c,e)}= H^i_{{\M}_1}(K_R)_{(c,e)} = 0.$$
Then, from the Mayer-Vietoris exact sequence, we get that for any
$i$, $e>0$, $c \geq e \overline f$,
$$H^i_{R_+}(R)_{(c,e)}= H^i_{{\M}_2}(R)_{(c,e)} ,$$
$$H^i_{R_+}(K_R)_{(c,e)}= H^i_{{\M}_2}(K_R)_{(c,e)} ,$$
and so $(ii)$ and $(iii)$ are equivalent. To get the vanishing of
the local cohomology modules with respect to the maximal ideal
$\M$, it is just enough to take  $c> \max \, \{ a_*^1(R),
a_*^1(K_R) \}$. By Theorem \ref{S8} there exists $\alpha >0$ such
that $a_*(I^e) < de+\alpha$, $a_*(K^e) < de+\alpha$, so by taking
$c \geq de+ \alpha$  we have $H^i_{{\M}_1}(R)_{(c,e)}= H^i_\fm
(I^e)_c =0$ and $H^i_{{\M}_1}(K_R)_{(c,e)}= H^i_{\fm}(K^e)_c = 0$.
$\B$

\medskip
Also note that the last proposition holds if we replace the
condition for all $c \geq ef$ and $e >0$ by the following one: for
all $c \geq de+ \alpha$ and $e>0$. As a direct consequence, we
obtain:

\begin{cor}
\label{UU83} Assume that $R$ has some Cohen-Macaulay diagonal. If
$a_*^2(R) \leq 0$ and $a_*^2(K_R) \leq 0$, there exists $\alpha$
such that $\Kk$ is Cohen-Macaulay for all $c \geq de + \alpha$ and
$e>0$.
\end{cor}

{\pf} It is a direct consequence of Proposition \ref{UU82} by
noting that for any $i$ and $e>0$ we have $H^i_{{\M}_2}(R)_{(c,e)}
= 0$, $H^i_{{\M}_2}(K_R)_{(c,e)} = 0$  by Proposition \ref{W11}.
$\B$

\begin{rem}
\label{UU84} {\rm The converse of the last corollary is not true.
Let us take the homogeneous ideal $I = (x^7, y^7, x^6y+x^2y^5)$ in
the polynomial ring $A=k[x, y]$. We have $\fm^{14} \subset I$, so
$(I^e)_c = A_c $ for any $c \geq 14 e$, and then $\Kk = k [x^a y^b
\mid a+b=c]$ is Cohen-Macaulay for all these $c, e$. But $a_*^2(R)
= 4>0$ by \cite[Example 3.13]{HM}. }
\end{rem}

\chapter*{$\;\;\;$}

\newpage

\medskip

\chapter{Gorenstein coordinate rings of blow-up schemes}
\typeout{Gorenstein coordinate rings of blow-up schemes}

\bigskip
Let $Y= {\rm Proj}(A)$ be a closed subscheme of $\Bbb P_k ^
{n-1}$, and let $X$ be the blow-up of $Y$ along ${\cal
I}=\widetilde{I}$, where $I$ is a homogeneous ideal in $A$. If $I$
is generated by forms of degree $\leq d$, we have already shown
that for any $c \geq de+1$ the ring $\Kk$ is the homogeneous
coordinate ring of a projective embedding of $X$ in $\Bbb
P^{N-1}_k$, where $N= \dim_k (I^e)_c$. In this chapter we are
interested in the (quasi-) Gorenstein property of the rings $\Kk
$. The results work in the case of the blow-up of the projective
space $\Bbb P^{n-1}_k$.

\medskip
If the Rees algebra is Cohen-Macaulay and the associated graded
ring is Gorenstein we will determine exactly for which pairs
$(c,e)$ the ring $\Kk$ is quasi-Gorenstein and, in particular, we
will obtain that there is just a finite set of diagonals with this
property.
 This result can be applied to several families of ideals. In
particular, to any complete intersection ideal, extending in this
way \cite[Corollary 4.7]{CHTV}, and to the ideal generated by the
maximal minors of a generic matrix.

\medskip
   After that, we show that there are always at most a finite number of
rings $\Kk$ which are quasi-Gorenstein and we give upper bounds
for such diagonals whenever $R_A(I)$ is Cohen-Macaulay.

\medskip
     Finally, we prove that under some restrictions the existence of
a diagonal $(c,e)$ such that $\Kk$ is quasi-Gorenstein forces the
associated graded ring to be Gorenstein.

\medskip
     At the end of the chapter we apply our results to
Room surfaces. These surfaces are obtained by blowing-up $\Bbb
P^2_k$ along ${d+1 \choose 2}$ points, $d \geq 2$, which do not
lie in any curve of degree $d-1$ and then embedding in $\Bbb
P^{2d+2}_k$. We will show that the only Room surface which is
Gorenstein is the del Pezzo sestic surface in $\Bbb P^6$, so
recovering that well known result (see \cite[Example 1]{GG}).

\bigskip
\section{The case of ideals whose form ring is Gorenstein}
\label{K} \markboth{CHAPTER IV. GORENSTEIN BLOW-UP SCHEMES}{IDEALS
WHOSE FORM RING IS GORENSTEIN}

\medskip
    Throughout all the chapter we shall use the following notations.
$A=k[X_1,...,X_n]$ will denote the usual polynomial ring with
coefficients in a field $k$, and $I \subset A$ a homogeneous ideal
minimally generated by forms $f_1,...,f_r$ of degrees $d_1
\leq...\leq d_r=d$. We put $u=\sum_{j=1}^r d_j$. Let
$S=k[X_1,...,X_n,Y_1,...,Y_r]$  be the polynomial ring with the
grading obtained by setting $\deg X_i = (1,0)$ for  $ i=1,...,n$,
$\deg Y_j = (d_j,1)$ for $  j=1,...,r$; so that $R= R_A(I)$ and
$G=G_A(I)$ can be seen in a natural way as bigraded $S$-modules.
We will assume $n \geq r \geq 2$.

\medskip
Notice that any diagonal $S_\Delta$ of the polynomial ring $S$ is
Cohen-Macaulay by Corollary \ref{W666}. We begin this section by
showing that, on the contrary, $S_\Delta$ is Gorenstein only for a
finite number of diagonals. Furthermore, we may determine them.

\begin{prop}
\label{K1}
 ${S_\Delta}$ is Gorenstein if and only if
  $\frac{r}{e}  =\frac{n+u}{c} = l \in {\Bbb Z}$.
 Then $a(S_\Delta)= -l$.
\end{prop}

{\pf} Let $T= S_\Delta = \bigoplus_{s \geq 0} U_s$, where $U_s$ is
the $k$-vector space generated by the monomials $X^{\alpha_1}_1
\cdots X^{\alpha_n}_n Y^{\beta_1}_1 \cdots Y^{\beta_r}_r$ with
$\alpha_i, \beta_j \geq 0$ satisfying the equations

 $$(\star) \hspace{5mm} \cases{ \hspace{3mm}
 \sum_{i=1}^n {\alpha_i} + \sum_{j=1}^r {d_j \beta_j} =cs \cr
 \vspace{1mm} \cr
 \hspace{3mm}\sum_{j=1}^r {\beta_j}=es \cr} $$

\noindent By Corollary \ref{W3}, we have
$$H^{n+r-1}_m (T) \cong H^{n+r}_{\cal M}(S)_\Delta \cong
(\bigoplus_{\alpha<0, \beta<0} k X^{\alpha} Y^{\beta}
\,)_\Delta.$$ Therefore, $K_T= {\uHom}_k (H^{n+r-1}_{m}(T), k) =
\bigoplus_{s \geq 1} V_s$ with $V_s$ the $k$-vector space
generated by the monomials $X^{\alpha_1}_1 \cdots X^{\alpha_n}_n
Y^{\beta_1}_1 \cdots Y^{\beta_r}_r$, and $\alpha_i>0, \beta_j>0 $
which satisfy ($\star$). Since $T$ is Cohen-Macaulay, $T$ is
Gorenstein if and only if $K_T \cong T(a(T))$.

Assume first that $\frac{r}{e}  =\frac{n+u}{c} = l \in {\Bbb Z}$.
Then, the multiplication by $X_1 \cdots X_n Y_1 \cdots Y_r \in
T_l$ induces an isomorphism $T \cong K_T(l)$ and so $T$ is
Gorenstein with $a(T) = -l$.

 To prove the converse set $({\bf \alpha},{\bf \beta}) =
({\alpha_1},...,{\alpha_n},{\beta_1},...,{\beta_r})$ with
$\alpha_i, \beta_j >0$ and assume the contrary. This means that
$({\bf 1},{\bf 1})$ is not a solution of ($\star$) for any $s$. On
the other hand, the set of solutions of $(\star)$ for some $s$ is
partially ordered by means of $({\bf \alpha},{\bf \beta})$ $\leq$
$({\bf \gamma},{\bf \rho})$ $\Longleftrightarrow$ $\alpha_i \leq
\gamma_i$, $\beta_j \leq \rho_j$, $\forall i,j $. Then, one can
easily check that for any $i=1, \dots, n$, $j=1, \dots, r$ there
exists a solution of $(\star)$ for some $s$ such that $\alpha_i =
\beta_j =1$. This implies the existence of at least two minimal
solutions, and so $T$ is not Gorenstein. $\B$

\medskip
\begin{rem}
\label{K12} {\rm Note that the number of minimal elements in the
set of solutions of the system $(\star)$ coincides with the type
of $S_ {\Delta}$. It is not difficult to see that if $S_{\Delta}$
is not Gorenstein, then its type is $\geq r$.}
\end{rem}

\medskip
\begin{rem}
\label{K13} {\rm Throughout the chapter we assume $r \geq 2$. In
the case that $r=1$  we have that $I$ is a principal ideal, and
$R_A(I) \cong S=k[{\x}, Y]$. Then, it is easy to check that
$S_\Delta$ is always Cohen-Macaulay, and $S_\Delta$ is Gorenstein
if and only if $\Delta=(n+d,1)$.}
\end{rem}

\medskip
The last proposition leads to the question of whether there exist
diagonals $(c,e)$ such that  ${\Kk}$ be quasi-Gorenstein, and how
we can determine them. Our answer will be partially based on the
following proposition which links the diagonal of the canonical
module of ${\Rr}$ to the canonical module of the diagonal of
${\Rr}$. It is stated and proved for complete intersection ideals
in \cite[Proposition 4.5]{CHTV}, but in fact the same statement
and proof are valid in general. We include the proof for
completeness.

\begin{prop}
\label{K2}
         $K_{R_\Delta} = (K_{R})_{\Delta}$.
\end{prop}

{\pf} Let us consider a presentation of $R$ as $S$-module
$$ 0 \to C \to S \to R \to 0 ,$$
which leads to the bigraded exact sequence of local cohomology
modules
$$ 0 \to H^{n+1}_{\cal M} (R)  \to H^{n+2}_{\cal M}(C) \to
H^{n+2}_{\cal M}(S) \to 0.$$

 Similarly, we get the graded exact sequence
$$ 0 \to H^{n}_{m} (R_{\Delta})  \to H^{n+1}_{m} (C_{\Delta}) \to
H^{n+1}_{m}(S_\Delta) \to 0.$$

On the other hand, by Corollary \ref{W3} we have a commutative
diagram

$$\begin{array}{ccccccccc}
 0 & \to & H^{n+1}_{\cal M}(R)_{\Delta} & \to &
H^{n+2}_{\cal M}(C)_{\Delta} & \to & H^{n+2}_{\cal M}(S)_{\Delta}
&
\to & 0 \\
 & & \varphi^{n}_{R} \uparrow & & \varphi^{n+1}_C \uparrow & &
\varphi^{n+1}_S \uparrow & & \\
0 & \to & H^{n}_{m} (R_{\Delta}) &  \to & H^{n+1}_{m} (C_{\Delta})
& \to & H^{n+1}_{m}(S_\Delta) &  \to & 0   \end{array} $$

\medskip
\noindent where $\varphi^{n+1}_C$, $\varphi^{n+1}_S$ are
isomorphisms, and so $\varphi^{n}_{R}$ is also an isomorphism.
Therefore $H^{n}_{m} ({\Rdd}) \cong H^{n+1}_{\cal M} (R)_{\Delta}$
and we get

\vspace{3mm}

\hspace{10mm} $K_{R_{\Delta}} = \uHom_k (H^n_{m}(R_{\Delta}), k )=
\uHom_k (H^{n+1}_{\cal M}(R)_{\Delta}, k )=$

\vspace{3mm}

\hspace{19mm} $=\uHom_k (H^{n+1}_{\cal M}(R),k)_{\Delta} =
(K_R)_\Delta.$ $\B$

\medskip
\begin{rem}
\label{K21}
 {\rm The hypothesis $n \geq r \geq 2$ fixed before is only used in
this chapter to prove Proposition \ref{K2}, and of course its
applications. Nevertheless, the isomorphism $K_{R_\Delta}=
(K_{R})_{\Delta}$ is also valid if $n, r \geq 2$, $I$ is
equigenerated and $R$ is  Cohen-Macaulay. To prove this, assume $r
>n$ (otherwise we may apply Proposition \ref{K2}). Let us consider
the bigraded minimal free resolution of $R$ over $S$
$$0\to {D_{r-1}}\to \dots \to {D_1} \to {D_0=S} \to R \to 0,$$

\noindent with $D_p =\bigoplus_{(a,b) \in {\Omega}_p} S(a,b)$.
Since $R$ is Cohen-Macaulay, we have $-b \leq r-1$ and $-a \geq
-bd$ for all $(a,b) \in \Omega_R$ by Lemma \ref{S13}. On the other
hand, recall from Corollary \ref{W6} that $H^r_{{\cal
M}_2}(S(a,b))_{(cs,es)} \not = 0$ if and only if
 $\frac{bd-a}{c-ed} \leq s \leq \frac{-b-r}{e}$, so we get
$H^r_{{\cal M}_2} (D_p)_{\Delta} = 0$ for all $p$. Since $r>n$, we
also have that $H^i_{{\cal M}_1} (D_p)_{\Delta} = H^i_{{\cal M}_2}
(D_p)_{\Delta} = 0$  for all $i>n$. Then, by Proposition \ref{W61}
and Proposition \ref{W2} we have that $\varphi^{n+1}_C$ is an
isomorphism, and the same proof as in Proposition  \ref{K2} shows
that $K_{R_{\Delta}} = (K_R)_{\Delta}$.

This means that all the results we are going to prove in this
chapter are also valid if $n, r \geq 2$, $I$ is equigenerated and
$\Rr$ is Cohen-Macaulay.}
\end{rem}

\medskip
  In view of Proposition \ref{K2} any information on the bigraded
structure of $K_{R}$ will be of interest. Let $B$ be a
$d$-dimensional local ring, $d \geq 1$, which has a canonical
module $K_B$ and $I \subset B$ an ideal of positive height such
that $R_B(I)$ is Cohen-Macaulay. In \cite[Theorem 2.2]{TVZ} it is
given a description of $K_{R_B(I)}$ in terms of a filtration of
submodules of $K_B$. Assume now that $B = \bigoplus_{n \geq 0}B_n$
is a positively graded ring of positive dimension over a local
ring $B_0$, which has a canonical module $K_B$. Let $I \subset B$
be a homogeneous ideal of positive height. Then, the Rees algebra
$R_B(I)$ has a bigraded structure by means of
 $[R_B(I)]_{(i,j)} = (I^j)_i t^j$ for all $i, j \geq 0$. We
also have a bigraded structure on the form ring by means of
$[G_B(I)]_{(i,j)} = (I^j)_i / (I^{j+1})_i$ for all $i,j \geq 0$.

\medskip
Then, the proof of \cite[Theorem 2.2]{TVZ} may be "bigraded" and
we thus obtain a description of the bigraded structure of
$K_{R_B(I)}$. Namely, we get:

\begin{thm}
\label{K3} With the notation above assume that $R_B(I)$ is
Cohen-Macaulay. Then, there exists a homogeneous filtration $\{K_m
\}_{m \geq 0}$ of $K_B$ and isomorphisms of bigraded modules such
that
 $$ K_{R_B(I)} \cong \bigoplus_{(l,m),\, m \geq 1} [K_m]_l \, ,$$
 $$ K_{G_B(I)} \cong \bigoplus_{(l,m),\, m \geq 1}
[K_{m-1}]_l/[K_m]_l.$$
 \end{thm}

\medskip
Several other results of \cite{TVZ} may also be "bigraded". In
particular \cite[Lemma 4.1]{TVZ} which makes precise when the
canonical module of the Rees algebra has the expected form. Recall
that $K_{R_B(I)}$ has the expected form if
$$K_{R_B(I)} \cong Bt \oplus Bt^2
\oplus \cdots \oplus Bt^l \oplus It^{l+1} \oplus I^2t^{l+2} \oplus
\cdots \,,$$
 for some $l \geq 0$. This definition was introduced by
J. Herzog, A. Simis and W. Vasconcelos in \cite{HSV2}. We still
use the same notation and again omit the proof.

\begin{cor}
\label{K4}
        Assume $R_B(I)$ is Cohen-Macaulay and $G_B(I)$ is
quasi-Gorenstein. Set ${\bf a}(G_B(I)) = (-b, -a)$. Then $K_B
\cong B(-b)$ and
$$K_{R_B(I)} \cong \bigoplus_{(l,m),\, m \geq 1} [I^{m-a+1}]_{l-b},$$
\noindent where $I^n = B$ if $n \leq 0$.
\end{cor}

\medskip
 Note that $-a$ coincides with the usual a-invariant of $G_B(I)$.
By Ikeda-Trung's criterion \cite{IT} it is always negative if
$R_B(I)$ is Cohen-Macaulay, and it has been calculated in many
cases (see for instance  \cite{HRZ}, \cite{GH}). As for $b$, it is
clear that under the hypothesis of Corollary \ref{K4} we get $-b=
a(B)$. It is then also easy to compute the bigraded a-invariant of
$R_B(I)$. Namely, we get that if $a = 1$ then ${\bf a} (R_B(I)) =
(-d_1 + a(B), -1)$, and if $a
>1$ then ${\bf a}(R_B(I)) = (a(B), -1).$

\begin{rem}
\label{K45}
   {\rm     Assume that $B = A = k[{\x}]$ and $I$ is a complete
intersection ideal. Then, the Eagon-Northcott complex provides a
${\Bbb Z^2}$-graded minimal free resolution of $\Rr$. Following
the proof of Yoshino \cite{Y} it is possible to see that
$$K_{\Rr} = J ( (r-2)d_1-n,-1)$$
\noindent with $J= (f_1^{r-2}, f_1^{r-2}\;t,...,
f_1^{r-2}\;t^{r-2}) {\Rr}$.

Observe that in this case $a(G_A(I)) = (-n,-r)$ and by Corollary
\ref{K4}
$$K_{R_A(I)} = \bigoplus_{(l,m), \, m \geq 1} [I^{m-r+1}]_{l-n}.$$

A straightforward computation shows that, in fact, multiplication
by $f_1^{r-2}$ provides an explicit isomorphism
$$ \bigoplus_{(l,m), \, m \geq 1} [I^{m-r+1}]_{l-n}
\cong
 J ( (r-2)d_1-n,-1).$$}
 \end{rem}

\bigskip
Let us now assume that $I \subset A = k[{\x}]$ is a homogeneous
ideal whose form ring is Gorenstein. We are now ready to prove the
main result of this section determining the quasi-Gorenstein
diagonals of $\Rr$. Namely,

\begin{thm}
\label{K5} Assume $\h (I) \geq 2$, $\dim (A/I) >0$, and $G_A(I)$
is Gorenstein. Set $a=-a^2(G_A(I))$. Then $\Kk$ is
quasi-Gorenstein if and only if $\frac{n}{c}= \frac{a-1}{e} =l_0
\in \Bbb Z$. In this case, $a(\Kk) = -l_0$.
\end{thm}

{\pf} Note that the Rees algebra $R$ is Cohen-Macaulay by using a
result of Lipman \cite[Theorem 5]{L}. Now, by applying Corollary
\ref{K4} we have that $b=-a(A)=n$ and $K_{R} = \bigoplus_{(l,m), m
\geq 1} [I^{m-a+1}]_{l-n}$, so that by Proposition \ref{K2} we get
$K_{R_\Delta} = (K_R)_{\Delta} = \bigoplus_{l \geq 1}
[I^{el-a+1}]_{cl-n}$. Let $l_0 = \min \{ l \in \Bbb Z \mid l \geq
\frac{n}{c} \}$, $s= a-1-el_0$. We shall now distinguish three
cases.

If $s=0$, then the first non-zero component of $K_{R_\Delta}$ is
$[K_{R_\Delta}]_{l_0}= [I^{e l_0 -a +1}]_{c l_0 -n} = A_{c l_0
-n}$, so that if $\Rdd$ is quasi-Gorenstein $c l_0 - n =0$ and we
get that
 $l_0$ = $\frac{n}{c}$ =$\frac{a-1}{e}$
  and $a(R_{\Delta})= -l_0$. Conversely, if
 $l_0 = \frac{n}{c} =\frac{a-1}{e}$
 then
$[K_{R_\Delta}]_{l_0 + m}$ =$[I^{e l_0 -a+1+em}]_{c l_0 + cm - n}
= [I^{em}]_{cm} = [{\Rdd}]_m$
 for all $m$ and so $\Rdd$ is
quasi-Gorenstein.

If $s < 0$, let $l_1 = \min \{l \mid el - a + 1 > 0,
 cl - n \geq d_1 (el - a +1) \}$. Then $l_1 \geq l_0$ and the first
 non-zero component of $K_{R_\Delta}$ is
 $[K_{R_\Delta}]_{l_1} = [I^{e l_1 - a +1}]_{c l_1 - n}$. In particular,
 $a(\Rdd) = -l_1$. Assume $\Rdd$ is quasi-Gorenstein. Then
 $K_{R_\Delta} \cong \Rdd (-l_1)$ and so $[K_{R_\Delta}]_{l_1} \cong k$.
This implies that $c l_1 - n = d_1 ( e l_1 -a +1)$ :  If
 $c l_1 - n - d_1 ( e l_1 -a +1) =r > 0$ we may choose two lineary
 independent forms $g, h \in A_r$ such that $g f_1^{el_1-a+1},
h f_1^{el_1-a+1} \in [I^{e l_1 - a + 1}]_{c l_1 -n} \cong k$,
which is a contradiction. From the isomorphism one gets that
$K_{R_\Delta}$ is generated by $f_1^ {e l_1-a+1}$ as
$\Rdd$-module. Now let $f_j \not \in \rad (f_1)$ (it exists
because $\h(I) \geq 2$), and choose $m$ such that $m (c - d_j e)
> d_j - d_1$ and there exists $f \in A_{d_1 + cm -d_j(em+1)}$ such
that $(f, f_1)=1$. Then $f_1^{e l_1-a} f_j^{em +1} f \in [I^{e l_1
-a +1 +em}]_{d_1(e l_1-a+1) + cm} = f_1^{e l_1-a+1}
[I^{em}]_{cm}$, and we get $f_j^{em+1} f \in (f_1)$ which is a
contradiction.

If $s>0$, the first non-zero component of $K_{R_\Delta}$ is $[I^{e
l_0-a+1}]_{c l_0-n} = A_{c l_0-n}$, so if  $\Rdd$ is
quasi-Gorenstein we  get $c l_0-n = 0$. Furthermore, for all $m
\geq 1$ we have $[K_{R_\Delta}]_{l_0 +m}= [I^{-s+em}]_{c l_0-n+cm}
= [I^{-s+em}]_{cm} \cong [I^{em}]_{cm}$. Since $s>0$ and
$[I^{em}]_{cm} \subset [I^{-s+em}]_{cm}$ this isomorphism is
possible if and only if $[I^{em}]_{cm} = [I^{-s+em}]_{cm}$. Now
choose $X_i$  such that $X_i \not \in \rad (I)$ (it always exists
because $\dim (A/I) >0$) and $m$ with $em-s \geq 1$. For any $j$
consider $F_j= X_i^{\alpha_j} f_j^{em-s}$ where $\alpha_j = cm
-d_j(em-s) \geq 1$, and assume $[I^{em}]_{cm} = [I^{-s+em}]_{cm}$.
Then $F_j \in [I^{em-s}]_{cm}$ and so $ X_i^{\alpha_j} f_j^{em-s}
\in I^{em}$. Now let $f_1^{c_1} \dots f_r^{c_r}$ such that $c_1 +
\dots + c_r \geq r (em-s)$. This implies that there exists $l$
with $c_l \geq em - s$ and so $X_i^{\alpha_1} f_1^{c_1} \dots
f_r^{c_r} = X_i^{\alpha_1} f_l^{em-s} f_1^{c_1} \dots
f_l^{c_l-em+s} \dots f_r^{c_r} \in I^{c_1+ \dots + c_r +s}$, since
$\alpha_1 \geq \alpha_i$ for all $i$. Thus we get $X_i^{\alpha}
I^h \subset I^{h+s}$ for $h >>0$, which implies that $X_i^{\alpha}
\in I^s \subset I$ since $R_A(I)$ is Cohen-Macaulay (see
\cite[Lemma 4.3]{TVZ}). But this contradicts $X_i \not \in \rad
(I)$ and so $R_{\Delta}$ is not quasi-Gorenstein. $\B$

\medskip
The remaining cases $\h(I) =1,\, n$ in the above theorem are
studied separately in the following remarks.

\medskip
\begin{rem}
\label{K6}
 {\rm   If $\h(I)=1$ then $\Kk$ is never quasi-Gorenstein. In fact,
by \cite[Proposition 4.6]{TVZ}, $a^2(G_A(I)) = -1$ and so $a=1$.
Following the same proof as in Theorem \ref{K5} we have that $s=
-e l_0 <0$. On the other hand, since $\h(I) = 1$ we may write $I =
gJ$, with $\h(J) \geq 2$, $J= (\overline f_1, \dots, \overline
f_r)$ and $f_j = \overline f_j g$ for all $j$. The same argument
as in Theorem \ref{K5} for the case $s<0$ but taking $\overline
f_j \not  \in \rad (\overline f_1)$ and $f \in
A_{d_1+cm-d_j(em+1)}$ such that $(f, \overline f_1) =1$ leads to
$\overline f_j^{em+1} f \in (\overline f_1)$, which is a
contradiction. }
\end{rem}

\begin{rem}
\label{K7} {\rm   If $\dim (A/I) = 0$, then the condition
$\frac{n}{c}= \frac{a-1}{e}= l_0 \in {\Bbb Z}$ is sufficient but
not necessary for $\Kk$ to be quasi-Gorenstein. For instance, let
$A= k[X_1, X_2, X_3]$ and $I= (X_1, X_2, X_3 )$. Note that $n=3
\geq r =3 \geq 2$, $G$ is Gorenstein and $a=3$. Then, by Corollary
\ref{K4} we have that $K_R = \bigoplus_{(l,m), m \geq 1}
[I^{m-2}]_{l-3}$. According to Proposition \ref{K2}, by taking the
$(3,1)$-diagonal we have
 $$K_{R_{\Delta}} =
\bigoplus_{l \geq 1} [I^{l-2}]_{3(l-1)} = \bigoplus_{l \geq 1}
A_{3(l-1)}= (\bigoplus_{l \geq 0} A_{3l})(-1) = R_{\Delta}(-1),$$
and so $R_{\Delta} = k[I_3]$ is quasi-Gorenstein. In this case,
$\frac{n}{c}$ =$1 \not = 2$ =$\frac{a-1}{e}$.}
\end{rem}

\medskip
As a consequence of Theorem \ref{K5} we obtain the following
result for the case of complete intersection ideals. It
generalizes \cite[Corollary 4.7]{CHTV} where the case of ideals
generated by two elements was considered.

\begin{cor}
\label{K8} Let $I \subset k[{\x}]$ be a complete intersection
ideal minimally generated by $r$ forms of degrees $d_1 \leq \dots
\leq d_r=d$, with $ r<n$. Then for $ c \geq de+1$,
 $\Kk$ is Gorenstein if and only if
 $\frac{n}{c}= \frac{r-1}{e} =l_0 \in \Bbb Z$.
 In this case, $a(\Kk) = -l_0$.

\end{cor}

{\pf} Since $a^2(G_A(I)) = -r$, by Theorem \ref{K5} we have that
$\Kk$ is quasi-Gorenstein if and only if $\frac{n}{c}$=
$\frac{r-1}{e}$ =$l_0 \in \Bbb Z$. But then $u + (e-1)d -n \leq rd
+ed -d -n= (r-1)d +de -n = e \frac{n}{c} d + de -n = n
(\frac{ed-c}{c}) +de \leq de <c$, and according to Proposition
\ref{S4}, $\Kk$ is also Cohen-Macaulay and so Gorenstein. $\B$

\medskip
We may also study the ideals generated by the maximal minors of a
generic matrix. We thank A. Conca for suggesting to consider this
case.

\begin{ex}
\label{K9} {\rm Let ${\bf X}= (X_{ij})$ be a generic matrix, with
$ 1 \leq i \leq n$, $1 \leq j \leq m$ and $m \leq n$. Let us
consider $I \subset A= k[{\bf X}]$ the ideal generated by the
maximal minors of ${\bf X}$, where $k$ is a field. In this case,
the Rees algebra $R_A(I)$ is Cohen-Macaulay and the form ring
$G_A(I)$ is Gorenstein \cite[Theorem 3.5]{EH}. Moreover, it has
been proved by A. Conca (personal communication) that all the
diagonals of $R_A(I)$ are Cohen-Macaulay (see also Example
5.2.23).

Now we want to study the Gorenstein property of these diagonals.
Note that $I$ is an equigenerated ideal whose Rees algebra is
Cohen-Macaulay, so we can apply Theorem \ref{K5} thanks to Remark
\ref{K21}. Since $I$ is generically a complete intersection, we
have that $a^2 (G_A(I)) = - \h(I) = - (n-m+1)$. We shall
distinguish two cases.

   If $ m < n$, then $\h (I) \geq 2$, and we get that $\Kk$ is
Gorenstein if and only if $\frac{nm}{c}= \frac{n-m}{e} \in \Bbb
Z$. So there exists always at least one diagonal which is
Gorenstein by taking
 $c= nm, e= n-m$.

If $m = n$, note that $I$ is a principal ideal and so the Rees
algebra is isomorphic to a polynomial ring. Then the only diagonal
which is Gorenstein occurs when $c= n(n+1), e=1$ by Remark
\ref{K13}.}
\end{ex}

\bigskip
\section{Restrictions to the existence of
Gorenstein diagonals. Applications.} \label{L} \markboth{CHAPTER
IV. GORENSTEIN BLOW-UP SCHEMES}{RESTRICTIONS TO THE EXISTENCE.
APPLICATIONS}

\medskip
   In Section 4.1 we have proved that under the assumptions of Theorem
\ref{K5} there is just a finite set of diagonals $(c,e)$ such that
$\Kk$ is quasi-Gorenstein. Our next result shows that this holds
in general.

\begin{prop}
\label{L1}
  There exist at most a finite number of diagonals $(c,e)$ such
that $\Kk$ is quasi-Gorenstein.
\end{prop}

{\pf} Let $w_1,...,w_m \in K_R$ be a homogeneous system of
generators of $K_R$ as $R$-module with $\deg w_i =(\alpha_i,
\beta_i)$ for all $i$, and so $K_R = \sum_{i=1}^m R w_i$. Note
that since $R$ is a domain $K_R$ is torsion free. For any diagonal
$\Delta = (c, e)$ we then have by Proposition \ref{K2} that for
all $l \geq 1$
$$[K_{R_\Delta}]_l = \sum_{i=1, \dots, m, el- \beta_i \geq 0} [I^{el -
\beta_i}]_{cl- \alpha_i} w_i.$$
 If $\Rdd$ is quasi-Gorenstein there exists an integer $l$ such that
$[K_{R_\Delta}]_l \cong k$ and so $[I^{el- \beta_i}]_{cl-\alpha_i}
\not =0$ for some $i$ ($\star$). We shall distinguish two cases.

Assume first that $I$ is an equigenerated ideal in degree $d$.
Then condition ($\star$) implies that $el-\beta_i=0$ and
$cl-\alpha_i \geq 0$ or $el-\beta_i>0$ and $cl-\alpha_i \geq d(el-
\beta_i)$. If $el-\beta_i =0$, then $k \cong [K_{R_\Delta}]_l
\supset A_{cl-\alpha_i} w_i$ and since $K_R$ is torsion-free we
get $cl-\alpha_i = 0$. Hence $(c,e)$ satisfies
  $\frac{\beta_i}{e} = \frac{\alpha_i}{c} = l \in \Bbb Z$ and the
statement holds. If $el-\beta_i >0$ then $k \cong [K_{R_\Delta}]_l
\supset [I^{el-\beta_i}]_{cl-\alpha_i} w_i$ which is impossible
since $K_R$ is torsion free and $cl-\alpha_i \geq d (el-\beta_i)
$.

Assume now that $I$ is not equigenerated. Condition ($\star$)
implies that $el-\beta_i= 0$ and $cl-\alpha_i \geq 0$ or
$el-\beta_i> 0$ and $cl-\alpha_i \geq d_1(el-\beta_i)$. In the
first case we may proceed as before to get the statement. In the
second case we have that $k \cong [K_{R_\Delta}]_l \supset
[I^{el-\beta_i}]_{cl-\alpha_i} w_i$ and so  $cl-\alpha_i =
d_1(el-\beta_i)$ and $d_1 < d_2$. Then
 $\alpha_i - d_1 \beta_i = cl -d_1 el \geq c-d_1 e \geq (d-d_1)e$
since $l \geq 1$  and $c \geq de+1 >de$. Thus we obtain the
inequality $e \leq \frac{\alpha_i - d_1 \beta_i}{d-d_1}$ and for
each $e$, we have $c \leq d_1 e + \alpha_i -d_1 \beta_i$. In any
case, these inequalties hold for at most a finite number of
diagonals and so we get the result. $\B$

\medskip
 For a real number $x$, let us denote by $\lceil x \rceil = \min \, \{m
\in \Bbb Z \mid m \geq x \}$. If the Rees algebra $\Rr$ is
Cohen-Macaulay we can also give bounds for the diagonals $(c,e)$
such that $\Kk$ is quasi-Gorenstein.

\begin{prop}
\label{L2} Assume that $\h (I) \geq 2$ and $\Rr$ is
Cohen-Macaulay. Let $a= -a^2(G_A(I))$. If $\Kk$ is
quasi-Gorenstein, then $e \leq a-1$ and $c \leq n$. Moreover, if
$\dim (A/I) >0$ then $\lceil \frac{a}{e} \rceil-1 = \frac{n}{c}
\in \Bbb Z$. In particular, if $a=1$ there are no diagonals
$(c,e)$ such that $\Kk$ is quasi-Gorenstein.
\end{prop}

{\pf}  By Theorem \ref{K3}, there exists a homogeneous filtration
$\{K_m \}_{m \geq 0}$ of $K_A$ such that
 $ K_{R} \cong \bigoplus_{ m \geq 1}  K_m$ and
 $ K_{G} \cong \bigoplus_{ m \geq 1}  K_{m-1}/K_m$. Bigrading the proof
of \cite[Corollary 2.5]{TVZ}, we have that $K_m= \Hom_A(I,
K_{m+1})$ for every $m \geq 0$. Note that $K_A$ may be viewed as
an ideal of $A$. Assume that $R_{\Delta}$ is quasi-Gorenstein.
Then there is an integer $l_0$ such that $[K_{R_\Delta}]_{l_0}
\cong k$. By Proposition \ref{K2} we may find an element $f \in
[K_{e l_0}]_{c l_0} =[K_R]_{(c l_0, e l_0)}$, $f \not =0$,
$K_{R_\Delta} = {\Rdd} f$.

\noindent CLAIM: $K_{e l_0} = A f.$

To prove the claim we first show that for any $g \in K_{e l_0}$,
$g \not =0$, then $g$ has degree $\geq c l_0$. Assume the
contrary: deg $g = k <c l_0$. Then $[Ag]_{c l_0} = A_{c l_0 -k} g
\subset [K_{e l_0}]_{c l_0} \cong k$. But since $c l_0 - k > 0$,
$\dim_k A_{c l_0 -k} >1$, so we get a contradiction.

Now let $g  \in K_{e l_0}$. If $\deg g = c l_0$, then $g \in Af$
because $[K_{e l_0}]_{c l_0} = k f$. Let us assume that deg $g = k
> c l_0$. Then, for each $l>0$, $[I^{el}]_{cl} f +
 [I^{el}]_{c(l_0 +l)-k} g \subset [K_{e(l_0 + l)}]_{c(l_0+l)}
 \cong [I^{el}]_{cl}$ as $k$-vector spaces, and so
 $[I^{el}]_{c(l_0+l)-k} g \subset [I^{el}]_{cl} f$. Now let
 $I^{el} =(F_1, \dots, F_t)$ where $F_i$ is a homogeneous
 polynomial of degree $\leq del$ for all $i$,
 and set $\alpha = c (l_0+l)-k- \deg F_i$.
  Note that for $l >> 0$, $\alpha \geq c(l_0 +l)-k-del = (c-de)l +c
l_0-k >0$ and we can find $h \in A_{\alpha}$ such that $(h,f) =1$.
Then $hgF_i \in [I^{el}]_{c(l_0+l)-k} g \subset [I^{el}]_{cl} f
\subset Af$ and we have that $g F_i \in Af$ for all $i$. Thus
$I^{el} g \subset (f)$ and writing $g = d \overline g$, $f = d
\overline f$ with $(\overline f, \overline g)= 1$ we get $I^{el}
\overline g \subset    A \overline f$. If $g \not \in Af$, then
$\overline f \not \in k$ and so $I^{el} \subset (\overline f)$
which is absurd because $\h (I) \geq 2$.

Now, as grade$(I) \geq 2$ we have $K_m = K_{e l_0}$ for all $m
\leq e l_0$, which implies that $K_A = K_{e l_0}$ and so $ c \leq
c l_0 = n$. Furthermore, $e \leq e l_0 \leq \min \{ m \mid K_m
\varsubsetneq
 K_{m-1} \} -1 = a-1.$

Finally assume that $\dim (A /I) >0$. We shall distinguish two
cases. If $e=1$ we have that $K_{l_0+1} \varsubsetneq K_{l_0 }$:
If not, then $I_c \cong [K_{l_0+1}]_{c(l_0 +1)} = [Af]_{c(l_0
+1)}\cong A_c$ which is absurd if $\dim (A/I) >0$. Therefore $a =
l_0 +1 = \frac{n}{c} +1$. If $e >1$, let
 $\widetilde  \Delta =(c,1)$ and
  $\widetilde R =R_A(I^e)$.
   Note that
  $\widetilde R_
 {\widetilde  \Delta} = \Rdd$ is quasi-Gorenstein.
 Applying the case before we obtain that $-a(G_A(I^e))= \frac{n}{c} +
1$. By \cite[Proposition 2.6]{HRZ}, $a(G_A(I^e)) = [\frac{-a}{e}]=
- \lceil \frac{a}{e} \rceil$ and so $\lceil \frac{a}{e} \rceil -1
= \frac{n}{c} = l_0 \in \Bbb Z$. $\B$

\medskip
Our next result shows that in some cases the existence of a
diagonal $(c,e)$ such that $\Kk$ is quasi-Gorenstein forces the
form ring to be Gorenstein. It may be seen as a converse of
Theorem \ref{K5} for those cases.

\begin{thm}
\label{L3}
  Assume that $\Rr$ is Cohen-Macaulay, $\h(I) \geq 2$,
$l(I)<n$ and $I$ is equigenerated. If there exists a diagonal
$(c,e)$ such that $\Kk$ is quasi-Gorenstein then $G_A(I)$ is
Gorenstein.
\end{thm}

{\pf}
 Let $\Delta= (c,e)$. Assume first that
$e=1$. We have seen in the proof of Proposition \ref{L2} that
there exists a homogeneous filtration $\{K_m \}_{m \geq 0}$ of
$K_A$ such that
 $ K_{R} \cong \bigoplus_{ m \geq 1}  K_m$ and
 $ K_{G} \cong \bigoplus_{ m \geq 1}  K_{m-1}/K_m$, and an integer
$l_0 = -a(\Rdd)$ such that $K_0= \dots =K_{l_0} =Af$, with $f \in
K_R$ and deg $f=cl_0$. It is then clear that for all $m \geq  0$,
$I^m f \subset K_{l_0 +m}$ and so $[I^m]_{cm} f \subset
[K_{l_0+m}]_{c(l_0+m)} \cong [I^{m}]_{cm}$ since $\Rdd$ is
quasi-Gorenstein. This implies that $[K_{l_0+m}]_{c( l_0+m)}=
[I^{m}]_{cm}f$.

We want to show that $K_{l_0+m} = I^m f$ for all $m \geq 0$.
Suppose that there exists $m_0$ such that
 $I^{m_0} f \varsubsetneq K_{l_0+m_0}$.
Then let $g \in K_{l_0+m_0} $, $g \not \in I^{m_0}f$ be a
homogeneous element of degree $k$. Note that from the inclusion
$K_{l_0+m_0} \subset K_{l_0}= Af$ we also have
 $g = f \overline g$ with  $\overline g \not \in I^{m_0}$.

    If $k \geq c(l_0+m_0)$ then for any $m > m_0$ we have $I^m f +
I^{m-m_0} g \subset K_{l_0+m}$ and so $[I^m]_{cm} f +
[I^{m-m_0}]_{c(l_0+m)-k}g \subset [K_{l_0+m}]_{c(l_0+m)} \cong
[I^m]_{cm}$. Hence $[I^{m-m_0}]_{c(l_0+m)-k} g \subset [I^m]_{cm}
f$ and we get that $[I^{m-m_0}]_{c(l_0+m)-k} \overline g \subset
[I^m]_{cm}$. Let $\lambda = c(l_0+m)-k-d(m-m_0)= (c-d)m +cl_0+
dm_0-k$. For $m>>0$ we have that $\lambda >0$. Then, if
$A_{\lambda} \overline g \subset I^{m_0}$ we would have that
$\overline g \in (I^{m_0})^*$, the saturation of  $I^{m_0}$. Since
$G_A(I)$ is Cohen-Macaulay, we have that the inequality of Burch
becomes an equality by \cite[Proposition 3.3]{EH}, that is,
 $$ \inf_{j \geq 0} \{ {\rm depth}(A/I^j) \} = \dim A - l(I).$$
Since $l(I) < n$, we then get $\depth \, A/I^{m_0} >0$, and so
$(I^{m_0})^* =I^{m_0}$. Hence $\overline g \in I^{m_0}$, which is
a contradiction. We may conclude that there exist $\lambda >0$, $h
\in A_{\lambda}$ such that $ \overline g  h \not \in I^{m_0}$. On
the other hand, $\overline g h [I^{m-m_0}]_{d(m-m_0)} \subset
\overline g [I^{m-m_0}]_{c(l_0+m)-k} \subset [I^m]_{cm}$. Since
$I$ is equigenerated we get $\overline g h I^{m-m_0} \subset I^m$.
Therefore, $\overline g h \in (I^m: I^{m-m_0})= I^{m_0}$ because
$R$ is Cohen-Macaulay. This is  a contradiction.

    If $k < c(l_0+m_0)$, let us write $k= c(l_0+m_0)-s$ with $s>0$.
Then $A_s g \subset [K_{l_0+m_0}]_{c(l_0+m_0)}=[I^{m_0}]_{c
m_0}f$, and $\overline g \in (I^{m_0})^* = I^{m_0}$ which, as
before, is a contradiction.

    Hence we have proved that $K_{l_0+m}= I^m f$ for all $m \geq 0$, so
$$K_R= f(At \oplus \cdots \oplus A t^{l_0} \oplus I t^{l_0+1}
\oplus \cdots),$$ i.e. $K_R$ has the expected form. By
\cite[Theorem 4.2]{TVZ} this implies that both $R_A(I^{l_0})$ and
$G_A(I)$ are Gorenstein.

Finally assume $e>1$, and denote by $\widetilde  \Delta =(c,1)$
and $\widetilde R =R_A(I^e)$. Then $\widetilde R_ {\widetilde
\Delta} = \Rdd$ is quasi-Gorenstein and so there exists $l_0$ such
that $R_A(I^{el_0})$ is Gorenstein. By \cite[Theorem 4.2]{TVZ}
this implies again that $G_A(I)$ is Gorenstein. $\B$

\medskip
\begin{ex} (Room surfaces)
\label{L4}
  {\rm  Let $k$ be an algebraically closed field.
Set $t= { d+1 \choose 2}$, with $d \geq 2$. Let $P_1, \dots, P_t$
be a set of $t$ distinct points in $\Bbb P^2_k$ which do not lie
on a curve of degree $d-1$. We assume further that there is not a
subset of $d$ points on a line if $d \geq 3$. We are going to
study the rational projective surfaces which arise as embeddings
of blowing-ups of ${\Bbb P}_k^2$ at this set of points via the
linear system $I_{d+1}$.

Let  $I$ be the ideal defining the set of points $\{ P_1, \dots,
P_t \}$. Since the points are not on a curve of degree $d-1$, we
have that $I$ is generated by forms in degree $d$ \cite{GG}, and
so for any $c \geq d+1$ the linear system $I_c$ gives a projective
embedding of the blow-up. For $c=d+1$ the surface obtained is
called Room surface.

Assume $d \geq 3$. Since there are not $d$ points on a line, we
also have that the rational map defined by the linear system $I_d$
give an embedding of the blow-up in the projective space $\Bbb
P^{d}_k$ (see \cite{GG}), and the resulting surface is called
White surface. A. Gimigliano proved that White surfaces have the
defining ideal given by the $3 \times 3$ minors of a $ 3 \times d$
matrix of linear forms, and it has a linear minimal graded free
resolution which comes from the Eagon-Northcott complex
\cite[Proposition 1.1]{G}. By applying Theorem \ref{B77}, we
obtain $a(k[I_d]) = -1$ and so by Lemma \ref{ED9} the reduction
number of $I$ is $r(I) = a(k[I_d]) + l(I) = -1+3 = 2$. Moreover,
the analytic deviation of $I$ is ${\rm ad}(I) = l(I) - \h(I) = 1$
and $I$ is generically a complete intersection ideal. So according
to \cite{GN} we may conclude that $G_A(I)$ is Cohen-Macaulay and
hence $a^2(G_A(I)) = r(I) - \h(I) -1 = -1$ by \cite[Proposition
2.4]{GH}. Therefore, $R_A(I)$ is also Cohen-Macaulay by using
Ikeda-Trung's criterion. From Proposition \ref{L2} we get that
there are not diagonals $(c,e)$ such that $\Kk$ is Gorenstein,
that is, there are not Gorenstein embeddings for the blow-up. In
particular, $k[I_{d+1}]$ is not Gorenstein for $d \geq 3$.

If $d=2$, by choosing the points to be [1:0:0], [0:1:0] and
[0:0:1], we have $I= (X_1 X_2, X_1 X_3, X_2 X_3)$. Notice that $I$
has $\mu (I) = 3 =\h I+1$ and $A/I$ is Cohen-Macaulay. Moreover,
$\mu (I_{\frak p}) \leq \h(\frak p)$ for any prime ideal $\frak p
\supset I$. Then $G_A(I)$ is Gorenstein with and $a^2(G_A(I)) = -
\h(I) = -2$ by \cite{HRZ}. Now, according to Theorem \ref{K5},
$\Kk$ is quasi-Gorenstein if and only if $\frac{3}{c} =
\frac{1}{e} \in \Bbb Z$. Hence $(3,1)$ is the only diagonal with
the quasi-Gorenstein property. This embedding  corresponds to the
del Pezzo sestic surface in ${\Bbb P}^6$.}
\end{ex}

\medskip
With more generality, we may consider the blow-up of $\Bbb P^2_k$
at a set of $t$ arbitrary distinct points.

\begin{ex}
\label{L5} {\rm  Let $k$ be an algebraically closed field. Let
$P_1, \dots, P_t$ be a set of $t$ distinct points in $\Bbb P^2_k$,
and let $I ={\cal P}_1 \cap \dots \cap {\cal P}_t$, where ${\cal
P}_i \subset A=k[X_1, X_2, X_3]$ is the defining ideal of $P_i$.
Now we consider the surfaces which arise as embeddings of the
blow-up of $\Bbb P^2_k$ at these points via the linear systems
$(I^e)_c$. We want to study the Gorenstein property of the rings
$\Kk$.

Set $d = {\rm reg}(I)$. We will assume that $P_1, \dots, P_t$ do
not lie on a curve of degree $d-1$ and that there is not a subset
of $d$ points on a line. Then, $I$ is generated by forms in degree
$d$ and $I_d$ defines a projective embedding of the blow-up
\cite[Theorem A]{GG}.

On the other hand, observe that $R_{A_{\fp}}(I_{\fp})$ is
Cohen-Macaulay for all $\fp \in \Proj (A)$. Let ${\cal L}=
\widetilde I {\cal O}_X$, ${\cal M} = \pi^*{\cal O}_{\Bbb
P^2}(1)$, where $\pi: X \longrightarrow \Bbb P^2_k$ is the blow-up
of $\Bbb P^2_k$ along $\cal I$. Then we have that $R^j \pi_* {\cal
L}^e = 0$ for all $e \geq 0$, $j >0$ and $\pi_* {\cal L}^e =
\widetilde I^e$ for all $e \geq 0$ by Corollary \ref{G02}.
Therefore, for any $s \geq 0$, ${\mit \Gamma} (X, {\cal L}^s
\otimes {\M}^{sd}) = {\mit \Gamma} ({\Bbb P^2}, \widetilde
I^s(sd)) = (I^s)^*_{sd}$ and $H^{i}(X, {\cal L}^s \otimes
{\M}^{sd}) = H^{i}({\Bbb P^2}, \widetilde I^s(sd)) =
H^{i+1}_{\fm}(I^s)_{sd}$ for $i \geq 1$. By \cite[Theorem 1.1 and
Corollary 1.4]{GGP}, we have
 $a_*(I^s) < \reg (I^s) \leq sd$ and $(I^{s})^*_{sd}= (I^{s})_{sd}$.
Then, by Remark \ref{G1} we get $H^i_{\fm}(k[I_d])_s=0$ for all $s
\geq 0$. Furthermore, recall that the fiber cone $F$ of $I$
coincides with $k[I_d]$ because $I$ is generated in degree $d$, so
we have $a_*(F) \leq -1$, and then $r(I) \leq \max_{i \leq 3} \{
a_i(F)+i \} \leq 2$ by Lemma \ref{ED9} . The analytic deviation of
$I$ is ${\rm ad}(I) = l(I) - \h (I) = 1$ and $I$ is generically a
complete intersection ideal, so we may conclude by \cite{GN} that
$G_A(I)$ is Cohen-Macaulay  and hence by \cite[Proposition
2.4]{GH}
$$a^2(G_A(I)) \leq r(I) - \h (I)-1 \leq  -1.$$
So $R_A(I)$ is also Cohen-Macaulay by Ikeda-Trung's criterion.

If $r(I)=2$ then $a^2(G_A(I))=-1$ by \cite[Proposition 2.4]{GH}.
Then, according to Proposition \ref{L2}, there are not
quasi-Gorenstein diagonals $\Kk$.

Otherwise, $r(I) \leq 1$. By \cite[Theorem 1.3]{GN}, the case
$r(I)=1$ is not possible, so $r(I)=0$ and then $G_A(I)$ is
Gorenstein. Furthermore, $a^2(G_A(I))= -\h (I) = - 2$ by
\cite[Proposition 2.4]{GH}. Therefore, by Theorem \ref{K5}, $\Kk$
is quasi-Gorenstein if and only if $\frac{3}{c}= \frac{1}{e} \in
\Bbb Z$. So $k[I_3]$ is the only quasi-Gorenstein diagonal. }
\end{ex}

\medskip
\begin{rem}
{\rm All throughout this chapter we have treated the case where
$A=k[{\x}]$ is the polynomial ring and $I$ is a homogeneous ideal
in $A$ satisfying $r \leq n$ or the assumptions in Remark
\ref{K21}. This set up was used to study the relationship between
the canonical module of the Rees algebra $R$ and the canonical
modules of its diagonals $R_\Delta$. Now let $A$ be an arbitrary
standard $k$-algebra and let $I$ be a homogeneous ideal in $A$
generated by $r$ forms in degree $\leq d$. Set $\nn = \dim A$. For
any $c \geq de+1$, from the Mayer-Vietoris sequence (see
Proposition \ref{W2}) we have a graded monomorphism
$$\psi: (K_R)_\Delta \rightarrow K_{R_\Delta} $$
such that
\begin{enumerate}
\item If $l(I)< \nn$ or $r< \nn$ or $I$ equigenerated, then
 $\psi_s$ is an isomorphism for any $s>0$.
\item Assume $l(I)< \nn$ or $a_2^*(R)< 0$. If $\h (I) \geq 2$,
$H^{\nn}_\fm (A)_0 = 0$, $a(A) < c$, then $\psi_s$ is an
isomorphism for any $s \leq 0$.
\item Assume that $R$ is Cohen-Macaulay. Then

           $\psi $ isomorphism  $\iff$
       $\cases{H^{\nn}_\fm (A)_0=0 \cr
       H^{\nn}_\fm (I^{es})_{cs}=0 &  for $s>0$ \cr
       H^{\nn}_{{\cal M}_2}(R)_{(cs,es)} = 0& for $s \in \Bbb Z$ \cr}$

\end{enumerate}

In the cases where $\psi$ is an isomorphism, some of the results
of the chapter can be extended. For instance, if $R$ is
Cohen-Macaulay and $G$ is quasi-Gorenstein, for any $c,e$ such
that $\psi $ is an isomorphism and $\frac{n}{c}= \frac{a-1}{e} \in
\Bbb Z$ the ring $\Kk$ is quasi-Gorenstein.}
\end{rem}

\chapter*{$\;\;\;$}

\newpage

\medskip
\chapter{The a-invariants of the powers of an ideal}
\typeout{The a-invariants of the powers of an ideal}

\bigskip

Our aim in this chapter is to study in more detail the bigraded
$a$-invariant and the bigraded regularity of any finitely
generated bigraded $S$-module $L$, for $S= k[{\x},{\y}]$ the
polynomial ring with deg$(X_i)= (1, 0)$, deg$(Y_j)= (0, 1)$.

\medskip
In Section 5.1 we will give a new description of the
$a_*$-invariant ${\bf a}_*(L)$ of $L$ and the regularity ${\rm \bf
reg} (L)$ of $L$ by means of the $a_*$-invariants and the
regularities of the graded $S_1$-modules $L^e$ and the graded
$S_2$-modules $L_e$.

\medskip
This result is used in Section 5.2 to study the behaviour of the
$a_*$-invariant of the powers of a homogeneous ideal in the
polynomial ring. In particular, we will bound it for several
families of ideals such as equimultiple ideals and strongly
Cohen-Macaulay ideals. Those results will be then applied to
determine Cohen-Macaulay diagonals of their Rees algebras.

\medskip
The last section is devoted to study the regularity of homogeneous
ideals $I$ in the polynomial ring $S$. First, we will provide a
bigraded version of the well-known Bayer-Stillman's Theorem
characterizing the regularity of $I$ in terms of generic forms.
After that, similarly to the graded case, we define the bigraded
generic initial ideal ${\rm \bf gin} \, I$ of $I$ and we establish
its basic properties. In the graded case, a classical result due
to D. Bayer and M. Stillman states the existence of an order such
that for any homogeneous ideal $I$ it holds $\reg I = \reg ({\rm
gin} \, I)$. We will show that the analogous bigraded statement is
not true. We finish the chapter by explaining how these results
can be used to study the Koszulness of the diagonals $\Kk$.

\bigskip
\section{The a-invariant of a standard bigraded algebra}
\label{BB} \markboth{CHAPTER IV. a-INVARIANTS OF THE POWERS OF AN
IDEAL}{a-INVARIANT OF A STANDARD BIGRADED ALGEBRA}
\bigskip

Let $S=k[{\x},{\y}]$ be the polynomial ring over a field $k$ in
$n+r$ variables with deg$(X_i)= (1, 0)$, deg$(Y_j)= (0, 1)$, and
let us distinguish two bigraded subalgebras: $S_1=k[{\x}]$, $S_2=
k[{\y}]$, with homogeneous maximal ideals ${\frak m}_1= (X_1, ...,
X_n)$, ${\frak m}_2 = ({\y}) $ respectively. Given $e \in \Bbb Z$
and a bigraded $S$-module $L$, recall that we may define the
graded $S_1$-module $L^e= \bigoplus_{i \in \Bbb Z} L_{(i,e)}$ and
the graded $S_2$-module $L_e= \bigoplus_{j \in \Bbb Z} L_{(e,j)}$.

\medskip
The first result shows how to compute the bigraded $a_*$-invariant
of any finitely generated bigraded $S$-module $L$ by means of the
$a_*$-invariants of the graded $S_1$-modules $L^e$ and the graded
$S_2$-modules $L_e$. Namely,

\medskip
\begin{thm}
\label{BB1} Let $L$ be a finitely generated bigraded $S$-module.
Then :
\begin{enumerate}
\item $a_*^1(L) = \max_e \{a_*(L^e) \} =
                  \max_e \{a_*(L^e) \mid e \leq a_*^2(L) + r \}$.
\item $a_*^2(L) = \max_e \{a_*(L_e) \} =
                  \max_e \{a_*(L_e) \mid e \leq a_*^1(L) + n \}$.
\end{enumerate}
\end{thm}

{\pf} Let us consider
$$ 0 \to D_t \to \dots \to D_1 \to D_0  \to L \to 0$$
the minimal bigraded free resolution of $L$ over $S$, where $D_p =
\bigoplus_{(a,b) \in {\Omega}_p} S(a,b)$. We have $a_*^1 (L) =
\max \,\{\, -a \mid (a,b) \in  {\Omega}_L \} - n$ by Theorem
\ref{B77}.

Let us denote by $\underline \beta = (\beta_1, \ldots, \beta_r)
\in \Bbb N^r$ and $\mid \underline \beta \mid = \beta_1 + \cdots +
\beta_r$. By applying the functor $(\,\,)^e$ to the resolution
note that

\vspace{1.5mm}

\hspace{10mm} $S(a,b)^e = \bigoplus_{i \in \Bbb Z} S(a,b)_{(i,e)}
               = \bigoplus_{i \in \Bbb Z} S_{(a+i,b+e)}  $

\vspace{1.5mm}

\hspace{24mm} $ = \bigoplus_{i \in \Bbb Z}
      \bigoplus_{\mid \underline \beta \mid = b+e}
         [S_1]_{a+i} Y_1^{\beta_1} \dots Y_r^{\beta_r} $

\vspace{1.5mm}

\hspace{24mm} $= S_1(a)^{\rho_{ab}^e}$

\medskip
\noindent for certain ${\rho_{ab}^e} \in \Bbb Z$
(${\rho_{ab}^e}=0$ if $b+e<0$). In this way, we have obtained a
graded free resolution of $L^e$ over $S_1$
$$ 0 \to D_t^e \to \dots \to D_1^e \to D_0^e \to L^e   \to 0 \, ,$$
with $D_p^e = \bigoplus_{(a,b) \in {\Omega}_p} S_1(a)^
{\rho_{ab}^e}$. The minimal graded free resolution of $L^e$ may be
obtained by picking out some terms \cite[Exercise 20.1]{E}.
Therefore,
$$a_* (L^e) \leq \max \{-a \mid (a,b) \in  \Omega_L \}-n = a_*^1(L).$$

Now let $ \alpha = \max \, \{\,-a \mid (a,b) \in \Omega_L \}$. Let
$p$ be the first place in the resolution of $L$ with a shift of
the form $(-\alpha, b)$, and let $\beta$ be one of these $-b$'s.
We are done if we prove that $-\alpha$ is a shift which appears in
the place $p$ of the minimal graded free resolution of
$L^{\beta}$. Note that it is enough to show that
 $$ \Tor^S_p (S / {\fm_1} S, L)_{(\alpha,\beta)} =
 \Tor^{S_1}_p (k, L^\beta)_{\alpha}  \not =0.$$
Let us consider
 $$ D_{p+1} \stackrel{\psi_{p+1}}{\longrightarrow} D_p
\stackrel{\psi_{p}}{\longrightarrow} D_{p-1}$$ the differential
maps appearing in the resolution of $L$. Tensorazing by $S/
{\fm_1} S$, we have the sequence
 $$ D_{p+1}/ \fm_1 D_{p+1} \stackrel{\overline
\psi_{p+1}}{\longrightarrow}
  D_p/ \fm_1 D_p \stackrel{\overline \psi_{p}}{\longrightarrow}
   D_{p-1} / \fm_1 D_{p-1} \; .$$
Now let us take $v \in D_p$ one of the elements of the homogeneous
basis of $D_p$ as free $S$-module with $\deg (v) = (\alpha,
\beta)$. If $w_1, \dots, w_s$ is the homogeneous basis of
$D_{p-1}$, we can write
$$ \psi_p (v) = \sum_{j=1}^{s} \lambda_j w_j ,$$
with $\lambda_j \in {\cal M}$ homogeneous. Set $\deg (w_j) =
(\alpha_j, \beta_j)$. By taking into account the way we have
choosen $\alpha$ and $p$, we have $\alpha > \alpha_j$ for any $j$.
Therefore the first component of the degree of $\lambda_j$ is
positive, so $\lambda_j \in \fm_1 S$. We conclude $ \overline
\psi_p (v) = 0$, that is, $v \in \Ker \, \overline \psi_p$.
Furthermore, notice that
 $v \not \in {\rm Im} \, \overline \psi_{p+1}$ because Im $\overline
\psi_{p+1} \subset {\cal M} (D_p / \fm_1 D_p)$. So $v \in \Tor^S_p
(S / \fm_1 S, L)_{(\alpha,\beta)}$, $v \not = 0$. By symmetry, we
get $(ii)$. $\Box$

\medskip
Next we are going to consider the bigraded regularity of a
finitely generated bigraded $S$-module $L$. Assume that
$$0 \to D_t \to \dots \to D_1 \to D_0  \to L \to 0,$$
with $D_p =\bigoplus_{(a,b) \in {\Omega}_p} S(a,b)$, is the
minimal bigraded free resolution of $L$ over $S$. The bigraded
regularity of $L$ is defined by ${\rm \bf reg}(L) = (\reg_1 L,
\reg_2 L)$, where
$$\reg_1 L= \max_p \, \{-a-p : (a,b) \in {\Omega}_p \}$$
$$\reg_2 L= \max_p \, \{-b-p : (a,b) \in {\Omega}_p \}.$$

\medskip
Let $A=k[{\x}]$ be the polynomial ring with the usual grading. For
any finitely generated graded $A$-module $L$, it is well known
that
$$\reg (L) = \max_{p \geq 0} \,\{ \, t_p(L)-p \}=
 \max_{p \geq 0} \, \{ \, a_p(L)+ p\} .$$
This equality does not hold in the bigraded case. For instance,
let us consider $f_1, \dots, f_r \in A$ a regular sequence of
forms in degree $d$, and $I=(f_1, \dots, f_r)$. Let
$S=k[{\x},{\y}]$ be the polynomial bigraded by setting
$\deg(X_i)=(1,0)$, $\deg(Y_j)=(d,1)$, and let $R$ be the Rees
algebra of $I$. Since $R$ is Cohen-Macaulay, we immediately get
$a^2_{n+1}(R)=-1$, $a^2_{i}(R)=0$ for $i \not = n+1$. Furthermore,
the Eagon-Northcott complex gives the bigraded minimal free
resolution of $R$ over $S$:
$$0 \to D_{r-1} \to \dots \to D_0=S \to R \to 0,$$
with $D_p = \bigoplus_{m=1}^p S (-(p+1)d, -m)^{r \choose p+1}$ for
$p \geq 1$. Therefore,
$$\max_{p \geq 0} \,\{ \, t_p^2 (R)-p \}= 0 $$
$$\max_{p \geq 0} \,\{ \, a_p^2(L)+ p\} = n,$$
which are different.

\medskip
The following result shows that the regularity of $L$ can also be
described by means of the regularity of the graded $S_1$-modules
$L^e$ and the graded $S_2$-modules $L_e$. Namely,

\begin{thm}
\label{BB2} Let $L$ be a finitely generated bigraded $S$-module.
Then :
\begin{enumerate}
\item $\reg_1 (L) = \max_e \{\reg (L^e) \}
      = \max_e \{\reg (L^e) \mid e \leq a_*^2(L) + r \}$.
\item $\reg_2 (L) = \max_e \{\reg (L_e) \}
      = \max_e \{\reg (L_e) \mid e \leq a_*^1(L) + n \}$.
\end{enumerate}
\end{thm}

{\pf} The proof follows the same lines as Theorem \ref{BB1}. By
applying the functor $(\;)^e$ to the minimal bigraded free
resolution of $L$ over $S$, we obtain a graded free resolution of
$L^e$ over $S_1$
$$ 0 \to D_t^e \to \dots \to D_1^e \to D_0^e \to L^e   \to 0 \, ,$$
with $D_p^e = \bigoplus_{(a,b) \in {\Omega}_p} S_1(a)^
{\rho_{ab}^e}$. Since the minimal graded free resolution of $L^e$
is then obtained by picking out some terms, we have
$$\reg(L^e) \leq \max_p \{-a -p \mid (a,b) \in  \Omega_p \} =
\reg_1 L.$$ Hence $\max_e \{\reg (L^e) \} \leq \reg_1 L$. To prove
the equality, let us take  $(a, b) \in \Omega_p$  such that
$\reg_1 L = -a-p$, and set $ \alpha =-a $, $\beta = -b$. We are
done if we prove $a \in \Omega_{p, L^\beta}$, that is, $a$ is a
shift which appears in the place $p$ of the minimal graded free
resolution of $L^{\beta}$. So we want to show that
 $$ \Tor^S_p (S / {\fm_1} S, L)_{(\alpha,\beta)} =
 \Tor^{S_1}_p (k, L^\beta)_{\alpha}  \not =0.$$
Let us consider
$$ D_{p+1} \stackrel{\psi_{p+1}}{\longrightarrow} D_p
\stackrel{\psi_{p}}{\longrightarrow} D_{p-1}$$ the differential
maps appearing in the resolution of $L$. Tensorazing by $S/
{\fm_1} S$, we have the sequence
 $$ D_{p+1}/ \fm_1 D_{p+1} \stackrel{\overline
\psi_{p+1}}{\longrightarrow}
  D_p/ \fm_1 D_p \stackrel{\overline \psi_{p}}{\longrightarrow}
   D_{p-1} / \fm_1 D_{p-1} .$$
Now let $v \in D_p$ be an element of the homogeneous basis of
$D_p$ as free $S$-module with $\deg (v) = (\alpha, \beta)$. If
$w_1, \dots, w_s$ is the homogeneous basis of $D_{p-1}$, we can
write
$$ \psi_p (v) = \sum_{j=1}^{s} \lambda_j w_j ,$$
with $\lambda_j \in {\cal M}$ homogeneous. Set $\deg (w_j) =
(\alpha_j, \beta_j)$. Since $\alpha -p \geq \alpha_j-(p-1)$ for
any $j$, we have that $\alpha > \alpha_j$, and so the first
component of the degree of $\lambda_j$ is positive. Therefore
$\lambda_j \in \fm_1 S$, and we can conclude $ \overline \psi_p
(v) = 0$, that is, $v \in \Ker \, \overline \psi_p$. It is clear
that
 $v \not \in {\rm Im} \, \overline \psi_{p+1}$ because Im $\overline
\psi_{p+1} \subset {\cal M} (D_p / \fm_1 D_p)$. So $v \in \Tor^S_p
(S / \fm_1 S, L)_{(\alpha,\beta)}$, $v \not = 0$. We get $(ii)$ by
symmetry. $\Box$

\bigskip
\section{The a-invariants of the powers of an ideal}
\label{Yy} \markboth{CHAPTER IV. a-INVARIANTS OF THE POWERS OF AN
IDEAL}{a-INVARIANTS OF THE POWERS OF AN IDEAL}

\medskip
  Let $A= k[{\x}]$ be the usual polynomial ring over a field $k$, and
let $I$ be a homogeneous ideal in $A$. Recently, the question of
how the regularity changes with the powers of $I$ has been studied
by many authors. I. Swanson in \cite{S} proved that there exists
an integer $B$ such that $\reg (I^e) \leq Be$ for all $e$. The
problem is then to make $B$ explicit.

\medskip
For ideals such that dim$(A/I)= 1$, A. Geramita, A. Gimigliano and
Y. Pitteloud \cite{GGP} and K. Chandler \cite{C} had shown that
$\reg(I^e) \leq \reg(I) \, e$; and this bound also holds for
Borel-fixed monomial ideals by using the Eliahou-Kervaire
resolution \cite{EK}.

\medskip
Another kind of bound is given by R. Sj\"{o}gren \cite{Sj}: If $I$ is
an ideal generated by forms in degree $\leq d$ with dim$(A/I) \leq
1$, then $\reg (I^e) < (n-1)de$. Also A. Bertram, L. Ein and R.
Lazarsfeld \cite{BEL} gave a bound for the regularity of the
powers of an ideal in terms of the degrees of its generators. More
explicitly, if $I$ is the ideal of a smooth complex subvariety $X$
in ${\Bbb P}_{\Bbb C}^{n-1}$ of codimension $c$ and $I$ is
generated by forms in degrees $d_1 \geq d_2 \geq \dots \geq d_r$ ,
then
$$ H^i ({\Bbb P}_{\Bbb C}^{n-1}, {\cal I}^e (k)) = 0, \,\,
\,\forall i \geq 1,\, \forall k \geq ed_1 + d_2 + \dots +d_c
-(n-1).$$ This result has been improved by A. Bertram \cite{B} for
some determinantal varieties.

\medskip
   Recently, work by S.D. Cutkosky, J. Herzog and N.V.
Trung \cite{CHT}, V. Kodiyalam \cite{Ko2} and O. Lavila--Vidal
(see Theorem \ref{S8}) provides by different methods bounds for
arbitrary graded ideals by means of the degrees of the generators
similar to the ones given in \cite{Sj} and \cite{BEL}. Namely, if
$I$ is a graded ideal generated by forms in degree $\leq d$, then
there exists $\beta$ such that
$$\reg(I^e) \leq de + \beta , \; \forall e.$$

\medskip
We are also interested in the behaviour of the $a_*$-invariant of
the powers of $I$, which can be used to apply the criteria seen in
Chapter 3 for the Cohen-Macaulayness of the diagonals. We have
already proved in Theorem \ref{S8} the existence of an integer
$\alpha$ such that $a_*(I^e) \leq de+ \alpha$ for all $e$. Our
first purpose will be to find for any graded ideal an explicit
$\alpha$. Furthermore, for equigenerated ideals we will compute
the best $\alpha$ we can take in terms of an appropiate
$a$-invariant of the Rees algebra. After that, these results will
be applied to give bounds for the $a_*$-invariant of the powers of
several families of ideals such as equimultiple ideals and
strongly Cohen-Macaulay ideals. Finally, we will use those bounds
to study the Cohen-Macaulay property of the diagonals of the Rees
algebra.

\medskip
Let $k$ be a field, $A$ a standard noetherian graded $k$-algebra,
$\nn = \dim A$. Then $A$ has a presentation $A=k[{\x}]/K=
k[{\xx}]$, where $K$ is a homogeneous ideal and each $X_i$ has
degree 1. Let $I$ be a homogeneous ideal in $A$ generated by forms
of degree $\leq d$. From Theorem \ref{S8}, there exists $\alpha$
such that
$$a_*(I^e) \leq de+ \alpha, \; \forall e.$$
Now let us assume that $I$ is generated by forms in degree $d$. By
defining $\varphi (p, q) = (p-dq, q)$, we have that $R^\varphi$ is
a standard bigraded $k$-algebra with $[R^\varphi]_{(p, q)} =
R_{(p+dq, q)}$. The next result precises the best $\alpha$ we can
take.

\begin{thm}
\label{YY2} Let $I$ be a homogeneous ideal of $A$ generated by
forms in degree $d$. Set $l = l(I)$. Then
\begin{enumerate}
\item $a_*^1 (R^\varphi) = \max_{e} \{\,a_*(I^e)-de \} =
 \max \{\,a_*(I^e)-de  \mid  e \leq a_*^2(R) + l\}.$
\item $\reg_1 (R^\varphi) = \max_{e} \{\, \reg(I^e)-de \} =
 \max \{\, \reg (I^e)-de  \mid  e \leq a_*^2(R) + l\}.$
\end{enumerate}
\end{thm}

{\pf} We may assume that $k$ is infinite (tensorazing by $k(T)$).
Then there exists a minimal reduction $J$ of $I$ generated by $l$
forms in degree $d$. By considering the polynomial ring
$S=k[{\x},Y_1, \dots,Y_l]$, we have a natural epimorhism $S
\rightarrow R_A(J)$. Then $R_A(J)$ is a finitely generated
bigraded $S$-module, and so $R=R_A(I)$ because it is a finitely
generated $R_A(J)$-module. Note that $S^\varphi$ is standard and
$R^\varphi$ is a finitely generated bigraded $S^\varphi$-module,
so according to Theorem \ref{BB1}
$$a_*^1(R^\varphi) = \max_e \{a_*([R^\varphi]^e) \} =
              \max_e \{a_*([R^\varphi]^e) \mid e \leq a_*^2(R^\varphi) +
l\}.$$ First, observe that $a_*^2(R^\varphi) = a_*^2(R)$ by Lemma
\ref{B1}. Moreover, for each $e \geq 0$ we have $[R^\varphi]^e =
\bigoplus_i (I^e)_{i+de}= I^e (de)$, so $a_*((I^e)^\psi) =
a_*(I^e)-de$. The proof of $(ii)$ follows the same lines. $\Box$

\begin{rem}
\label{YY22} {\rm Let $I$ be a homogeneous ideal in $A$ generated
by forms in degree $d$. By repeating the previous arguments for
the form ring, we also get
\begin{enumerate}
\item $a_*^1 (G^\varphi) = \max_{e} \{\,a_*(I^e/I^{e+1})-de \} $

\hspace{12mm} $=\max_{e} \{\,a_*(I^e/I^{e+1})-de  \mid  e \leq
a_*^2(G)
               + l\}.$
\item $\reg_1 (G^\varphi) = \max_{e} \{\, \reg (I^e/I^{e+1})-de \} $

\hspace{16mm} $= \max_{e} \{\, \reg (I^e/I^{e+1})-de  \mid  e \leq
a_*^2(G) + l\}.$
\end{enumerate}
}
\end{rem}

\begin{ex}
\label{YY4} {\rm Let $I \subset A= k[X_1, X_2, X_3, X_4]$ be the
defining ideal of the twisted cubic in $\Bbb P^3_k$, that is,
$$I= (X_1X_4-X_2X_3, X_2^2-X_1X_3, X_3^2-X_2X_4).$$
It is well known that $I$ is the ideal of the Veronese embedding
of $\Bbb P^1_k$ in $\Bbb P^3_k$ :
$$ \begin{array}{cccc}
 \Bbb P^1_k & \stackrel{\mu}{\longrightarrow} & \Bbb P^3_k  \\
(u:v) & \longmapsto & (u^3:u^2v:uv^2:v^3) \;.
\end{array}$$
$I$ is licci because it is linked to $J=(X_1, X_2)$ by the regular
sequence $\underline \alpha = X_2^2-X_1X_3, X_3^2-X_2X_4$
\cite[Example 2.3]{Ul},  so $I$ is a strongly Cohen-Macaulay ideal
\cite[Theorem 1.14]{Hu2}. Since $I$ is a prime ideal, we easily
get $\mu(I_\fp) \leq \h (\fp)$ for any prime ideal $\fp \supseteq
I$. Therefore $R_A(I)$ is Cohen-Macaulay by \cite[Theorem
2.6]{HSV1}, so in particular $a_*^2(R_A(I))= -1$. On the other
hand, $I$ is an ideal generated by forms of degree 2 with $l(I) =
\mu(I)=3$. By using CoCoa \cite{CNR}, we have that the minimal
graded free resolutions of $I$ and $I^2$ are:
$$ 0 \to A(-3)^2 \to A(-2)^3 \to I \to 0 \,\,,$$
$$ 0 \to A(-6) \to A(-5)^6 \to A(-4)^6 \to I ^2 \to 0 \;,$$
so according to Theorem \ref{B77} we have $a_*(I)= -1$, $a_*(I^2)
=2$. By Theorem \ref{YY2} we get
$$a_*(I^e) \leq 2(e-1), \; \forall e.$$
Furthermore, notice that since $\reg(I)= 2$, $\reg (I^2) =4$ we
also get
$$\reg(I^e) \leq 2e, \; \forall e.$$
Therefore, we have that $I^e$ has a linear resolution for any $e
\geq 1$. This has already been proved by A. Conca \cite{C1} by
different methods. }
\end{ex}

\begin{rem}
\label{YY3} {\rm Let $S=k[{\x},{\y}]$ be the polynomial ring
bigraded by setting $\deg (X_i)=(1,0)$, $\deg (Y_j)=(d_j, 1)$,
with $d_1, \dots, d_r \in \Bbb Z_{\geq 0}$, and $u= \sum_{j=1}^r
d_j$. For a finitely generated bigraded $S$-module $L$, let us
consider the minimal bigraded free resolution of $L$ over $S$
$$ 0 \to D_t \to \dots \to D_1 \to D_0 \to L \to 0,$$
where $D_p =\bigoplus_{(a,b) \in {\Omega}_{p}} S(a,b)$. By
applying the functor $(\,\,)^e$, note that
$$S(a,b)^e = \bigoplus_{\mid \underline \beta \mid = b+e}
S_1(a-d_1 \beta_1- \dots -d_r \beta_r),$$ where $\underline \beta
= (\beta_1, \ldots, \beta_r) \in \Bbb N^r$ and $\mid \underline
\beta \mid = \beta_1 + \cdots + \beta_r$. So we get a graded free
resolution of $L^e$ over $S_1$
$$ 0 \to D_t^e \to \dots \to D_1^e \to D_0^e \to L^e   \to 0 \, ,$$
with $D_p^e = \bigoplus_{(a,b) \in {\Omega}_{p}}
      \bigoplus_{\mid \underline \beta \mid = b+e}
           S_1 (a-d_1 \beta_1- \dots -d_r \beta_r)$.
The minimal graded free resolution is then obtained by picking out
some terms. Therefore, for any $i \leq n$ we have that
$$a_i (L^e) \leq
 \max \{d_1 \beta_1+ \dots +d_r \beta_r-a \mid (a,b) \in
{\Omega}_{n-i}, \mid \underline \beta \mid = b+e \}-n  $$

\hspace{12mm}  $\leq de -n + \max \{db-a \mid (a,b) \in
{\Omega}_{n-i} \}.$

\vspace{2mm}

\noindent Therefore, $a_*(L^e) \leq de -n + \max \{db-a \mid (a,b)
\in {\Omega}_{L} \} \leq d (e - {\rm indeg}_2 L) + a^1_*(L)+ u.$
In particular,  for any homogeneous ideal $I$ of $A$ we have
$$a_* (I^e) \leq de + a^1_*(R) + u.$$
 }
\end{rem}

\bigskip
\noindent
\subsection{Explicit bounds for some families of ideals}
\label{Bb}

\bigskip
The next purpose is to get explicit bounds for the $a_*$-invariant
of the powers of an ideal, and we will focus our attention to the
case of ideals in the polynomial ring. Throughout the rest of this
section, $A=k[{\x}]$ will denote the usual polynomial ring in $n$
variables over a field $k$ and $I$ will be a homogeneous ideal in
$A$. First of all, for equigenerated ideals whose Rees algebra is
Cohen-Macaulay we have

\begin{prop}
\label{Bb111}
 Let $I$ be a homogeneous ideal generated by forms in degree $d$ whose
Rees algebra is Cohen-Macaulay. Set $l= l(I)$. Then
$$ -n + d(-a^2(G)-1) \leq
 \max_{e \geq 0} \{ a_* (I^e) -de \}
 \leq -n + d (l-1).$$
\end{prop}

{\pf} As in the proof of Theorem \ref{YY2}, we may assume that $k$
is infinite. By considering then the polynomial ring
$S=k[{\x},Y_1, \dots,Y_l]$, we have that $R_A(I)$ is a finitely
generated bigraded $S$-module in a natural way. Then let
$$ 0 \to D_t \to \dots \to D_1 \to D_0   \to R \to 0$$
be the minimal bigraded free resolution of the Rees algebra $R$
over $S$, with $D_p =\bigoplus_{(a,b) \in {\Omega}_p} S(a,b)$. The
shifts $(a,b) \in {\Omega}_R$, $(a,b) \not = (0,0)$, satisfy $b
\leq -1$ and $-a \leq dl+n+a^1(R) \leq dl$ by Lemma \ref{S1}.
Therefore, we have $a_* (I^e) \leq de + d(l-1) -n$ by Remark
\ref{YY3}.

Let $R_{++} = \bigoplus_{(i,j), j>0} R_{(i,j)}$. From the bigraded
exact sequences
$$0 \to R_{++} \to R \to A \to 0\,,$$
$$0 \to R_{++}(0,1) \to R \to G \to 0 \,,$$
we get the following exact sequences of local cohomology
$$ 0 \to H^{n}_{\fm}(A)_{(i,j)} \to
 H^{n+1}_{\cal M}(R_{++})_{(i,j)} \to
 H^{n+1}_{\cal M}(R)_{(i,j)} \to  0 \;\;\; (\star),$$
$$ 0 \to H^{n}_{\cal M}(G)_{(i,j)} \to
 H^{n+1}_{\cal M}(R_{++})_{(i,j+1)} \to
 H^{n+1}_{\cal M}(R)_{(i,j)} \to  0 \;\;\; (\star \star).$$
Since $a^2(R) = -1$, from the above exact sequences we have
$a^2(G) \leq -1$. If $a^2(G) =-1$, the lower bound is obvious by
considering $e=0$. So we may assume $a^2(G) < -1$, and by Theorem
\ref{YY2} we must prove $a^1(R^\varphi) \geq -n-d(a^2(G)+1)$. The
local cohomology modules behave well under a change of grading by
Lemma \ref{B1}, hence we have
$$H^{n+1}_{\cal M}(R^\varphi)_{(p,q)} =
H^{n+1}_{\cal M}(R)^\varphi_{(p,q)} = H^{n+1}_{\cal
M}(R)_{(p+dq,q)},$$ so $a^1(R^\varphi) = \max \, \{\,p \mid
\exists \; q \; { s.t.} \, H^{n+1}_{\cal M}(R)_{(p+dq,q)} \not = 0
\}$. Since $H^{n}_{\fm}(A)_{(-n,0)} \not = 0$ we have
$H^{n+1}_{\cal M}(R_{++})_{(-n,0)} \not =0$ from the exact
sequence $(\star)$. As $a^2(G) < -1$, from the second exact
sequence $(\star \star)$ we get $H^{n+1}_{\cal M}(R)_{(-n,-1)}
\not =0$, and by using once more $(\star)$ we have $H^{n+1}_{\cal
M}(R_{++})_{(-n,-1)} \not =0$. Note that we can repeat this
procedure while the second component of the degree be greater than
$a^2(G)$, and finally we get $H^{n+1}_{\cal M}(R)_{(-n,a^2(G)+1)}
\not =0$. In particular, $a^1(R^\varphi) \geq -n-d(a^2(G)+1)$.
$\B$

\medskip
\begin{rem}
\label{Bb1111} {\rm Let $I$ be an ideal generated by forms of
degree $d$ in a general standard graded noetherian $k$-algebra
$A$. By setting $l=l(I)$, one can similarly prove that if the Rees
algebra is Cohen-Macaulay then
$$ a(A) + d(-a^2(G)-1) \leq  \max \{ a_* (I^e) -de \}  \leq a(A)+ dl.$$
}
\end{rem}

\medskip
For non-equigenerated ideals, we can also give an upper bound. A
similar result for the regularity was already proved in
\cite[Corollary 2.6]{CHT}.

\begin{rem}
\label{Bb1} {\rm
 Let $I$ be a homogeneous ideal generated by forms $f_1, \dots, f_r$ in
degrees $d_1 \leq \dots \leq d_r= d$ whose Rees algebra is
Cohen-Macaulay. Set $u = \sum_{i=1}^r d_i$. Then $a_* (I^e) \leq
d(e-1) + u -n$. }
\end{rem}

{\pf} By considering the bigraded minimal free resolution of $R$
over $S=k[{\x},Y_1, \dots,Y_r]$, we have that any shift $(a,b) \in
{\Omega}_p$ with $p \geq 1$ satisfies $b \leq -1$ and $-a \leq u$
by Lemma \ref{S13} and $\Omega_0$ only contains the shift $(0,0)$.
Therefore, $a_* (I^e) \leq d(e-1) + u -n$ by Remark \ref{YY3}.
$\B$

\medskip
     The $a_*$-invariant of the powers of an ideal can be computed for
complete intersection ideals (see the proof of Proposition
\ref{S4}), and then we have that the inequalities in Proposition
\ref{Bb111} and Remark \ref{Bb1} are sharp. Next we are going to
compute explicitly $\max_{e \geq 0} \{a_*(I^e)-de \} =
a^1(R^\varphi)$ for several families of ideals. First we consider
the case of equimultiple ideals.

\begin{prop}
\label{Bb2} Let $I$ be an equimultiple ideal generated in degree
$d$. Set $h= \h (I) \geq 1$. If the Rees algebra is
Cohen-Macaulay,
\begin{enumerate}
\item  $a (I^e/I^{e+1}) = de + a(A/I)$. In particular,
$a^1(G^\varphi) = a(A/I)$.
\item   $a_{n-h+1} (I^e) = d (e-1) + a(A/I) $. In particular,
$a^1(R^\varphi) = a(A/I)-d$.
\end{enumerate}
\end{prop}

{\pf} We may assume that $k$ is infinite. Then there exist $f_1,
\dots, f_{n-h} \in A$ of degree 1 such that $\overline f_1, \dots,
\overline f_{n-h} \in A/I$ is a homogeneous system of parameters.
Denoting by $f^*$ the initial form of $f \in A$ in $G$, let us
consider $\overline G= G/(f_1^*,\dots, f_{n-h}^*)$. Since $\rad
((f_1^*,\dots, f_{n-h}^*) G) = \rad (\fm G)$, we have that a
system of parameters of $F_\fm(I)$ is also a system of parameters
of $\overline G$. As $F_\fm(I)$ is a bigraded $k$-algebra
generated by forms in degree $(d,1)$, there exist $F_1,\dots, F_h
\in I$ of degree $d$ such that $\overline F_1, \dots, \overline
F_h $ is a system of parameters of $F_\fm (I)$. Then $f_1^*,
\dots, f_{n-h}^*,F_1^*,\dots, F_h^*$ is a homogeneous system of
parameters of $G$, and so algebraically independent over $k$.
Therefore, there is a finite extension
$$ T = k[U_1, \dots, U_h, V_1, \dots, V_{n-h}]  \rightarrow G,$$
where $T$ is a polynomial ring with $\dg(U_i) = (d,1)$ , $\dg(V_j)
= (1,0)$ for $i= 1, \dots, h$, $ j=1, \dots, n-h$. Since $G$ is
Cohen-Macaulay and $T$ is regular, we have that $G$ is a free
$T$-module, that is, $ G = \bigoplus_{(a,b) \in \Lambda} T (a,
b)$, where $\Lambda \subset \Bbb Z^2$ is a finite set. Let us
denote by $T_1= k[V_1, \dots, V_{n-h}]$, $m= (V_1, \dots,
V_{n-h})$, and for a given $\underline \beta = (\beta_1, \ldots,
\beta_h) \in \Bbb N^h$, let $\mid \underline \beta \mid = \beta_1
+ \cdots + \beta_h$. Note that

\vspace{2mm}

$\hspace{10mm} T(a,b)^e =\bigoplus_i T(a,b)_{(i,e)} = \bigoplus_i
T_{(a+i,b+e)} $

\vspace{2mm}

$\hspace{25mm} = \bigoplus_i \bigoplus_ {\mid \underline \beta
\mid = b+e}
 [T_1]_{a+i-d(b+e)}  U_1^{\beta_1} \dots U_h^{\beta_h} $

\vspace{2mm}

$\hspace{25mm} =T_1(a-db-de)^{\rho_b^e},$

\vspace{2mm}

\noindent with $\rho_b^e \in \Bbb N$ and $\rho_b^e = 0$ if $b+e
<0$. Therefore,
$$ I^e/ I^{e+1} =  G^e  =
\bigoplus_{(a,b) \in \Lambda} T_1(a-db-de)^{\rho_b^e},$$ \noindent
and by taking local cohomology
$$H^{n-h}_{\frak m} (I^e/ I^{e+1}) =
 \bigoplus_{(a,b) \in \Lambda}
 H^{n-h}_{m}(T_1(a-db-de))^{\rho_b^e}.$$
Hence $a_*(I^e/I^{e+1})= a (I^e/I^{e+1}) = \max \, \{ -(n-h) -a +
db +de : -b \leq e \} $. In particular,  $a (A/I) = \max \{ -(n-h)
-a + db : b= 0\}$, and so we get $a (I^e/I^{e+1}) \geq de + a(A/I)
$ for all $e$. On the other hand, since the modules $I^e/I^{e+1}$
are $A/I$-modules of maximal dimension, we have an epimorhism
$\bigoplus A/I (-de) \to I^e /I^{e+1}$ and we may deduce that $a
(I^e/I^{e+1}) \leq de + a(A/I)$ for all $e$. To get $(ii)$, it is
just enough to consider the short exact sequences
$$ 0 \to I^{e+1} \to I^e \to I^e/I^{e+1} \to 0 \,\,,$$
and then the result follows from $(i)$ by induction on $e$. $\B$

\medskip
   Next we study equigenerated ideals whose form ring is Gorenstein.
In this case, we prove that the lower bound given in Proposition
\ref{Bb111} is sharp.

\begin{prop}
\label{Bb5} Let $I$ be a homogeneous ideal equigenerated in degree
$d$ whose form ring is Gorenstein. Set $l = l(I)$. Then
\begin{enumerate}
\item $\max_{e \geq 0} \{ a_*(I^e) -de \}=  d(- a^2(G) -1) -n$.
\item For $e> a^2(G) -a(F)$,
$\depth (A/I^e) = n-l$ and
 $a_*(I^e)= a_{n-l} (A/I^e) =d(e-a^2(G)-1)-n$.
\end{enumerate}
\end{prop}

{\pf}  We may assume that the field $k$ is infinite. Since $I$ is
generated by forms in degree $d$, there exists a minimal reduction
$J$ of $I$ generated by forms $g_1, \dots, g_l$ of degree $d$. By
considering $S=k[{\x},Y_1, \dots,Y_l]$ bigraded by setting $\dg
(X_i)= (1, 0)$, $\dg(Y_j)= (d, 1)$, we have a bigraded epimorphism
$S \rightarrow R_A(J)$. Suppose that $I^{m+1}= J I^m$. Then
$R_A(I)$ is finitely generated over $R_A(J)$ by the generators of
$A$, $I$, \dots, $I^m$; so in particular by homogeneous elements
in degree $(di,i)$ for $i = 0, \dots, m$. Then we have an
epimorphism $ F  \rightarrow R_A(I) $, where $F$ is a finite free
$S$-module with a basis of elements in degrees $(di,i)$ for $i =
0, \dots, m$.

 Let us consider the minimal bigraded free resolution of $G$ over $S$
$$ 0 \to D_l \to \dots \to D_1 \to D_0  \to G \to 0,$$
where $D_p = \bigoplus _{(a,b) \in {\Omega}_p} S(a,b)$. From
Remark \ref{YY3}, for all $e \geq 0$ we have
$$a_* (I^e/I^{e+1}) \leq de -n + \max \{db-a \mid (a,b) \in {\Omega}_G \}.$$
Assume we prove that the maximum is accomplished for a shift
$(a,b) \in {\Omega}_l$. Denoting by $(\;\;)^* = \uHom_S(\;\; ,
K_S)$, then
$$ 0 \to D_0^* \to \dots  \to D_{l}^*  \to \uExt^l_S(G, K_S) = K_G \to 0$$
is the minimal bigraded free resolution of the canonical module
$K_G$ of $G$ over $S$.
Since $G$ is Gorenstein, according to Corollary \ref{K4} there is
a bigraded isomorphism
$$K_G \cong G(-n, a^2(G)). $$
Now the shifts $(a,b) \in {\Omega}_{l}$ are of the type $(di,
i-a^2(G))$ for certain integers $i$, so we get $a_* (I^e/I^{e+1})
\leq de -n + \max_i \{ d(i-a^2(G)) -di \} = d(e-a^2(G))-n$. From
Remark \ref{Bb111} we have $a^1(G^\varphi) = \max
\,\{a_*(I^e/I^{e+1})-de \} \leq -n-d a^2(G)$, and $a^1(R^\varphi)
\leq -n-d a^2(G)-d $. Observe that Proposition \ref{Bb111} gives
the other inequality, so we obtain $ a^1(R^\varphi) =
d(-a^2(G)-1)-n$.

    Now let us prove that the maximum for the differences $db-a$ is
taken for $(a,b) \in {\Omega}_l$. Let us consider $\phi : \Bbb Z^2
\to \Bbb Z$ defined by  $\phi (i,j) = i-dj$. Note that
 $X_i \in S^{\phi}$ has degree 1 and $Y_j \in S^{\phi}$ has degree 0, and
$$ 0 \to D_{l}^{\phi} \to \dots \to D_0^{\phi} \to G^{\phi} \to 0$$
is a graded free resolution of $G^{\phi}$ over $S^{\phi}$. Since
$S(a,b)^{\phi} = S^{\phi}(a-db)$, it is clear that
$$\min \{db-a \mid (a,b) \in {\Omega}_{p+1} \} \geq
 \min \{db-a \mid (a,b) \in {\Omega}_{p} \}.$$
Applying the same argument for the resolution of $K_G$, one gets
that
 $$\max \{db-a \mid (a,b) \in \Omega_{p+1} \} \geq
 \max \{db-a \mid (a,b) \in {\Omega}_{p} \},$$
so $(i)$ is already proved.

 To prove the rest of the statement, we will use
that $\dpp_A (I^e/I^{e+1}) = l$ if and only if $e \geq a^2(G)-
a(F)$ (see Proposition \ref{GG1}). By applying the functor
$(\;)^e$ to the resolution of $G$ over $S$ we obtain a free
resolution of $I^e/I^{e+1}$ as $A$-module, and the shifts
appearing in the place $p$ of these resolutions of the type
$a-db-de$, with $(a,b) \in \Omega_p$. Therefore, $t_p
(I^e/I^{e+1}) \leq de-da^2(G)$ for any $p$, and for $e \geq
a^2(G)- a(F)$ we have $t_l (I^e/I^{e+1}) = de-da^2(G)$. Now, since
$\dpp_A I^e \leq l-1$ by Proposition \ref{GG1}, from the short
exact sequences
$$0 \to I^{e+1} \to I^e \to I^e/I^{e+1} \to 0,$$
we get $t_{l-1}(I^e) \geq d(e-a^2(G)-1)$ for $e>a^2(G)-a(F)$. On
the other hand, we have $t_{l-1}(I^e) \leq t_*(I^e) = a_*(I^e)+n
\leq d(e-a^2(G)-1)$, so we get the equality. We finally obtain
$a_{n-l}(A/I^e) =a_{n-l+1}(I^e) = d(e-a^2(G)-1) -n$ for
$e>a^2(G)-a(F)$ by Theorem \ref{B77}. $\B$

\begin{ex}
\label{Bb81} {\rm Let ${\bf X}=(X_{ij})$ be a $d \times n$ generic
matrix, with $1 \leq i \leq d$, $1 \leq j \leq n$ and $d \leq n$,
and let $A= k[{\bf X}]$ be the polynomial ring in the entries of
${\bf X}$. Let $I= I_d ({\bf X})$ be the ideal generated by the
maximal minors of ${\bf X}$. We are going to apply to this example
the different bounds we have found.

\begin{itemize}
\item The Rees algebra of $I$ is Cohen-Macaulay, so by applying
Remark \ref{Bb1} we get
$$a_*(I^e) \leq d(e-1)+ {n \choose d}d-nd.$$
(a similar bound for the regularity of the powers has been given
in \cite[Example 2.7]{CHT}).
\item Note that $\F= k[I_d]$ is the coordinate ring of the
Grassmannian $G(d,n)$, so we have that the analytic spread of $I$
is $l(I)= d(n-d) +1$. Therefore by Proposition \ref{Bb111} we get
the bound
$$a_*(I^e) \leq de + d^2(n-d)-nd.$$
\item Since $G_A(I)$ is Gorenstein with $a^2(G_A(I)) = - \h (I)
=-(n-d+1)$, by using Proposition \ref{Bb5} we get the better bound
$$a_*(I^e) \leq d(e+n-d) -nd = de-d^2.$$ Furthermore, the a-invariant of $F$
is $-n$ by \cite[Corollary 1.4]{BH2}. So we also obtain
$a_{d^2}(I^e) = de - d^2$ for any $e >  d-1$.
\end{itemize}

   K. Akin et al. \cite{ABW} have constructed resolutions for the powers of
$I$, in particular showing that all the powers of $I$ have linear
resolutions. Note that this fact also allows to prove the last
bound: According to Proposition \ref{GG1} and Theorem \ref{B77},
for any $e > a^2(G)-a(F)$ we have  $a_*(I^e) = a_{nd-l+1}(I^e)=
t_{l-1}(I^e)-nd = (de+l-1)-nd = de-d^2$. }
\end{ex}

\medskip
We may also use Proposition \ref{Bb5} to study the $a_*$-invariant
of the powers of a strongly Cohen-Macaulay ideal. Let $I$ be an
ideal of $A$, and let $\underline{f}=f_1, \dots, f_r$ be a system
of generators of $I$. Recall that $I$ is a strongly Cohen-Macaulay
ideal if for any $p \geq 0$ the Koszul homology $H_p (K
(\underline{f}))$ is a Cohen-Macaulay $A/I$-module.

\begin{cor}
\label{Bb8} Let $I$ be a strongly Cohen-Macaulay ideal generated
in degree $d$ such that $\mu(I_{\frak p}) \leq \h (\frak p)$ for
any prime ideal $\frak p \supseteq I$. Let $h=\h (I)$, $l= l(I) $.
Then
\begin{enumerate}
\item $a_*(I^e) \leq d(e +h -1) -n$, $\forall e$.
\item For $ e > l-h$, $\depth (A/I^e) = n-l$ and
$a_*(I^e)=a_{n-l}(A/I^e) = d(e +h -1) -n.$
 \end{enumerate}
\end{cor}

{\pf} In this situation, $G_A(I)$ is Gorenstein and $I$ is an
ideal of linear type by \cite[Theorem 2.6]{HSV1}, so $a(F_{\frak
m}(I)) = -l$. Furthermore, according to \cite[Proposition
2.5]{HRZ} we have $a^2(G_A(I)) = -h$. Then the result follows from
Proposition \ref{Bb5}.$\B$

\begin{ex}
\label{Bb9} {\rm  Let
 $I \subset  A= k[X_1, X_2, X_3, X_4]$ be the defining ideal of the
twisted cubic in $\Bbb P^3_k$. From Example \ref{YY4}, recall that
$I$ is a strongly Cohen-Macaulay ideal generated in degree $2$
with ht$(I)= 2$, $l(I)=\mu(I)= 3$. Now, by Corollary \ref{Bb8},
for any $e>1$ we have that depth$(A/I^e)=1$, $a_1(A/I^e) = 2e-2$
and $a_2(A/I^e) \leq 2e-2$. In the case $e=1$, since $I$ is linked
to $J=(X_1, X_2)$ by the regular sequence $\underline \alpha =
X_2^2-X_1X_3, X_3^2-X_2X_4$, we have that $A/I$ is Cohen-Macaulay
and there is a graded isomorphism
 $K_{A/I} \cong J/(\underline \alpha)$. In particular,  $a(A/I) = -1$.
}
\end{ex}

\medskip
In trying to extend the bounds in Proposition \ref{Bb5} to the
non-equigenerated case many difficulties appear. Next we will use
approximation complexes to do this for strongly Cohen-Macaulay
ideals.

\medskip
 Let $I$ be a homogeneous ideal in the polynomial ring $A= k[X_1, \dots,
X_n]$ and $\underline f= f_1, \dots, f_r$ a homogeneous system of
generators of $I$, with $d_i =$deg$(f_i)$, and let us consider the
graded Koszul complex $K(\underline{f})$ of $A$ with respect to
$\underline f$. Denote by $S = A[Y_1, \dots, Y_r]$ with the
bigrading deg$(X_i) = (1,0)$, deg$(Y_j) = (d_j, 1)$. Then the
approximation
 complex of $I$ is
 $${\cal M}(\underline{f}) :
  0 \to {\cal M}_r \to \dots \to {\cal M}_1 \to {\cal M}_0 \to 0,$$
\noindent with  ${\cal M}_p = H_p (K (\underline{f})) \otimes_A
S(0, -p)$, and the differential maps are homogeneous. Assume that
$I$ is a strongly Cohen-Macaulay ideal with $\h (I) \geq 1$ such
that for any prime ideal $\frak p \supseteq I$, $\mu (I_{\frak p})
\leq \h (\frak p)$. Then ${\cal M}(\underline{f})$ is exact and
provides a resolution of ${\rm Sym}_A(I/I^2) \cong  G_A(I)$
\cite[Theorem 2.6]{HSV1}. We will use it to get a bound for the
a-invariants of the powers of these ideals.

\begin{prop}
\label{Bb10}
 Let $I$ be a strongly Cohen-Macaulay ideal such that for any prime ideal
  $\frak p \supseteq I$, $\mu (I_{\frak p}) \leq \h (\frak p)$. Assume
that $I$ is minimally generated by forms $f_1, \dots, f_r$ of
degree $d=d_1 \geq \dots \geq d_r$, and set $h = \h (I)$, $t=
r-h$. Then:
\begin{itemize}
\item[(i)]  $a(H_m(\underline{f})) \leq -n + d_1 + \dots + d_{h+m}$,
for all $m \leq t$.
\item[(ii)]  If $1 \leq e \leq t$, $\depth(A/I^e) \geq n-h-e+1 $ and
 for any $0 \leq m \leq e-1$,

\hspace{15 mm}$a_{n-h+1-m}(I^e) \leq -n + d_1 + \dots + d_{h+m} +
d(e-m-1).$

\item[(iii)]  If $e > t$, $\depth (A/I^e) = n-r $ and
 for any $0 \leq m \leq t$,

\hspace{15 mm}$a_{n-h+1-m}(I^e) \leq -n + d_1 + \dots + d_{h+m} +
d(e-m-1).$
\end{itemize}
\end{prop}

\medskip
{\pf} Recall that $H_p = H_p (K(\underline f)) = 0$ for all $p>t$.
So the resolution of $G$ given by the approximation complex is
 $$0 \to {\cal M}_t \to \dots \to {\cal M}_0 \to G \to 0,$$
\noindent with  ${\cal M}_p = H_p \otimes_A S(0, -p)$. Let us
denote by $\underline \beta = (\beta_1, \ldots, \beta_r) \in \Bbb
N^r$, and $\mid \underline \beta \mid = \beta_1 + \cdots +
\beta_r$. Applying the functor $(\,)^e$ to the modules of this
resolution:

\vspace{3mm}

\hspace{5mm} $G^e = \bigoplus_i G_{(i,e)} = I^e/I^{e+1},$

\vspace{2mm}

\hspace{5mm} ${\cal M}_p^e = \bigoplus_i H_p[Y_1, \dots, Y_r]_{(i,
e-p)}$

\hspace{14mm}$= \cases{ 0 & if $e<p$ \cr
  \bigoplus_{\mid \underline \beta \mid = e-p}
   H_p (-d_1\beta_1 - \dots - d_r \beta_r) & if $e \geq p$. \cr}$

\vspace{3mm} So we get graded exact sequences
$$0 \to {\cal M}_q^e \to \dots \to {\cal M}_0^e \to I^e/I^{e+1} \to 0,$$
\noindent with $q = \min \{ e, t \}$. Since $I$ is a strongly
Cohen-Macaulay ideal, we have that ${\cal M}_p^e$ is a  maximal CM
$A/I$-module for any $p \leq q$. By taking short exact sequences,
we obtain that if $e<t$, depth$(I^e/I^{e+1}) \geq n-h-e$ and if $e
\geq t$, depth$(I^e/I^{e+1}) \geq n-h-t = n-r$. On the other hand,
by Proposition \ref{GG1} we also have that depth$(I^e/I^{e+1})=
n-r$  if and only if $e \geq t$. Furthermore, we get
$a_{n-h-m}(I^e/I^{e+1}) \leq a({\cal M}_m^e) = a(H_m) + d(e-m)$
for all $0 \leq m \leq \min \{e, t \}$. From the exact sequences
$$0 \to I^e /I^{e+1} \to A/ I^{e+1} \to A/I^{e} \to 0, $$

\noindent we have now

\noindent \vspace{3mm} (a) depth$(A/I^e) \geq n-h-e+1$ if $1 \leq
e \leq t$

\noindent \vspace{2mm} (b) depth$(A/I^e) = n-r$ if $e >t$

\noindent \vspace{2mm} (c) $a_{n-h-m}(A/I^e) \leq \max_{0 \leq j
\leq e-1} \{ a_{n-h-m}(I^j/I^{j+1}) \} \leq a (H_m) + d(e-m-1)$,

\noindent \vspace{2mm} for $0 \leq m \leq \min \{e-1, t \}$.

\noindent So, if we prove the bound for the a-invariant of the
Koszul homology we have finished. Let us assume that among the
forms $f_1, \dots, f_r$ we can choose a regular sequence of length
$h$. Let $g_1=f_{j_1},\dots, g_h=f_{j_h}$ be this sequence, and
$g_1, \dots g_r$ the minimal system of generators of $I$.

   Let us consider the morphism from $A$ to $A/(g_1)$. By
\cite[Lemma 1.1]{Hu1}, there is a graded exact sequence
 $$0 \to H_m(I;A) \to H_m(I/(g_1);A/(g_1)) \to H_{m-1}(I;A)(-\deg \,g_1)
 \to 0$$
\noindent for all $m \geq 1$; where $H_m(I/(g_1);A/(g_1))$ denotes
the Koszul homology of the elements $0,\overline g_2, \dots,
\overline g_r \in A/(g_1)$. From this exact sequence, we have in
particular $a(H_m (I; A)) \leq a(H_m(I/(g_1);A/(g_1))$.
   Denote by "--" the morphism from $A$ to $\overline A = A/(g_1, \dots,
g_{h-1})$. Repeating $h-1$ times the previous procedure, we get
$a(H_m (I; A)) \leq a(H_m (\overline I; \overline A))$ for all $m
\geq 1$. But now $\overline I$ is a height one ideal in the CM
ring $\overline A$. Let us denote the Koszul complex of $\overline
I$ by $\overline K = K(\overline I ; \overline A)$, and the
differential from $\overline K_{m+1}$ to $\overline K_m$ by
$d_{m+1}$. Set $\overline Z_m = {\rm Ker} (d_m)$, $\overline B_m =
{\rm Im} (d_{m+1})$, $\overline H_m = H_m (\overline I; \overline
A)$.
 Then there are exact sequences
$$0 \to
\overline B_m \to \overline Z_m  \to \overline H_m \to 0,$$
$$0 \to
\overline Z_{m+1} \to \overline K_{m+1} \to \overline B_m \to 0.$$
By \cite[Lemma 1.6]{Hu2} $\overline H_m $ are CM modules for all
$m$, and then by \cite[Lemma 1.8]{Hu2}, $\overline Z_m$ and
$\overline B_m$ are maximal CM modules for $\overline A$. The
exact sequences imply now $a(\overline H_m) \leq a(\overline B_m)
\leq a(\overline K_{m+1})$.
 Denoting by $\delta_i =$deg$(g_i)$,

\vspace{3mm}

$a(\overline K_{m+1}) =
 a(\overline A)+
\max \{ \delta_{i_1} + \dots + \delta_{i_{m+1}} \mid h \leq i_1 <
\dots < i_{m+1} \leq r \} $

\vspace{2mm}

\hspace{2mm}$= a(A)+ \delta_1 + \dots +\delta_{h-1} + \max \{
\delta_{i_1} + \dots + \delta_{i_{m+1}} \mid h \leq i_1 < \dots <
i_{m+1} \leq r \} $

\vspace{2mm}

\hspace{2mm}$\leq -n + d_1 + \dots + d_{h+m}$.

\vspace{3mm} \noindent So we are done if we prove the following
lemma. $\B$

\begin{lem}
\label{Bb12} Let $A= k[t_1, \dots, t_s]$ be a CM graded algebra
over an infinite field $k$, with $\deg(t_i) =1$. For any
homogeneous ideal $I$, there exists a minimal homogeneous system
of generators $g_1, \dots g_r$ of $I$ such that $g_1, \dots g_h$
is a maximal regular sequence in $I$.
\end{lem}

{\pf}  Set $r =\mu(I)$, $h = \h (I)$. Let $f_1, \dots, f_r$ be a
minimal homogeneous system of generators of $I$, with
 $d_i =\deg(f_i)$, $d_1 \leq \dots \leq d_r=d$. We are going to prove
the statement by induction on $h$. If $h= 0$ there is nothing to
prove. Assume $h \geq 1$. Then $I \not \subset z(A) = \bigcup
_{\frak p \in {\rm Ass}(A)} \frak p$, and so $I_d \not \subset
\frak p$, $\forall \frak p \in \Ass (A)$ (otherwise, we would have
$f_i^d \in \frak p$ for all $i$, and then $I \subset \frak p$).
Since $k$ is infinite, we get $I_d \not \subset \bigcup _{\frak p
\in {\rm Ass} (A)} \frak p \cap A_d$, and so there exists $g \in
I_d$ such that $g \not \in \frak p$ for all $\frak p \in \Ass(A)$.
Note that $I_d$ is a $k$-vector space generated by the forms
 $f_{j_1}, \dots, f_{j_i}$ in degree $d$ and the forms
$M f_j$, with $d_j < d$ and $M$ a monomial in $t_1, \dots, t_s$ of
degree $d-d_j$. Now we can write
 $$g = \lambda_1 f_{j_1}+ \dots+ \lambda_i f_{j_i}
+ \sum \mu_{jM} M f_j,$$ with $\lambda_1, \dots, \lambda_i,
\mu_{jM} \in k$. If there exists $p$ such that $\lambda_p \not =
0$, we can replace $f_{j_p}$ by $g$ in the minimal system of
generators of $I$. Otherwise, we have an element $g$ in the ideal
generated by the forms in $I$ of degree $d'=d_{r-1}$ with the
property that $g \not \in \frak p$, $\forall \frak p \in \Ass(A)$.
We can repeat the arguments for $I_{d'}$, and finally we will
replace one of the forms $f_j$ by $g$. By considering $\overline
A= A/(g)$, the ideal $\overline I = I/(g)$ has $\mu (\overline I)
= r-1$, $\h(\overline I) = h-1$. Then we get the result by
induction. $\B$

\bigskip
\subsection{Applications to the study of the diagonals of the Rees
algebra} \label{T}

\bigskip
Next we are going to apply the results about the $a$-invariant of
the powers of an ideal to study the Cohen-Macaulay property of the
diagonals of the Rees algebra. If the Rees algebra is
Cohen-Macaulay, according to Theorem \ref{S3} there exists $\alpha
\in \Bbb Z$ such that $\Kk$ is Cohen-Macaulay for any $c
>de+\alpha$ and $e>0$. For equigenerated ideals, we obtained
$\alpha=d(l-1)$ as upper bound. The following result precises the
best $\alpha$.

\begin{prop}
\label{T1}
  Let $I$ be an ideal in $A= k[{\x}]$ generated by
forms in degree $d$ whose Rees algebra is Cohen-Macaulay. Set $l =
l(I)$. For $\alpha \geq 0$, the following are equivalent
\begin{itemize}
 \item[(i)]   For all $c > de + \alpha$, $\Kk $ is CM.
 \item[(ii)]   $a_i (I^e) \leq d e +\alpha$, $\forall i$, $\forall e$.
 \item[(iii)]   $a_i (I^e) \leq d e +\alpha$, $\forall i$, $\forall e
  \leq  l-1$ .
 \item[(iv)]   $H^{n +1}_{\cal M} (R_A(I))_{(p,q)} = 0$, \,
  $\forall p >dq + \alpha $,
   that is, $ \alpha \geq a^1(R^\varphi)$.
 \item[(v)]   The minimal bigraded free resolution of $R_A(I)$ is
  good for any diagonal $\Delta=(c,e)$ such that $c >de + \alpha .$
\end{itemize}
\end{prop}

{\pf} If $\Kk$ is CM for $c > de + \alpha$ then we have
$H^i_{\fm}(I^e)_c =0$ for any $i<n$ and $c > de + \alpha$ by
Proposition \ref{S611}, so $a_* (I^e) \leq d e +\alpha$,  $\forall
e$. The converse follows similarly, and we get the equivalence
between $(i)$ and $(ii)$.

Since the Rees algebra $R$ of $I$ is Cohen-Macaulay, we have
$a_*^2(R)=a^2(R)=-1$. Then, conditions $(ii)$ and $(iii)$ are
equivalent to $a^1(R^\varphi)= a_*^1(R^\varphi) \leq \alpha$ by
Theorem \ref{YY2}.

Finally, we want to prove the equivalence to $(v)$. If the
resolution of $R$ is good for diagonals $\Delta=(c,e)$ such that
$c > de + \alpha$, then we have $H^i_{m}(\Kk) = H^{i+1}_{\cal
M}(R)_\Delta=0$ for $i<n$ by Corollary \ref{W64}, so $\Kk$ is
Cohen-Macaulay for any $c > de + \alpha$ and we obtain $(i)$. Now
assume that $a^1(R^\varphi) \leq \alpha$, and let us consider the
minimal bigraded free resolution of $R$ over $S$
$$ 0 \to D_t \to \dots \to D_1 \to D_0 \to R \to 0 ,$$
with $D_p =\bigoplus_{(a,b) \in {\Omega}_p} S(a,b)$. By applying
the functor $(\,\,)^\varphi$, we have that
$$ 0 \to D_{t}^{\varphi} \to \dots \to D_0^{\varphi} \to R^{\varphi}
\to 0$$ is the bigraded minimal free resolution of $R^{\varphi}$
over $S^{\varphi}$, with $D_p^{\varphi} = \bigoplus_{(a,b) \in
\Omega_p}  S^{\varphi}(a-db,b)$. Therefore, according to Theorem
\ref{B77}, for any $(a,b) \in \Omega_R$ we have $db-a-n \leq
\alpha$. Then the sets $X^\Delta$, $Y^\Delta$ introduced in Remark
\ref{W10} are empty for diagonals $\Delta=(c,e)$ with $c > de +
\alpha$, so the resolution is good for these $\Delta$. $\B$

\medskip
    If the form ring is Gorenstein, we can express this
criterion by means of the second a-invariant of the form ring.
Namely,

\begin{cor}
\label{T2}
 Let $I$ be an ideal in $A= k[{\x}]$ generated by
forms of degree $d$ whose form ring is Gorenstein. For $\alpha
\geq 0$, the following are equivalent
\begin{itemize}
   \item[(i)]   For all $c > de+\alpha$, $\Kk $ is CM.
   \item[(ii)]   $\alpha \geq d(-a^2(G)-1)-n$.
\end{itemize}
\end{cor}

{\pf} By Proposition \ref{Bb5}, $a^1(R^\varphi) = d(-a^2(G)-1)-n$.
Then the result follows from Proposition \ref{T1}. $\B$

\medskip
Let $I$ be an equigenerated ideal in $A$. If the Rees algebra is
Cohen-Macaulay, it can happen that some of its diagonals are not
Cohen-Macaulay. Now, by taking $\alpha=0$ in Proposition \ref{T1}
we have a criterion to decide when all the diagonals of a
Cohen-Macaulay Rees algebra are Cohen-Macaulay.

\begin{cor}
\label{T3}
        Let $I$ be an ideal in $A= k[{\x}]$
generated by forms in degree $d$ whose Rees algebra is
Cohen-Macaulay. Set $l = l(I)$. Then the following are equivalent
\begin{itemize}
   \item[(i)]   For all $c \geq de+1$, $\Kk $ is CM.
   \item[(ii)]   $a_i (I^e) \leq d e $, $\forall i$, $\forall e \leq
 l-1$ .
   \item[(iii)]   $H^{n +1}_{\cal M} (R_A(I))_{(p,q)} = 0$, \,
$\forall p >dq .$
 \item[(iv)]   The minimal bigraded free resolution of $R_A(I)$ is
  good for any $\Delta$.
\end{itemize}
Assuming that $G_A(I)$ is Gorenstein, these conditions are also
equivalent to
\begin{itemize}
  \item[(v)] $-a^2(G) \leq \frac{n}{d}+1.$
\end{itemize}
\end{cor}

\begin{ex}
\label{T4} {\rm  We may recover Corollary \ref{S14} as an easy
application of Corollary \ref{T3}. Let $\{L_{ij} \}$  be a set of
$d \times (d+1)$ homogeneous linear forms in a polynomial ring
$A=k[{\x}]$, and let $M$ be the matrix $(L_{ij})$. Let $I_t(M)$ be
the ideal generated by the $t \times t$ minors of $M$ and assume
that ht$(I_t(M)) \geq d-t+2$ for $1 \leq t \leq d$. Set $I =
I_d(M)$. The ideal $I$ is generated by $d+1$ forms of degree $d$,
and we have a presentation of the Rees algebra of the form
        $$R_A(I)= k[{\x},Y_1, \dots, Y_{d+1}]/(\phi_1, \dots, \phi_d),$$
with deg$(Y_j) =(d,1)$, deg$(\phi_i)=(d+1,1)$, such that $\phi_1,
\dots, \phi_d$ is a regular sequence. Therefore $R_A(I)$ is
Gorenstein, and so $a^2(G_A(I))= -2$. Since $d \leq n-1$, we have
that $-a^2(G) \leq \frac{n}{d}+1$. Therefore $\Kk$ is
Cohen-Macaulay for any $c \geq de+1$. }
\end{ex}

\medskip
Next, we are going to use the bounds of the $a$-invariants of the
families of ideals considered in Subsection 5.3.1 to study the
Cohen-Macaulayness of the diagonals of their Rees algebras. First
we recall a well-known result about the vanishing of the graded
pieces of the local cohomology modules.

\begin{lem}
\label{T55}
  Let $A$ be a standard noetherian graded $k$-algebra with
graded maximal ideal $\fm$. Let $L$ be a finitely generated graded
$A$-module with $d = \dim L >0$. Then
$$H^d_{\fm}(L)_j \not =0, \; \forall j \leq a(L).$$
\end{lem}

{\pf} Since $d>0$, we can assume $H^0_{\fm} (L) =0$ because
otherwise by considering $\overline L = L / H^0_{\fm} (L)$ we have
$H^0_{\fm} (\overline L) =0$ and $H^d_{\fm} (\overline L) =
H^d_{\fm} (L)$. We may also assume that the field $k$ is infinite.
Then there exists $x \in A_1$ such that $x \not \in z_A(L)$, and
the exact sequence
$$0 \to L(-1) \stackrel{\cdot x} \longrightarrow L \to L/xL \to  0 \;$$
induces the graded exact sequence of local cohomology modules
$$H^{d-1}_{\fm}(L/xL) \to H^d_{\fm}(L)(-1) \to H^d_{\fm}(L) \to  0.$$
From this exact sequence, we have that $H^d_{\fm}(L)_s =0$ implies
$H^d_{\fm}(L)_j=0$ for $j \geq s$, so we are done. $\B$

\begin{prop}
\label{T5}
   Let $I$ be an equimultiple ideal in $A$ of height $h>1$ generated by
forms in degree $d$ whose Rees algebra is Cohen-Macaulay. For any
$c \geq de +1$, $\Kk$ is Cohen-Macaulay if and only if $c> d(e-1)
+ a(A/I)$.
\end{prop}

{\pf} We have proved in Proposition \ref{Bb2} that $a^1(R^\varphi)
=a(A/I)-d$. Therefore, $\Kk$ is Cohen-Macaulay for any $c> d (e-1)
+ a(A/I)$ by Proposition \ref{T1}.

 On the other hand, since $a_{n-h}(I^e/I^{e+1}) =de+a(A/I)$ by
Proposition \ref{Bb2}, we have $H^{n-h}_{\frak m} (I^e/I^{e+1})_ s
\not = 0$ for all $s \leq de +a(A/I)$ according to Lemma
\ref{T55}. By considering the short exact sequences
 $$0 \to I^{e+1} \to I^e \to I^e/I^{e+1} \to 0,$$
and by induction on $e$, we get  $H^{n-h+1}_{\frak m}(I^e)_ s \not
= 0$ for all $s \leq d (e-1)+ a(A/I)$. Now, if $\Kk$ is
Cohen-Macaulay then $H^i_{\fm}(I^{es})_{cs} =0$ for $i<n$ and
$s>0$ by Proposition \ref{S611}. In particular, it holds
$H^{n-h+1}_{\frak m}(I^e)_ c =0$, and so $c>d(e-1)+a(A/I)$. This
proves the converse. $\B$

\begin{prop}
\label{T6} Let $I$ be a strongly Cohen-Macaulay ideal such that
$\mu(I_{\frak p}) \leq \h(\frak p)$ for any prime ideal $\frak p
\supseteq I$. Assume that $I$ is minimally generated by $r$ forms
of degree $d=d_1 \geq  \dots \geq d_r$, and let $h = \h (I)$. For
$c >  d(e-1)+d_1 + \dots + d_h -n $, $\Kk$ is Cohen-Macaulay.
\end{prop}

{\pf} According to Corollary \ref{S61}, for a given $c \geq de+1$
we have that  $\Kk$ is Cohen-Macaulay if and only if $H^i_{\frak
m}(I^{es})_{cs} = 0$, for $i < n$, $s>0$, and $H^i_{\frak
m}(I^{es-h+1})_{cs-n} = 0$, for $ 1 < i \leq n$, $s>0$.

From Proposition \ref{Bb10}, note that $a_*(I^e) \leq (e-1) d+d_1+
\dots + d_h-n$. Therefore, to get the vanishing of the cohomology
modules it suffices to see that $cs > (es-1)d+d_1+ \dots+d_h-n$
and $cs-n > (es-h)d+d_1+ \dots+d_h-n$  for any $s \geq 1$. The
first condition is equivalent to $(c-de)s > d_1+ \dots +d_h -d-n$
for $s \geq 1$, that is, $c-de> d_1+ \dots+d_h-d-n$. The second
one is equivalent to $(c-de)s > d_1+ \dots +d_h -dh$ for $s \geq
1$, that is, $c-de> d_1+ \dots+d_h-dh$; and this always holds
because $d_1+\dots+d_h-dh \leq 0$.$\B$

\medskip
To finish, let us consider the case that the Rees algebra has
rational singularities. Then all the diagonals $\Kk$ have rational
singularities by \cite{Bo}, so in particular the Rees algebra and
its diagonals are Cohen-Macaulay. By Proposition  \ref{S611}, we
get immediately

\begin{prop}
\label{T7} Let $I$ be a homogeneous ideal in $A = k[{\x}]$
generated by forms of degree $\leq d $, where $k$ is a field with
$char k =0$. If $R_A(I)$  has rational singularities, then
$a_*(I^e) \leq de$ for all $e.$
\end{prop}

\begin{ex}
\label{T8} {\rm  Let ${\bf X}=(X_{ij})$ be an $m \times n$ matrix
of indeterminates, with $1 \leq i \leq m$, $1 \leq j \leq n$ and
$m \leq n$. Let $A = k[{\bf X}]$ be the polynomial ring with
variables the entries in ${\bf X}$, where $k$ is a field with char
$k = 0$ and $k = \overline k$. Let $I = I_d({\bf X})$ be the ideal
generated by the $d$-minors of ${\bf X}$, $1 < d < m$.

By \cite[Theorem 3.2]{Br}, $R_A(I)$ has rational singularities. So
we have $a_* (I^e) \leq de$, for all $e$. This also holds for
ideals generated by minors of symmetric generic matrices and
ideals generated by pfaffians of alternating generic matrices by
\cite[Remark 3.4]{Br}. The defining ideals of the varieties
considered by A. Bertram in \cite{B} are:

\begin{itemize}
\item[(a)]
 $I_2({\bf X})$, with ${\bf X}$ a generic matrix, for the defining
ideals of the products $\Bbb P^r_k \times \Bbb P^s_k$.
\item[(b)]
  $I_2({\bf X})$, with ${\bf X}$ a generic symmetric matrix, for the
defining ideals of quadratic Veronese embeddings of $\Bbb P^r_k$.
\item[(c)]
 $Pf_2({\bf X})$, the ideal generated by the pfaffians of a generic
alternating matrix, for the defining ideal of the Pl\"{u}cker
embedding of $G(2, r)$.
\end{itemize}

In these cases we get $a_* (I^e) \leq 2 e $, for all $e$. A.
Bertram gets the following bounds for $M=\max_{i \geq 2} \{a_i
(I^e)\}$ :

\begin{itemize}
\item[(a)] $M \leq 2e-4$.
\item[(b)] $M \leq 2e -3$.
\item[(c)] $M \leq 2e -6$.
\end{itemize}
}
\end{ex}

\bigskip
%
%
\section{Bayer--Stillman Theorem}
\label{AA} \markboth{CHAPTER IV. a-INVARIANTS OF THE POWERS OF AN
IDEAL}{BAYER--STILLMAN THEOREM}

\medskip
Let $S=k[{\x},{\y}]$ be the polynomial ring in $n+r$ variables
with the bigrading given by deg$(X_i)= (1, 0)$, deg$(Y_j)= (0,
1)$, so that $S$ is a standard bigraded $k$-algebra. For a
homogeneous ideal $I$ in $S$, we have already defined in Section
5.1 the bigraded regularity ${\rm \bf reg} (I)$ of $I$. The aim of
this section is to give a new description of the regularity of $I$
analogous to the one given by D. Bayer and M. Stillman \cite{BS1}
in the graded case. To this end, we are going to prove several
technical lemmas which are the bigraded version of the ones in
\cite{BS1}. To state them, we need to introduce the ${\it
saturations}$ of $I$ with respect to the variables $\underline X$
and $\underline Y$ ($I^{*1}$ and $I^{*2}$), and the generic forms
for $I$ with respect to $\underline X$ and $\underline Y$.
Furthermore, Theorem \ref{BB1} and Theorem \ref{BB2} will be
needed to prove some of these lemmas. We will include all the
proofs for the completeness.

\medskip
For a given finitely generated bigraded $S$-module $L$, we say
that $L$ is $(m, \cdot)$-regular if $\reg_1 L \leq m$. Similarly,
$L$ is $(\cdot, m)$-regular if $\reg_2 L \leq m$. Denote by
${\M}_1$ and ${\M}_2$ the ideals of $S$ generated by ${\frak m}_1=
(X_1, ..., X_n)$ and ${\fm}_2 = ({\y}) $ respectively. Then we
have

\medskip
\begin{prop}
\label{AA11} Let $L$ be a finitely generated bigraded $S$-module.
Then the following are equivalent:
\begin{enumerate}
   \item $L$ is $(m, \cdot)$-regular.
   \item $H^i_{{\cal M}_1} (L)_\pq =0$ for all $i$, $q$, $p \geq m-i+1$.
\end{enumerate}
\end{prop}

{\pf} By Theorem \ref{BB2}, $L$ is $(m, \cdot)$-regular if and
only if $\reg (L^q) \leq m$ for any $q$, that is,
$H^i_{{\fm}_1}(L^q)_p=0$ for all $i$, $q$ and $p \geq m-i+1$. Now
the result follows from Proposition \ref{W11}. $\B$

\medskip
Given a homogeneous ideal $I$, let us define $I^{*1}$ and $I^{*2}$
to be the homogeneous ideals
$$I^{*1} = \{ f \in S : \exists \; k \; s.t. \;{{\cal M}_1}^k f \subset
I \},$$
$$I^{*2} = \{ f \in S : \exists \; k \; s.t. \;{{\cal M}_2}^k f \subset
I \}.$$

\begin{lem}
\label{AA2} Assume $k$ infinite, and let $s = \max \{i \mid
H^i_{{\cal M}_1}(S/I) \not =0 \}$. Then,
\begin{enumerate}
\item If $s=0$, $I^{*1}=S$.
\item If $s>0$, there exists $h \in S_{(1,0)}$ such that
$(I^{*1} : h) = I^{*1}$.
\end{enumerate}
\end{lem}

{\pf} First note that

\vspace{1.5mm}

\hspace{10mm} $H^0_{{\cal M}_1}(S/I)= \{ \, \overline{f} \in S/I
\mid \exists \;k \;{s.t.} \; {{\cal M}_1}^k \;  \overline{f} = 0
\}=$

\vspace{1.5mm}

\hspace{30mm} $= \{\, f \in S \mid \exists \; k \;{s.t.}\;{{\cal
M}_1}^k \; f \subset I \} /I = I^{*1}/I.$

\vspace{1.5mm} \noindent If $s=0$, we have $ H^i_{{\fm}_1}
((S/I)^q)_p= H^i_{{\cal M}_1}(S/I)_{(p,q)} =0$ for any $p$, $q$,
$i>0$, so $(S/I)^q $ has dimension $0$ as graded $S_1$-module, and
then $H^0_{{\fm}_1} ((S/I)^q)= (S/I)^q$. Therefore, $H^0_{{\cal
M}_1}(S/I)=S/I$, and we get $I^{*1}= S$.

If $s>0$, note that $I^{*1} \not = S$ because otherwise
$H^0_{{\cal M}_1}(S/I) =S/I$ and then we would have $H^i_{{\cal
M}_1}(S/I) =0$ for all $i>0$. Now consider $\overline S=
S/I^{*1}$, and denote by $T = {\overline S} ^0 = S_1 /(I^{*1})^0$,
${\overline {{\fm}_1}} = {\fm}_1 T$. We have that $H^0_{{\cal
M}_1}(\overline S)= (I^{*1})^{*1}/I^{*1} = 0$, and so
$H^0_{{\fm}_1}({\overline S}^q)=0$ for all $q$. Therefore,
$\overline {{\fm}_1} \not \in \bigcup_q \Ass_{T} ({\overline
S}^q)$. On the other hand, according to \cite[Proposition
23.6]{HIO} we have that $\bigcup_q \Ass_{T}({\overline S}^q)$ is a
finite set. Since $k$ is infinite, we can find $h \in S_1$ of
degree 1 such that $\overline h \not \in z_{T}({\overline S}^q)$
for all $q$. Then $h \in S_{(1,0)}$ satisfies that $h \not \in
z_{S_1}(S/I^{*1})$. Therefore, $(I^{*1} : h) = I^{*1}$. $\B$

\medskip
From now on in this section we will assume that the field $k$ is
infinite. Let $s= \max \, \{i \mid H^i_{{\cal M}_1}(S/I) \not =0
\}$. If $s>0$, $h \in S_{(1,0)}$ is generic for $I$ if $h \not \in
z_{S_1}(S/I^{*1})$, that is, $(I^{*1}:h) =I^{*1}$. If $s=0$, we
say that any $h \in S_{(1,0)}$ is generic for $I$. Given $j>0$,
we define $U^1_j(I)$ to be the set
$$\{ (h_1, \dots, h_j) \in S_{(1,0)}^j
\mid h_i \, {\rm is} \; {\rm generic} \; {\rm for} \; (I, h_1,
\dots, h_{i-1}), 1 \leq i \leq j \}.$$

\begin{lem}
\label{AA3} Let $h \in S_{(1,0)}$. The following are equivalent:
\begin{enumerate}
\item $(I:h)_\pq = I_\pq$ for $p \geq m$.
\item $h$ is generic for $I$ and $(I^{*1})_\pq =I_\pq$ for $p \geq m$.
\end{enumerate}
\end{lem}

{\pf} First, let us notice that for $p > a_*^1(S/I)$ we have
$$(I^{*1}/I)_\pq =H^0_{{\cal M}_1}(S/I)_\pq = H^0_{{\fm}_1}
((S/I)^q)_p = 0 $$ by Theorem \ref{BB1}. Therefore, for $p$ large
enough it holds $(I^{*1})_\pq = I_\pq$, $\forall q$.

Now let us assume that $(i)$ holds, and let $f \in I^{*1}$ be a
homogeneous element not in $I$ such that deg$_1 f$ is maximum.
Then $hf \in I^{*1}$ has $\deg_1(hf)> \deg_1 f$, so $hf \in I$.
Hence $\deg_1 f < m$, and $(I^{*1})_\pq = I_\pq$ for any $p \geq
m$. To show that $h$ is generic for $I$, we may assume that $s=
\max \{i \mid H^i_{{\cal M}_1}(S/I) \not =0 \} >0$ (if not, any
element in $S_{(1,0)}$ is generic for $I$). Then we want to prove
$h \not \in z_{S_1}(S/I^{*1})$. Otherwise, there exists a
homogeneous element $f \not \in I^{*1}$ such that $hf \in I^{*1}$.
By Lemma \ref{AA2}, there exists $g \in S_{(1,0)}$ such that $g
\not \in z_{S_1}(S/I^{*1})$. Then, for any $s \geq 0$ we have $g^s
f \not \in I^{*1} $ and $h g^s f \in I^{*1}$, so $(I^{*1}:h)_\pq
\not = (I^{*1})_\pq$ for all $p \gg 0$. But note that for any $p
\geq m$,
$$(I^{*1}:h)_\pq = (I:h)_\pq = I_\pq = (I^{*1})_\pq,$$
and we get a contradiction.

Now assuming $(ii)$, we have that for $p \geq m$
$$(I:h)_\pq =(I^{*1}:h)_\pq = (I^{*1})_\pq =I_\pq. \; \B$$

\begin{lem}
\label{AA4} Let $h \in S_{(1,0)}$ be generic for $I$. The
following are equivalent:
\begin{enumerate}
\item $I$ is $(m,\cdot)$-regular.
\item $(I,h)$ is $(m,\cdot)$-regular and $(I^{*1})_\pq =I_\pq$ for all
$p \geq m$.
\end{enumerate}
\end{lem}

{\pf} If $I$ is $(m,\cdot)$-regular, then $S/I$ is
$(m-1,\cdot)$-regular. Then for any $i$, $q$ and $p \geq m-i$,  we
have $ H^i_{{\cal M}_1}(S/I)_\pq = 0$ by Proposition \ref{AA11}.
In particular, for $p \geq m$
$$ 0 =  H^0_{{\cal M}_1}(S/I)_\pq =
(I^{*1}/I)_\pq,$$ and so $(I^{*1})_\pq = I_\pq$ for $p \geq m$.

Let us consider $Q:= (I:h)/I$. In the assumptions of $(i)$ or
$(ii)$, observe that for $p \geq m$, $(I:h)_\pq = (I^{*1}:h)_\pq =
(I^{*1})_\pq = I_\pq$, so $Q_\pq = 0$ for $p \geq m$. Therefore
$H^i_{{\cal M}_1}(Q)=0$ for all $i>0$ and $H^0_{{\cal M}_1}(Q)=Q$.
From the bigraded exact sequence
       $$0 \to I \to (I:h) \to Q \to 0,$$
the long exact sequence of local cohomology gives
$$H^i_{{\cal M}_1} (I)_\pq \cong H^i_{{\cal M}_1}((I:h))_\pq \; ,
\forall i, p \geq m-i+1.$$

Assume first $(i)$. We have already shown that $(I^{*1})_\pq =
I_\pq$ for all $p \geq m$. Since $I$ is $(m,\cdot)$-regular, we
have $H^i_{{\cal M}_1}((I:h))_\pq= 0$ for all $i$, $p \geq m-i+1$.
By considering the exact sequence
$$0 \to I\cap (h)= (I:h)h = (I:h)(-1,0) \to I \oplus (h) \to (I,h) \to
0,  $$ we get $H^i_{{\cal M}_1}((I,h))_\pq= 0$ for all $i$, $p
\geq m-i+1$, so $(I,h)$ is $(m,\cdot)$-regular.

   Now by assuming $(ii)$, from the previous exact sequence we obtain
that $H^i_{{\cal M}_1}((I:h))_{(p-1,q)} \cong H^i_{{\cal M}_1}
(I)_\pq$ for $p \geq m-i+2$. For $p \geq m-i+1$, we then have that
$H^i_{{\cal M}_1} (I)_\pq \cong H^i_{{\cal M}_1}((I:h))_\pq \cong
H^i_{{\cal M}_1} (I)_{(p+1,q)}$. Therefore $H^i_{{\cal
M}_1}(I)_{(p,q)}= 0$ for $p \geq m-i+1$, so $I$ is
$(m,\cdot)$-regular. $\B$

\begin{lem}
\label{AA5} Let $I$ be an ideal generated by forms in $\deg_1 \leq
m$ and $h \in S_{(1,0)}$. If $(I,h)$ is $(m,\cdot)$-regular, then
$(I:h)$ is generated by forms in $\deg_1 \leq m$.
\end{lem}

{\pf} Let $f_1, \dots, f_u, h f_{u+1}, \dots, h f_v$ be a minimal
system of homogeneous generators for $I$, where $f_1, \dots,
f_u,h$ is a minimal system of generators for $(I,h)$. If $f \in
(I:h)$, then
$$hf= g_1 f_1+ \dots + g_u f_u + h(g_{u+1} f_{u+1} + \dots + g_v f_v),$$
for $g_1, \dots, g_v \in S$. Thus
$$(f - g_{u+1} f_{u+1} - \dots - g_v f_v)h -
 g_1 f_1- \dots - g_u f_u = 0 .$$
The first map in the bigraded minimal free resolution of $(I,h)$
is
 $$\begin{array}{ccc}
  Se \oplus Se_1 \oplus \dots \oplus Se_u  & \longrightarrow & (I,h) \\
                   e   & \longmapsto &  h  \\
                   e_j & \longmapsto & f_j
 \end{array}$$
and we have that
$$(f - g_{u+1} f_{u+1} - \dots - g_v f_v) e -
 g_1 e_1- \dots - g_u e_u  $$
is a first syzygy of $(I, h)$. Conversely, if $le+l_1e_1+ \dots
+l_ue_u$ is a first syzygy of $(I,h)$ then $l h+l_1f_1+
\dots+l_uf_u=0$, so $lh \in (f_1, \dots, f_u) \subset I$, and $l
\in (I:h)$. Because $(I, h)$ is $(m, \cdot)$-regular, each first
syzygy of $(I,h)$ can be expressed in terms of syzygies of $(I,h)$
in $\deg_1 \leq m+1$.
 Then
$$(f - g_{u+1} f_{u+1} - \dots - g_v f_v) e -
 g_1 e_1- \dots - g_u e_u = \sum_i \lambda_i
(\gamma_i e + \gamma_{i1} e_1 + \dots+ \gamma_{iu} e_u ),$$
 with
$\deg_1(\gamma_i e + \gamma_{i1} e_1 + \dots+ \gamma_{iu} e_u )
\leq m+1$. So
$$f = g_{u+1} f_{u+1} + \dots + g_v f_v + \sum_i \lambda_i \gamma_i,$$
with $\gamma_i \in (I:h)$, $\deg_1 \gamma_i \leq m$. Since
$f_{u+1}, \dots, f_v$ also belong to $(I:h)$ and have $\deg_1 \leq
m$, we finally obtain that $(I:h)$ can be generated by elements in
$\deg_1 \leq m$. $\B$

\medskip
  We are now ready to prove a bigraded version of the
Bayer-Stillman's Theorem characterizing the regularity of a
homogeneous ideal in terms of generic forms.

\begin{thm}
\label{AA6} Let $I$ be a homogeneous ideal in $S$ generated by
forms in $\deg_1 \leq m$. Then the following are equivalent:
\begin{enumerate}
\item $I$ is $(m,\cdot)$-regular.
\item There exist $h_1, \dots, h_j \in S_{(1,0)}$ for some $j \geq 0$
such that
$$((I,h_1, \dots,h_{i-1}):h_i)_{(m,q)} =
 (I,h_1,\dots,h_{i-1})_{(m,q)} , \; \forall q, \;1 \leq i \leq j.$$
$$(I,h_1,\dots,h_j)_{(m,q)} =S_{(m,q)}, \; \forall q.$$
\item Let $s = \max \{i \mid H^i_{{\cal M}_1}(S/I) \not = 0 \}$.
 For all $(h_1, \dots, h_s) \in U_s^1(I)$, $p \geq m$,
$$((I,h_1, \dots,h_{i-1}):h_i)_{(p,q)} =
 (I,h_1,\dots,h_{i-1})_{(p,q)} , \; \forall q, \; 1 \leq i \leq s.$$
$$(I,h_1,\dots,h_s)_{(p,q)} =S_{(p,q)}, \; \forall q.$$

\end{enumerate}
\end{thm}

{\pf} Note that $(iii) \Rightarrow (ii)$ is obvious. Now we are
going to show that $(ii) \Rightarrow (i)$ by induction on $j$. If
$j=0$, we have that $I_\mq = S_\mq$ for all $q$, so $I_\pq =
S_\pq$ for all $q$ and $p \geq m$. Therefore,
$$H^i_{{\cal M}_1}(S/I)= \cases{0 & if $i>0$ \cr
                                S/I & if $i=0$ \cr}.$$
In particular, we have that $H^i_{{\cal M}_1}(S/I)_\pq=0$ for all
$i$, $q$ and $p \geq m-i$, so $I$ is $(m,\cdot)$-regular. If
$j>0$, we have that $(I,h_1)$ is $(m,\cdot)$-regular by the
induction hypothesis. Since $I$ is generated by forms in $\deg_1
\leq m$, we have that $(I:h_1)$ is generated by forms in $\deg_1
\leq m$ by Lemma \ref{AA5}. As $(I:h_1)_\mq =I_\mq$, we then
conclude $(I:h_1)_\pq =I_\pq$ for all $p \geq m$. According to
Lemma \ref{AA3}, we have that $h_1$ is generic for $I$ and
$(I^{*1})_\pq = I_\pq$ for all $p \geq m$. Then $I$ is
$(m,\cdot)$-regular by Lemma \ref{AA4}.

  Now let us prove $(i) \Rightarrow (iii)$ by induction on $s$.
If $s=0$, since $I$ is $(m,\cdot)$-regular we have $H^0_{{\cal
M}_1}(S/I)_\pq = (I^{*1}/I)_\pq =0$ for $p \geq m$, and $I^{*1}=S$
by Lemma \ref{AA2}. Therefore, $I_\pq = (I^{*1})_\pq =S_\pq$ for
$p \geq m$. Assume now $s>0$. Since $I$ is $(m,\cdot)$-regular and
$h_1$ is generic for $I$, we get that $(I,h_1)$ is
$(m,\cdot)$-regular and $(I^{*1})_\pq =I_\pq$ for all $p \geq m$
by Lemma \ref{AA4}. As $(h_2, \dots, h_s) \in U_{s-1}^1((I,h_1))$,
by the induction assumption it is just enough to show $(I:h_1)_\pq
=I_\pq$ for $p \geq m$, which is satisfied by Lemma \ref{AA3}.
$\B$

\medskip
We are going to use this criterion to compute the bigraded
regularity of the generic initial ideal of a homogeneous ideal $I$
in $S$. Let $S=k[{\x},{\y}]$ be the polynomial ring over an
infinite field $k$ with the bigrading given by deg$(X_i)= (1, 0)$,
deg$(Y_j)= (0, 1)$. Let $<$ be an order on the monomials of $S$.
Let us denote by ${\cal G}= {\cal G}_1 \times {\cal G}_2$, with
${\cal G}_1=GL(n,k)$, ${\cal G}_2= GL(r,k)$. Given an element
$g=(f,h) \in {\cal G}$, where $f=(f_{ij})_{1 \leq i,j \leq n}$ and
$h=(h_{ij})_{1 \leq i,j \leq r}$, $g$ acts on $S$ by acting on the
variables in the following way
$$ X_j \longmapsto \sum_{i=1}^n f_{ij} X_i
\hspace{5mm} , \hspace{5mm}
 Y_j \longmapsto \sum_{i=1}^r h_{ij} Y_i.$$
We will denote by ${\cal B} = {\cal B}_1 \times {\cal B}_2$, where
${\cal B}_1, {\cal B}_2$ are the Borel subgroups of ${\cal G}_1,
{\cal G}_2$ consisting of upper triangular matrices, and by ${\cal
B}' = {\cal B}'_1 \times {\cal B}'_2$, where ${\cal B}'_1, {\cal
B}'_2$ are the Borel subgroups of ${\cal G}_1, {\cal G}_2$
consisting of lower triangular matrices. We will denote by ${\cal
U} = {\cal U}_1 \times {\cal U}_2$, where ${\cal U}_1, {\cal U}_2$
are the unipotent matrices. By bigrading the proof of
\cite[Theorem 15.18]{E}, we get

\begin{thm} (Galligo, Bayer--Stillman)
\label{bB1} Let $I \subset S$ be a homogeneous ideal. There exists
a non-empty Zariski open $ U = {\cal B}' U \subset{\cal G}$, $U
\cap {\cal U} \not = Id$, and a monomial ideal $J$ such that
$${\rm in}(gI) = J , \; \forall  g \in U .$$
\end{thm}

\medskip
We call $J$ the (bi)graded generic initial ideal of $I$, written
$J= {\bf gin}(I)$. Given a homogeneous ideal $I \subset S$, we say
that $I$ is Borel-fix if $g I = I$ for any $g \in {\cal B}$. It
was proved that the generic initial ideal of a graded ideal is
Borel fix. By bigrading the proof of \cite[Theorem 15.20]{E}, we
easily obtain that the generic initial ideal is Borel-fix.

\begin{thm}
\label{bB2} Let $I \subset S$ be a homogeneous ideal. For any $g
\in {\cal B}$,
 $$g ({\bf gin}(I))= {\bf gin}(I).$$
\end{thm}

\medskip
Let $p \geq 0$. Given $s,t \in \Bbb N$, we define $s <_p  t  \iff
{t \choose s} \not  \equiv 0 \; (mod \, p)$. We also can give an
equivalent characterization of the Borel-fix bihomogeneous ideals
analogous to the one in the graded case \cite[Theorem 15.23]{E}.
Namely,

\begin{thm}
\label{bB3} Let $I$ be a homogeneous ideal of $S$. Let $p = char k
\geq 0$. Then
\begin{enumerate}
\item $I$ is diagonal-fix iff $I$ is monomial.
\item $I$ is Borel-fix iff $I$ is generated by
monomials $m$ such that satisfy the following conditions

 \begin{itemize}
   \item[--] If $m$ is divisible by $X_j^t$ but by no
    higher power of $X_j$, then $( X_i / X_j) ^s m \in I$, $\forall
   i<j$, $s<_p t$.
   \item[--] If $m$ is divisible by $Y_j^t$ but by no higher
    power of $Y_j$, then $( Y_i / Y_j) ^s m \in I$, $\forall i<j$, $s<_p t$.
  \end{itemize}

\end{enumerate}

\end{thm}

\medskip
For a homogeneous ideal $I$, let us denote by $\delta_1(I)$ the
maximum first component of the degree in a minimal system of
generators of $I$. In a similar way, we may define $\delta_2(I)$.
Then we have

\begin{prop}
\label{bB4} Let $I \subset S$ be a Borel-fix ideal. If $char k=0$,
then
$$\reg_1(I) = \delta_1(I),$$
$$\reg_2(I) = \delta_2(I).$$
\end{prop}

{\pf} Set $m = \delta_1(I)$. From the definition of the regularity
it is clear that $\reg_1(I) \geq m$. According to Theorem
\ref{AA6}, to prove the equality it is enough to show that for $p
\geq m$, $i \leq n$, we have
$$((I,X_n, \dots,X_{i+1}):X_i)_\pq =(I,X_n,\dots,X_{i+1})_\pq.$$

Let $f \in ((I,X_n, \dots,X_{i+1}):X_i)$  be a monomial with
$\deg_1 f \geq m$. If there exists $k \geq i+1$ such that $X_k
\vert f$, we immediately have $f \in (I,X_n,\dots,X_{i+1})$.
Otherwise,
$$X_i f = X^{\alpha} Y^{\beta} (X^A Y^B),$$
where $X^A Y^B \in I$, $\deg_1 X^A \leq m$, $\deg_1 X^\alpha \geq
1$. If $X_i \vert X^\alpha$, we then easily get $f \in I$. If not,
by taking $k \leq i$ such that $X_k \vert X^\alpha$, we can write
$$f = \frac{X^{\alpha}}{X_k} Y^{\beta} (\frac{X_k}{X_i} X^A Y^B).$$
Since $I$ is Borel-fix,  we have that $\frac{X_k}{X_i} X^A Y^B \in
I$ by Theorem \ref{bB3}, and so $f \in I \subset
(I,X_n,\dots,X_{i+1})$. $\B$

\medskip
This result has been proved independently by A. Aramova et al.
\cite{ADK} by other methods.

\medskip
In the graded case, it was proved by D. Bayer and M. Stillman
\cite{BS1} that there exists an order in $A = k[{\x}]$ (the
reverse lexicographic order) with the property that $\reg (I)=
\reg ({\rm gin} I)$ for any homogeneous ideal $I$ in $A$. We may
wonder if the analogous bigraded result also holds, that is, if
there exists an order in the polynomial ring $S$ such that ${\rm
\bf reg}(I) = {\rm \bf reg}({\rm \bf gin} \, I)$ for any
homogeneous ideal $I$. We show that this is not true by giving a
homogeneous ideal in $S$ such that ${\rm \bf reg}(I) \not = {\rm
\bf reg}({\rm \bf gin} \,I )$ for any order on $S$.

\begin{ex}
\label{bB5} {\rm Let us consider the polynomial ring $S=k[X_1,
X_2, Y_1, Y_2]$, with $\deg (X_1)= \deg (X_2)=(1,0)$, $\deg (Y_1)=
\deg (Y_2)=(0,1)$. Let $>$ be a term order in $S$, that is, an
order satisfying
\begin{enumerate}
\item $X_1>X_2$, $Y_1>Y_2$.
\item For monomials $m$, $m_1$, $m_2$ in $S$, if $m_1> m_2$ then
$m m_1> m m_2$.
\end{enumerate}
Let $I$ be the homogeneous ideal in $S$ generated by the forms
$f_1=X_1 Y_1$ and $f_2= X_1Y_2+ X_2 Y_1$ in degree $(1,1)$. Note
that $f_1, f_2$ is a regular sequence, so the Koszul complex of
these forms provides the minimal bigraded free resolution of $I$:
$$0 \to S(-2,-2) \to S(-1,-1)^2 \to I \to 0.$$
Then the regularity of $I$ is ${\rm \bf reg} (I) =(1,1)$. Note
that $X_1 Y_1> X_1 Y_2$, $X_2 Y_1 >X_2 Y_2$. Therefore, if we want
to define an order on the monomials of $S$ we only must decide if
$X_1 Y_2 >X_2 Y_1$ or $X_1 Y_2 <X_2 Y_1$ in degree $(1,1)$. Assume
first that $X_1 Y_2 >X_2 Y_1$. Recall that $g \in GL(2,k) \times
GL(2,k)$, with
$$g= (A,B) =\bigg ( \left(\begin{array}{cc}
                      a & c \\
                      b & d
                      \end{array} \right),
          \left(\begin{array}{cc}
                     \alpha & \gamma \\
                     \beta & \rho
                     \end{array} \right) \bigg),$$
operates in $S$ by means of
   $$  \begin{array}{lll}
     X_1    & \longmapsto  &   a X_1+ b X_2  \\
     X_2    & \longmapsto  &   c X_1+ d X_2  \\
     Y_1    & \longmapsto  &   \alpha Y_1+ \beta Y_2  \\
     Y_2    & \longmapsto  &   \gamma Y_1+ \rho Y_2
      \end{array} $$

Since $\dim_k {\rm \bf gin}(I)_{(i,j)}=\dim_k I_{(i,j)}$ for any
$(i,j)$, we have that ${\rm \bf gin}(I)_{(i,j)}=0$ for $(i,j) \in
\{(0,0),(1,0), (0,1) \}$. In degree $(1,1)$, the forms $f_1, f_2$
are a $k$-basis of $I_{(1,1)}$. By computing $g (f_1 \wedge f_2)$,
we get

\vspace{1.5mm}

\hspace{15mm} $g ( f_1 \wedge f_2) = a^2 (\alpha  \rho - \beta
\gamma) X_1Y_1 \wedge X_1Y_2 + \dots $

\vspace{1.5mm}

\noindent so ${\rm \bf gin}(I)_{(1,1)}$ is the $k$-vector space
generated by $X_1 Y_1$, $X_1 Y_2$. If ${\rm \bf gin}(I) = (X_1
Y_1, X_1 Y_2)$, then $\dim_k {\rm \bf gin}(I)_{(1,2)}=3$ because
$X_1 Y_1^2$, $X_1 Y_1 Y_2$, $X_1 Y_2^2$ is a $k$-basis, which is a
contradiction because $\dim_k I_{(1,2)} \geq 4$. Therefore, ${\rm
\bf gin}(I)$ has minimal generators with $\deg_1 \geq 2$ or
$\deg_2 \geq 2$, so $\reg_1 ({\rm \bf gin} \,I) \geq 2$ or $\reg_2
({\rm \bf gin}\, I )\geq 2$. In the case $X_1 Y_2 < X_2 Y_1$, it
can be proved that ${\rm \bf reg} ({\rm \bf gin}\, I) \not =
(1,1)$ by similar arguments. Therefore, we get
 ${\rm \bf reg}(I) \not = {\rm \bf reg}({\rm \bf gin} \,I )$ for any
order in $S$. }
\end{ex}

\medskip
Finally, note that these results can be applied to study the
Koszul property of the diagonals of a bigraded standard
$k$-algebra. By using \cite[Theorem 18]{ERT}, and following the
same lines as \cite[Theorem 2]{ERT} in the graded case, it can be
proved that that for a homogeneous ideal $I$ of $S$,
$(S/I)_\Delta$ has a Gr\"{o}bner basis of quadrics for $c \gg 0$, $e
\gg 0$ (see also \cite{ADK}).

\chapter*{$\;\;\;$}

\newpage

\medskip
\chapter{Asymptotic behaviour of the powers of an ideal}
\typeout{Asypmtotic behaviour of the powers of an ideal}

\medskip
Let $A= k[{\x}]$ be a polynomial ring over a field $k$, and let
$I$ be a homogeneous ideal in $A$. In this chapter we are
concerned with the asymptotic behaviour of the powers of $I$. We
will use the bigraded structure of the Rees algebra to get
information about the Hilbert polynomials, the Hilbert series and
the graded minimal free resolutions of the powers of $I$.

\medskip
In Section 6.1 we show that the Hilbert polynomials of the powers
of the ideal $I$ have a uniform behaviour. In particular, the
Hilbert polynomials of a finite set of these powers allow to
compute the Hilbert polynomials of its Rees algebra and its form
ring, without needing an explicit presentation of these algebras.
In Section 6.2, similar results are stated for the Hilbert series
of the powers of $I$.

\medskip
The last section begins by studying the projective dimension of
the powers of $I$. The approach to this question by means of the
bigrading of the Rees algebra allows to recover some classical
results as the constant asymptotic value for the projective
dimension, as well as to determine the powers of the ideal which
take the asymptotic value whenever the form ring is Gorenstein.
After that, we study the graded minimal free resolutions of the
powers of an ideal. In the equigenerated case, it is proved that
the shifts are given by linear functions asymptotically and the
graded Betti numbers of these resolutions are given by polynomials
asymptotically. This result is then applied to guess the
resolutions of the powers of some families of ideals from a finite
set of these resolutions.

\bigskip
\section{Hilbert polynomial of the powers of an ideal}
\label{D}
\medskip
\label{FF} \markboth{CHAPTER VI. ASYMPTOTIC BEHAVIOUR}{HILBERT
POLYNOMIAL OF THE POWERS OF AN IDEAL}

\medskip
First of all, let us recall some standard definitions and
notations referred to the Hilbert polynomial (see for instance
\cite{BH}). Let $A= k[{\x}]$ be the polynomial ring over a field
$k$, and let $M$ be a finitely generated graded $A$-module. The
numerical function

\vspace{2mm}

\hspace{29mm} $H(M, \;) : \Bbb Z \longrightarrow \Bbb Z $

\hspace{46mm} $ j \mapsto \dim_k M_j$

\vspace{2mm} \noindent is the Hilbert function of $M$. Denoting by
$d = \dim M$, there exists an unique polynomial $P_M(s) \in \Bbb
Q[s]$, of degree $d-1$, for which $H(M,j) =P_M(j)$ for all $j \gg
0$. We can write
$$P_M(s) = \sum_{k=0}^{d-1} (-1)^{d-1-k} e_{d-1-k} {s+k \choose k},$$
with $e_0, \dots, e_{d-1} \in \Bbb Z$. $P_M(s)$ is called the
Hilbert polynomial of $M$.

\medskip
  Our first result shows the uniform behaviour of the Hilbert polynomial
of the powers of any homogeneous ideal in a polynomial ring.

\begin{thm}
\label{DD1} Let $I$ be a homogeneous ideal in $A$. Set $c=
a_*^2(R_A(I))$, $h = \h (I)$. Then there are polynomials $e_0(j),
\dots, e_{n-h-1}(j)$ with integer values such that for all $j \geq
c+1$
$$P_{A/I^j}(s) = \sum_{k=0}^{n-h-1} (-1)^{n-h-1-k}e_{n-h-1-k}(j) {s+k
\choose k}.$$ Furthermore, $\deg e_{n-h-1-k}(j) \leq n-k-1$ for
all k.

\end{thm}

{\pf} Assume that $I$ is generated by forms $f_1, \dots, f_r$ in
degrees $d_1 \leq \dots \leq d_r=d$ respectively. Then the Rees
algebra $R=R_A(I)$ of $I$ can be endowed with the bigrading given
by $R_{(i,j)} =(I^j)_i$, so that $R$ is a bigraded $S$-module, for
$S=k[{\x},{\y}]$ the polynomial ring with deg$(X_i)= (1, 0)$,
deg$(Y_j)= (d_j, 1)$. Since $R$ is a domain, it has relevant
dimension $n+1$. Then by Proposition \ref{F0} and Proposition
\ref{F2} there exists a polynomial $P_R(s,t)$ of total degree
$n-1$ such that for all $(i, j)$
 $$ \dim_k R_{(i,j)} - P_R (i,j)= \sum_q (-1)^q \dim_k
H^q_{R_+}(R)_{(i,j)},$$ where $R_+$ is the ideal generated by the
products $X_i f_jt$ for $1 \leq i \leq n$, $1 \leq j \leq r$. By
taking in $S$ the homogeneous ideals ${\cal M}_1 = ({\x}) S$ and
${\cal M}_2=({\y}) S$, the Mayer-Vietoris long exact sequence
gives then
$$ \cdots \to H^q_{{\cal M}_1}(R) \oplus H^q_{{\cal M}_2}(R) \to
H^q_{S_+}(R) \to H^{q+1}_{\cal M}(R) \to \cdots $$ Notice that for
$j>c$ we have $H^q_{\cal M}(R)_{(i,j)} = 0$ for all $i, q$. Then,
by Proposition \ref{W11} we also get $H^q_{{\cal M}_2}(R)_{(i,j)}
= 0$ for $j > c$, for all $i,q$. Furthermore, $H^q_{{\cal
M}_1}(R)_{(i,j)} =H^q_{\frak m}(I^j)_i$ for any $j \geq 0$ by
Proposition \ref{W11}. Therefore, for any $j >c$ there exists an
integer $i_0=a_*(I^j)$ (depending on $j$) such that
$H^q_{R_+}(R)_{(i,j)} = H^q_{S_+}(R)_{(i,j)}= H^q_{\frak m}(I^j)_i
= 0$ for all $q$ and $i >i_0$. Hence
 $P(i,j) = \dim_k R_{(i,j)} = \dim_k (I^j)_i$ for any $j>c$ and
$i>i_0$.

   Now, by defining $P_j(s) = {n+s-1 \choose n-1}- P (s,j)$,
for $j > c$, $s \gg 0$ we have that $P_j(s) = \dim_k (A/I^j)_s$.
Hence $P_j(s)$ is the Hilbert polynomial of $A/I^j$. Furthermore,
we can write

\vspace{2mm}

\hspace{25mm} $P_j(s) = {n+s-1 \choose n-1}- P (s,j) $

\vspace{3mm}

\hspace{35mm} $={n+s-1 \choose n-1} - \sum_{l+m \leq n-1} a_{lm}
{s-dj \choose l}{j \choose m} $

\vspace{3mm}

\hspace{35mm} $= \sum_{k=0}^{n-1} b_k(j)
  {s+k \choose k},$

\vspace{2mm}

\noindent with $b_{k}(j)$ polynomials in $j$. Since $\deg P_j(s) =
n-h-1$ for any $j >c$, we have $b_{k}(j)=0$ for $k \geq n-h$ and
$j>c$, so $b_{k}(j) \equiv 0$ for $k \geq n-h$. Then we may write
$$P_j(s) = \sum_{k=0}^{n-h-1} (-1)^{n-h-1-k}
 e_{n-h-1-k}(j) {s+k \choose k},$$
for $j>c$. Moreover, since $P_R(s,t)$ has total degree $n-1$ and
$P_R(s,t) = {n+s-1 \choose n-1} - P_t (s)$, we easily obtain that
$\deg e_{n-h-1-k}(j) \leq n-k-1$. $\B$

\begin{rem}
\label{DD35} {\rm We have seen that $\deg e_{n-h-1-k}(j) \leq
n-k-1$ for all $k$, so in particular the polynomial $e_0(j)$ which
gives the multiplicity of $A/I^j$ has degree $\leq h$. By Nagata's
formula,
$$e_0(j)= e(A/I^{j}) = \sum_{\frak p \in {\rm Assh}(A/I)} {\rm
length}(\,{A_{\frak p}}/I_{\frak p}^{j} ) e(A/\frak p),$$ with
Assh$(A/I) = \{ \, \frak p \in \Ass (A/I) \mid \dim A/\frak p =
\dim A/I \}$. Note that for all those $\frak p$, we have that
$\dim A_{\frak p}= h$ and then
$${\rm length}(\, A_{\frak p}/I_{\frak p}^{j}) =
e(IA_{\frak p},A_{\frak p}) {h+j \choose j} + {\rm polynomial}
\;{\rm in} \;j\; {\rm of} \; {\rm degree}\; {\rm lower}\;{\rm
than} \; h .$$ Therefore $e_0(j)$ has degree $h$, so let us write
$$e_0(j)= \lambda_h {j \choose h} +
{\rm polynomial} \;{\rm in} \;j\; {\rm of} \; {\rm degree}\; {\rm
lower}\;{\rm than} \; h .$$ We can give an upper bound for the
leading coefficient $\lambda_h$. According to \cite[Corollary
3.8]{HS}, we have
$$ e_0(j) \leq {\reg (I^j) +h -1 \choose h }. $$
Assume that $I$ is generated by forms in degree $\leq d$. Then
there exists a positive integer $\alpha$ such that $\reg(I^j) \leq
dj + \alpha$ by Theorem \ref{S8}, and so $\lambda_h \leq d^h$. In
Proposition \ref{DD8} we will show that $\lambda_h$ and, more
generally, the leading coefficients of the polynomials
$e_{n-h-1-k}(j)$ play an important role in the mixed
multiplicities of the Rees algebra and the form ring. }
\end{rem}

\medskip
Now let us consider a homogeneous ideal $I$ generated by forms in
degree $d$. Let us take the Hilbert polynomial $P_R(s,t)$ of its
Rees algebra with the usual bigrading, and let us write
$$P_R(s + dt ,t) = \sum_{k+m \leq n-1} a_{km} {s \choose k}{t \choose m}.$$
 Following \cite{HHRT}, we call $e_i(R)= a_{i,n-1-i}$ the mixed
multiplicity of $R$ of type $i$ for $i= 0, \dots, n-1$. According
to Proposition \ref{F0} we have $e_i(R) \geq 0$, and then $e(R) =
\sum_{i=0}^{n-1}e_i(R)$ by \cite[Theorem 4.3]{HHRT}. Next we are
going to study the multiplicity of the Rees ring and to relate it
to the multiplicity of the form ring. First, we need to compute
the relevant dimension of the form ring.

\begin{lem}
\label{DD75} Let $I$ be a homogeneous ideal in $A$ generated by
forms in degree $d$. Then the relevant dimension of $G$ is $n$ if
and only if $I$ is not $\frak m$-primary.
\end{lem}

{\pf} If $I$ is $\frak m$-primary, then $G_+ \subset P$ for any
homogeneous prime $P \in \Spec (G)$ because $G_{(1,0)}$ is
nilpotent, and so $\rdim G =1 < \dim G =n $. If $I$ is not
$\fm$-primary, for $k = 0, \dots, n-h-1$ let us write
$$e_{n-h-1-k}(j)= \lambda_{n-k-1}{j \choose n-k-1}+ {\rm polynomial
\; in \;} j  {\rm \; of \; lower \;degree}. $$ Then for $s \gg 0$,
$t \gg 0$, we have

\vspace{3mm}

\hspace{2mm} $P_G(s+dt,t) = P_{A/I^{t+1}}(s+dt)-P_{A/I^t}(s+dt) $

\vspace{4mm}

\hspace{10mm} $= \sum_{k=0}^{n-h-1} (-1)^{n-h-1-k}
 (e_{n-h-1-k}(t+1)- e_{n-h-1-k}(t))   {s+dt+k \choose k} $

\vspace{4mm}

\hspace{10mm} $= \sum_{k=0}^{n-h-1} (-1)^{n-h-1-k}
 \lambda_{n-k-1}
 ({t+1 \choose n-k-1}-{t \choose n-k-1})  {s+dt+k \choose k} \;\; +$

\vspace{2mm}

\hspace{14mm}  $+ \;$ polynomial in $s,t$ of lower total degree

\vspace{4mm}

\hspace{10mm} $= \sum_{k=0}^{n-h-1} (-1)^{n-h-1-k} \lambda_{n-k-1}
 {t \choose n-k-2}  \frac{(s+dt)^k}{k!} \;\; +$

\vspace{2mm}

\hspace{14mm}  $+ \;$ polynomial in $s,t$ of lower total degree

\vspace{4mm}

\hspace{10mm} $= \sum_{i=0}^{n-h-1} \bigg[ \sum_{k = i}^{n-h-1}
(-1)^{n-h-1-k}  \lambda_{n-k-1} d^{k-i}
 {n-2-i \choose k-i} \bigg]
{s \choose i} {t \choose n-2-i} \;\; +$

\vspace{2mm}

\hspace{14mm}  $+ \;$ polynomial in $s,t$ of lower total degree

\vspace{3mm}

  In particular,  $\lambda_h$ is the coefficient of
${s \choose n-h-1} {t \choose h-1}$ which is not zero by Remark
\ref{DD35}. So the total degree of the Hilbert polynomial of the
form ring is $n-2$, and then the relevant dimension of $G$ is $n$
by Proposition \ref{F0}. $\B$

\medskip
If $I$ is a homogeneous ideal generated by forms in degree $d$
which is not $\frak m$-primary, let us consider the Hilbert
polynomial of its form ring
$$P_G(s + dt ,t) = \sum_{k+m \leq n-2} b_{km} {s \choose k}{t \choose m}.$$
We call $e_i(G)= b_{i,n-2-i} \geq 0$ the mixed multiplicity of $G$
of type $i$ for $i=0, \dots, n-2$. Then $e(G) = \sum_{i=0}^{n-2}
e_i(G)$ again by \cite[Theorem 4.3]{HHRT}.

\medskip
 Now we can give the mixed multiplicities of the Rees algebra
and the form ring of an equigenerated ideal by means of the
leading coefficients of the polynomials  $e_{n-h-1-k}(j)$ given by
Theorem \ref{DD1}, and to relate the mixed multiplicities of both
rings.

\begin{prop}
\label{DD8}
 Let $I$ be a homogeneous ideal generated in degree $d$ which is not
$\frak m$-primary. Set $h = \h (I)$, $l= l(I)$. For each $k$, let
us write
 $$e_{n-h-1-k}(j)= \lambda_{n-k-1}{j \choose n-k-1}+
{\rm polynomial \; in \;} j {\rm \; of \; lower \; degree}.$$ Then
\begin{itemize}
\item[(i)] $ e_i(G) = 0$ if $i \geq n-h$ or $i \leq n-l-2$.
For each $n-l-2 <i <n-h$, we have
 $$e_i(G)=
\sum_{k=i}^{n-h-1} (-1)^{n-h-1-k} \lambda_{n-k-1} d^{k-i} {n-2-i
\choose k-i}.$$
\item[(ii)]
$ e_i(R) = 0$ if $i \leq n-l-1$.
 For each $i > n-l-1$, we have
     $$ e_i(R)   = \cases
{ d^{n-1-i}  \hspace{10mm}  {\rm if} \; i \geq n-h  & \cr
 d^{n-1-i} - \sum_{k=i}^{n-h-1}(-1)^{n-h-1-k} \lambda_{n-k-1} d^{k-i}
{n-1-i \choose k-i} & {\rm otherwise} \cr}
$$
\item[(iii)]
$e_i(G) = d e_{i+1}(R) - e_i(R)$, for $i=0, \dots, n-2$. In
particular, we have $e_i(R) \leq  d e_{i+1}(R)$, for $i=0, \dots,
n-2$. Furthermore,
$$e(G)= \cases {(d-1) e(R)+1 & if $l \leq n-1$ \cr
(d-1) e(R)+1- de_0(R) & if $l =n$ \cr} $$
\end{itemize}
\end{prop}

{\pf} Let us fix $j >a_*^2(G)$. Then we have that for $ s \gg 0$,
$$ \dim_k \bigg( \frac{I^j}{I^{j+1}} \bigg)_{s+dj} = P_G(s+dj,j) =
 \sum_{k+m \leq n-2} b_{km} {s \choose k}{j \choose m},$$
so $P_G(s+dj,j)$ is the Hilbert polynomial of the $A/I$-module
$I^j/I^{j+1}$ for large $j$. Hence $b_{km}=0$ for any $k \geq
n-h$, so in particular $e_i(G)=0$ for $i \geq n-h$.

Let us fix now $i >a_*^1(G^\varphi)$. Then we have that $H^q
_{\M}(G)_{(i+dj,j)} = H^q _{{\M}_1}(G)_{(i+dj,j)} = 0$ for all $q$
and $j$ by Remark \ref{YY22}. From the Mayer-Vietoris long exact
sequence, we have that for $t \gg 0$,
$$ \dim_k \bigg( \frac{I^t}{I^{t+1}} \bigg)_{i+dt} = P_G(i+dt,t) =
 \sum_{k+m \leq n-2} b_{km} {i \choose k}{t \choose m}.$$
Therefore, we have that $P_G(i+dt,t)$ is the Hilbert polynomial of
the $F_{\fm}(I)=k[I_d]$-module $E_i = \bigoplus_{j \geq 0}
({I^j}/{I^{j+1}})_{i+dj}$. Hence $b_{km}=0$ for $m \geq l$. Then
the first part of $(i)$ is proved, and for the rest it suffices to
notice that for $s,t \gg0$, $P_G(s+dt,t) =
P_{A/I^{t+1}}(s+dt)-P_{A/I^t}(s+dt).$

To get $(ii)$ and $(iii)$, it is just enough to take into account
that for $s,t \gg 0$ we have
$$P_R(s+dt,t) = P_{A}(s+dt)-P_{A/I^t}(s+dt),$$
$$P_G(s+dt,t) = P_R(s+dt,t) - P_R(s+dt, t+1).\; \; \B$$

\medskip
\begin{rem}
\label{DD88} {\rm Let $(A, \fm, k)$ be a local ring, and $I \not =
A$ an ideal. Set $n = \dim A$, $l=l(I)$, $h= \h (I)$. Let us
denote by $\overline G = G_\fm(G_I(A))$ bigraded by means of
$$\overline G_{(i,j)} = \frac{\fm^i I^j+ I^{j+1}} {\fm^{i+1} I^j+
I^{j+1}}.$$ Since $\overline G$ is a standard bigraded
$k$-algebra, we may consider its Hilbert polynomial
$$P_{\overline G} (s,t) = \sum_{k+m \leq n-2} c_{km} {s \choose k}{t
\choose m},$$ and let us denote by $c_i(\overline G)=
c_{i-1,n-1-i}$ for $1 \leq i \leq n-1$. R. Achilles and M.
Manaresi \cite{AM} show that $ c_i(\overline G) = 0$ if $i> \dim
A/I$ or $i< n-l$.

   This definition and the results proved in \cite{AM} can be extended
to the graded case, and then we get that $ c_i(\overline G) = 0$
if $i> n-h$ or $i< n-l$. For a homogenous ideal $I$ in $A=k[{\x}]$
generated by forms of degree $d$, note that
$$\overline G_{(i,j)} = \bigg( \frac{I^j}{I^{j+1}} \bigg)_{i+dj} =
G_{(i+dj,j)},$$ so $e_i(G)= c_{i+1}(\overline G)$, and we get part
of $(i)$ of the previous proposition. In fact, the idea for
proving this part of $(i)$ is similar to \cite{AM}.}
\end{rem}

\begin{rem}
\label{DD888} {\rm  If $I$ is a complete intersection ideal then
$l=h$, and from Lemma \ref{DD8} $(ii)$ we have $ e_i(R) = 0$ if $i
\leq n-l-1$, and $ e_i(R) = d^{n-1-i}$ if $i \geq n-l$. This
result was proved in \cite[Theorem 3.6]{STV}.}
\end{rem}

\begin{cor}
\label{DD9} For a homogeneous ideal $I$ generated by forms in
degree $d$ of height $h$, we have
\begin{itemize}
\item[(i)]  $e(R) \geq 1+ d+ \dots + d^{h-1}$.
\item[(ii)] If $I$ is equimultiple,
$e(R) = 1 +d+ \dots+ d^{h-1}$. Assume further that $I$ is not a
$\frak m$-primary ideal. Then $e(G) = e_{n-h-1}(G) = \lambda_h=
d^h$.
\item[(iii)]  Assume that $I$ is not a $\frak m$-primary ideal and
$A/I^j$ is Cohen-Macaulay for $j \gg 0$ (or Buchsbaum). Then
$e_{n-h-1}(G) = d^h$, so $e(G) \geq d^h$.
\end{itemize}

\end{cor}

{\pf} $(i)$ and $(ii)$ are trivial. To get $(iii)$, according to
\cite[Proposition 2.3]{HRTZ} for $j \gg 0$ we have that $e(A/I^j)
\geq {dj+h-2 \choose h}$, and so $\lambda_h \geq d^h$.
Furthermore, $\lambda_h \leq d^h$ because $e_{n-h-1}(R)= d^{h}
-\lambda_{h} \geq 0$ by Lemma \ref{DD8}. We conclude
$e_{n-h-1}(G)= \lambda_h = d^h$, and so $e(G) \geq
e_{n-h-1}(G)=d^h$. $\B$

\medskip
Notice that as a consequence of Theorem \ref{DD1} we have that
with the Hilbert polynomials of a finite set of the powers of an
ideal we can compute the Hilbert polynomials of its Rees algebra
and its form ring, without needing an explicit presentation of
these bigraded algebras. For equigenerated ideals, we may also
compute the multiplicities of the Rees algebra and the form ring.
Namely,

\medskip
\begin{cor}
\label{DD2}
 Let $I$ be a homogeneous ideal in $A$. Set  $c=a_*^2(R_A(I))$,
$h = \h (I)$.
 Then the Hilbert polynomials of $I^j$ for $c + 1 \leq j \leq c+n$
determine
\begin{itemize}
\item[(i)] The polynomials $e_{n-h-1-k}(j)$ for $k = 0, \dots, n-h-1$.
\item[(ii)] The Hilbert polynomials of $A/I^j$ for $j> c + n$.
\item[(iii)] The Hilbert polynomial of the Rees algebra of $I$ and the
Hilbert polynomial of the form ring of $I$.
\item[(iv)] If $I$ is equigenerated and not $\frak m$-primary, the
 multiplicity of the Rees algebra of $I$ and the multiplicity of the
form ring of I.
\end{itemize}
\end{cor}

\medskip
We describe all these computations by means of an explicit
example.

\begin{ex}
\label{DD6}
 {\rm  Let us consider $I \subset A= k[X_1, X_2, X_3, X_4]$ the defining
ideal of the twisted cubic in $\Bbb P^3_k$. Recall from Example
\ref{YY4} that the Rees algebra of $I$ is Cohen-Macaulay, and so
$a_*^2(R_A(I))= -1$. Moreover, $I$ is an ideal of height $2$
generated by forms in degree $2$. Then, according to Corollary
\ref{DD2} we can get the Hilbert polynomials of $A/I^j$ for $j >3$
from the Hilbert polynomials of $I$, $I^2$ and $I^3$. By using
CoCoa, we have

\vspace{1mm}

\hspace{35mm} $P_{A/I}(s) \,= 3 s +1$

\vspace{2mm}

\hspace{35mm} $P_{A/I^2}(s) = 9 s -7$

\vspace{2mm}

\hspace{35mm} $P_{A/I^3}(s) = 18 s -34$

\vspace{2mm} By imposing $e_0 (0) = 0$, $e_0 (1)=3$ and $e_0 (2)
=9$, we get the multiplicity function $e_0(t) = \frac{3}{2} t
(t+1)$. Similarly, one gets $e_1(t) = \frac{5}{3} t
(t+1)(t-\frac{2}{5})$. Then the Hilbert polynomial of $A/I^j$ is
$$P_{A/I^j}(s) = e_0(j) {s+1 \choose 1} - e_1(j)$$
and the Hilbert polynomial of the Rees algebra $R$ of $I$ is
$$P_R(s,t) = {s+3 \choose 3} -e_0 (t) {s+1 \choose 1} + e_1(t).$$
In this case $\lambda_2 = 3$, $\lambda_3 = 10$, so by Proposition
\ref{DD8} we have $e_3(R)=1$, $e_2(R)=2$, $e_1(R)=1$, $e_0(R)=0$,
and $e_2(G)=0$, $e_1(G)=3$, $e_0(G)=2$. Therefore, the
multiplicity of the Rees algebra is $e(R)= \sum_{i=0}^3 e_i(R)= 4$
and the multiplicity of the form ring is
 $e(G)= \sum_{i=0}^2 e_i(G)= 5$. }
\end{ex}

\bigskip
\section{Hilbert series of the powers of an ideal}
\label{EE} \markboth{CHAPTER VI. ASYMPTOTIC BEHAVIOUR}{HILBERT
SERIES OF THE POWERS OF AN IDEAL}
\medskip
Let $A = k[{\x}]$ be the polynomial ring over a field $k$ in $n$
variables. For any finitely generated graded $A$-module $M$,
recall that the Hilbert series of $M$ is defined as
$$H_M(s) = \sum_{j \in \Bbb Z} H(M,j) s^j=
\sum_{j \in \Bbb Z} \dim_k M_j s^j \in \Bbb Z [[s]].$$

\medskip
Following A. Conca and G. Valla \cite{CV}, for a given class $\cal
C$ of homogeneous ideals in $A$, we say that $\cal C$ has rigid
powers if for any ideals $I,J$ in $\cal C$ such that $H_{A/I}(s) =
H_{A/J}(s)$ then $H_{A/I^j}(s) = H_{A/J^j}(s)$ for all $j$. For
example, the class of complete intersection ideals has rigid
powers. The class of the homogeneous ideals in $A$ which are
Cohen-Macaulay of codimension 2 and the class of the homogeneous
ideals in $A$ which are Gorenstein of codimension 3 do not have
rigid powers, but their subclasses consisting of the ideals of
linear type have this property as it has been proved in \cite{CV}.

\medskip
 Our first aim in this section is to show that for an equigenerated
ideal $I$ we can compute the Hilbert series of  $A/I^j$ for $j
\geq 1$ from a finite set among these Hilbert series, and so we
can also compute the bigraded Hilbert series of its Rees algebra.
This fact will be a direct consequence of the noetherian property
of the Rees algebra, and the finite set of Hilbert series will be
found thanks to the bounds for the shifts of the bigraded minimal
free resolution of the Rees algebra given by Theorem \ref{B77}. In
particular, we will have that if the Hilbert series of the powers
of two ideals $I,J$ coincide for certain exponents then all the
Hilbert series of the powers of $I$ and $J$ must coincide.

\begin{thm}
\label{EE2} Let $I$ be an equigenerated homogeneous ideal. Set $l
= l(I)$, $c =a_*^2(R_A(I))$. The Hilbert series of $I^j$ for $c +
1 \leq j \leq c + l$ determine the Hilbert series of $I^j$ for $j
> c + l$.
\end{thm}

\medskip
{\pf}
 Let us assume that $I$ is generated by forms in degree $d$.
Then we have that $R=R_A(I)$ is a finitely generated bigraded
$S$-module in a natural way, for $S=k[{\x}, Y_1, \dots, Y_l]$ the
polynomial ring with deg$(X_i)= (1, 0)$, deg$(Y_j)= (d, 1)$. Let
$H_R(s,t)$ be the bigraded Hilbert series of $R$, that is
  $$H_R(s,t) = \sum_{i,j} \dim_k  R_{(i,j)} s^i t^j =
   \sum_{i,j} \dim_k (I^j)_i \; s^i t^j =
    \sum_j H_{I^j}(s) t^j. $$

By considering the bigraded minimal free resolution of $R$ as
$S$-module
$$0 \to D_t \to \dots \to D_0 \to R \to 0,$$
$$D_p =\bigoplus_{(a,b) \in \Omega_p} S(a,b),$$
we can write
$$H_R(s,t) = \frac{Q(s,t)}{(1-s)^n (1-s^dt)^l} \;\;,$$
with $Q(s,t)= \sum_{p=0}^t (-1)^p \sum_{(a,b) \in \Omega_p}
s^{-a}t^{-b} \in \Bbb Z[s,t]$. Now let us fix $\alpha \in \Bbb
Z_{\geq 0}$. Denoting by $\beta_p^j = \dim_k \Tor^A_p(k,
I^j)_{\alpha + dj}$, then

\vspace{3mm}

\hspace{20mm} $\sum_p (-1)^p \beta_p^j =
 [(1-s)^n H_{I^j}(s) \;]_{\deg s=\alpha+dj} $

\vspace{3mm}

\hspace{40mm} $ =[ (1-s)^n H_R (s,t) \; ] _{\begin{array}{c}
\scriptsize{\dg s=\alpha+ dj} \\
\scriptsize{\dg t=j}
\end{array}} $

\vspace{3mm}

\hspace{40 mm} $ = \biggl [ \, \frac{Q(s,t)}{(1-s^dt)^l} \; \biggr
] \,_{\begin{array}{c}
\scriptsize {\dg s=\alpha+ dj} \\
\scriptsize  {\dg t=j}
\end{array}} $

\vspace{3mm}

 Let us write $Q(s,t) = \sum_k m_k s^{\alpha + dk}t^k + \overline Q
(s,t)$, with $\overline Q (s,t)$ containing all the monomials of
the type $s^{\beta + dk}t^k $ for any $\beta \not = \alpha$ and
any $k$. The pairs $(-\alpha- dk, -k)$ are shifts in the bigraded
minimal free resolution of $R$ as $S$-module, so $k \leq t_*^2(R)
= a_*^2(R) +l = k_0$ by Theorem \ref{B77}. Then we have that for
any $j \geq k_0$,

\vspace{4mm}

$\sum_p (-1)^p \beta_p^j = \bigg[ \bigg(\sum_{v=0}^j {v+l-1
\choose l-1} s^{dv} t^v \bigg) (\sum_{k=0}^{k_0} m_k s^{\alpha +
dk}t^k) \;  \, \bigg] _{\begin{array}{c}
\footnotesize {\dg s=\alpha+ dj} \\
\footnotesize {\dg t=j}
\end{array}} $

\vspace{3mm}

\hspace{19mm} $= m_0 {j+l-1 \choose l-1} + \dots + m_{k_0}
{j-k_0+l-1 \choose l-1}$

\vspace{3mm}

\hspace{19mm} $= P_{\alpha} (j).$

\vspace{3mm}

It is easy to prove that this equality holds for $j \geq k_0-l+1 =
a_*^2(R) +1$. So we have found a polynomial $P_{\alpha}(j)$ of
degree $\leq l-1$ such that $P_{\alpha} (j) = \sum_p (-1)^p
\beta_p^j $ for any $j \geq a_*^2(R) +1$. Hence the Hilbert series
of the powers $I^j$ for $a_*^2(R)+1 \leq j \leq a_*^2(R) +l$ will
determine the Hilbert series of $I^j$ for any $j > a_*^2(R) +l$.
$\B$

\medskip
\begin{cor}
\label{EE3}
 Let $I$ be an equigenerated homogeneous ideal whose Rees algebra is
Cohen-Macaulay, and $l = l(I)$. Then the Hilbert series of $I^j$
for $j \leq l -1$ determine the bigraded Hilbert series of the
Rees algebra of $I$.
\end{cor}

\medskip
     Recent papers by A. M. Bigatti, A. Capani, G. Niesi and L. Robbiano
\cite{BCNR} and L. Robbiano and G. Valla \cite{RV} treat the
problem of computing the Hilbert series of the powers of a
homogeneous ideal $I$ in the polynomial ring $A=k[{\x}]$. The
strategy followed there to solve this problem is to compute the
Rees algebra $R_A(I)$ of $I$ and then a Gr\"{o}bner basis of it, from
which one can get easily the bigraded Hilbert series of the Rees
algebra, and so the Hilbert series of all the powers of $I$.
Notice that we can use Theorem \ref{EE2} to give another approach
to this problem: To get the Hilbert series of the powers of an
equigenerated ideal $I$ it suffices to compute the Hilbert series
of $l(I)$ of its powers. Next we apply this procedure to the
following example studied by A. Bigatti et al. \cite[Example
5.4]{BCNR}.

\begin{ex}
\label{EE22} {\rm  Let us consider the ideal $I$ generated by the
$2$ by $2$ minors of the generic symmetric $3$ by $3$ matrix
$$M =    \left( \begin{array}{c}
         X_1 \; \; X_2 \; \; X_3 \\
         X_2 \; \; X_4 \; \; X_5 \\
         X_3 \; \; X_5 \; \; X_6
         \end{array} \right) .$$
The Rees algebra of $I$ is Cohen-Macaulay, so $a_*^2(R)=-1$.
Therefore, the Hilbert series of $I^j$ for $j \leq 5$ will
determine the rest of the Hilbert series. By using CoCoa, we
obtain:

\vspace{3mm}

\hspace{8mm} $H_{I}(s)= \frac{6s^2-8s^3+3s^4}{(1-s)^6} ,$

\vspace{3mm}

\hspace{8mm}
$H_{I^2}(s)=\frac{21s^4-45s^5+38s^6-18s^7+6s^8-s^9}{(1-s)^6}\; ,$

\vspace{3mm}

\hspace{8mm}
$H_{I^3}(s)=\frac{56s^6-150s^7+165s^8-100s^9+36s^{10}-6s^{11}}{(1-s)^6},$

\vspace{3mm}

\hspace{8mm}
$H_{I^4}(s)=\frac{126s^8-385s^9+486s^{10}-330s^{11}+125s^{12}-21s^{13}}
{(1-s)^6},$

\vspace{3mm}

\hspace{8mm}
$H_{I^5}(s)=\frac{252s^{10}-840s^{11}+1155s^{12}-840s^{13}+330s^{14}+
56s^{15}}{(1-s)^6}. $

\vspace{3mm}

Then the polynomials $P_\alpha(j)$ defined in the proof of Theorem
\ref{EE2} are

\vspace{3mm}

\hspace{8mm} $P_{\alpha} (j) = 0  \;{\rm for} \; \alpha \not = 0,
\dots,5 \;,$

\vspace{3mm}

\hspace{8mm} $P_0 (j) = 1+ \frac{137}{60}j+\frac{15}{8}
j^2+\frac{17}{24} j^3+ \frac{1}{8}j^4 +\frac{1}{20}j^5 ,$

\vspace{3mm}

\hspace{8mm} $P_1 (j) =  - \frac{7}{4}j -\frac{83}{24}
j^2-\frac{53}{24} j^3- \frac{13}{24}j^4 -\frac{1}{24}j^5 ,$

\vspace{3mm}

\hspace{8mm} $P_2 (j) = -\frac{3}{2}j+\frac{13}{12}
j^2+\frac{29}{12} j^3+ \frac{11}{12}j^4 +\frac{1}{12}j^5 ,$

\vspace{3mm}

\hspace{8mm} $P_3 (j) = \frac{7}{6}j+\frac{3}{4} j^2-\frac{13}{12}
j^3- \frac{3}{4}j^4 -\frac{1}{12}j^5 ,$

\vspace{3mm}

\hspace{8mm} $P_4 (j) = -\frac{1}{4}j-\frac{7}{24}
j^2+\frac{5}{24} j^3+ \frac{7}{24}j^4 +\frac{1}{24}j^5 ,$

\vspace{3mm}

\hspace{8mm} $P_5 (j) = \frac{1}{20}j+\frac{1}{24}
j^2-\frac{1}{24} j^3- \frac{1}{24}j^4 -\frac{1}{120}j^5 ,$

\vspace{3mm}

\noindent and the Hilbert series of $I^j$ is
$$ \frac{P_0(j) s^{2j} + P_1(j) s^{2j+1} + P_2(j) s^{2j+2}
+ P_3(j) s^{2j+3} + P_4(j) s^{2j+4} + P_5(j) s^{2j+5}}
{(1-s)^6}\;.$$

}
\end{ex}

\medskip
  Next we also compute the Hilbert series of the powers of the ideal of
the twisted cubic in $\Bbb P^3_k$ studied in the previous section.

\begin{ex}
\label{EE6}
 {\rm  Let $I \subset A= k[X_1, X_2, X_3, X_4]$ be the defining ideal of the
twisted cubic in $\Bbb P^3_k$. From Example \ref{YY4}, let us
recall that $I$ is generated by quadrics with $l(I) = \mu(I)=3$
and $R_A(I)$ is Cohen-Macaulay. Therefore, according to Theorem
\ref{EE2}, we can get the Hilbert series of $I^j$ for $j >2$ from
the Hilbert series of $I$ and $I^2$. By using CoCoa, we have

\vspace{3mm}

\hspace{35mm} $H_I(s) = \frac{3 s^2-2s^3}{(1-s)^4} \;,$

\vspace{3mm}

\hspace{35mm} $H_{I^2}(s) = \frac{6 s^4- 6 s^5 + s^6}{(1-s)^4}.$

\vspace{3mm} Then the polynomials $P_\alpha (j)$ defined in the
proof of Theorem \ref{EE2} are

\vspace{3mm}

\hspace{35mm} $P_{\alpha} (j) = 0,   \;{\rm for} \; \alpha \not =
0,1,2 \;,$

\vspace{3mm}

\hspace{35mm} $P_0 (j) = \frac{1}{2} (j+1)(j+2) \;,$

\vspace{3mm}

\hspace{35mm} $P_1 (j) = - j(j+1) \;,$

\vspace{3mm}

\hspace{35mm} $P_2 (j) = \frac{1}{2} j(j-1) \;,$

\vspace{3mm} \noindent and the Hilbert series of $I^j$ is then
$$H_{I^j}(s) = \frac{P_0(j) s^{2j} + P_1(j) s^{2j+1} + P_2(j) s^{2j+2}}
{(1-s)^4}\;.$$ Now, the bigraded Hilbert series of its Rees
algebra is
$$H_R(s,t) = \sum_j H_{I^j}(s) t^j =
 \frac{1-2s^3 t + s^6 t^2}{(1-s)^4 (1-s^2t)^3}.$$
 }
\end{ex}

\medskip
\begin{rem}
\label{EE4} {\rm Similarly we can prove the following statement
for the Hilbert series of the form ring of an equigenerated ideal
$I$: If $l = l(I)$ and $e =a_*^2(G_A(I))$, then the Hilbert series
of $I^j/I^{j+1}$ for $e + 1 \leq j \leq e + l$ determine the
Hilbert series of $I^j/I^{j+1}$ for $j > e + l$. In fact, for any
$a \geq e$ the Hilbert series of $I^j/I^{j+1}$ for $a + 1 \leq j
\leq a + l$ determine the rest.   }
\end{rem}

\medskip
For any $m \geq 0$, let us define $C_m$ to be the class of
equigenerated homogeneous ideals $I$ in $A$ such that $
a_*^2(G_A(I))+l(I) \leq m$. Note that $C_0$ contains the class of
complete intersection ideals, and we have the chain $$C_0 \subset
C_1 \subset \cdots \subset C_m \subset C_{m+1} \subset \cdots$$ As
a corollary, we get

\begin{cor}
\label{EE5}
        Let $I, J \in C_m$ be such that
$$H_{I^j/I^{j+1}}(s) = H_{J^j/J^{j+1}}(s), \;
{\rm for} \;m-l+1 \leq j \leq m.$$
 Then $H_{I^j/I^{j+1}}(s)= H_{J^j/J^{j+1}}(s)$, $\forall j$. Therefore
$H_{I^j}(s)= H_{J^j}(s)$, for all $j$, and in particular $C_0$ has
rigid powers.

\end{cor}

\medskip
    For an arbitrary homogeneous ideal $I$ in $A$, we can also show
that a finite set of Hilbert series of the powers of $I$ determine
the rest. But in this case, the bound we get is worse.

\begin{prop}
\label{EE45}
 Let $I$ be a homogeneous ideal in $A$. Set $r = \mu (I)$,
$c =a_*^2(R_A(I))$. The Hilbert series of $I^j$ for $ j \leq c +r$
determine the Hilbert series of $I^j$ for $j > c + r$.
\end{prop}

 {\pf}
 Assume that $I$ is minimally generated by forms $f_1$,...,$f_r$
of degrees $d_1,\ldots, d_r$ respectively.
 Then let us consider the presentation of the
Rees algebra $R$ of $I$ as a quotient of the polynomial ring
$S=k[{\x},{\y}]$, with deg$(X_i)= (1, 0)$, deg$(Y_j)= (d_j, 1)$.
From the bigraded minimal free resolution of $R$ as $S$-module, we
have that there is a polynomial $Q(s,t) \in \Bbb Z[s,t]$ such that

$$H_R(s,t) = \frac{Q(s,t)}{(1-s)^n (1-s^{d_1}t)\dots (1-s^{d_r}t)}
\;\;.$$

\vspace{2mm}

According to Theorem \ref{B77} we can write $Q(s,t) =
\sum_{i=0}^{m} Q_i(s) t^i$, with $m= a_*^2(R)+r$. Since $H_R(s,t)
= \sum_{j \geq 0} H_{I^j}(s) t^j$, we have
$$Q(s,t)= (1-s)^n (1-s^{d_1}t) \dots (1-s^{d_r}t) (\sum_{j
\geq 0} H_{I^j}(s) t^j),$$ and then the result follows
immediately. $\B$

\bigskip
\section{Minimal graded free resolutions of the powers of
an ideal} \label{GG} \markboth{CHAPTER VI. ASYMPTOTIC
BEHAVIOUR}{MINIMAL RESOLUTIONS OF THE POWERS OF AN IDEAL}
\medskip

The results about the Hilbert series of the powers of a
homogeneous ideal imply in some particular cases the estability
(in a meaning which we will precise immediately) of the minimal
graded free resolutions of the powers of the ideal. For instance,

\begin{prop}
\label{GG50} Let $I$ be an ideal generated in degree $d$ with
$l(I)=2$ whose Rees algebra is Cohen-Macaulay. Then the minimal
graded free resolution of $I$ determines the minimal graded free
resolutions of all its powers. Namely, if the minimal graded free
resolution of $I$ is
$$ 0 \to A(-\alpha_1 -d)^{\beta_1}
 \oplus  \dots \oplus A(-\alpha_m -d)^{\beta_m} \to
A(-d)^{\beta} \to I \to 0,$$ then for any $j \geq 1$ the minimal
graded free resolution of $I^j$ is
$$ 0 \to A(-\alpha_1 -dj)^{\beta_1 j}
 \oplus  \dots \oplus A(-\alpha_m -dj)^{\beta_m j} \to
A(-dj)^{(\beta-1)j +1} \to I^j \to 0.$$
\end{prop}

{\pf} First, note that for any $j \geq 1$ we have that $\dpp_A I^j
\leq l(I)-2=1$ because $R$ is Cohen-Macaulay (see Proposition
\ref{GG1}). Therefore, we can conclude that $\dpp_AI^j = 1$ for
any $j$. On the other hand, since $R$ is Cohen-Macaulay, we have
$a^2_*(R)=-1$, so the Hilbert series of $I^j$ for $j \leq l-1 = 1$
determine the Hilbert series of $I^j$ for $j > 1$ according to
Theorem \ref{EE2}. The polynomials $P_\alpha(j)$ defined there are

\vspace{2mm}

\hspace{20mm} $P_\alpha(j)= 0, \, {\rm for} \; \alpha \not \in \{
0, \alpha_1, \dots, \alpha_m\},  $

\vspace{2mm}

\hspace{20mm} $P_0(j)= (\beta-1)j+1,$

\vspace{2mm}

\hspace{20mm} $P_{\alpha_i}(j)= -\beta_i \, j, \; {\rm for} \; i
\in \{0, \dots, m \}.$

\vspace{4mm} Then the Hilbert series of $I^j$ is \vspace{2mm}

\hspace{20mm} $H_{I^j}(s) = \frac{\sum_{\alpha} P_{\alpha}(j)
s^{\alpha+dj}}{(1-s)^n},$

\vspace{2mm} \noindent so the minimal graded free resolution of
$I^j$ must be
$$ 0 \to A(-\alpha_1 -dj)^{-P_{\alpha_1}(j)}
 \oplus  \dots \oplus A(-\alpha_m -dj)^{-P_{\alpha_m}(j)} \to
A(-dj)^{P_0 (j)} \to I^j \to 0.$$ $\B$

\medskip
This result leads to the question of when a finite number of
minimal graded free resolutions of the powers of $I$ determine the
rest (and, in this case, which set of resolutions determine the
others).

\medskip
Let us begin by studying the behaviour of the projective dimension
of the powers of $I$. It is well-known that these projective
dimensions are asymptotically constant (see \cite[Theorem
2]{Bro}), but not for which powers of the ideal the projective
dimension takes the asymptotic value. We will precise these powers
for ideals whose form ring is Gorenstein by considering the Koszul
homology of the Rees algebra $R$ of $I$ with respect to $\x$. This
also provides new proofs of well-known results as the Burch's
inequality or the constant asymptotic value for the depth.

\begin{prop}
\label{GG1} Let $I$ be a homogeneous ideal in $A$, and set
$l=l(I)$. Then:
\begin{enumerate}
\item ${\dpp}_A (I^j) \leq n - \depth_{(\frak m R)} (R)$ for all
$j$, and the equality holds for $j \gg 0$. So, $\inf_{j \geq 0} \{
\depth (A/I^j) \}= n - l-
 (\h (\frak m R)-\depth_{(\frak m R)}(R)).$
\item If $R$ is Cohen-Macaulay,
${\dpp}_A (I^j) \leq l-1$ for any $j $ and ${\dpp}_A (I^j) = l-1$
for $j \gg 0$. Furthermore,
 ${\dpp}_A (I^j) = l-1$ implies ${\dpp}_A (I^{j+1}) = l-1$.
\item[(iii)] If $G$ is Gorenstein,
${\dpp}_A I^j = l - 1$ if and only if $j > a^2(G) - a(F)$, and
${\dpp}_A I^j/I^{j+1} = l $ if and only if $j \geq a^2(G) - a(F)$.
\end{enumerate}

\end{prop}

\medskip
{\pf} Let us consider the Koszul complex $K.(\XX, R) = K.(\x, R)$
of the Rees algebra $R$ with respect to ${\underline X} =\x$. We
have the natural bigrading in the Rees algebra $R$ by means of
$R_{(i,j)} = (I^j)_i$, and then the modules $K_p(\XX, R)$ of the
Koszul complex are also bigraded in a natural way. Denoting by
$F=F_{\fm}(I)$,  we have that for any $p$ the Koszul homology
module $H_p = H_p^S (\XX; R)$ is a finitely generated bigraded
$F$-module. Moreover, since $K_p (\XX, R)_{(i,j)} = K_p(\XX,
I^j)_i$ we have
$$H_p^S({\underline X}, R)_{(i,j)} = H_p^A({\underline X}, I^j)_i =
\Tor^A_p (k, I^j)_i \,\,,$$ so the Koszul homology modules $H_p$
contain all the information about the graded minimal free
resolutions of the powers of $I$.

 Now set $s = n -$depth$_{(\frak m R)}(R)$. Recall that
$H_p$ is zero for any $p>s$, and so ${\dpp}_A (I^j) \leq s$ for
any $j$. Moreover, since $H_s$ is a $F$-module of dimension $l$
\cite[Remark 1.5]{Hu1} we can find for any $j \gg 0$ an integer
$i$ (depending on $j$) such that $[H_s]_{(i,j)} \not =0$.
Therefore we obtain that ${\dpp}_A (I^j) = s$ for $j \gg 0$. By
the graded Auslander-Buchsbaum formula, ${\dpp}_A(I^j) = n -
\depth \,A/I^j-1$, and noting that $\depth_{(\frak m R)}(R) \leq$
ht$(\frak m R) = n+1-l$ we get $(i)$.

 To prove $(ii)$, let us denote by $t =\depth_{(\frak m R)} (R) =
n+1-l$. We may assume that $k$ is infinite, and then there exists
a homogeneous regular sequence $b_1, \dots, b_t \in \frak m R$ of
degree $(1, 0)$. Then
$$H_s = \frac {(b_1, \dots, b_t): ({\rm X}_1, \dots, {\rm X}_n)}{(b_1,
\dots, b_t)} (t-n,0).$$ Note that $s = l-1$ because $R$ is CM.
Now, observe that ${\dpp}_A(I^j)= l-1$ implies that there exists
$i$ such that $[H_s]_{(i,j)} \not = 0$; so let us take $f \in
[H_s]_{(i,j)}, f \not = 0$. For a positively graded ring $A=
\bigoplus_{j \geq 0} A_j$, let us denote by $A_+= \bigoplus_{j >
0} A_j$, and in the following let us consider the fiber cone $F$
and the Rees algebra $R$ graded by means of $F_j = I^j/ \fm I^j $,
$R_j =I^j$. If $F_+ f = 0$, then $I^m f \subset (b_1, \dots, b_t)$
for any $m \geq 1$. So, denoting by $\overline R = R/(b_1, \dots,
b_t)$, we have that $\overline R_+ \subset {\rm Ann} (f) \subset
\frak p \in \Ass(\overline R)$ and so ht$(\overline R_+) = 0$. But
ht$(\overline R_+) = \dim \overline R - \dim  \overline R /
\overline R_+ = 1$. As a consequence, $(F_+)_1 f \not = 0$ and so
there exists $d$ such that $[H_s]_{(i+d,j+1)} \not = 0$ and we
have $(ii)$.

  Finally, we are going to determine the powers of $I$ whose
projective dimension is $l-1$ if $G$ is Gorenstein. To this end,
let us consider the Koszul homology modules of the form ring $G$
with respect to $\XX$, which will be also denoted by $H_p$. Set $s
= n -\depth_{(\frak m G)}(G) = l$, $t= \depth_{(\frak m G)} (G)$.
As before, $H_p$ is zero for $p>s$ and there exists a homogeneous
regular sequence $b_1, \dots, b_t \in \frak m G$ of degree $(1,
0)$, such that
$$H_l = \frac {(b_1, \dots, b_t): ({\rm X}_1, \dots, {\rm X}_n)}{(b_1,
\dots, b_t)} (t-n,0).$$ On the other hand, from the natural
bigraded epimorphism $ G \rightarrow G/ \fm G=F$, we can compute
the canonical module of the fiber cone $F$ by using Corollary
\ref{B3} :
 $$K_F = \uExt^{n-l}_G (F, K_G).$$
Since $G$ is Gorenstein, we have a bigraded isomorphism $K_G \cong
G(-n,a)$ with $a=a^2(G)$ by Corollary \ref{K4}. Therefore,

\vspace{3mm}

\hspace{15mm} $K_F=  \uExt^{n-l}_G (F, G) (-n, a) $

\vspace{3mm}

\hspace{22mm} $= \uHom_G (F, G/(b_1, \dots, b_t)) (-n+t, a)  $

\vspace{3mm}

\hspace{22mm} $=\frac {(b_1, \dots, b_t): (X_1, \dots, X_n)}{(b_1,
\dots, b_t)} (-n+t, a) $

\vspace{3mm}

\hspace{22mm} $= H_l (0, a).$

\vspace{3mm}

\noindent Now, observe that ${\dpp}_A(I^j / I^{j+1}) =l $ if and
only if there exists $i$ such that $[H_l]_{(i,j)} \not = 0$ if and
only if there exists $i$ such that $[K_F]_{(i, j-a)} \not = 0$,
that is, $j \geq a -a(F)$. From the exact sequences
$$0 \to I^{j+1} \to I^j \to I^j/I^{j+1} \to 0\,\,,$$
\noindent it is then easy to check that ${\dpp}_A(I^j) = l-1$ if
and only if $j > a - a(F)$, and so we are done.$\B$

\medskip
\begin{ex}
\label{GG11} {\rm Let $I$ be a strongly Cohen-Macaulay ideal such
that $\mu (I_{\frak p}) \leq \h (\frak p)$ for any prime ideal
$\frak p \supseteq I$.
 Set $l = l (I)$, $h =$ht$(I)$. Recall from Corollary \ref{Bb8}
that $G_A(I)$ is Gorenstein,  $a^2(G_A(I)) = -h$ and $a(F_{\frak
m}(I)) = -l$. So, by Proposition \ref{GG1} we have depth$(A/I^j) =
n-l$ if and only if $j > l-h$.}
\end{ex}

\begin{ex}
\label{GG12} {\rm  Let ${\bf X}= (X_{ij})$ be a generic matrix,
with $1 \leq i \leq d$ , $1 \leq j \leq n$ and $d \leq n$. Let $I
\subset A= k[{\bf X}]$ be
 the ideal generated by the maximal minors of ${\bf X}$.
Recall from Example \ref{Bb81} that the Rees algebra $R$ is
Cohen-Macaulay and the form ring $G$ is Gorenstein with $a^2(G) =
- \h (I) = - (n-d+1)$. Furthermore, $l(I) = d(n-d) + 1 $ and
$a(F)=-n$. Now by Proposition \ref{GG1}, we get that depth$(A/I^j)
= d^2 -1$ if and only if $j > d-1$. In the case $n = d+1$, this
was proved in \cite[Example 9.27]{BV}.}

\end{ex}

\begin{ex}
\label{GG88} {\rm  Let ${\bf X}= (X_{ij})$ be a generic
skew-symmetric matrix, with $1 \leq i < j \leq n$, and $n$ odd.
Let $I \subset A= k[{\bf X}]$ be the ideal generated by the
$(n-1)$-pfaffians of ${\bf X}$, where $k$ is a field. In this
case, the form ring $G$ is Gorenstein \cite{CD} and $l(I) = n$,
$a(F_\fm(I))=-n$ \cite{Hu3}. So $\depth(A/I^j)$ takes the
asymptotic value $\frac{n(n-1)}{2}-n$ for some $j \leq n$, and by
Proposition \ref{GG1} for all $j \geq n$.

    G. Boffi and R. S\'{a}nchez  \cite{BS2} have constructed a family
of complexes which give a resolution for all the powers $I^j$, for
$j \geq 1$, in particular proving that $\dpp_A(A/I^j)=n $ if and
only if $j \geq n-2$. Then Proposition \ref{GG1} shows that
$a^2(G_A(I))= -3$.}
\end{ex}

\medskip
Our next aim is to study the graded minimal free resolutions of
the powers of an equigenerated ideal by doing a deeper study of
the Koszul homology of the Rees algebra with respect to $\x$. The
general case will be studied later by different methods.

\bigskip
\subsection{Case study : Equigenerated ideals}

\medskip
First of all, we show that the shifts in the graded minimal free
resolutions of the powers of an equigenerated ideal are given by
linear functions asymptotically and the graded Betti numbers of
these resolutions are given by polynomials asymptotically.

\begin{prop}
\label{GG2}
   Let $I$ be an ideal generated by forms in degree $d$. Set $l =
l(I)$, $s = n - \depth_{(\frak m R)} (R)$. Then there is a finite
set of integers  $$\{\alpha_{pi} \mid 0 \leq p \leq s, 1 \leq i
\leq k_p \}$$ and polynomials of degree $\leq l-1$
$$\{Q_{\alpha_{pi}}(j)  \mid 0 \leq p \leq s, 1 \leq i \leq k_p \}$$
such that the graded minimal free resolution of $I^j$ for $j$
large enough is
$$ 0 \to D_s^j \to \dots \to D_0^j \to I^j \to 0 \,\, ,$$
\noindent with $D_p^j = \bigoplus_i A(-\alpha_{pi}
-dj)^{\beta_{pi}^j}$ and ${\beta_{pi}^j}= Q_{\alpha_{pi}} (j)$.

\end{prop}

{\pf} Let us consider again the Koszul homology of the Rees
algebra $R$ of $I$ with respect to $\XX = X_1,\dots, X_n$, and let
us denote by $F= \bigoplus_{j \geq 0} I^j/ \fm I^j$ the fiber cone
of $I$ and by $F_+ = \bigoplus_{j > 0} I^j/ \fm I^j$. For every $p
\leq s$, $H_p = H_p^S (\XX; R)$ is a finitely generated bigraded
$F$-module. Let $g$ be a homogeneous generator of $H_p$ with
deg$(g) = (a, b)$, and set $\alpha = a -d b$. If $F_+ \subset \rad
(\Ann(g))$, there exists $j$ such that $F_+^j g= 0$, and so $F_j g
= 0$ for all $j \gg 0$. Otherwise, there exists a homogeneous
element $f \in F$ of degree $d$ such that $f \not \in \rad (
\Ann(g))$. Then $f^j g \not = 0$ for all $j$, and so we have
$[H_p]_{(\alpha +dj, j)} \not =0$ for all $j \gg 0$. Let
$g_1,\dots,g_m$ be the homogeneous generators of $H_p$ with this
property, and set deg$(g_i)=(a_i, b_i)$, $\alpha_i = a_i -db_i$.
Then, for $j$ large enough we have that $[H_p]_{(a,j)} \not =0$
if, and only if, there exists $i \in \{1, \dots, m \}$ such that
$a= \alpha_i+dj$. Since $[H_p]_{(\alpha_i+dj,j)} = \Tor^A_p (k,
I^j)_{\alpha_i +dj}$, we obtain that $\alpha_i +dj$, for $1 \leq i
\leq m$, are the only shifts in the place $p$ of the graded
minimal free resolution of $I^j$ for $j \gg 0$.

For $\alpha \in \{\alpha_1, \dots, \alpha_m \}$, let us define
$H_p^ \alpha = \bigoplus_j [H_p]_{(\alpha + dj, j)}$. Notice that
$\dim H_p^\alpha \leq \dim H_p = l$ by \cite[Remark 1.5]{Hu1}.
Since $H_p^\alpha$ is a finitely generated graded $F$-module,
there exists a polynomial $Q_\alpha(j)$ of degree $\dim H_p^\alpha
-1 \leq l -1$ such that for $j$ large enough
$$Q_\alpha (j) = \dim_k[H_p^\alpha]_j =
\dim_k \Tor^A_p (k, I^j)_{\alpha +dj},$$ so $Q_\alpha(j)$ is the
Betti number of $I^j$ corresponding to $\alpha + dj$ in the place
$p$. $\B$

\begin{ex}
\label{GG22} {\rm Let $I$ be a Cohen-Macaulay homogeneous ideal of
codimension two in the polynomial ring $A= k[{\x}]$ such that:

(i) The entries of the Hilbert-Burch matrix of $I$ are linear
forms.

(ii) $I$ verifies $G_n$.

(iii) $\mu(I) \leq n$.

This example has been studied by A. Conca and G. Valla in
\cite{CV}. Set $r =\mu (I)$, $d= r-1$,
 $S = A[{\y}]$. Then
 $$R_A(I) \cong {\rm Sym}_A (I) \cong
S/(F_1,\dots,F_{r-1}),$$ with $F_1,\dots,F_{r-1}$ a regular
sequence of degree $(d, 1)$. So the Koszul complex of $S$ with
respect to $F_1,\dots,F_{r-1}$ gives the bigraded minimal free
resolution of $R_A(I)$. From this resolution one can get the
minimal graded free resolutions of $I^j$, for all $j \geq 0$.
Namely, $\dpp _A (I^j) = \min \{j, r-1\}$ and the minimal free
resolution of $I^j$ is
$$ 0 \to D_{r-1}^j \to \dots \to D_0^j \to I^j \to 0 \,\, ,$$

\noindent with $D_p^j = A (-p-dj)^{\beta_p^j}$,
 $\beta_p^j = {r-1 \choose p}{ r+j-p-1 \choose r-1}$.
 }

\end{ex}

\begin{rem}
\label{GG21} {\rm Let $(A, \frak m, k)$ be a noetherian local
ring, and let $I \subset A$ be an ideal. Set $l= l(I)$, $r=
\mu(I)$. Denote by $R$ the Rees algebra of $I$ graded by
$R_j=I^j$, and let $S=A[{\y}]$ be a polynomial ring over $A$ with
$\deg Y_j=1$ so that $R$ is a finitely generated graded
$S$-module. As above, we can prove that there are polynomials
$Q_p(j)$ of degree $\leq l-1$, the Hilbert polynomials of
$\Tor^S_p(S/\fm S, R)$, such that the minimal free resolutions of
$I^j$ for $j \gg 0$ are
$$\dots \to A^{Q_p(j)} \to \dots \to A^{Q_0(j)} \to I^j \to 0 .$$
\noindent This result was proved by V. Kodiyalam in \cite{Ko1}. If
$A$ is regular, let $\underline x$ be a regular sequence
generating $\fm$. Since $\Tor^S_p(S/\fm S,R) \cong H_p (\underline
x, R) $, the module $\Tor^A_ p(k, R)$ has dimension $l$ if it is
not zero. Therefore, the polynomial $Q_p(j)$ has degree $l-1$ if
$Q_p (j) \not =0$. This answers affirmatively \cite[Question
13]{Ko1} for any regular ring $A$.}
\end{rem}

\medskip
Observe that Proposition \ref{GG2} says that we can compute the
graded minimal free resolution of any power of an equigenerated
ideal from a finite set among these resolutions. Now we consider
the problem of determining this finite set of resolutions. To
begin with, let us study the asymptotic shifts of Proposition
\ref{GG2}.

\begin{lem}
\label{GG60} Let $I$ be a homogeneous ideal generated in degree
$d$ and let $R=R_A(I)$. Then
\begin{enumerate}
\item  For all $p$ and $i$, there exists $(a,b) \in \Omega_{p,R}$ such
that $\alpha_{pi} =db-a$.
\item For each $\alpha$, let
$$p =\min \, \{\, q \mid \exists \, b \, s.t. \, (\alpha + db, b) \in
\Omega_{q,R} \},$$ and let
$$b_0= \max \{ \, b \mid (\alpha + db, b) \in \Omega_{p,R} \}.$$
Then $\alpha + d b_0 \in \Omega_{p, I^{-b_0}}$, that is, $\alpha +
d b_0$ is a shift that appears in the graded minimal free
resolution of $I^{-b_0}$ at the place $p$.
\end{enumerate}
\end{lem}

{\pf} Let $0 \to D_m \to \dots \to D_0 \to R \to 0$ be the
bigraded minimal free resolution of $R$ over $S=k[{\x}, Y_1,
\dots, Y_l]$. By applying the functor $(\,\,)^j$ to this
resolution, we get a graded free resolution of $I^j$ over $A$
$$ 0 \to D_m^j \to \dots \to D_1^j \to D_0^j \to I^j   \to 0 \, ,$$
with $D_p^j = \bigoplus_{(a,b) \in \Omega_{p,R}}
A(a-db-dj)^{\rho_{ab}^j}$, for some ${\rho_{ab}^j} \in \Bbb Z$.
This resolution is the direct sum of the minimal graded free
resolution of $I^j$ and the trivial complex \cite[Exercise
20.1]{E}, so we obtain that for $j \gg 0$
$$\{ \alpha_{pi} +dj \}_{p,i} \subset
\{ db-a+dj \mid (a,b) \in \Omega_{p,R} \},$$ and so $(i)$ is
already shown.

Now let $\alpha$ be such that there exists $b$ with $(\alpha + db,
b) \in \Omega_R$. Let $p$ be the first integer such that $(\alpha
+ db, b) \in \Omega_{p,R}$, and let $b_0$ be the maximum of these
$b$'s. We must show that
 $$ \Tor^S_p (S / \fm S, R)_{(-\alpha -d b_0, -b_0)} =
      \Tor^A_p (k, I^{-b_0})_{-\alpha-d b_0} \not =0.$$
We will proceed as in Theorem \ref{BB1}: Let
 $$ D_{p+1} \stackrel{\psi_{p+1}}{\longrightarrow} D_p
\stackrel{\psi_{p}}{\longrightarrow} D_{p-1}$$ be the differential
maps appearing in the resolution of $R$. Tensorazing by $\otimes_S
S/ \frak m S$, we have the sequence
 $$ D_{p+1}/ \fm D_{p+1} \stackrel{\overline
\psi_{p+1}}{\longrightarrow}
  D_p/ \fm D_p \stackrel{\overline \psi_{p}}{\longrightarrow}
   D_{p-1} / \fm D_{p-1} .$$
Now let $v \in D_p$ be an element of the homogeneous basis of
$D_p$ as free $S$-module with $\deg (v) = (-\alpha-d b_0, -b_0)$.
If $w_1, \dots, w_s$ is the basis of $D_{p-1}$, we can write
$$ \psi_p (v) = \sum_{j=1}^{s} \lambda_j w_j ,$$
with $\lambda_j \in {\cal M}$ homogeneous. Set $\deg (w_j) =
(-\alpha_j - d b_j, -b_j)$. By looking at the degree of the
elements, we have that $\lambda_j$ must be zero for all $j$ such
that $-b_j > -b_0$. For the integers $j$ such that $-b_j = -b_0$,
we have that $\lambda_j \in \fm S$ necessarily. Finally, for $j$
such that $-b_j < -b_0$ we also have $\lambda_j \in \fm S$ because
$\alpha_j \not = \alpha$. We may conclude $ \overline \psi_p (v) =
0$, that is,  $v \in {\rm Ker} \, \overline \psi_p$. It is clear
that $v \not \in {\rm Im} \,\overline \psi_{p+1}$ because Im
$\overline \psi_{p+1} \subset {\cal M} D_p$. So $v \in \Tor^S_p (S
/ \fm S, L)_{(-\alpha - d b_0, -b_0)}$, $v \not =0$ and we are
done. $\Box$

\medskip
 As a consequence of this lemma we have that all the differences $a-db$
for $(a, b) \in \Omega_R$ appear in the minimal graded free
resolution of some power $I^j$ of $I$ for $j \leq a_*^2(R)+ l(I)$.
The problem is to distinguish which of these shifts will appear
asymptotically, and the place from where on the resolutions are
stable. We can solve this problem for ideals with a very
particular nice behaviour. For instance, we get

\begin{prop}
\label{GG6} Let $I$ be an equigenerated homogeneous ideal, and set
$b = a_*^2(R_A(I)) + l(I)$. If the graded minimal free resolutions
of $I, I^2, \dots, I^{b}$ are linear, then the graded minimal free
resolutions of $I^j$ are also linear for any $j$. Furthermore, the
minimal free resolutions of $I, I^2, \dots, I^b$ determine the
minimal graded free resolutions of $I^j$ for any $j$.
\end{prop}

{\pf} Assume that $I$ is generated by forms in degree $d$ and set
$s= \sup_{j=1, \dots, b} \{\dpp_A I^j \}$. According to Lemma
\ref{GG60}, we have that the shifts in $\Omega_R$ are of the type
$(p+db,b)$ with $0 \leq p \leq s$. Furthermore, there exists $b_0$
such that $(p+db_0,b_0) \in \Omega_{p,R}$, but for any $b$, $q<p$,
$(p+db,b) \not \in \Omega_{q,R}$. Therefore, $\Omega_{p,R}$ has
only shifts of the form $(a+db,b)$ for $0 \leq a \leq p$. Again by
Lemma \ref{GG60}, we get
$$\{\alpha_{pi}\}_i \subset \{ 0, \dots, p \}.$$
Finally, since $\min \, \{ \,-\beta : \beta \in \Omega_{p+1, I^j}
\} > \min \, \{ \,-\beta : \beta \in \Omega_{p, I^j} \},$ we have
that $\min \, \{ \,-\beta : \beta \in \Omega_{p, I^j} \} \geq
p+dj$. Therefore $I^j$ must have a linear minimal free resolution.

Moreover, by Theorem \ref{EE2} we also have that for $p=0, \dots,
s$ there exists a polynomial $Q_p(j)$ of degree $\leq l-1$ such
that
$$ Q_p (j)= \dim_k \Tor^A_p (k, I^j)_{p+dj},$$
for $j \geq a_*^2(R) +1$. So, if we know the minimal graded free
resolutions of $I^{b-l+1}, \dots, I^{b}$, we may determine the
polynomials $Q_p (j)$, and then the minimal graded free resolution
of $I^j$ for $j>b$.$\B$

\begin{rem}
\label{GG66} {\rm  The first part of Proposition \ref{GG6} can be
also obtained from Theorem \ref{YY2} $(ii)$.}
\end{rem}

\begin{rem}
\label{GG7} {\rm  Given an equigenerated homogeneous ideal $I$
with a linear minimal free resolution, it can happen that $I^2$
has a non linear minimal free resolution (see \cite[Remark
3]{C1}). We have shown in Proposition \ref{GG6} that if certain
powers of $I$ have linear resolution, then the rest of the powers
have this property too.}
\end{rem}

\medskip
We may apply this result to guess the minimal graded free
resolutions of the powers of the ideal defining the twisted cubic
in $\Bbb P^3_k$.

\begin{ex}
\label{GG99}
 {\rm  Let $I \subset A=k[X_1, X_2, X_3, X_4]$ be the defining ideal of
the twisted cubic in $\Bbb P^3_k$, and let us study the graded
minimal free resolutions of its powers. $I$ is generated by forms
in degree $2$ with $l(I)=3$ and $b =a_*^2(R)+l(I)= 2$. The minimal
resolutions of $I$ and $I^2$ (computed with CoCoa) are:
$$ 0 \to A(-3)^2 \to A(-2)^3 \to I \to 0 \,\,,$$
$$ 0 \to A(-6) \to A(-5)^6 \to A(-4)^6 \to I ^2 \to 0 \,\,.$$
Since these resolutions are linear, we have that the minimal
graded free resolutions of $I^j$ for $j>2$ are also linear by
Proposition \ref{GG6}, and we may compute them:
$$ 0 \to A(-2-2j)^{Q_2(j)} \to A(-1-2j)^{Q_1(j)} \to A(-2j)^{Q_0(j)}
\to I^j \to 0 \,\,,$$ with $Q_0(j) = \frac{1}{2}(j+1)(j+2)$,
$Q_1(j) = j (j+1)$ and $Q_2(j) = \frac{1}{2} j (j-1)$.}
\end{ex}

\medskip
Similarly, one can prove the following statement.

\medskip
\begin{prop}
\label{C66} Let $I$ be an ideal generated in degree $d$, and set
$b = a_*^2(R_A(I)) + l(I)$. Assume that there are integers
$\alpha_1, \dots, \alpha_s$ such that the graded minimal free
resolutions of $I, I^2, \dots, I^{b}$ take the form
$$0 \to D_s^j \to \dots \to D_1^j \to D_0^j \to I^j \to 0,$$
with $D_p^j = A(-\alpha_p-dj)^{\beta_p^j}$ and $\beta_p^j \geq 0$.
Then the graded minimal free resolutions of $I^j$ are of this type
too. Furthermore, the minimal graded free resolutions of $I, I^2,
\dots, I^b$ determine the minimal graded free resolutions of
$I^j$ for any $j$.
\end{prop}

\medskip
The following example does not belong to the family of ideals
considered in the previous propositions, but we can also guess the
asymptotic resolution of its powers.

\begin{ex}
\label{GG9} {\rm   Let $I = (X^7, Y^7, X^6 Y + X^2 Y^5) \subset A
= k[X,Y]$. Note that $I$ is a $\fm$-primary ideal generated by
forms of degree 7 with $l(I)=2$. Since $\dpp_A I^j =1$ for any $j
\geq 1$, we have that the shifts in the place $0$ and $1$ of the
resolution of $I^j$ can not coincide. Then, according to Theorem
\ref{EE2} we have that for any $\alpha \not =0$ there is a
polynomial $P_\alpha(j)$ of degree $\leq 1$ such that
$$P_\alpha(j) = \dim_k \Tor^A_1 (k,I^j)_{\alpha+dj},$$
for all $j \geq a_*^2(R)+1$.

This example was studied by S. Huckaba and T. Marley \cite[Example
3.13]{HM}. Denoting by $G_{++}= \bigoplus_{j>0} I^j/I^{j+1}$ and
by ${\underline a}_i (G) = \max \{j \mid H_{G_{++}}^i (G)_j \not =
0 \}$, it was proved that the form ring  $G$ has depth 0, and
 ${\underline a}_0 (G) < {\underline a}_1 (G) < {\underline a}_2 (G)=
4$. Now $a_*^2(G) = \max \{{\underline a}_i(G) : i=0,1,2 \} = 4$
according to \cite[Lemma 2.3]{H}, and then the short exact
sequences
$$0 \to R_{++} \to R \to A \to 0,$$
$$0 \to R_{++} (1) \to R \to G \to 0,$$
\noindent where $R_{++}= \bigoplus_{j>0} I^j$, show $a_*^2(R)= 4$.
The graded minimal free resolutions of $I^5$ and $I^6$ (computed
with CoCoa) are :
$$ 0 \to A(-37)^{15} \oplus A(-36)^5  \to A(-35)^{21} \to I^5 \to
0\,\,,$$
$$ 0 \to A(-44)^{15} \oplus A(-43)^{12}  \to A(-42)^{28} \to I^6 \to
0\,\,,$$ \noindent Then we may compute the polynomials
$P_\alpha(j)$, so the graded minimal free resolutions of $I^j$ for
$j \geq 5$ are:
$$ 0 \to A(-2-7j)^{15} \oplus A(-1-7j)^{7j-30}  \to A(-7j)^{7j-14} \to
I^j \to 0\,\,.$$

Furthermore, in this case we check that the bound can not be
improved because the resolution of $I^4$ is
$$ 0 \to A(-30)^{14} \to A(-28)^{15} \to I^4 \to 0\,\,.$$

}
\end{ex}

\medskip
\noindent ${\bf Open \; Question} \;$ Let $I$ be a homogeneous
ideal generated by forms in degree $d$. Denote by $l = l(I)$, $s =
n -\depth_{(\frak m R)} (R)$, $c= a_*^2(R)$. By Proposition
\ref{GG2}, there are integers
 $\{\alpha_{pi} \}$
 and polynomials
 $\{Q_{\alpha_{pi}}(j) \}$
 of degree $\leq l-1$
 such that for $j \gg 0$ the graded minimal free
resolution of $I^j$ is
$$ 0 \to D_s^j \to \dots \to D_0^j \to I^j \to 0 \,\, ,$$
\noindent with $D_p^j = \bigoplus_i A(-\alpha_{pi}
-dj)^{\beta_{pi}^j}$ and ${\beta_{pi}^j}= Q_{\alpha_{pi}} (j)$. In
some particular cases, we have shown that this holds for $j \geq
c+1$. The question is if this bound holds for any equigenerated
ideal $I$.

\bigskip
\subsection{General case}

\medskip
We may also study the minimal graded free resolutions of the
powers of any arbitrary homogeneous ideal in the polynomial ring
although in this case the asymptotic result is not so nice. Our
approach will be based on a detailed study of the proof of
\cite[Theorem 3.4]{CHT}. We need to introduce some notation.

\medskip
Let $I$ be a homogeneous ideal in $A$ generated by $r$ forms of
degrees $d_1, \dots, d_r$. Let us consider $S= k[{\x},{\y}]$ the
polynomial ring with $\deg X_i = (1, 0)$, $\deg Y_j = (d_j, 1)$,
and let $S_2= k [{\y}]$. For any finitely generated bigraded
module $L$ over $S_2$, let us define the set
     $$ \delta_L(j) = \{i : L_{(i,j)} \not = 0\}.$$
Given $\underline c=(c_1, \dots, c_r) \in \Bbb N^r$, let us denote
by $v (\underline c) = d_1 c_1 + \dots+ d_r c_r$ and by $\mid
\underline c \mid = c_1 + \dots + c_r$. Given a set $C \subset
\Bbb N^r$, $C + \Bbb N^r$ denotes the set of points of $\Bbb N^r$
of the type $\underline c + \underline c'$ with $\underline c \in
C$, $\underline c' \in \Bbb N^r$. Then we have:

\begin{lem}
\label{Gg1} Let $L$ be a finitely generated bigraded $S_2$-module.
Then there are pairs $(\alpha_{i}, \beta_{i}) \in \Bbb Z^2$ and
finite subsets $C_i$ of $\Bbb N^r$, $1 \leq i \leq m$, such that
for any $j$
$$\delta_L(j) = \bigcup_i \{ v(\underline c) + \alpha_i:
   \underline c \not \in C_i +\Bbb N^r,
   \mid \underline c \mid =j - \beta_i \} .$$
Therefore, $$\dim_k L_{(l,j)} = \sum_{i=1}^m \#
 \{\underline c \in \Bbb N^r : \underline c \not \in C_i + \Bbb N^r,
\mid \underline c \mid = j - \beta_i, v(\underline c) = l -
\alpha_i\}.$$
\end{lem}

{\pf} As said, the proof is based on \cite[Theorem 3.4]{CHT}.
Given any finitely generated bigraded $S_2$-module $L$, there
exists a sequence of bigraded submodules \noindent
$$0 = L_0 \subset L_1 \subset \dots \subset L_{m-1} \subset L_m = L$$
\noindent of $L$ such that $M_i = L_i / L_{i-1} \cong S_2/{\frak
p_i}(-\alpha_i,-\beta_i)$, $1 \leq i \leq m$, with $\frak p_i$
homogeneous prime ideals in $S_2$. Note that $\dd_L(j) = \bigcup_i
\dd_{M_i}(j) = \bigcup_i \dd_{S_2/{\frak p_i}}(j - \beta_i) +
\alpha_i $, and so we can assume that $L$ is cyclic.

 Now let $L = S_2/J$, with $J \subset S_2$ a homogeneous ideal.
By fixing a term order $<$ in $S_2$, then $L$ has a $k$-basis
consisting of the classes of the monomials which do not belong to
the initial ideal ${\rm in}(J)$ of $J$. So we get $ \dd_{S_2/J}(j)
= \dd_{S_2/{\rm in}(J)} (j)$, and we may assume $J$ is a monomial
ideal.

        Let us write $ J = (Y_1^{c_{11}}\cdots Y_r^{c_{1r}},
\dots, Y_1^{c_{p1}}\cdots Y_r^{c_{pr}})$, and
 $\underline c_i = (c_{i1}, \dots, c_{ir})$ for $1 \leq i \leq p$.
For any $\underline c \in \Bbb N^r$, note that $Y_1^{c_1} \cdots
Y_r^{c_r} \in J$ if and only if there exists $i$ such that
$\underline c = \underline c_i + \underline c'$, for some
$\underline c' \in \Bbb N^r$, i.e. $\underline c \in C + \Bbb N^r$
, where $C= \{\underline c_1,\dots, \underline c_p\}$. Therefore
$$\dd_L(j) =\{v(\underline c) : \underline c \not \in C+ \Bbb N^r,
 \mid \underline c \mid = j \},$$
and we are done.$\B$

\medskip
   Now we can show the asymptotic minimal graded free
resolution of the powers of an arbitrary homogeneous ideal $I$.

\begin{prop}
\label{Gg2}
   Let $I$ be a homogeneous ideal in the polynomial ring $A= k[{\x}]$
minimally generated by forms $f_1, \dots, f_r$ of degrees $d_1,
\dots,d_r$. Then there are pairs $(\alpha_{pi}, \beta_{pi}) \in
\Bbb Z^2$ and sets $C_{pi} \subset \Bbb N^r$, for $0 \leq p \leq
s, 0 \leq i \leq k_p$, such that for $j$ large enough the graded
minimal free resolution of $I^j$ is
$$ 0 \to D_s^j \to \dots \to D_0^j \to I^j \to 0 \,\, ,$$
\noindent with $D_p^j = \bigoplus_{i, \underline c}
A(-\alpha_{pi}- v(\underline c))$, for $\alpha_{pi}$ and
$\underline c$ such that
  $\underline c \not \in C_{pi} + \Bbb N^r$ and
 $\mid \underline c \mid = j -\beta_{pi}$.

\end{prop}

{\pf} Let us consider the Koszul homology modules $H_p
=H_p({\underline{X}}, R)$ of the Rees algebra $R$ with respect to
${\x}$. For any $p$, $H_p$ is a finitely generated bigraded
$S_2$-module with $[H_p]_{(i,j)} = \Tor^A_p (k, I^j)_i$. By Lemma
\ref{Gg1}, there exist $(\alpha_{pi}, \beta_{pi}) \in \Bbb Z^2$
and sets $C_{pi} \subset \Bbb N^r$ such that
 $$\delta_{H_p}(j) = \bigcup_i \{ v(\underline c) + \alpha_{pi}:
   \underline c \not \in C_{pi} +\Bbb N^r,
   \mid \underline c \mid =j - \beta_{pi} \},$$
so we get the statement. $\B$

\medskip
Similarly to the equigenerated case, we can also prove that the
rank of the modules of the graded minimal free resolution of the
powers of an ideal behaves as a polynomial.

\begin{prop}
\label{Gg3} Let $I$ be a homogeneous ideal of the polynomial ring
$A= k[{\x}]$ minimally generated by forms $f_1, \dots, f_r$ of
degrees $d_1, \dots,d_r$, and set $l= l(I)$. For any pair
$(\alpha, \beta)$ and set ${C}$ in the previous proposition, there
exists a polynomial $Q(j)$ of degree $\leq l-1$ such that for any
$j \gg 0$
  $$Q(j) = \# \{\underline c :
 \underline c \not \in C + \Bbb N^r, \mid \underline c \mid = j -\beta
\}\,\, .$$
\end{prop}

{\pf}
 Let $J$ be the homogeneous ideal of $S_2$ generated by the
monomials $Y_1^{c_1} \cdots Y_r^{c_r}$ with $\underline c=(c_1,
\dots, c_r) \in C$. Then, denoting by $Q(j)$ the Hilbert
polynomial of the graded $S_2$-module $L=S_2/J (-\beta)$, we have
that for $j \gg 0$

\vspace{2mm}

\hspace{15mm} $Q(j)= \dim_k \,(S_2/J)_{j-\beta}$

\vspace{2mm}

\hspace{25mm}  $= \# \{\underline c :  \underline c \not \in C +
\Bbb N^r, \mid \underline c \mid = j -\beta \}. $

\vspace{2mm}

On the other hand, since $H_p$ is a finitely generated graded
$F$-module of dimension $\leq l$, there exists a polynomial $P(j)$
of degree $\leq l-1$ such that $P(j)= \dim_k \, [H_p]^j = \dim_k
\, \Tor^A_ p (k,I^j)$ for $j \gg 0$. Since $Q(j) \leq P(j)$ for $j
\gg 0$, we immediately get that $\deg Q(j) \leq l-1$. $\B$

\medskip

\begin{rem}
\label{Gg5}
  {\rm  Given a homogeneous ideal $I$ in the polynomial ring $A$
generated by $r$ forms of degree $d$, we considered its Rees
algebra $R$ with a natural bigrading. By defining $S=
k[{\x},{\y}]$ the polynomial ring with the bigrading deg$X_i = (1,
0)$, deg$Y_j = (d, 1)$, $R$ is a bigraded finite $S$-module in a
natural way. Now let $E$ be any bigraded finitely generated
$S$-module, and let consider the graded $A$-modules $E^j= \oplus_i
E_{(i,j)}$. By taking $E$ instead of $R$, we can get analogous
results for the asymptotic behaviour of the $A$-modules $E^j$. In
particular, by considering $E$ to be the form ring $G$ of $I$, the
integral clousure of the Rees algebra $\overline R= \bigoplus_j
\overline{I^j}$ or the symmetric algebra $\Sym_A (I)$ of $I$ we
have the asymptotic behaviour of $I^j/I^{j+1}$, $\overline{I^j}$
and $\Sym_j(I)$.
         }
\end{rem}

\chapter*{$\;\;\;$}

\newpage

\medskip

\addcontentsline{toc}{chapter}{Bibliography}

\markboth{BIBLIOGRAPHY}{BIBLIOGRAPHY}


\end{document}